\documentclass[12pt]{amsart}
\usepackage{amsfonts}
\usepackage{amsmath}
\usepackage{amsxtra}
\usepackage{amssymb,latexsym}
\usepackage[mathcal]{eucal}
\usepackage{amscd}
\usepackage{tikz-cd}
\usepackage{hyperref}

\makeindex

\newtheorem{theo}{{\bfseries Theorem}}[section]
\newtheorem{prop}[theo]{{\bfseries Proposition}}
\newtheorem{lem}[theo]{{\bfseries Lemma}}
\newtheorem{cor}[theo]{{\bfseries Corollary}}
\newtheorem{df}[theo]{{\bfseries Definition}}

\newtheorem{ex}[theo]{{\bfseries Example}}

\def \N {\mathbb N}

\def \Z {\mathbb Z}
\def \R {\mathbb R}
\def \Q {\mathbb Q}
\def \I {\mathbb I}

\def \B {\mathcal B}

\def \F {\mathcal F}

\def \H {\mathcal H}
\def \J {\mathcal J}

\def \T {\mathcal T}
\def \G {\mathcal G}

\def \X {\mathcal X}

\def \a {\alpha }
\def \b {\beta}

\def \ep {\epsilon}
\def \om {\omega}
\def \d {\delta}
\def \g {\gamma}
\def \i {\iota}
\def \r {\rho}
\def \s {\sigma}
\def \t {\tau}

\usepackage{amssymb,latexsym}
\usepackage[mathcal]{eucal}
\usepackage{amscd}
\usepackage{amsmath}
\usepackage{graphicx}
\usepackage{amsfonts}
\usepackage{amscd}

\numberwithin{equation}{section}
\begin{document}
\begin{titlepage}
\large
\title{\bfseries Trees and Homogeneous LOTS}
\author{Ethan Akin and Karel Hrbacek}
 \vspace{.7cm}

\address{Mathematics Department \\
    The City College \\ 137 Street and Convent Avenue \\
       New York City, NY 10031, USA     }
\date{September, 2021}
\normalsize

\begin{abstract} We describe those complete linearly ordered topological spaces $X$ which are homogeneous (=CHLOTS).
That is, $X$ is order isomorphic with any nonempty open interval in $X$.
Using countable tail-like ordinals as indices,
we build towers of distinct CHLOTS. Using tree constructions we are able to extend the towers and to describe an inductive procedure
which yields every CHLOTS.
\end{abstract}


\vspace{1cm}
\end{titlepage}
\maketitle

\tableofcontents

\section{\textbf{Introduction}}

At first glance the Cantor Set $C \subset [0,1]$   does not
appear to be homogeneous.  Aside from the maximum and minimum
there is the countable family of endpoint pairs, each of which
forms a gap $x^{-} < x^{+}$ such that the intersection $C \cap
[x^{+},\infty) $ is clopen in $C$.  Then there is the uncountable
residuum of points whose very existence is not obvious until the
bijection from the set of zero/one sequences to $C$ is revealed.
This bijection is in fact a homeomorphism of $C$ with a
topological group from which topological homogeneity is clear as the
automorphism group $H(C)$ contains all the translations of the
group.

The original impression of non-homogeneity comes from the order structure on $C$
inherited from $\R$.  Indeed, if you restrict to
$H_{+}(C)$, the subgroup of order preserving automorphisms, then
there are five equivalence classes with respect to the action:
$\{ max \},\  \{ min \}$, the set of left endpoints $\{ x^{-} \}$, the
set of right endpoints $\{ x^{+} \}$, and the remaining residual
subset. If you allow order reversing automorphisms, using
$H_{\pm}(C)$ which contains $H_{+}(C)$ as a subgroup of index two,
then the $max$ and $min$ pair up and the left and right endpoints
are equivalent, leading to three classes.  Ignoring the
extrema
 we focus on the distinction between the gap pairs and
the rest.

There is an interesting construction called the
\emph{Alexandrov-Sorgenfrey Double Arrow} \index{Double Arrow}
\index{Alexandrov-Sorgenfrey Double Arrow}
\begin{equation}\label{eq1.1}
\R' \quad =_{def} \quad \R \times \{ -1, +1 \}
\hspace{2cm}
\end{equation}
in which we denote by $t^{-}$ the point $(t,-1)$  and similarly
$t^{+} = (t,+1)$.  On  $\R'$ we introduce the
lexicographic ordering and use the associated order topology.  For
every $t \in \R$  $t^{-} < t^{+}$ is a gap pair in
$\R'$ and so every point is either a left or a right
endpoint.  Every closed, bounded subset of $\R'$  is
compact and so the space is locally compact and $\sigma$-compact.
The family of clopen intervals $\mathcal{B} = \{ [t^{+},s^{-}] : t
< s \ \mbox{in} \ \R \} $ is uncountable and so $\R'$ is not metrizable although it is clearly separable.
Since $\mathcal{B}$ is a basis for the
topology the space is zero-dimensional. The space $\R'$
is a famous example in part because the subset $\{ t^{-} : t \in
\R \}$ is order isomorphic with $\R$ and so the
order topology is the usual one on $\R$, but the
subspace topology induced from $\R'$ is the
nonmetrizable, not locally compact topology on $\R$ with
basis the right-closed, left-open intervals.

We denote by $\bullet \R' \bullet$ the two point
compactification obtained by attaching a minimum $m$ and a maximum
$M$ to $\R'$.  We call this space the \emph{Fat Cantor
Set}\index{Fat Cantor Set}.  For every $t < s$  in $\R$ there is an order
preserving homeomorphism
\begin{equation}\label{eq1.2}
f : \bullet \R' \bullet \ \ \rightarrow \ \ [t^{+},s^{-}].
\hspace{2cm}
\end{equation}
Clearly, the group $H_{\pm}(\R')$ acts
transitively on $\R'$.

Our work began with an analogy question: The Cantor Set is to the real line as the Fat Cantor Set is to what?

 There should be a
linearly ordered space $X$ with the order topology, that is, a
\emph{LOTS}\index{LOTS},which is connected, which contains the Fat Cantor Set and
which is homogeneous in the sense that for all $a < b $ in $X$
there exists an order preserving homeomorphism
\begin{equation}\label{eq1.3}
f : X \rightarrow (a,b).  \hspace{2.5cm}
\end{equation}

Our first thought was to use $\R \times J$ with $J = [-1,+1]
\subset \R$ with the lexicographic ordering.  This LOTS
is connected but it is not homogeneous.  However, equipped with
the lexicographic ordering the countably infinite product
\begin{equation}\label{eq1.4}
\R_{\omega} \quad = \quad  \R \times J \times J
\times ...  \hspace{1.5cm}
\end{equation}
with $\om$ \index{$\om$} the first infinite ordinal.
This is a connected and homogeneous LOTS, which we call a CHLOTS.
Define for $t \in J$
\begin{equation}\label{eq1.5}
j(t^{-}) \  = \  (t,-1,-1,...) \qquad \mbox{and} \qquad j(t^{+}) \
= \ (t,+1,+1,...). \hspace{.5cm}
\end{equation}
Then $j : [(-1)^{+},(+1)^{-}] \rightarrow \R_{\omega} $
is an order preserving, topological embedding onto a closed
subset.  Thus, $\R_{\omega}$ naturally contains the Fat
Cantor Set.

Looking for other examples of  CHLOTS led us to look at products like (\ref{eq1.4}) but
indexed by more general countable ordinals than $\om$. Over each CHLOTS we constructed a tower of distinct CHLOTS, increasing in size in a
suitable sense. The tower is
indexed by the countable ordinals, i.e. the ordinals less than $\Omega$, the first uncountable ordinal.

In this way we re-discovered a construction begun by Arens \cite{Ar} and \cite{Ar2} and extended by Babcock \cite{B}, see also \cite{Tr}.
Their procedure was also re-discovered by others, see e.g. \cite{C}.

We had nearly completed the initial phase of our work in 2001, \cite{AH}, when we were directed to the paper of
Hart and van Mill \cite{HvM} whose work is complementary to ours,
and of course to that of Arens and Babcock.
Hart and van Mill construct an uncountable family of distinct CHLOTS no two of
which have comparable size but all of which are bigger than
$\R$ but smaller than $\R_{\omega}$.  So their
class of examples extends horizontally where ours proceeds
vertically.

We now apply  trees to the study of CHLOTS.

A \emph{tree}\index{tree} is a partially ordered set such that the set of the predecessors of each point is well-ordered by the
induced order and so is isomorphic to an ordinal.
The relation between trees and general LOTS is well-known, e.g. \cite{T} and \cite{BLR}.
By using trees we develop a number of constructions for building CHLOTS and
extend the tower over a CHLOTS to one indexed by $\Omega \times \Omega$.
In addition, we describe an inductive tree
construction from which every CHLOTS can be obtained.

We would like to thank Richard Wilson for some helpful discussions
as we began this work.\vspace{.5cm}

We now provide a brief sketch of what follows. \vspace{.5cm}

In Section 2 we introduce the elementary properties of LOTS and ordinals.

A subset $J$ of a LOTS is \emph{convex} when $a < c < b$ in $X$ and $a,b \in J$ implies $c \in J$. An open or closed interval is
convex. A LOTS $X$ is \emph{order dense}\index{LOTS!order dense}
when $a < b$ in $X$ implies that the interval $(a,b)$ is infinite and $X$ is \emph{complete}\index{LOTS!complete} when every bounded set has
a supremum and infimum. A LOTS is connected iff it is a complete and order dense LOTS. An order dense LOTS $X$ is contained as a dense subset
in an essentially unique connected LOTS $\hat{X}$ called its \emph{completion}\index{LOTS!completion}.

An order dense LOTS which is not complete contains holes.  There is a \emph{hole} between $a$ and $b$ if there is a clopen convex set which contains
$a$ and not $b$. A LOTS has \emph{dense holes}\index{LOTS!dense holes} when between any pair of distinct points there are holes.

We call a LOTS \emph{unbounded}\index{LOTS!unbounded} when it has neither a maximum nor a minimum and \emph{bounded}\index{LOTS!bounded} when it has both.

There is a rough notion of \emph{size} \index{size}for LOTS. In comparing two
LOTS $X$ and $X_{1}$ we say that $X$ \emph{injects into} $X_{1}$
if there exists a, not necessarily continuous,  injective order
map from $X$ into $X_{1}$.
We say that
$X_{1}$ is \emph{bigger}\index{LOTS!bigger} than $X$ if $X$ order injects into $X_{1}$ but
not vice-versa. For example, any CHLOTS which is not $\R$
itself is bigger than $\R$. If neither injects into the other we say that they are not comparable.

In outlining the use of ordinals we pay special attention to those which are tail-like.

An ordinal $\a$ is \emph{tail-like}\index{ordinal!tail-like} if $\b < \a$ implies $\b + \a = \a$, or,
equivalently, $\a = \om^{\gamma}$ for some ordinal $\gamma$.
For ordinal exponentiation, $\gamma$ countable implies $\om^{\gamma}$ is countable and so there are uncountably many countable
tail-like ordinals, indexed by the countable ordinals $\g$. We describe the Cantor Normal Form which writes any ordinal uniquely
as a  sum of a finite non-increasing sequence of tail-like ordinals.
 \vspace{.5cm}

In Section 3 we define various transitivity and
homogeneity properties for LOTS.

A LOTS $X$ is \emph{transitive}\index{LOTS!transitive} when $a, b \in X$ implies there exists an order isomorphism $f$ on $X$ with $f(a) = b$,
i.e. the group of order automorphisms acts transitively on $X$.  A LOTS $X$ is \emph{homogeneous}\index{LOTS!homogeneous}
when $X$ is order isomorphic with any
nonempty, open, convex subset of $X$ in which case we call it a \emph{HLOTS}\index{HLOTS}.
In particular, we call $X$  a \emph{CHLOTS}\index{CHLOTS} when it is a connected homogenous LOTS. If a HLOTS is not complete, then it has dense holes
and we call it an incomplete homogeneous LOTS,
an \emph{IHLOTS}\index{IHLOTS}. An IHLOTS is order-dense and its completion is a CHLOTS.

Any HLOTS is first countable and $\s$-bounded. A CHLOTS is locally compact and $\s$-compact. \vspace{.5cm}

In Section 4 we build over a given HLOTS $X$,  a tower  indexed by the
countable ordinals, i.e. by $\Omega$ \index{$\Omega$}, the first uncountable ordinal.

For a LOTS $X$  and  a positive ordinal $\a$ define
 $X^{\a}$ to be the set of maps from $\a$ to $X$, thought of as the lexicographically ordered product. In $X$
 we select a nontrivial closed interval $J $ and then define $X_{\a} = \{  s \in X^{\a} : s(i) \in J$ for all $i > 0 \}$.
 If $X$ is order dense then $X^{\a}$ and $X_{\a}$ are order dense. If, in addition, $X$ is complete, then $X_{\a}$ is complete
 and so is connected.

 Call a LOTS $X$ \emph{$\R$-bounded}\index{LOTS!$\R$-bounded} if there exists an order injection of
 $X$ into $\R_{\d}$ for some countable ordinal $\d$. \vspace{.5cm}

 \begin{theo}\label{introtheo1} Assume that $X$ is a HLOTS.

 \begin{enumerate}
 \item[(a)] If $\a$ is a countable, tail-like ordinal, then $X_{\a}$ is a HLOTS and so its completion $\widehat{X_{\a}}$ is a CHLOTS.
 In particular, if $X$ is a CHLOTS, then $X_{\a}$ is a CHLOTS for each countable, tail-like ordinal.

 \item[(b)] Assume $X$ is a CHLOTS, e.g. $X = \R$. If $\a$ and $\b$ are positive ordinals with $\a > \b$, then
 $X_{\a}$ is bigger than $X_{\b}$.  That is, there is an order injection from $X_{\b}$ into $X_{a}$ but no order injection in the
 other direction. In addition, $X_{\a}$ is not homeomorphic to $X_{\b}$.

\item[(c)] Assume $X$ is an $\R$-bounded IHLOTS.
If $\a$ and $\b$ are sufficiently large countable, tail-like  ordinals with $\a > \b$, then
 $\widehat{X_{\a}}$ is bigger than $\widehat{X_{\b}}$.

\end{enumerate}
\end{theo}\vspace{.5cm}

Part (a) is essentially the result of Arens \cite{Ar}, \cite{Ar2} and Babcock \cite{B}. In addition, they showed that for $X = \R$ the elements of
the tower at different heights are not isomorphic.  \vspace{.5cm}

In Section 5 we describe the definitions and elementary results for  trees.

Let $T$ be a tree. For a vertex $p \in T$ we let $A_p$ denote the set of predecessors of $p$. This is isomorphic to a unique ordinal $o(p)$, the
\emph{order} of $p$. We let $T_p$ consist of $p$ together with all its successors.
The immediate successors of $p$ are those $q \in T_p$ with $o(q) = o(p) + 1$; the set of all immediate successors of $p$ is denoted $S_{p}$.
For any subset $A$ of $T$ the \emph{height}\index{tree!height} $h(A)$
is the smallest ordinal greater than $o(p)$ for all $p \in A$.

A subset $T_1 \subset T$ is a \emph{subtree} if $p \in T_1$ implies that all the $T$ predecessors lie in $T_1$. So, for example, $T_p$ is
not a subtree if $o(p) > 0$. On the other hand, the \emph{truncation} $T^{\a} = \{ p \in T : o(p) < \a \}$ is a subtree.

Our trees are all assumed to be at least \emph{semi-normal}\index{tree!semi-normal} meaning:
\begin{itemize}
\item There is a unique root, $0 \in T$, with $o(0) = 0$.

\item For all $p \in T$, the set $S_p$ is either empty or contains at least two points.

\item If for $p, q \in T$, $A_p = A_q$ and $o(p) = o(q)$ is a limit ordinal, then $p = q$.
\end{itemize}

Notice that  $p, q \in T$, $A_p = A_q$ and $o(p) = o(q) = \b + 1$ iff $p, q \in S_r$ for some $r \in T$ with $o(r) = \b$.

A tree is \emph{normal}\index{tree!normal} if, in addition:
\begin{itemize}
\item  If $p \in T$ and $\a$ is an ordinal with $o(p) < \a < h(T)$, then there exists a successor $q$ of $p$ with $o(q) = \a$.
\end{itemize}

A \emph{branch}\index{tree!branch} is a maximal linearly ordered subset $x \subset T$. We denote by $X(T)$ the \emph{branch space}\index{tree!branch space},
i.e. the set of branches of $T$. We let $x_i \in T$ denote the element of $x$ with $o(x_i) = i$.

We call the tree $T$ \emph{$\Omega$-bounded}\index{tree!$\Omega$-bounded} when $h(x) < \Omega$ for every $x \in X(T)$. This implies $h(T) \le \Omega$,
but a tree of height $\Omega$ can be $\Omega$-bounded.

The tree $T$ is \emph{bi-ordered}\index{tree!bi-ordered} when the successor set $S_p$ has the structure of a LOTS for all $p \in T$.
In that case there is an induced
order on $X(T)$. For distinct branches $x, y$ we write $x < y$ when for some $\ep$, $x_{\ep} = y_{\ep}$ and $x_{\ep + 1} < y_{\ep +1}$. The
latter ordering is the LOTS ordering in $S_p$ with $p = x_{\ep} = y_{\ep}$. With the induced order $X(T)$ is a LOTS.
\begin{itemize}
\item $X(T)$ is order dense if either every $S_p$ is order dense, or else, every $S_p$ is unbounded and $h(T)$ is a limit ordinal.

\item $X(T)$ is complete if every $S_p$ is complete and $o(p) > 0$ implies $S_p$ is bounded.
\end{itemize}

A bijection $h : T_1 \to T_2$ is a \emph{tree isomorphism} when it preserves both orders. A tree isomorphism induces an order isomorphism between
the branch spaces.

A tree $T$ is \emph{reproductive}\index{tree!reproductive} if for all $p \in T$ there is an isomorphism from $T$ to $T_p$. A reproductive
tree has height a tail-like ordinal. For a reproductive tree, every $S_p$ is order isomorphic to the LOTS $S_0$.\vspace{.5cm}

\begin{theo}\label{introtheo2} If $T$ is an $\Omega$-bounded, reproductive tree with $S_0$ a HLOTS, then $X(T)$ is a HLOTS. \end{theo}
\vspace{.5cm}

In Section 6 a number of tree constructions are presented.

If $X$ is a LOTS and $\a$ is a positive ordinal, then the \emph{simple tree}\index{tree!simple} on $X, \a$ has vertices $X^i$
at level $i < \a$, and for $s \in X^i$ the predecessor at level $j < i$ is the restriction $s|j$. The simple tree has height $\a$ and
we can identify the branch space with $X^{\a}$.

For $p \in X^i, q \in X^j$ we define $p + q \in X^{i + j}$ by
\begin{equation}\label{eq1.6}
(p + q)(k) \ = \ \begin{cases} \ p(k)  \ \ \text{for} \ k < i, \\ q(k \setminus i) \ \text{for} \ i \le k < i+j. \end{cases}
\end{equation}
Here $k \setminus i = \{ \ell : i \le \ell < k \}$ is identified with the ordinal with which it is isomorphic so that, e.g.
$(i + j) \setminus i = j$.

A subtree $T$ of the simple tree on $X$ is called an \emph{additive tree}\index{tree!additive} if
\begin{itemize}
\item $p, q \in T \quad \Longleftrightarrow \quad p + q \in T.$
\end{itemize}
The map $a_p : T \to T_p$ given by $a_p(q) = p + q$ is then an isomorphism from $T$ to $T_p$ and so an additive tree is reproductive.
In particular, the height of an additive tree is tail-like.

For a CHLOTS $X$ let $\bullet X \bullet$ be the two point compactification of $X$ obtained by adding a minimum $m$ and a maximum $M$.
If $p \in X^{\a}$ we define $\hat p \in (\bullet X \bullet)^{\a + 1}$ by
 \begin{equation}\label{eq1.7}
 \hat p (0) = m, \quad \hat p (i) = sup \{ p(j) : j < i \} \ \ \text{for} \ 0 < i \le \a.
 \end{equation}

 We call $p$ \emph{sharply increasing}\index{sharply increasing} when for all $i < \a$, $\hat p (i) < p(i)$.
 We define the \emph{ order tree}\index{order tree} $T(X)$ to be the subtree of the simple tree on $X, \Omega$
 whose elements of order $\a$ are the bounded, sharply increasing elements of $X^{\a}$.  \vspace{.5cm}

\begin{theo}\label{introtheo3} If $X$ is a CHLOTS the order tree $T(X)$ is an $\Omega$-bounded reproductive tree of height $\Omega$ with
$S_0 \cong X$. The branch space of $T(X)$ is an IHLOTS with completion a CHLOTS. \end{theo}\vspace{.5cm}

We will denote by $a(X)$ the completion of the branch space of $T(X)$.\vspace{.5cm}

In Section 7 we build the Double Tower over a CHLOTS $X$.

Inductively, we let $a_0(X) = X$ and $a_{\a + 1}(X) = a(a_{\a}(X))$. If $\a$ is a countable limit ordinal we define the CHLOTS $a_{\a}(X)$
as an inverse limit of  $\{ a_{\b}(X) : \b < \a \}$. For $(\a,\b) \in \Omega \times \Omega$ we obtain the CHLOTS $(a_{\a}(X))_{\b}$.
\vspace{.5cm}

\begin{theo}\label{introtheo4} If $(\a',\b') > (\a,\b)$ in the lexicographical ordering on $\Omega \times \Omega$, then
$(a_{\a'}(X))_{\b'}$ is bigger than $(a_{\a}(X))_{\b}$. \end{theo} \vspace{.5cm}

In particular, $a(X)$ is bigger than $X_{\b}$ for any countable ordinal $\b$.  It follows that $a(X)$ is not $\R$-bounded.\vspace{.5cm}

In Section 8 for any given CHLOTS $X$ we construct an $\Omega$-bounded, additive subtree of the simple tree on $\Z, \Omega$ the completion of whose
branch space is isomorphic to $X$. From it we obtain the following.

 \begin{theo}\label{introtheo5} If $X$ is a LOTS, then the following are equivalent.
 \begin{itemize}
 \item $X$ is a CHLOTS.

 \item There exists an $\Omega$-bounded, additive subtree of the simple tree on $\Z, \Omega$ the completion of whose
branch space is isomorphic to $X$.

 \item There exists an IHLOTS $D$, dense in $\R$ and an $\Omega$-bounded, additive subtree of the simple tree on $D, \Omega$ the completion of whose
branch space is isomorphic to $X$.

 \item There exists  an $\Omega$-bounded, reproductive tree with $S_0$ a HLOTS  the completion of whose
branch space is isomorphic to $X$.
\end{itemize} \end{theo}

In addition, we show that a CHLOTS $X$ is $\R$-bounded if and only if there exists a tree of countable height the completion of whose
branch space is isomorphic to $X$.
\vspace{.5cm}

In Section 9 we describe the Hart-van Mill results.

A subset $Y \subset \R$ is a \emph{Bernstein subset}\index{Bernstein subset} (a \emph{B-set}) if it meets every Cantor subset of $\R$.
It is a \emph{Bi-Bernstein subset}\index{Bi-Bernstein subset} (a BB-set) when, in addition, its complement is a Bernstein subset.

Let $\G$ denote the group of positive, rational, affine transformations on $\R$, i.e. maps of the form $t \mapsto at + b$ with $a,b \in \Q$
and $a > 0$. We observe that nonempty subset of $\R$ which is $\G$ invariant is a HLOTS.

Let ${\mathbf c} = 2^{\aleph_0}$ denote the cardinality of $\R$.

  \begin{df}\label{introdf6} A Hart-van Mill collection $\H$ with base $V$ is a set of cardinality ${\mathbf c}$ which satisfies the following
  conditions.
  \begin{itemize}
  \item $\H \cup \{ V \}$ is a collection of pairwise disjoint BB-subsets of $\R$ each of which is  $\G$ invariant.

  \item $\Q \subset V$.

  \item $Y \in \H$ implies $-Y \in \H$ with $-Y = \{ -x : x \in Y \}$ distinct from $Y$.

  \item Assume $Y \in \H$ and $f$ is an order automorphism of $\R$. If $f(Y) \setminus Y$ has cardinality ${\mathbf c}$, then
  $f(Y) \cap V$ has cardinality ${\mathbf c}$.
  \end{itemize} \end{df}

  For any $\J \subset \H$, let $X(\J)$ denote the complement of the union of the elements of $\J$. Because it contains $V$ and so $\Q$, and
  because it is $\G$ invariant, $X(\J)$ is a HLOTS.

  Hart and van Mill show that such collections exist and prove the following.\vspace{.5cm}

  \begin{theo}\label{introtheo7} For a Hart-van Mill collection, let $\J, \J_1$ be distinct subsets of $ \H$ with associated
  HLOTS $X, X_1$. For countable tail-like ordinals $\a, \a_1$
 the CHLOTS $\widehat{X_{\a}}$ and $\widehat{(X_1)_{\a_1}}$ are not order isomorphic.
  If $\a = \a_1 = \om$, then the two do not even have the same size. \end{theo} \vspace{.5cm}

  It follows from Theorem \ref{introtheo1} that for any subset $\J$ of $\H$, we obtain a tower $\widehat{X(\J)_{\a}}$, indexed by
  the countable tail-like ordinals $\a$. Each tower consists of CHLOTS and is nondecreasing in size. For the $2^{{\mathbf c}}$ distinct
  subsets, Theorem \ref{introtheo7} implies that the towers are distinct, i.e. no two contain any pairwise isomorphic elements.
    \vspace{.5cm}

In Section 10 we conclude with some results on complete, perfect, zero-dimensional LOTS,
that is, complete LOTS with no isolated points for which the clopen
intervals form a basis. Thus, we return to generalizations of the Cantor Set which motivated our original inquiry.

If $T$ is a normal tree of type $2 = \{ 0, 1 \}$ and of height a limit ordinal, then the branch space is a compact, perfect, zero-dimensional
LOTS which is first countable if the tree is $\Omega$-bounded. We let $\bar 0$ and $\bar 1$ denote the minimum and maximum branches
in $X(T)$, so that $\bar 0_i = 0$, or $\bar 1_i = 1$, for every successor ordinal $i < h(\bar 0)$, resp. $i < h(\bar 1)$.

We focus on those zero-dimensional LOTS which satisfy the \emph{clopen interval condition} \index{clopen interval condition} which
holds when any two clopen intervals are isomorphic. In particular, we prove the following.

\begin{theo}\label{introtheo8} If $T$ is an additive tree of type $2 = \{ 0, 1 \}$ with height $\a$ a  tail-like ordinal,
such that $h(\bar 0) = h(\bar 1) = \a$, then the branch space $X(T)$ is a compact, perfect, zero-dimensional LOTS which satisfies
the clopen interval condition.  In addition, if $\a$ is countable, then $X(T)$ is topologically homogeneous. \end{theo}  \vspace{.5cm}

In particular, from this we recover results of Maurice \cite{M}, \cite{M2} on the product space $2^{\a}$ for $\a$ a countable tail-like
ordinal.

\vspace{1cm}

\section{\textbf{LOTS and Ordinals}}

\subsection{The Category of LOTS}

For a totally ordered set $X$ we will use the usual notation
$(a,b)$ for the \emph{open interval} and $[a,b]$ for the
\emph{closed interval} with \emph{endpoints} $a \leq b$ in $X$ and
we will write $(a,+\infty)$ and $(-\infty,b)$ for unbounded open
intervals. An interval is called \emph{proper} when it
contains more than one point. A subset $A$ of $X$ is
\emph{bounded} \index{subset!bounded}when $A \subset [a,b]$ for some $a,b \in X$. So $X$
itself is bounded iff it has a maximum and a minimum (hereafter
$max$ and $min$). Somewhat abusively, we will call $X$
\emph{unbounded}\index{unbounded}\index{LOTS!unbounded} when it has neither $max$ nor $min$.

As usual, we will let $\R, \Q,\Z, \N$ stand for the set of reals, rationals, integers and non-negative integers, respectively. In particular, $0 \in \N$.
They are all equipped with the usual order.

A linearly ordered topological space (hereafter a LOTS) is a
totally ordered set equipped with the order topology\index{LOTS}.  That is,
the set of open intervals is a base for the topology.  The
topology is Hausdorff and the order and topology properties are
closely related.

A LOTS $X$ is called \emph{order complete}\index{LOTS!order complete} \index{LOTS!complete}\index{complete}
(hereafter \emph{complete})when every bounded subset $A$ has a
supremum and an infimum (denoted $sup \  A$ and $inf \  A$), or,
equivalently, when every closed bounded interval $[a,b]$ is
compact. Thus, a complete LOTS is locally compact. In particular,  $X$ is compact iff it is complete and
bounded. On the other hand, local compactness is not sufficient for completeness. For example, $X = \R \setminus \Z$ is locally
compact, but not complete.

A LOTS $X$ is called \emph{order dense}\index{LOTS!order dense} \index{order dense}
when between any two
points of $X$ there lie other points of $X$, or, equivalently,
when every nonempty open interval in $X$ is infinite. $X$ is connected iff
it is complete and order-dense, in which case, every subinterval
is connected. If $X$ is not order-dense, then there exists a
\emph{gap pair}\index{gap pair}, $a < b$ in $X$ with $(a,b) = \emptyset$.  The point
$a$ is then called the \emph{left endpoint}\index{gap pair!endpoints}\index{endpoint!left} and $b$ is called the
\emph{right endpoint}\index{endpoint!right} of the pair.  By convention the $max$ of
$X$, if it exists, is a left endpoint and the $min$ is a right
endpoint.  Thus, $a \in X$ is a left (or right) endpoint iff the
closed interval $(-\infty,a]$ (resp. $[a,+\infty)$) is open.  A
point $a$ is
\emph{isolated}, i.e. $\{ a \}$ is clopen, iff it is both a left and a right endpoint.

\begin{lem}\label{lem2.0} A subset $A$ of a LOTS $X$ is closed iff  for all $B \subset X$,
$a =  sup  (A \cap B)$ or $a =  inf  (A \cap B)$ implies $a \in A$. \end{lem}

\begin{proof} Assume the conditions hold and that $a$ is a  point of the closure of $A$. Let
 $a_1 =  sup  (A \cap (-\infty,a)) $ and $a_2 =  inf (A \cap (a,\infty))$. If $a = a_1$ or $a = a_2$, then $a \in A$ by hypothesis.
 If $a$ is neither of these, then $a \in (a_1,a_2)$ and  $ A \cap (a_1,a_2) \setminus \{ a \} = \emptyset $. So $a \in \overline{A} $
 implies $a \in A$.

 The converse is clear.

\end{proof} \vspace{.5cm}

A subset $A$ of a LOTS $X$ is  \emph{convex}\index{convex}\index{subset!convex}  if $a <
c < b$ and $a, b \in A$ imply $c \in A$.
Intervals are convex subsets and if $X$ is complete then every
convex subset is an interval.  If a convex set $J$ contains at least three points  $a <
c < b$, then the nonempty open interval $(a, b)$ is a subset of $J$ and so $J$ has a nonempty
interior.

A \emph{Dedekind cut}\index{Dedekind cut} in $X$ is a
partition of $X$ by a pair of nonempty disjoint sets
$(A_{1},A_{2})$ such that for all $a < b$ in $X$, $b \in A_{1}$
implies $a \in A_{1}$ and so $a \in A_{2}$ implies $b
\in A_{2}$.  A \emph{hole}\index{hole} between $a$ and $b$ is a
Dedekind cut $(A_{1},A_{2})$ with $a \in A_{1}$, $b \in
A_{2}$ such that $A_{1}$ and $A_{2}$ are clopen.  For example, if
$a < b$ is a gap pair then $((-\infty,a],[b,+\infty))$ is a hole between $a$ and $b$.  On the other
hand, while the LOTS $\Q$ of rational numbers is order dense, every
irrational number creates a hole in $\Q$.  We say that a LOTS $X$
has \emph{dense holes}\index{dense holes} if there is a hole between every pair
$a < b$ in $X$.  Thus $\Q$ and the Cantor Set have dense
holes.

\begin{lem}\label{lem2.0a} Let $X_1$ be a subset of a LOTS $X$.
If $X$ is order dense and the subset $X_1$ is dense in a convex subset of $X$, then, regarded as a
LOTS in its own right, $X_1$ is order dense. \end{lem}

\begin{proof} If $X_1$ is dense in $J$, a convex subset of $X$, and $a < b$ in $X_1$, then because $X$ is order dense
the interval $(a,b) $ in $X$ is infinite.  Since $J$ is convex, $(a,b)$ is a subset of $J$. Between any two points of $J$
there are points of $X_1$ because $X_1$ is dense in $J$. Hence, $(a,b) \cap X_1$ is infinite.

\end{proof}

The \emph{reverse}\index{reverse} \index{LOTS!reverse} of a LOTS $X$, denoted $X^{*}$, is the set
equipped with the reverse order.  Clearly, the intervals, and so
the topology, for $X$ and $X^{*}$ are the same.

A function $f: X_{1} \rightarrow X_{2}$ between LOTS is an
\emph{order map}\index{order map} if it is order preserving, i.e. $a \leq
b$ implies $f(a) \leq f(b)$, while $f$ is called an
\emph{order* map}\index{order* map} if it is order reversing, i.e. $f: X^{*}_{1}
\rightarrow X_{2}$ is an order map or, equivalently, if $f: X_{1}
\rightarrow X^{*}_{2}$ is an order map.  An injective (or
surjective) order map is called an \emph{order injection}\index{order injection} (resp.
an \emph{order surjection}\index{order surjection}).  A bijective order map is called an
\emph{order isomorphism}\index{order isomorphism} or just an \emph{isomorphism}.  This is the isomorphism concept for the
category of LOTS with order maps.  We say that order isomorphic
LOTS $X_{1}, X_{2}$ have the same \emph{order type}\index{same order type} and we write
$X_{1} \cong X_{2}$\index{$X_{1} \cong X_{2}$}.

While an order isomorphism is a homeomorphism, an order map need
not be continuous.  In particular, if $X_{1}$ is a subset of a
LOTS $X$, then with the induced order $X_{1}$ is itself a LOTS.
However, the LOTS topology on $X_{1}$ need not be the topology induced
from $X$ and the inclusion map might not be continuous. For
example, consider $X_{1} = [-\infty,0]\cup(1,\infty)$ in $X = \R$. Clearly, the disconnected subset $X_{1}$ is order isomorphic
to $\R$ itself.

We will call a map $f: X_{1} \rightarrow X_{2}$ an \emph{order
embedding}\index{order embedding} if it is an order map which is a topological embedding,
i.e. $f: X_{1} \rightarrow f(X_{1})$ is a homeomorphism with the
topology on $f(X_{1})$ induced from $X_{2}$.

\begin{prop}\label{prop2.1} Let $f: X_{1} \rightarrow X_{2}$ be an order map.
\begin{enumerate}
\item[(a)] Assume $f$ is surjective.
\begin{itemize}
\item[(i)] If each point-inverse is closed,
then $f$ is continuous. If each point-inverse is compact, then the map $f$ is closed as well as continuous.
If each point-inverse is compact and $X_1$ is complete, then the map $f$ is topologically proper, i.e. the preimage
of every compact subset of $X_{2}$ is a compact subset of $X_{1}$.

\item[(ii)] If $X_{2}$ is order dense, then $f$ is continuous.

\item[(iii)] If $X_{2}$ is unbounded, then $X_{1}$ is unbounded and the preimage of every
bounded subset of $X_{2}$ is a bounded subset of $X_{1}$.

\item[(iv)] If $X_{1}$ is
complete, $X_{2}$ is unbounded and $f$ is continuous, then $f$ is
topologically proper.
\end{itemize}

\item[(b)]Assume $f$ is injective.
If $f$ is continuous, then it is an
order embedding.  This occurs if one of the following holds.
\begin{itemize}
\item[(i)] The image $f(X_{1})$ is convex in $X_{2}$.

\item[(ii)]$X_{2}$ is complete and $f(X_{1})$ is closed in $X_{2}$.

\item[(iii)] $X_{2}$ is order dense and $f(X_{1})$ is dense in $X_{2}$.
\end{itemize}

\item[(c)]Let $A$ be a subset of $X_{1}$.  If $A$ is bounded in $X_{1}$
then the image $f(A)$ is bounded in $X_{2}$.  If $f$ is continuous and
$x =  inf  A$ (or $=  sup  A$) then $f(x) =  inf  f(A)$ (resp. $=  sup  f(A)$) in $X_{2}$. Conversely, if for every bounded
subset $A$ of $X_1$  $x =  inf  A$ (or $=  sup  A$) implies $f(x) =  inf  f(A)$ (resp. $=  sup  f(A)$) in $X_{2}$, then $f$ is continuous.

\item[(d)] If $f$ is surjective, then there exists a map $g: X_2 \to X_1$ such that
$f \circ g = 1_{X_2}$. Any such map $g$ is a (not usually continuous) order injection.

\item[(e)] If $X_1$ is order dense and $f$ is injective on a dense subset $D$ of $X_1$, then $f$ is injective.
\end{enumerate}
\end{prop}

\begin{proof} (a): If $f(a_{1}) = a_{2}$ and
$f(b_{1}) =  b_{2}$ then
\begin{equation}\label{eq2.1}
f^{-1}((a_{2},b_{2})) = (a_{1},b_{1}) \setminus (f^{-1}(a_{2})\cup
f^{-1}(b_{2}))
\end{equation}
which is open if $f$ has closed point inverses.

Now with $f$ continuous assume that $A$ is a subset of $X_1$ with $y$ a limit point of $f(A)$ not in $f(A)$.
By replacing $A$ by $A \cap f^{-1}((-\infty,y])$ or by $A \cap f^{-1}([y,\infty))$ (which are closed when $A$ is)
we may assume that $y =  sup  f(A)$ or $=  inf  f(A)$. Assume the first. Since $f^{-1}(y)$ is compact it has an infimum
which we denote $x$. So  $a \in A$ implies $f(a) < y = f(x)$ and so $a < x$. For any $z < x$, $f(z) < f(x) = y$ because $x =  inf  f^{-1}(y)$.
Since $y =  sup f(A)$ there exists $a_z \in A$ such that $f(z) < f(a_z) < f(x)$ and so $z < a_z < x$. This means that
$x =  sup  A$. Since $x \not\in A$, it follows that $A$ is not closed. Contrapositively,
when all point inverses are compact, $A$ closed
implies that $f(A)$ contains all its limit points and so is closed.

Now assume that all point inverses are compact and $X_1$ is complete.  If $B$ is a compact subset of $X_2$, then it has a supremum $y$.
Since $f^{-1}(y)$ is compact, it has a supremum $x$.  It is clear that $x =  sup  f^{-1}(B)$. Similarly, $f^{-1}(B)$ has an infimum.
Since $f$ is continuous, $f^{-1}(B)$ is closed as well as bounded and so is compact by completeness of $X_1$.

(ii): If $f(x) \in
(a_{2},b_{2})$ and $X_{2}$ is order dense then there exist
$a_{3},b_{3}$ in $X_{1}$ such that
\begin{equation}\label{eq2.2}
\begin{split}
 a_{2} < f(a_{3}) < f(x) < f(b_{3}) < b_{2}\\
\mbox{and so}\hspace{1in} \\
x \in (a_{3},b_{3}) \subset f^{-1}((a_{2},b_{2})).
\end{split}
\end{equation}
Thus, the latter is a neighborhood of $x$ in $X$.

(iii), (iv): If $M$ is an upper bound for $A \subset X_{1}$ then $f(M)$ is an
upper bound for $f(A)$.  In particular, if $M =  max  X_{1}$ then
$f(M) =  max  X_{2}$ since $f$ is surjective.  So if $X_{2}$ is
unbounded then $X_{1}$ is.  Furthermore, if $B$ is bounded above
in $X_{2}$, then because $X_{2}$ has no $max$ and $f$ is
surjective, there exists $a \in X_{1}$ such that $y < f(a)$ for
all $y \in B$. Hence, $x < a$ for all $x \in f^{-1}(B)$. If  $B$
is compact, then it is closed and bounded in $X_{2}$.  If $f$ is
continuous, then $f^{-1}(B)$ is closed as well as bounded and so is
compact if $X_{1}$ is complete.

(b): Let $A = f(X_{1}) \subset X_{2}$.  The order injection $f$ is
an order isomorphism of $X_{1}$ with $A$ regarded as LOTS and so
is a homeomorphism.  The problem concerns the comparison between
the order topology on $A$ and the topology induced from $X_{2}$.
If $a,b \in A$ then the interval $(a,b)$ in $A$ is the
intersection of $A$ with the corresponding interval in $X_{2}$.
Hence the topology on $A$ is included in the topology induced from
$X_{2}$ and the two topologies agree exactly when the inclusion is
continuous.

If $a \in A, x \in X_{2}\setminus A $ and $a < x$ then
$(-\infty,x) \cap A $ is a neighborhood in the induced topology.
If $x$ is an upper bound for $A$ then the intersection is $A$.
Otherwise, we require a point $\tilde{a} \in A$ with $a <
\tilde{a}$ such that $(-\infty,x) \cap A \supset
(-\infty,\tilde{a}) \cap A$.  Such a point exists iff the
following condition holds:
\begin{equation}\label{eq2.3}
\begin{split}
a \in A,\  x \in X_{2},\  a < x,\  (a,x] \cap A = \emptyset ,\  [x,\infty) \cap A \not= \emptyset  \\
\Longrightarrow  \hspace{2in}\\
\exists b \in A\ \mbox{such that}\   x \leq b\ \mbox{and}\   (x,b)
\cap A = \emptyset.
\end{split}
\end{equation}

This condition and its analogue for the reverse orders are those
required for the two topologies to agree.

If $A$ is convex then $(a,x] \cap A = \emptyset $ implies
$[x,\infty) \cap A = \emptyset $.
If $X_{2}$ is order dense and
$A$ is dense in $X_{2}$ then $a < x $ implies $(a,x) \cap A \not=
\emptyset$.  So the conditions hold vacuously in these cases.

If
$X_{2}$ is complete and $A$ is closed in $X_{2}$, then $b =  inf [x,\infty) \cap A $
is a point of $A$ and $(x,b) \cap A =
\emptyset $.  So the conditions hold in this case as well.

(c): If $x$ is a lower bound for $A$, then $f(x)$ is a lower bound
for $f(A)$.  If $f(x) < a_{2}$ and $a_{2}$ is a lower bound for
$f(A)$,  then $f(x)$ is not in the closure of $f(A)$ and so by
continuity $x$ is not $inf \ A$.

For the converse, if $A$ is a closed subset of $X_2$ and for $B \subset X_1$ $x =  inf (B \cap f^{-1}(A))$ or $x = sup (B \cap f^{-1}(A))$,
then, by assumption, $f(x) =  inf f(B \cap f^{-1}(A))$ or $f(x) =  sup f(B \cap f^{-1}(A))$ and so $f(x)$ is in the closed set $A$. Hence,
$x \in f^{-1}(A)$ and so $f^{-1}(A)$ is closed by Lemma \ref{lem2.0}.

(d):  We can define $g(x)$ by
choosing any element of $f^{-1}(x)$. Such choices exactly define
the functions $g$ such that  $f \circ g = 1_{X}$.  In that case,
if $g(x_{1}) \leq g(x)$ then $ x_{1} = f(g(x_{1})) \leq
f(g(x_{2})) = x_{2} $.  Contrapositively, $x_{1} > x_{2} $ implies
$g(x_{1}) > g(x_{2})$.

(e): If $y_1 < y_2$ in $X_1$, then there exist $x_1, x_2$ in the dense subset with $y_1 < x_1 < x_2 < y_2$. Then
$f(y_1) \leq f(x_1) < f(x_2) \leq f(y_2)$.  Thus, $f$ is injective on $X_1$.

\end{proof} \vspace{.5cm}

\begin{cor}\label{cor2.1a} If $f : X_1 \to X$ is an order map with $X$ order dense and $f(X_1)$ is dense in a convex subset of $X$,
then $f$ is continuous. \end{cor}

\begin{proof} Assume $J$ is a convex subset of $X$ and $f(X_1)$ is dense in $J$.  By Lemma \ref{lem2.0a}, regarded as LOTS in their own
right, both $J$ and $f(X)$ are order dense. By Proposition \ref{prop2.1}(a)(ii) the order surjection $f : X_1 \to f(X_1)$ is continuous.  The inclusions
$f(X_1) \to J$ and $J \to X$ are continuous by (b)(iii) and (b)(i), respectively. Hence, the composition $f: X_1 \to X$ is continuous.

\end{proof} \vspace{.5cm}

If $I$ is a LOTS and $\{ X_{i} : i \in I \}$ is a family of
nonempty LOTS indexed by $I$ (a \emph{LOTS indexed family} \index{LOTS indexed family}),
then we define the \emph{order space
sum}\index{order space sum}:

\begin{equation}\label{eq2.4}
\Sigma _{i \in I} X_{i} \ \ =  \bigcup _{i \in I} \{ i \} \times
X_{i}
\end{equation}
with $(i,x) < (j,y)$ if $i < j $ or $i = j$ and $x < y$ in
$X_{i}$. If $I = \{ 0,1 \}$ then we write $X_{0} + X_{1}$ for the
sum.

If $\{ X_{i} : i \in I \}$ and $\{ Y_{j} : j \in J \}$ are LOTS
indexed families of nonempty LOTS and $f:I \rightarrow J$ and $\{
g_{i} : X_{i} \rightarrow Y_{f(i)} \}$ are order maps with $f$
injective, then the \emph{sum order map} \index{order map!sum} is the map obtained by \emph{putting together the
family} $\{ g_{i}\}$ :
\begin{equation}\label{eq2.5}
\begin{split}
g = \Sigma_f  g_i : \Sigma_{i \in I} X_{i} \rightarrow \Sigma _{j \in J} Y_j \\
g(i,x) = (f(i),g_{i}(x)).
\end{split}
\end{equation}
This is clearly an order  map which is
injective/surjective/bijective if $f$ and each $g_{i}$ satisfies
the corresponding property.

The projection map
\begin{equation}\label{eq2.6}
\begin{split}
\pi : \Sigma _{i \in I} X_{i} \rightarrow I \\
\pi (i,x) = i
\end{split}
\end{equation}
can be thought of as the special case of putting together the
surjections from $X_{i}$ to the singleton LOTS $\{i\}$. By
Proposition \ref{prop2.1}(a), $\pi$ is continuous if $I$ is order dense, or, more generally, when each $X_i$ is a closed subset of
the sum.

\begin{prop}\label{prop2.2} Let $\{ X_{i} : i \in I \}$ be a LOTS indexed family of
nonempty LOTS and let $X = \Sigma _{i \in I} X_{i}$.

\begin{enumerate}
\item[(a)] A pair $(i,x) < (j,y)$ is a gap pair in $X$ iff either $i = j$
and $x < y$ is a gap pair in $X_{i}$ or $i < j$ is a gap pair in $I$ and
$x =  max  X_{i}$ and $y =  min  X_{j}$.  So if $X$ is order dense, then each
$X_{i}$ is order dense. Conversely, assume that each
$X_{i}$ is order dense. If $I$ is also order dense or if each
$X_{i}$ is unbounded, then $X$ is order dense.

\item[(b)] If $I$ is complete and each $X_{i}$ is compact then $X$ is complete.

\end{enumerate}
\end{prop}

\begin{proof} (a): Obvious.

(b): If $A \subset X$ is bounded then $\pi (A)$ in $I$ is by
Proposition \ref{prop2.1}(c).  Let $i =  inf \ \pi (A)$.  Because $X_{i}$ is
compact, $ \pi^{-1}(i) \cap A$ is bounded in $X_{i}$. Let $x$ be
its $inf$ in $X_{i}$. If $\pi^{-1}(i) \cap A = \emptyset$ let $x =
 max \ X_{i}$. Clearly, $(i,x) =  inf \ A$.

\end{proof} \vspace{.5cm}

When it is identified with  $\pi^{-1}(i)= \{i\}\times X_{i},
X_{i}$ is a convex subset of $\Sigma _{i \in I} X_{i}$.   By Proposition \ref{prop2.1}(b) the inclusion of $X_i$ into the sum is
an embedding.

On the
other hand, let $\{ X_{i}: i \in I\}$ be a family of nonempty convex subsets
which partition a LOTS $X$, i.e. an $I$ indexed \emph{convex partition}\index{convex partition} of $X$.
Observe that if $A$ and $B$ are disjoint convex subsets of $X$, then $x_1 < y_1$ for some
pair $x_1 \in A, y_1 \in B$ implies $x < y$ for all $x \in A, y \in B$. Hence, $I$ is a LOTS with ordering uniquely
defined by $i < j$ when $x < y$ for all
 $x \in X_{i},y \in X_{j}$. A convex partition of $X$ indexed by a LOTS $I$
is equivalent to an order surjection $\pi : X \to I$.

 Clearly,we have an order
isomorphism:
\begin{equation}\label{eq2.7}
\begin{split}
\Sigma _{i \in I} X_{i} \ \cong \ X \\
(i,x) \ \mapsto \ x,
\end{split}
\end{equation}
which we can regard as an identification.

This will allow us to
put together a family of order isomorphisms between elements of
convex partitions, to obtain an order isomorphism between the partitioned
spaces. \vspace{1cm}

\subsection{Ordinal Constructions}

Of special interest are the ordinals. As usual we let $0 =
\emptyset$ and define the ordinal $\alpha$ to be the set of
ordinals smaller than $\alpha$, with the ordering by set
inclusion. The \emph{successor}\ $\alpha + 1$ of $\alpha$ is
$\alpha \cup \{ \alpha \}$.  If $A \subset \alpha$ then $inf \ A = \bigcap \ A$
is the first element of $A$ and $sup \ A =  \bigcup \ A$. Thus, any ordinal is a complete LOTS.

Any well-ordered set has the order type of an ordinal.  If
$A \subset \alpha $ then there is a unique order isomorphism of
$A$ onto an ordinal $\beta \leq \alpha$. We will usually identify
the subset $A$ with the ordinal $\b$ whose order type is that of $A$.

We let $\om$ denote the first infinite ordinal and $\Omega$ denote the first uncountable ordinal. We identify the ordinal $\om$
with the set $\N$ by letting $n$ label the $n^{th}$ (finite) ordinal.  Thus $n = \{0, 1, \dots, n-1 \}$.

For ordinal results we follow Rosenstein \cite{R} and Jech \cite{J}.
The arithmetic of ordinals is defined inductively so that $\a + \b, \a \cdot \b$ and $\a^{\b}$ are
continuous in the $\beta$ variable.\index{ordinal!sum}\index{ordinal!product}
\begin{align}\label{eq2.8}
\begin{split}
\alpha + 0 = \alpha \quad &\mbox{and} \quad  \alpha + (\beta +1) = (\alpha + \beta) + 1. \\
\alpha \cdot 0 = 0 \quad &\mbox{and} \quad \alpha \cdot (\beta + 1) = (\alpha \cdot \beta) + \alpha. \\
\alpha^{0} = 1 \quad &\mbox{and} \quad \alpha^{(\beta + 1)} =
(\alpha^{\beta}) \cdot \alpha.
\end{split}
\end{align}
In particular, if $\alpha$ and $\beta$ are countable ordinals then
the results of all of these operations are countable ordinals. We
will use the usual order type sloppiness, writing $A+B$ for
well-ordered sets $A$ and $B$ to mean the ordinal which is the sum
of the ordinals having order types $A$ and $B$.

Notice that the
ordinal sum is a special case of the two term order space sum
defined above. Also the product $\a \cdot \b$ is isomorphic to the order space sum of the $\b$ indexed family of
copies of $\a$.

Ordinal addition and multiplication are associative and by induction on $\b_2$ the following arithmetic identities hold

\begin{align}\label{eq2.8b}
\begin{split}
\a \cdot (\b_1 + \b_2) &= \a \cdot \b_1 + \a \cdot \b_2, \\
\a^{\b_1 + \b_2} &= \a^{\b_1}  \cdot \a^{\b_2}, \hspace{1cm}\\
(\a^{\b_1})^{\b_2} &= \a^{\b_1  \cdot \b_2}.\hspace{1cm}\
\end{split}\end{align}

In addition, we have
\begin{align}\label{eq2.8a}
\begin{split}
0 < \a, \b \quad \ &\Longrightarrow \ \a < \a + 1 \leq \a + \b, \hspace{1.5cm}\\
0 < \a \ \text{and} \ 1 < \b \ &\Longrightarrow \ \a < \a + \a = \a \cdot 2 \leq \a \cdot \b, \\
1< \a, \b \quad  \ &\Longrightarrow \ \a < \a \cdot 2 \leq \a \cdot \a = \a^{2} \leq \a^{\b},\\
1 < \a \ \text{and} \ \b_1 < \b_2 \ &\Longrightarrow \ \a^{\b_1} < a^{\b_1} \cdot \a = a^{\b_1 +1} \leq a^{\b_2}.
\end{split}\end{align}

If $\alpha$ is an ordinal and $\beta < \alpha$ then the
\emph{tail}\index{tail}\index{ordinal!tail}
\begin{equation}\label{eq2.9}
\begin{split}
\alpha \setminus \beta = \{ i : \beta \leq i < \alpha \} \subset
\alpha, \hspace{1cm}\\ \mbox{so that}\quad  \beta + (\alpha
\setminus \beta) = \alpha. \hspace{1cm}
\end{split}
\end{equation}
As usual, we identify the subset $\a \setminus \b$ with the ordinal having the same order type.

An ordinal $\alpha$ is called \emph{tail-like}\index{tail-like}\index{ordinal!tail-like} if it is positive and all of
the tails of $\alpha$ have order type $\alpha$, i.e. $\a \setminus \b = \a$ for all $\b < \a$. Observe that
if  $\beta < \alpha$, then $\a \setminus \b = \a$ is
equivalent to \ $\beta + \alpha = \alpha$. Thus, $\a$ is tail-like iff \  $\beta_{1} + \beta_{2} < \alpha$ for all
$\beta_{1},\beta_{2} < \alpha$.

We recall from \cite{R} Theorem 3.46 the \emph {Cantor Normal Form Theorem}\index{Cantor Normal Form}.\index{Cantor Normal Form}

\begin{prop}\label{propCNF} An ordinal $\a$ is tail-like iff
$\a = \omega^{\beta}$ for some ordinal $\b$.

Any positive ordinal
$\alpha$ can be written uniquely as the sum
\begin{equation}\label{eq2.11}
\begin{split}
\alpha = \omega^{\beta_{1}} + ... +\omega^{\beta_{N}} \\
\mbox{with} \quad  \beta_{1} \geq ...\geq \beta_{N}
\end{split}
\end{equation}
\end{prop}

\begin{proof} First observe that if $\g < \b < \a$ and $\a \setminus \b = \a$, then
\begin{equation}\label{eq2.9a}
 \a = \a \setminus \b \leq \a \setminus \gamma \leq \a.
\end{equation}

Clearly, $1 = \omega^{0}$ is
tail-like. It is the only tail-like ordinal which is not a limit ordinal.

Inductively,
we have, for $\epsilon < \beta$
and $N < \omega$:
\begin{equation}\label{eq2.10}
\begin{split}
\omega^{\epsilon} + \omega^{\beta} = \omega^{\epsilon}\cdot
(1 + \omega^{\beta \setminus \epsilon}) =  \omega^{\epsilon}\cdot
\omega^{\beta \setminus \epsilon} = \omega^{\beta}, \\
\omega^{\beta} \cdot N + \omega^{\beta+1}=\omega^{\beta}\cdot (N +
\omega) = \omega^{\beta}\cdot \omega  = \omega^{\beta + 1}
\end{split}
\end{equation}
Thus, $\om^{\b} \setminus \om^{\ep} = \om^{\b}$ and $\om^{\b+1} \setminus \omega^{\beta} \cdot N = \om^{\b+1} $.
It then follows from (\ref{eq2.9a}) that
 $\omega^{\beta}$ \  is tail-like for any ordinal
$\beta$.

Continuity in $\b$ implies that the set $\{ \b : \om^{\b} \leq \a \}$ is closed and so we can choose
$\om^{\b_1}$ to be the largest ordinal of this form less than or equal to $\a$.  We see that $\b_1$ is the unique
ordinal such that $\om^{\b_1} \leq \a < \om^{\b_1 + 1}$.

We proceed by induction on $\b_1$.
Observe that if $\a$ satisfies (\ref{eq2.11}), then $\om^{\b_1} \leq \a \le  \om^{\b_1} \cdot N   < \om^{\b_1 + 1}$.
Thus, $\b_1$ is uniquely determined by this inequality.

If $\a = \om^{\b_1}$, then we have Cantor Normal Form Theorem for $\a$ with $N = 1$.

If $\a > \om^{\b_1}$, then since $\a < \om^{\b_1+1}$, there exists $k \in \om$ such that $\a < \om^{\b_1}\cdot (k+2)$. The minimum
such value
$k = k(\a)$ is uniquely determined by $\a$. We have $\a \geq \om^{\b_1}\cdot (k+1)$ and
$\gamma = \a \setminus \om^{\b_1}\cdot (k+1) < \om^{\b_1}$.
Let $\b_i = \b_1$ for $1 \leq i \leq k+1$.

Applying the induction hypothesis, let $\gamma =  \omega^{\beta_{k+2}} + ... +\omega^{\beta_{N}}$ be the unique Cantor Normal Form for $\gamma$.
We have $\omega^{\beta_{k+2}} \leq \gamma < \om^{\b_1} = \om^{\b_{k+1}}$ and so by (\ref{eq2.8a}) $\b_{k+2} < \b_{k+1}$. Summing the two
decompositions we obtain the unique normal form for $\a$.

Finally, observe that if $N > 1$, then by (\ref{eq2.8a}) $\om^{\b_1} < \a$. Clearly $\a \setminus \om^{\b_1}$ is equal to
$\omega^{\beta_{2}} + ... +\omega^{\beta_{N}}$ and by uniqueness of the Cantor Normal Form this does not equal $\a$. Hence, $\a$ is
not tail-like.

It follows that the only tail-like ordinals are of the form $\om^{\gamma}$.

\end{proof} \vspace{.5cm}

\begin{cor}\label{corCNF}An ordinal $\a$ is a limit ordinal iff
$\a = \omega \cdot \beta$ for some positive ordinal $\b$.

Any infinite ordinal $\a$ can be written uniquely as $\a = \b + k$ with
$\b$ a limit ordinal and $k < \om$. \end{cor}

\begin{proof} If $\omega^{\beta_{1}} + ... +\omega^{\beta_{N}}$ is Cantor Normal Form for $\a$, then
$\omega^{\beta_{1}} + ... +\omega^{\beta_{N}} + \om^0$ is Cantor Normal Form for $\a + 1$. So we can uniquely
write $\a$ as
\begin{equation}\label{eq2.10a}
\om \cdot (\omega^{\beta_{1}- 1} + ... +\omega^{\beta_{N-k} - 1}) + k
\end{equation}
with $\beta_{N-k} > 0$ and $\b_i = 0$ for $N-k < i \le N$. For $\a$ an infinite ordinal, e.g. a limit ordinal, $\b_1 > 0$.
It is a limit ordinal iff $k = 0$.

\end{proof}\vspace{.5cm}

A \emph{cardinal}\index{cardinal} $\aleph$ is the ordinal which is minimum among the ordinals of that cardinality. That is,
if $\a < \aleph$, then the cardinality of $\a$ is strictly less than that of $\aleph$. By mapping $\b + k$ to $\b + 2k$ or to
$\b + 2k + 1$ for $\b$ any limit ordinal less than $\aleph$, we see that an infinite cardinal can be written as the disjoint union of
two sets of the same cardinality.   It follows that any infinite cardinal $\aleph$ is
tail-like.

The elements of $\om$ are the finite cardinals, and $\om, \Omega$ are the infinite cardinals $\aleph_0$ and $\aleph_1$, respectively.\vspace{.5cm}

If $\alpha$ is a positive ordinal and $\{X_{i}: i \in \alpha\}$ is
an $\alpha$ indexed family of nonempty LOTS then we define the
\emph{order space product} to be the set $\Pi_{i \in \alpha} X_{i}$
with the lexicographic ordering \index{lexicographic ordering}.  That is, for $x \not= y$ in the
product
\begin{equation}\label{eq2.12}
x < y \quad \Longleftrightarrow \quad x_{\beta} < y_{\beta}\ \
\mbox{with} \ \beta = min \{ j : x_{j} \not= y_{j} \}.
\end{equation}
If $\alpha = 2$, we write $X_{0} \times X_{1}$ for the product.

When $X_{i} = X$ for all $i$ then we obtain $X^{\alpha}$, the
space of functions from $\alpha$ to $X$, as a LOTS.  Observe that
the LOTS topology is \emph{not} the product topology.

If $X_{i} = X$ for all $i \in I$ in an $I$ indexed family $\{
X_{i} : i \in I \}$ then the order sum $\Sigma_{i \in I} X_{i}$ is
the product $I \times X$.

If $Y$ is a LOTS then $(\bigcup _{i \in I} \{ i \} \times X_{i}) \times Y = \bigcup _{i \in I} (\{ i \} \times X_{i}) \times Y)$ implies
\begin{equation}\label{eq2.12a}
(\Sigma _{i \in I} X_{i}) \times Y \ \ \cong \ \ \Sigma _{i \in I} (X_{i} \times Y).
\end{equation}
Furthermore, if $\{ X_{i} : i \in I \}$ is an arbitrary family of subsets of a LOTS $X$ and $Y$ is a LOTS, then
$\{ X_{i} \times Y : i \in I \}$ is a family of subsets of $X \times Y$ and regarded as LOTS in their own right, we have
\begin{equation}\label{eq2.12b}
(\bigcup _{i \in I} X_{i}) \times Y \ \ \cong \ \ \bigcup _{i \in I} (X_{i} \times Y).
\end{equation}

In contrast with the sum, the ordinal conventions for products and
powers given by (\ref{eq2.8}) disagree with these new definitions.

It follows by induction on $\b$ using
(\ref{eq2.12a}) and (\ref{eq2.12b}) that the order space product
$\b \times \a$ is the
ordinal product $\a \cdot \b$.

The order space power $2^{\omega}$ is the
uncountable space of all zero/one valued sequences.  The ordinal
$2^{\omega}$ is the limit of the finite ordinals $2^{N}$ and so is
just $\omega$.

If $\{X_{i}: i \in \alpha\}$ is an ordinal indexed family of
nonempty LOTS and $0 < \beta \leq \alpha$ then we can write
$\Pi_{\beta}$ for the subproduct $\Pi_{i \in \beta}X_{i}$ and if
$0 < \epsilon \leq \beta \leq \alpha$ then we denote by
\begin{equation}\label{eq2.13}
\pi_{\ep}^{\b} :\Pi_{\beta}  \longrightarrow
\Pi_{\epsilon}
\end{equation}
the projection map obtained by forgetting the coordinates in
$\beta \setminus \epsilon$. Since $\epsilon$ is an initial segment
of $\beta$ it is clear that $ \pi_{\ep}^{\b}$ is an order
surjection.  However, it need not be continuous.  The first
coordinate projection $\omega \times \omega \rightarrow \omega$ \
does not have closed point inverses.

\begin{prop}\label{prop2.3} Let $X = \Pi_{i \in \alpha}X_{i}$ be an order space product of LOTS.
\begin{enumerate}
\item[(a)]If each $X_{i}$ is order dense, then $X$ is order dense.  In that
case, each projection $ \pi_{\ep}^{\b}$ for $0 < \epsilon \leq \beta \leq \alpha$ is continuous.
\item[(b)] If each $X_{i}$ for $i > 0$ is bounded, then each projection
            $ \pi_{\ep}^{\b}$ is continuous.  If, in addition, each $X_{i}$ is complete, then $X$ is complete.
\end{enumerate}
\end{prop}

\begin{proof} (a) If $x < y$ in $X$, then with
$\beta = min \{ j: x_{j} \not= y_{j} \}$ \ we can choose
$z_{\beta} \in X_{\beta}$ such that $x_{\beta} < z_{\beta} <
y_{\beta}$.  Define $z_{j} = x_{j}$ for all $j \not= \beta$.  Then
$x < z < y$ in $X$. The projections are continuous by Proposition
\ref{prop2.1}(a)(ii).

(b)  Given $a < b$ in $\Pi_{\epsilon}$ with $\epsilon > 0$ define
$a+$ and $b-$ in $\Pi_{\beta}$ by:
\begin{equation}\label{eq2.14}
 (a+)_{i} = \begin{cases} a_{i} & i \in \epsilon \\ max \ X_{i} & i \in \beta \setminus \epsilon \end{cases}
\end{equation}
\begin{displaymath}
 (b-)_{i} = \begin{cases} b_{i} & i \in \epsilon \\ min \ X_{i} & i \in \beta \setminus \epsilon \end{cases}
\end{displaymath}
Note that $i \in \beta \setminus \epsilon$ is positive and so
$X_{i}$ is bounded. Clearly,
\begin{equation}\label{eq2.15}
\begin{split}
 (\pi_{\ep}^{\b})^{-1}((a,b)) = (a+,b-),\\
  (\pi_{\ep}^{\b})^{-1}([a,b]) = [a-,b+].
\end{split}
\end{equation}
Hence, $\pi_{\ep}^{\b}$ is continuous.

If $A \subset X$ is contained in $[a,b]$ then by replacing $X_{0}$
by $[a_{0},b_{0}]$ we can assume that every $X_{i}$ is
bounded. We prove by induction on $\beta$ that if every
$X_{i}$ is compact, then $A \subset \Pi_{\beta}$ has an $inf$.
With similar results for $sup$, completeness follows.

If $\beta = 1$ then $\Pi_{\beta} = X_{0}$ which is compact.

Now assume that the result holds for all $\epsilon < \beta$ and
let $x^{\epsilon} = inf \ \pi_{\ep}^{\b}(A)$ in
$\Pi_{\epsilon}$.  By Proposition \ref{prop2.1}(c) it is clear that $\delta
\leq \epsilon$ implies
\begin{equation}\label{eq2.16}
\pi_{\d}^{\ep}(x^{\epsilon}) = x^{\delta}.
\end{equation}

\textbf{Case 1:} If $\beta$ is a limit ordinal then define $x^{\beta}$ so
that
\begin{equation}\label{eq2.17}
x^{\beta}_{i} = x^{\epsilon}_{i} \qquad \mbox{for} \  i < \epsilon
< \beta,
\end{equation}
which is well-defined by (\ref{eq2.16}).  Clearly, $x^{\beta}= inf \ A$ in
this case. \vspace{.25cm}

\textbf{Case 2:} If $\beta = \epsilon + 1$ then $\Pi_{\beta}$ is the order
space sum of copies of the compact LOTS $X_{\epsilon}$ indexed by
the points of $\Pi_{\epsilon}$.  By induction hypothesis and
Proposition \ref{prop2.2}(b)  $\Pi_{\beta}$ is complete.  Hence, $inf \ A$
exists in this case as well.

\end{proof} \vspace{.5cm}

We will also use the inverse limit construction under special
assumptions which we will refer to as a \emph{special inverse
system}\index{special inverse system}\index{inverse system!special}\index{inverse system}.  We begin with  an ordinal indexed family
$\{ X_{i} :   i \in \alpha \}$ of connected LOTS.   For $i
< j<\alpha$ we are given order surjections $p^{j}_{i} : X_{j}
\rightarrow X_{i}$. such that
\begin{equation}\label{eq2.18}
 p^{i}_{k} \circ p^{j}_{i} = p^{j}_{k} \qquad \mbox{for} \  k < i < j < \alpha
\end{equation}
We assume that each $p^{j}_{i}$ has compact point inverses and so each $p^{j}_{i}$ is continuous and topologically proper by
Proposition \ref{prop2.1}(a). If each $X_i$ is unbounded, i.e. has no $max$ or $min$, then we call the system \emph{unbounded}
\index{inverse system!unbounded} and Proposition \ref{prop2.1}(a)
implies directly that each $p^{j}_{i}$ is continuous and topologically proper, i.e. the compact point inverses condition is automatically satisfied.

We then define the \emph{inverse limit}\index{inverse limit}
\begin{equation}\label{eq2.19}
\overleftarrow{Lim} \{X_{i}\} \ = \  \{ x \in \Pi_{i
\in \alpha} X_{i} : x_{i} = p^{j}_{i}( x_{j}) \  \mbox{for all} \
i < j < \alpha \}.
\end{equation}
We denote by $p_{j} :\overleftarrow{Lim} \{X_{i}\} \rightarrow
X_{j}$ the projection to the $j$ coordinate.

If $x<y$ in $\overleftarrow{Lim}$ then with $\ep =  min \{ j :
x_{j} \not= y_{j} \}, \ x_{\ep} < y_{\ep}$ and $x_{i} = y_{i}$
for all $i < \ep$. On the other hand, if $j > \ep$ then $x_{j}
< y_{j}$ because $p_{\ep}^{j}$ is an order map. It follows that
each $p_{j}$ is an order map and we have
\begin{equation}\label{eq2.19x}
x < y \ \Longleftrightarrow \ \exists \ep < \a \ \text{such that} \begin{cases} x_i = y_i \ \ \text{for all} \ \ i < \ep, \\
x_i < y_i \ \ \text{for all} \ \ i \ge \ep. \end{cases}
\end{equation}
\vspace{.5cm}

\begin{prop}\label{prop2.4} If $(\{ X_{i} : i \in \alpha \},\{ p^{j}_{i} : i < j < \alpha \})$
is a special inverse system, then the inverse limit $\overleftarrow{Lim} \{X_{i}\}$
is a connected LOTS and each projection $p_{j}$ is a
continuous, topologically proper order surjection. If the system is unbounded, then the inverse
limit is unbounded. If each $X_i$ is compact, then the inverse limit is compact.

\end{prop}

\begin{proof} Given $a \in X_{j}$ we construct
by the usual compactness argument $x \in \overleftarrow{Lim}
\{X_{i}\}$ such that $p_{j}(x) = a$. We use the topological product topology on $Z = \prod_i X_i$.
For $j \leq i < \a$ let
\begin{equation}\label{eq2.19a}
\begin{split}
Q_i = \{ z \in Z: z_k \in (p_i^k)^{-1}(z_i) \ \text{for} \ k > i, \\
z_k = p_k^i(x_i) \ \text{for} \ k \leq i,  \ \text{ and with} \ z_j = a \}
\end{split}
\end{equation}
Because each
$p^{i}_{k}$ is a continuous, topologically proper map $\{ Q_i \}$ is a filterbase of nonempty compact sets in the
topological product. Hence, the intersection, which is $(p_j)^{-1}(a)$ is nonempty.

By Proposition \ref{prop2.1}(a) each $p_{j}$ is continuous and
$\overleftarrow{Lim}$ is unbounded if some $X_j$ is unbounded.

If $x<y$ in $\overleftarrow{Lim}$, then let $\ep$ satisfy (\ref{eq2.19x}).  Because $X_{\ep} $ is order
dense, we can choose $z_{\ep} \in X_{\ep}$ such that
$x_{\ep} < z_{\ep} < y_{\ep}$.  Choose $z \in
\overleftarrow{Lim}$ such that $p_{\ep}(z) = z_{\ep}$. Because
$p_{i}^{\ep}$ is order preserving, $x_{i} = z_{i} = y_{i}$ for
all $i < \ep$.  Hence, $x < z < y$ in  $\overleftarrow{Lim}$ and
so the inverse limit is order dense.

If $A \subset \overleftarrow{Lim}$, let $A_{j} = p_{j}(A)$.
Clearly, $i < j$ implies $A_{i} = p^{j}_{i}(A_{j})$.  If $A$ is
bounded in  $\overleftarrow{Lim}$ then by Proposition \ref{prop2.1}(c) each
$A_{j}$ is bounded in $X_{j}$. If each $A_j$ is bounded and if $x_{j} =  inf  A_{j}$ then
$x_{i} = p^{j}_{i}(x_{j})$.  This defines a point $x$ of
$\overleftarrow{Lim}$ which is $ inf  A$. With a similar argument
for the $sup$ we see that $\overleftarrow{Lim}$ is complete and so
is connected.

If $C \subset A_j$ is bounded, then let $A = (p_j)^{-1}(C)$ and $A_i = p_i(A)$ for all $i$.
In particular, since $p_j$ is surjective, $A_j = C$, $A_i = p^j_i(C)$ for $i < j$ and $A_i = (p^i_j)^{-1}(C)$ for $i > j$.
Thus, each $A_j$ is bounded because every $p^i_j$ is a topologically proper, continuous map. It follows that
Hence, $A = (p_j)^{-1}(C)$ is bounded as well as closed and so is compact by completeness. Thus,
each $p_{j}$ is topologically proper.

It follows that if $X_j$ is compact, then the inverse limit ($= (p_j)^{-1}(X_j)$) is compact.

\end{proof} \vspace{.5cm}

\noindent{\bfseries Remarks.}  While it is easy to show that
$\overleftarrow{Lim}\{X_{i}\}$ is a closed subset of the order
space product $\Pi_{i \in \alpha} X_{i}$, the topology on the
inverse limit is not the relative topology induced from that product.  For example, if $X_{0} = X_{1} =
\R$ and $p^{1}_{0}$ is the identity map then the inverse
limit is the diagonal $\{ (t,t) : t \in \R \}$ which is a
discrete subset of the order space product $\R\times
\R$. The topology on the inverse limit is instead the relative topology induced from the topological product space.

Absent the assumption of connectedness, the inverse limit projections need not be
continuous. Define the $\om$ indexed inverse limit system by
\begin{equation}\label{eq2.19b}
X_n = \{1, \dots, n \} \cup \{\om \} \ \text{ and} \ p_{n}^{n+1}(k) =
\begin{cases} k \ \text{for} \ k \leq n, k = \om, \\ n \ \text{for} \ k = n+1.\end{cases}
\end{equation}
The inverse limit is isomorphic to $\om + 1$ with $(p_1)^{-1}(1) = \om$.
 \vspace{.5cm}

Any LOTS $X$ can be regarded as a subset of a smallest complete
LOTS $\hat{X}$ called its \emph{completion}\index{completion}\index{$\hat X$}, see \cite{R}
Theorem 2.32.  We will only need $\hat{X}$ in the case when $X$ is
order dense.

First, assume that $X$ is unbounded.  In that case, define
$\hat{X}$ to be the set of open subsets $A$ of $X$ such that $(A,X
\setminus A)$ is a Dedekind cut in $X$, i.e. $A$ is a proper open
subset of $X$ and $x \in A$ implies $(-\infty,x) \subset A$.
Order $\hat{X}$ by inclusion and identify $x \in X$ with
$(-\infty,x) \in \hat{X}$.  Since $X$ is unbounded and $A, X
\setminus A$ are nonempty, $\hat{X}$ is unbounded.  The $sup$ of a
bounded subset of $\hat{X}$ is its union and the $inf$ is the
interior of its intersection.  If $A_{2} < A_{1}$ in $\hat{X}$,
then because $A_1$ is open, there exist  $x < y $ in  $ A_{1} \setminus A_{2}$.
So $A_{2} < (-\infty,y) < A_{1}$ because $x \in (-\infty,y) \setminus A_2$ and $y \in A_1 \setminus (-\infty,y)$.
 It follows that $\hat{X}$ is order dense, and so is connected, and
that $X$ is dense in $\hat{X}$.

If $X$ has a $max$ or $min$ then we obtain $\hat{X}$ by removing
such endpoints, completing, and then reattaching them.  If $M$ is
$ max  X$ then we regard $M$ as $ max  \hat{X}$ and if $m = min X$ then we regard $m $ as $ min \hat{X}$.

Let $Y$ be a complete LOTS and $f : X \rightarrow Y$ be an order
map with $X$ order dense.  Define $\hat{f} : \hat{X} \rightarrow Y$ so that for $\hat x \in \hat X$
\begin{equation}\label{eq2.20}
\hat{f}(\hat x) \quad = \quad  sup \{ f(x) : x \in X \cap (-\infty, \hat x]  \}.
\end{equation}
Clearly, $\hat f$ is an order map which extends $f$.

\begin{prop}\label{prop2.4a} Assume  $Y$ is a complete LOTS and $f : X \rightarrow Y$ is an order
map with $X$ order dense.
 \begin{itemize}
 \item[(a)] If $f$ is injective, then $\hat{f}$ is injective.

 \item[(b)] If $f(X)$ is dense in a connected subset of $Y$, then the map $f$ and its extension $\hat{f}$ are continuous order maps
  with $\hat{f}(\hat{X})$ connected.

  \item[(c)] Assume $Y$ is connected and $f(X)$ is dense in $Y$. If $Y$ is unbounded, then $\hat{f}$ is surjective and $X$ is unbounded.
 More generally, if
 \begin{equation}\label{eq2.20a}
 \begin{split}
 max \ Y \  \text{exists} \ \Longrightarrow \  max \ X \ \text{exists and} \\
 min \ Y \  \text{exists} \ \Longrightarrow \  min \ X \ \text{exists},
\end{split}
\end{equation}
then $\hat{f}$ is surjective.
\end{itemize}
In particular, if $X$ is a
dense subset of a connected LOTS $Y$, then $X$ is order dense.
If $Y$ is unbounded or (\ref{eq2.20a}) holds, then $\hat{X} \cong Y$.
\end{prop}

\begin{proof} (a): This follows from Proposition \ref{prop2.1}(e).

(b): Let $Y_0$ be a  connected subset of $Y$ in which $f(X)$ is dense. Since the closure of a connected set is connected, we may assume
that $Y_0$ is closed and so $Y_0 = \overline{f(X)}$. Clearly, $f(X) \subset \hat f(\hat X) \subset Y_0$. It then follows from Corollary \ref{cor2.1a}
that $f : X \to Y_0$ and $\hat f : \hat X \to Y_0$ are continuous maps. By Proposition \ref{prop2.1}(b)(ii)
the inclusion of $Y_0$ into $Y$ is continuous.

Composing we see that $f$ and $\hat f $ are continuous. Since $\hat X$ is connected, its image is connected.

(c): By (b) $\hat f$ is continuous and the subset $\hat{f}(\hat{X})$ is  connected and dense in $Y$. If $Y$ is
unbounded, $\hat{f}(\hat{X})= Y$ and so $\hat{f}$ is surjective. Also $X$ is unbounded by Proposition \ref{prop2.1}(c).
If $ max  X$ exists, then by Proposition \ref{prop2.1}(c) again $\hat{f}( max X) =  max Y$ and
similarly for the minimum. Thus, (\ref{eq2.20a}) implies that $\hat{f}(\hat{X}) \subset Y$ is  dense and connected
and contains $ max Y$ or $ min Y$ when they exist. Hence, $\hat{f}(\hat{X}) = Y$.

In particular, if $f$ is injective and (\ref{eq2.20a}) holds, then from (a) it follows that
$\hat{f}$ is bijective and so is an order isomorphism.

\end{proof}\vspace{.5cm}

A bounded, convex subset $J$ in a complete space is an interval
with endpoints $ inf  J$ and $ sup  J$.  With $X \subset \hat{X}$
as above, a bounded, convex subset $J \subset X$ is equal to
$\hat{J} \cap X$ where $z \in \hat{J}$ iff there exist $x_{1},
x_{2} \in J$ such that $x_{1} \leq z \leq x_{2}$.  It follows that
a convex set in $X$ is the intersection of an interval of
$\hat{X}$ with $X$. \vspace{1cm}

\subsection{Countability Conditions}

A topological space $X$ is \emph{separable} if it has a countable
dense subset.  $X$ satisfies the \emph{countable chain condition}\index{countable chain condition (c.c.c)},
hereafter denoted \emph{c.c.c}, if any collection of pairwise
disjoint, nonempty open subsets is countable. A LOTS $X$ satisfies
c.c.c. if any collection of pairwise disjoint, nonempty open
subintervals is countable because the subintervals form a base for
the topology.

A subset $A$ of a LOTS $X$ is \emph{cofinal}\index{cofinal}\index{subset!cofinal} if for any $x \in X$
there exists $a \in A$ such that $ x \leq a$.  If $M =  max  X$
exists then $A$ is cofinal iff $M \in A$. If $X$ has no $max$ then
any dense subset of $X$ is cofinal.
$A$ is  \emph{coinitial} \index{coinitial}\index{subset!coinitial}if it is cofinal for the reverse $X^{*}$.
  $A$ is  \emph{$\pm$cofinal} \index{$\pm$cofinal}\index{subset!$\pm$ cofinal}if it
is both cofinal and coinitial.  $X$ is called \emph{$\sigma$-bounded} \index{$\sigma$-bounded}\index{subset!$\sigma$-bounded}
if it admits  countable $\pm$cofinal
subsets.  If $X$ is complete then it is $\sigma$-bounded iff it is
$\sigma$-compact.

A point $x$ in a LOTS $X$ has a countable neighborhood base iff
the interval $(-\infty, x)$ has a countable cofinal subset and the
interval $(x,+\infty)$ has a countable coinitial subset.  Thus,
$X$ is first countable iff for every point $x \in X, x$ is either
a right endpoint or the limit of an increasing sequence, and $x$
is either a left endpoint or the limit of a decreasing sequence.
Equivalently,  $X$ is first countable iff every bounded interval
in $X$ is a $\sigma$-bounded LOTS in its own right.

We say that a LOTS $X$ has \emph{countable type}\index{countable type}\index{LOTS!countable type} if every convex
subset of $X$ is a $\sigma$-bounded LOTS.  It clearly suffices
that every open, convex subset of $X$ be $\sigma$-bounded.

\begin{prop}\label{prop2.5} Let $X$ be a LOTS.
\begin{enumerate}
\item[(a)] $X$ satisfies c.c.c. iff there does not exist an order injection
into $X$ of a LOTS with uncountably many isolated points.
\item[(b)] $X$ is of countable type iff there does not exist an injective
order map or order* map of the first uncountable ordinal $\Omega$ into $X$.
\item[(c)] If $X$ satisfies c.c.c., then it has only countably many isolated points.
\item[(d)] If $X$ is of countable type, then it is first countable and
$\sigma$-bounded. If $X$ is complete, first countable and $\sigma$-bounded then it is of countable type.
\item[(e)] If $X$ is separable, then it satisfies c.c.c.  If $X$ satisfies c.c.c., then it is of countable type.
\item[(f)] Let $f: X_{1} \rightarrow X$ be an order injection.  If $X$ is
separable, satisfies c.c.c. or is of countable type, then $X_{1}$ satisfies the corresponding property.
\item[(g)] Let $f: X \rightarrow X_{1}$ be an order surjection.  If $X$ is
separable, satisfies c.c.c. or is of countable type, then $X_{1}$ satisfies the corresponding property.
\item[(h)] Let $f: X \rightarrow X_{1}$ be an order map with image $f(X)$
dense in $X_{1}$.  If $X$ is of countable type, then $X_{1}$ is. If $X_{1}$
has only countably many isolated points and $X$ is separable or satisfies
c.c.c., then $X_{1}$ satisfies the corresponding property.  If $f$ is continuous
and $X$ is separable or satisfies c.c.c., then $X_{1}$ satisfies the corresponding property.
\end{enumerate}
\end{prop}

\begin{proof}(a): Suppose $f:A \rightarrow X$ is
an order injection and with $I \subset A$ defined to be the set of
isolated points in $A$, excluding the $max$ and $min$ if any,
suppose that $I$ is uncountable.  For each $x \in I$ there are
unique $x-,x+ \in A$ such that $x- < x < x+$ and $x$ is the only
point in the interval $(x-,x+)$.  By Zorn's Lemma we can choose
$\tilde{I} \subset I$ maximal with respect to the property that $x
\in \tilde{I} \Rightarrow  x-,x+ \not\in \tilde{I}$.  By
maximality $x \in I$ implies that either $x, x-$ or $x+ \in
\tilde{I}$.  Hence, $\tilde{I}$ is uncountable. If $x_{1} < x_{2}$
in $\tilde{I}$, then $x_{1}+ \leq x_{2}$ and the inequality is
strict because $x_{1}+ \not\in \tilde{I}$.  Hence, $x_{1}+ \leq
x_{2}-$.  It follows that $\{(f(x-),f(x+)) : i \in \tilde{I} \}$
is an uncountable family of pairwise disjoint, nonempty open
intervals in $X$.  Thus, $X$ does not satisfy c.c.c.

Now assume that $X$ does not satisfy c.c.c.  If $X$ has
uncountably many isolated points, then the identity on $X$ is the
required order injection.  So we can assume that $X$ itself has
only countably many isolated points and that $\{ J_{i}\}$  is an
uncountable family of pairwise disjoint open subintervals in $X$.
If $J_{i}$ is finite, then it consists of isolated points and there
are only countably many such.  Thus, we can assume that each
$J_{i}$ is infinite.  For each $i$ choose points $x_{i}- < x_{i} <
x_{i}+ $ in $J_{i}$.  The collection of all these points is a
subset $A$ of $X$ and the inclusion $f : A \rightarrow X$  is a,
not necessarily continuous, order injection of LOTS.  In $A$ each
$x_{i}$ is an isolated point and so $A$ has uncountably many
isolated points.

(b):  If $f : \Omega \rightarrow X$ is an order injection, then $J =
\{ x : x < f(i)$ for some $i \in \Omega \}$ is a nonempty, open
convex subset of $X$ and if $A$ is a countable subset of $J$, then
we can choose for each $a \in A, i(a) \in \Omega $ such that $a <
f(i(a))$. Let $j \in \Omega$ with $j > sup \ \{i(a) : a \in A \}$.
Then $f(j) > a$ for all $a \in A$ and $f(j) < f(j+1)$ so that $
f(j) \in J$.  Thus, $A$ is not cofinal in $J$.

Conversely, if $J$ is a nonempty, open convex subset of $X$ with
no countable cofinal subset, then we can construct an order
injection $f : \Omega \rightarrow J$ inductively by choosing $f(j)
\in J$ larger than all the $f(i)$ previously chosen for $i < j$.

Similarly, every nonempty, open convex subset of $X$ has countable
coinitial subsets iff there does not exist an order injection of
the reverse $\Omega^{*}$ into $X$.

(c): This follows from (a) or directly because the isolated points
form a pairwise disjoint collection of open singletons.

(d): If $X$ is of countable type, then $X$ is a convex subset and so
is $\sigma$-bounded.  Every interval is $\sigma$-bounded and so
$X$ is first countable.  Conversely, if $X$ is first countable and
complete then every bounded convex set is an interval and so is
$\sigma$-bounded by first countability.  If, in addition, $X$ is
$\sigma$-bounded, then every interval is $\sigma$-bounded and so
$X$ is of countable type.

(e): Disjoint nonempty open sets contain distinct elements of any
dense subset.  Hence separability implies c.c.c.  $\Omega$ and
$\Omega^{*}$ contain uncountably many isolated points and so
c.c.c. implies countable type by (a) and (b).

(f):  If $g :A \rightarrow X_{1}$ is an order injection, then $f
\circ g : A \rightarrow X$ is an order injection.  By (a) and (b)
$X$ is not c.c.c. or of countable type if $X_{1}$ is not.

Now assume that $D$ is a countable dense subset of $X$. For every
pair $d_{1} < d_{2}$ in $D$ with $f(X_{1}) \cap (d_{1},d_{2})
\not= \emptyset$ choose a point $x \in X$ such that $d_{1} < f(x) <
d_{2}$.  The collection of such points is a countable subset
$D_{1}$ of $X_1$.  Since $X$ is separable it satisfies c.c.c. and
so $X_{1}$ satisfies c.c.c. by what we have already shown.  Thus,
the set $D_{2}$ of isolated points of $X_{1}$ is countable by (c).
We show that the countable set $D_{1} \cup D_{2}$ is dense in
$X_{1}$.  Let $a < b$ in $X_{1}$.  If the interval $(a,b) \subset X_1$ is
finite but nonempty, then it consists of isolated points and so
meets $D_{2}$.  If it is infinite then we can choose $c_{1} <
c_{2} < c_{3}$ in $(a,b)$.  The open subintervals
$(f(a),f(c_{2}))$ and $(f(c_{2}),f(b))$ contain $f(c_{1})$ and
$f(c_{3})$ respectively and so they meet $D$.  That is, there exist
$d_{1},d_{2} \in D$ such that $f(a) < d_{1} < f(c_{2}) < d_{2} <
f(b). $  By definition there exists $x \in D_{1}$ such that $f(x)
\in (d_{1},d_{2})$.  Hence, $a < x < b$ and so $(a,b)$ meets
$D_{1}$.

(g): Choose for each $x \in X_{1}$ a point $g(x) \in f^{-1}(x)
\subset X$, which is nonempty since $f$ is surjective.  We obtain
an order injection $g: X_{1} \rightarrow X$ and so the results
follow from (f) applied to $g$.

(h):  Suppose that $\tilde{g} : \Omega \rightarrow X$ is an order
injection.  We will define, by induction, an order injection $g :
\Omega \rightarrow X_{1}$ such that
\begin{equation}\label{eq2.21}
\tilde{g}(i) < f \circ g(i) < \tilde{g}(i+2) \quad \mbox{for all}\
i \in \Omega.
\end{equation}
For $i= 0$ or $i$ a limit ordinal, $\tilde{g}(i+1)$ in the interval
$(\tilde{g}(i),\tilde{g}(i+2))$ implies that the dense set
$f(X_{1})$ meets the interval and so we can choose $g(i)$ so that
$f(g(i))$ is in the interval.  If $i = \beta +1$ and $g$ has been
defined through $\beta$, then
\begin{equation}\label{eq2.22}
f \circ g(\beta ) < \tilde{g}(\beta +2) = \tilde{g}(i+1)
\end{equation}
implies that $\tilde{g}(i+1)$ lies in the interval between $max \
(\tilde{g}((i),f \circ g(\beta ))$ and $\tilde{g}(i+2)$.  Choose
$g(i)$ so that $f(g(i))$ lies in this interval.  With a similar
argument for $\Omega^{*}$ we use (b) to see that if $X$ is not of
countable type then $X_{1}$ is not of countable type.

Now assume that $X_{1}$ has only countably many isolated points
and let $D_{2}$ denote the set of isolated points of $X_{1}$
together with its $max$ and $min$ if any.  So $D_{2}$ is countable
by assumption.

If $X_{1}$ does not satisfy c.c.c., then there exists an
uncountable family $\{ J_{i} \}$ of pairwise disjoint, nonempty
open intervals in $X_{1}$.  Since $D_{2}$ is countable we can
assume that $J_{i} \cap D_{2} = \emptyset$ for all $i$ and so each
$J_{i}$ is infinite.  In $J_{i}$ choose $y_{1} < ...< y_{4}$.
Because $f(X)$ is dense there exist $x_{1},x_{2},x_{3}$ in $X$
such that
\begin{equation}\label{eq2.23}
y_{1} < f(x_{1}) < y_{2} < f(x_{2}) < y_{3}< f(x_{3}) <  y_{4}.
\end{equation}
Hence, $(x_{1},x_{3})$ is a nonempty open interval contained in
$f^{-1}(J_{i})$. Thus, we have constructed a pairwise
disjoint,uncountable family of nonempty open intervals in $X$.
Thus, $X$ does not satisfy c.c.c.

Now assume that $X$ contains a countable dense subset $D$.  I
claim that $D_{2} \cup f(D)$ is dense in $X_{1}$, which will prove
that $X_{1}$ is separable.

If $J$ is a nonempty open subinterval of $X_{1}$ and $J$ is finite,
then $J$ meets $D_{2}$.  Choose $y_{1} < ...< y_{4}$ in $J$ as
before and get $x_{1},x_{2},x_{3}$ in $X$ which satisfies (\ref{eq2.23}).
The interval $(x_{1},x_{3})$ in $X$ is nonempty and so meets $D$.
Thus, $J$ meets $f(D)$.

If $f$ is continuous and $D$ is dense in $X$, then $f(D)$ is dense
in $f(X)$ and hence in $X_{1}$.  If $\{ O_{1} \}$ is a pairwise
disjoint family of nonempty open subsets of $ X_{1}$ then $
\{f^{-1}(O_{i}\}$ is such a family in $X$.  So if $f$ is
continuous with a dense image then $X_{1}$ is separable or c.c.c.
if $X$ is. This does not even require that $f$ be order
preserving.

\end{proof}\vspace{.5cm}

\noindent{\bfseries Remark.} Here is an example which illustrates
why the extra condition is required in part (h).  The order space
product $ X = \R \times \{ -1,+1\} $ is separable but the
points $(t,0)$ in $X_{1} = \R \times \{-1,0,+1\}$ are all
isolated.  Mapping $(t,-1)$ to $(t,0)$ and $(t,+1)$ to $(t,+1)$
for all $t \in \R$ we obtain a discontinuous order
injection $f : X \rightarrow X_{1}$ with a dense image.
\vspace{.5cm}

\begin{cor}\label{lem2.6} If $f :\Omega \rightarrow X$ is an order map or an order* map and
$X$ is of countable type, then $f$ is eventually constant.  That is, there
exists $\alpha \in \Omega$ such that $f(i) = f(\alpha)$ for all $i \in \Omega$ with $i \geq \alpha$.
\end{cor}

\begin{proof} If $f$ is not eventually constant,
then we can define an order map $q: \Omega \rightarrow \Omega$ so
that $f \circ g$ is injective.  We consider the case when $f$ is
order preserving.  Define $q(0) = 0$. If $q$ has been defined for
all $j < i$, then let $i^{*} =  sup  \{q(i) : i < j \}$ and define
$q(i) =  min   \{ k \in \Omega : f(k) > f(i^{*}) \}$.  This set is
nonempty because $f$ is not eventually constant. Clearly, $ f
\circ g$ is strictly increasing and so is an order injection.
Since $X$ is of countable type, this contradicts Proposition \ref{prop2.5}(b).

\end{proof}\vspace{.5cm}

\begin{prop}\label{prop2.7} Let $\{X_{i} : i \in I \}$ be a LOTS indexed family of nonempty LOTS.
\begin{enumerate}
\item[(a)] The order space sum $\Sigma _{i \in I} X_{i}$ is of countable type iff
$I$ and each $X_{i}$ is of countable type.
\item[(b)] Assume $I = \alpha$ is a positive ordinal.  If $\alpha$ is countable, then
the order space product $\Pi_{i \in I} X_{i}$ is of countable type iff each $X_{i}$
is of countable type.  If each $X_{i}$ is nontrivial and $\Pi_{i \in I} X_{i}$ is of
countable type then $\alpha$ is countable.
\item[(c)] Assume $I = \alpha$ is a positive ordinal and that
 $(\{X_{i} : i \in \alpha \}, \{p_{i}^{j} : i < j < \alpha \})$ is a
 special inverse family.  If $\alpha$ is countable and each $X_{i}$ is
 of countable type, then the inverse limit $\overleftarrow{Lim} \{X_{i}\}$
 is of countable type.  If $\overleftarrow{Lim} \{X_{i}\}$ is of countable
 type, then each $ X_{i}$ is of countable type.
\end{enumerate}
\end{prop}

\begin{proof} (a):  Each $X_{i}$ is a convex
subset of $\Sigma$ and the sum admits an order surjection $\pi :
\Sigma \rightarrow I$.  So if the sum is of countable type, I and
each $X_{i}$ is of countable type by Proposition \ref{prop2.5}(f),(g).

Conversely, assume that $I$ and each $X_{i}$ is of countable type.
Suppose that $f: \Omega \rightarrow \Sigma$ is an order map.  By
Corollary \ref{lem2.6} the order map $\pi \circ f$ is eventually constant, i.e.
for some $\beta \in \Omega$, $\pi \circ f$ takes the constant
value $i \in I$ on the tail $\Omega \setminus \beta$.  Then on the
tail $\Omega \setminus \beta$, $f$ itself maps into $X_{i}$,
regarded as a subset of the sum. Because $\Omega \setminus \beta
\cong \Omega$ Corollary \ref{lem2.6} applied to $X_{i}$ implies that $f $ is
eventually constant. Applying a similar argument to the reverse
$\Omega^{*}$ it follows from Proposition \ref{prop2.5}(b) that $\Sigma$ is of
countable type.

(b):  Assume that $\alpha$ is a countable ordinal.  We prove the
result by induction on $\alpha$.  If $\alpha = 1$, then the product
$\Pi$ is $X_{0}$ and so the product is of countable type iff
$X_{0}$ is.

Now assume, inductively, that the equivalence holds for all $\beta
< \alpha$.\vspace{.25cm}

\textbf{Case 1:} If $\alpha = \beta + 1$, then $\Pi_{i \in \alpha} X_{i}$ is
order isomorphic to the order sum indexed by $\Pi_{i \in \beta}
X_{i}$ with each term a copy of $X_{\beta}$.  So by part
(a)  $\Pi_{i \in \alpha} X_{i}$ is of countable type iff $X_{\beta}$
and $\Pi_{i \in \beta} X_{i}$ are.  By the inductive hypothesis
this holds iff $X_{i}$ is of countable type for all $ i \in \alpha
= \beta \cup \{ \beta \} $.\vspace{.25cm}

\textbf{Case 2:}  Assume now that $\alpha $ is a limit ordinal.

If $\Pi_{i \in \alpha} X_{i}$ is of countable type, then by
Proposition \ref{prop2.5}(g) applied to the order surjections
$\pi^{\beta}_{\alpha} :\Pi_{i \in \alpha} X_{i} \rightarrow\Pi_{i
\in \beta} X_{i}$ the latter are all of countable type. By
inductive hypothesis this implies that $X_{i}$ is of countable
type for all $ i < \beta < \alpha $ and so for all $i \in \alpha $
since $ \alpha $ is a limit ordinal.

On the other hand, if each $X_{i}$ is of countable type and $f:
\Omega \rightarrow \Pi_{i \in \alpha} X_{i}$ is an order map, then
for each $\beta < \alpha$ we can apply the inductive hypothesis
and Corollary \ref{lem2.6} to see that each $\pi^{\beta}_{\alpha} \circ f :
\Omega \rightarrow \Pi_{i \in \beta} X_{i}$ is eventually
constant, i.e. there exists $\epsilon(\beta)$ such that
$\pi^{\beta}_{\alpha} \circ f$ is constant on the tail $\Omega
\setminus \epsilon(\beta)$. Because $\alpha $ is a countable
ordinal $\epsilon (\alpha) =   sup  \{\epsilon(\beta) : \beta
< \alpha \}$ is a countable ordinal.  Clearly, $f$ is constant on
the tail $\Omega \setminus \epsilon (\alpha)$.  Together with a
similar argument for $\Omega^{*}$, this shows that $ \Pi_{i \in
\alpha} X_{i} $ is of countable type.\vspace{.25cm}

This completes the induction for the case when $\alpha$ is assumed
to be countable.

Finally, suppose that each $X_{i}$ is nontrivial.  There exist $a,
b \in \Pi_{i \in \alpha} X_{i}$ such that $a_{i} < b_{i}$ for all
$i \in \alpha$.  Define $f : \alpha + 1 \rightarrow \Pi_{i \in
\alpha} X_{i}$ by
\begin{equation}\label{eq2.24}
f(\beta)_{i} = \begin{cases} b_{i} \quad i< \beta \\ a_{i} \quad
\beta \leq i < \alpha.
\end{cases}
\end{equation}
$f$ is an order injection and if $\alpha$ is uncountable then it
restricts to an order injection of $\Omega$ into $\Pi$.

(c): If $\overleftarrow{Lim}$ is of countable type, then by
Proposition \ref{prop2.5}(g) applied to the surjection $p_{j}:
\overleftarrow{Lim} \rightarrow X_{j}$, each $X_{j}$ is of
countable type.

Now assume that $\alpha$ is countable and each $X_{j}$ is of
countable type.  If $f : \Omega \rightarrow \overleftarrow{Lim}$
is an order map or an order* map, then by Corollary \ref{lem2.6} each
$p_{j}\circ f$ is eventually constant.  That is, there exist
$\epsilon(j) \in \Omega$ and $x_{j} \in X_{j}$ such that
$p_{j}f(\beta) = x_{j}$ for all $\beta \in \Omega\setminus
\epsilon(j)$.  Since $i < j$ implies $p_{i}^{j} \circ p_{j} =
p_{i}$, it follows that $p_{i}^{j}(x_{j}) = x_{i}$.  Hence, there
is a unique $x \in \overleftarrow{Lim}$ such that $p_{j}(x) =
x_{j}$ for all $j \in \Omega$.  Since $\alpha$ is countable there
exists $\epsilon \in \Omega$ such that $\epsilon > \epsilon(j)$
for all $j \in \alpha$.  Clearly, $f(\beta) = x$ for all $\beta
\in \Omega \setminus \epsilon$.  By Proposition \ref{prop2.5}(b),
$\overleftarrow{Lim}$ is of countable type.

\end{proof}\vspace{.5cm}

\begin{prop}\label{prop2.7a} Let $\{X_{i} : i \in \a \}$ be an ordinal  indexed family of nonempty LOTS with $\a$ countable.
If each $X_i$ is first countable, then $X = \Pi_{i \in \alpha} X_{i}$ is first countable and order dense. \end{prop}

\begin{proof} For $x \in X$ with $x \not= max X$, we show that $x$ is the limit of a decreasing sequence in $X$.

Let $K = \{ i < \a : x_i \not= max X_i \}$. Since $x$ is not a maximum point, $K \subset \a$ is nonempty and so is isomorphic to
an ordinal $\g$. \vspace{.25cm}

{\bfseries Case 1:} $ \g = \b + 1$. If $x_{\b}$ a left endpoint of a gap with associated right end-point $a$, then with $z$
$z_i = x_i$ for $i \not= \b$ and $z_{\b} = a$, the pair $x < z$ is a gap pair.\vspace{.25cm}

{\bfseries Case 2:} $ \g = \b + 1$. If $x_{\b}$ not the left endpoint of a gap, then
we can choose $\{ a^n \}$ be a sequence in $X_{\b}$ decreasing and with limit $x_{\b}$.
Define $y^n_i = x_i$ for $i \not= \b$ and $y^n_{\b} = a^n$. \vspace{.25cm}

{\bfseries Case 3:} $ \g$ is a limit ordinal. Choose $\b_n$ an increasing sequence in $K$ converging to $\g$. For each $\b_n$ choose
$a^n \in X_{\b_n}$ with $x_{\b_n} < a^n$. Define $y^n_i = x_i$ for $i \not= \b_n$ and $y^n_{\b_n} = a^n$.\vspace{.25cm}

In  cases 2 and 3, $\{ y^n \}$ is a decreasing sequence in $X$ converging to $x$.

Similarly, for  $x \not= min X$, we can construct an increasing sequence in $X$ converging to $x$ or
find a left endpoint for a gap pair with $x$ on the right.

\end{proof} \vspace{.5cm}

\begin{prop}\label{prop2.8} Let $X$ be an unbounded, order dense LOTS. Let $\hat{X}$
be the completion of $X$ and let $a < b$ be points in $X$.
\begin{enumerate}
\item[(a)] If $A$ is a countable LOTS, then there exists an order injection
$f: A \rightarrow (a,b)$.  If $X$ is countable and $A$ is order dense, then $f$ can be chosen to be an order isomorphism.
\item[(b)]  There exists a continuous order map
$f : [a,b] \rightarrow I$ where $I$ is the unit interval in $\R$ and with $f((a,b))$ dense in $I$.
If $X$ is complete, then $f$ is surjective. If $X$ is
complete and separable, then $f$ can be chosen to be an isomorphism.
\item[(c)]  If $\alpha$ is a positive ordinal and $f : \alpha + 1 \rightarrow X$ is
an order embedding with $f(0) = a$ and $f(\alpha) = b$, then $\{[f(i),f(i+1)) : i \in \alpha \}$ is a convex partition of $[a,b)$ and so
we can identify
\begin{equation}\label{eq2.25}
[a,b) \cong \Sigma_{i \in \a} \  [f(i),f(i+1))
\end{equation}
expressing $[a,b)$ as the order sum of an $\alpha$ indexed family
of intervals.
\item[(d)] If $X$ is first countable and $\alpha$ is a positive, countable ordinal,
then there exists an order embedding $f : \alpha + 1 \rightarrow X$ with $f(0) = a$ and $f(\alpha) = b$.
\item[(e)] The following conditions are equivalent.
\begin{itemize}
\item[(i)] $X$ is of countable type.
\item[(ii)] $\hat{X}$ is of countable type.
\item[(iii)] $\hat{X}$ is first countable and $\sigma$-bounded.
\item[(iv)] There does not exist $f: \Omega \rightarrow \hat{X}$ a continuous
injective map which is either order preserving or order reversing.
\end{itemize}
\end{enumerate}
\end{prop}

\begin{proof} (a): This is a standard inductive
argument using a counting of the points of $A$.  If $A$ is order dense and $X$, too, is
countable, then one counts $X$ as well and proceeds back and forth
between $A$ and $X$ to build the isomorphism.

(b):  With $A$ the set of rationals in $I$, use (a) to get an order
injection $g : A \rightarrow [a,b]$ such that $g(0) = a$ and $g(1)
= b$.  Define for $x \in [a,b]$
\begin{equation}\label{eq2.26}
f(x) \ = \  sup  \  g^{-1}([a,x])\   = \  inf  \ g^{-1}([x,b]).
\end{equation}

By Proposition \ref{prop2.1}(a), $f$ is continuous
with $f(g(c)) = c$ for $c \in A$.

If $X$ is complete, then it is connected and so the connected dense image is $I$.

If $X$ is separable, choose $g$ with a dense image in $[a,b]$, then $f$ is injective by Proposition \ref{prop2.1} (e).

(c): If $x \in [a,b)$, then $b = f(\alpha) > x$.  Let $\beta =  min \{ j \in \alpha + 1 : f(j) > x \}$.  Since $a = f(0) \leq x$,
$\beta $ is positive and by continuity of $f$, $\beta$ is not a
limit ordinal.  Hence, $\beta = i+1$ for some $i \in \alpha$ and
$x \in [f(i), f(i+1))$.

(d): We construct the embedding $f$ by induction on $\alpha$.  If
$\alpha = 1$, let $f(0) = a $ and $f(1) = b$.  Now assume the result
is true for all $\beta < \alpha$.\vspace{.25cm}

\textbf{Case 1:}  If $\alpha = \beta + 1$, then choose $\tilde{b}$ so that
$a < \tilde{b} < b$. By inductive hypothesis there exists an order
embedding $\tilde{f} : \beta +1 \rightarrow [a,\tilde{b}]$ with
$\tilde{f}(0) = a$ and $\tilde{f}(\beta) = \tilde{b}$. Extend the
definition by $f(\alpha) = b$ to get $f$.\vspace{.25cm}

\textbf{Case 2:} If $\alpha$ is a limit ordinal, then because it is
countable, there exists an increasing cofinal sequence $\{
\beta_{n} \}$ in $\alpha$.  Because $X$ is first countable and
order dense there exists an increasing sequence $\{ x_{n} \}$ in
$(a,b)$ with limit $b$.  By inductive hypothesis there exists an
order embedding of the interval $( \beta_{n}, \beta_{n+1}]$ in
$\alpha$ to $(x_{n}, x_{n+1}]$ with $\beta_{n+1}$ mapped to
$x_{n+1}$.  Put these together and map $0$ to $a$ and $\alpha$ to
$b$ to get $f$.\vspace{.25cm}

(e): (1)$\Leftrightarrow$(ii) by Proposition \ref{prop2.5}(f)(h) since the
inclusion of $X$ into $\hat{X}$ is continuous by Proposition \ref{prop2.1}(a).

(ii)$\Leftrightarrow$(iii) by Proposition \ref{prop2.5}(d).

(ii)$\Rightarrow$(iv) by Proposition \ref{prop2.5}(b).

(iv)$\Rightarrow$(ii) If $\hat{X}$ is not of countable type, then by
Proposition \ref{prop2.5}(b) there exists an injective map $\tilde{f} : \Omega
\rightarrow \hat{X}$ which is either order preserving or order
reversing.  Without loss of generality assume that $\tilde{f}$ is
an order map.  Define $f : \Omega \rightarrow \hat X$  by
\begin{equation}\label{eq2.27}
f(\beta) = \begin{cases} sup \ \{ \tilde{f}(i) : i < \beta \}
\quad \beta \  \mbox{is a limit ordinal} \\ \tilde{f}(\beta)
\hspace{8em} \mbox{otherwise}. \end{cases}
\end{equation}
It is easy to check that $f$ is a continuous, injective order map.

\end{proof} \vspace{1cm}

\section{\textbf{Complete Homogeneous LOTS}}

\subsection{Doubly Transitive and Homogeneous LOTS}

If a group $G$ acts on a set $S$, then we say that $s_{1}$ and $s_{2}$
in $S$ are \emph{$G$ equivalent}\index{$G$ equivalent points} if $g(s_{1}) = s_{2}$ for some $g \in
G$.  We say that $G$ \emph{acts transitively}\index{$G$ acts transitively}  on $S$ when all
points are $G$ equivalent.

For a topological space $X$ we let $H(X)$\index{$H(X)$} denote the automorphism
group of $X$, i.e. the group of homeomorphisms from $X$ to itself.
We call $X$ \emph{topologically homogeneous}\index{homogeneous} if $H(X)$ acts
transitively on $X$.  If $X$ is a LOTS, then any order automorphism
of $X$, i.e. order preserving bijection of $X$ to itself, is a
homeomorphism.  We denote by $H_{+}(X)$\index{$H_{+}(X)$} the subgroup of order
automorphisms and by $H_{\pm}(X)$\index{$H_{\pm}(X)$} the subgroup of bijections which
either preserve or reverse order.  If $X$ admits an order
reversing homeomorphism then we call $X$ a \emph{symmetric} LOTS\index{LOTS!symmetric}
in which case $H_{+}(X)$ is a subgroup of $H_{\pm}(X)$ of index 2.
Otherwise, $H_{+}(X) = H_{\pm}(X)$.

We call a LOTS $X$ \emph{$\pm$transitive}\index{LOTS!$\pm$transitive} if  $H_{\pm}(X)$ acts
transitively on $X$ and \emph{transitive}\index{LOTS!transitive} if  $H_{+}(X)$ acts
transitively on $X$.  $X$ is \emph{doubly transitive}\index{LOTS!doubly transitive} if
$H_{+}(X)$ acts transitively, via the diagonal action, on the set
$\{(x_{1},x_{2}) : x_{1} < x_{2} \} \subset X \times X$. Since any order
automorphism of $X$ is also an order automorphism of the reverse
$X^{*}$, i.e.$H_{+}(X) = H_{+}(X^{*})$, the reverse LOTS $X^*$ is $\pm$transitive, transitive or doubly transitive  when
the corresponding property holds for $X$.

\begin{lem}\label{lem3.1} If $X_{1} $ and $X_{2}$ are LOTS with $X_{1}$ connected and
$ f: X_{1} \rightarrow X_{2} $ is a continuous  map then the image $f(X_{1})$
is a convex subset of $X_{2}$.  If, in addition $f$ is injective, then it is
either order preserving or order reversing.
\end{lem}

\begin{proof} The image of $f$ is connected by
continuity.  Any connected subset of a LOTS is convex.

Assume now that $f$ is injective.  If $f(a) < f(c) <
f(b)$, then the image of the open interval between $a$ and
$b$ is connected and so contains $f(c)$.  Since $f$ is
injective, $c$ must
therefore lie in the interval. That is, either
$a < c < b$ or $a > c > b$.  Thus, on each
triple of points in $X_{1}$ $f$ either preserves or reverses
order.  If $f$ preserves the order of some pair $a,b$ in $X_{1}$
then it preserves the order of every triple which includes $a,b$
and so of every pair which includes either $a$ or $b$.  So it
preserves every triple which includes $a$ and so preserves every
pair.  The remaining possibility is that $f$ reverses every pair.

\end{proof}\vspace{.5cm}

\begin{prop}\label{prop3.2}Let $X$ be a LOTS with at least three points.
\begin{enumerate}
\item[(a)] If $X$ is connected, then $H(X) = H_{\pm}(X)$ and any order reversing
homeomorphism has a unique fixed point.  In addition, $X$ is transitive if it is topologically homogeneous.
\item[(b)] If $X$ is $\pm$transitive, then it is unbounded.  If $X$ is transitive, then it is $\pm$transitive.
\item[(c)]  $X$ is transitive iff for all $a,b \in X$
\begin{equation}\label{eq3.1}
(a,\infty) \cong (b,\infty) \qquad \mbox{and} \qquad (-\infty,a)
\cong (-\infty,b).
\end{equation}
If $X$ is transitive, then it is either discrete, i.e. every point
is isolated, or it is order dense.  If $X$ is transitive and
complete, then $X$ is first countable and it is either order
isomorphic to $\Z$, the LOTS of integers, or it is
connected.
\item[(d)]  The following conditions are equivalent.
\begin{itemize}
\item[(i)]  $X$ is doubly transitive.
\item[(ii)]  $X$ and every nonempty open subinterval of $X$ are transitive.
\item[(iii)]  $X$ has no $min$ and $(a,\infty)$ is transitive for every $a \in X$.
\item[(iv)]  $X$ has no $max$ and $(-\infty,a)$ is transitive for every $a \in X$.
\item[(v)]  $X$ is transitive and $(a,\infty)$ is transitive for some $a \in X$.
\item[(vi)]  $X$ is transitive and $(-\infty,a)$ is transitive for some $a \in X$.
\item[(vii)]  $X$ is unbounded and any two nontrivial, closed, bounded subintervals of $X$ are order isomorphic.
\item[(viii)]  For every positive integer $n$, $H_{+}(X)$ acts transitively on
\begin{displaymath}
\{ (x_{1},...,x_{n}) : x_{1} < ...< x_{n} \} \subset X^{n} .
\end{displaymath}
\end{itemize}
\item[(e)] If $X$ is doubly transitive, then it is order dense and unbounded and every nonempty convex open subset
is doubly transitive.

\end{enumerate}
\end{prop}

\begin{proof} (a): $H(X) = H_{\pm}(X)$ by Lemma
\ref{lem3.1}.  If $f$ is order reversing on $X$ and $ y = f(x) > x $, then
$f(y) < f(x) = y $.  Since $f$ is not the identity $\{ x : f(x) > x \}$
and $\{x : f(x) < x \}$ are disjoint nonempty open subsets of
$X$.  Because $X$ is connected their union is a proper subset and
so some fixed point $e$ exists.  If $x > e$, then $f(x) < f(e) = e$
and so $x$ is not a fixed point.  Similarly, if $x < e$.  Thus,
$e$ is the unique fixed point.

If $H(X)$ acts transitively, then $H(X) = H_{\pm}(X) $ implies $X$
is $\pm$transitive.  If $H_{+} = H_{\pm}$, then $X$ is transitive.
On the other hand, if $X$ is symmetric, let $f_{0}$ be an order
reversing homeomorphism with fixed point $e$.  If $x \in X$, then
there exists $f \in H_{\pm}$ such that $f(x) = e$.  Also,
$f_{0} \circ f(x) = e$ and either $f$ or $f_{0} \circ f$ is in
$H_{+}(X)$.  Thus, every $x \in X$ is $H_{+}(X)$ equivalent to $e$
and so $X$ is transitive.

(b):  If $X$ has a $max$, then the $max$ is fixed by any element of
$H_{+}(X)$ and is mapped to $min$ by any order reversing
isomorphism.  Hence, if $X$ is $\pm$transitive and has at least
three points, then it is unbounded.

Since $H_{+} \subset H_{\pm}$, transitivity implies
$\pm$transitivity.

(c):  If $f \in H_{+}(X)$ maps $a$ to $b$, then it restricts to an
order isomorphism of $(a,\infty)$ with $(b,\infty)$ and of
$(-\infty,a)$ with $(-\infty,b)$.  Conversely, we can put together
isomorphisms $(a,\infty) \cong (b,\infty)$ and $(-\infty,a) \cong
(-\infty,b)$ to get an automorphism which maps $a$ to $b$.

Now assume that $X$ is transitive but not order dense.  There
exists a gap pair with $a < b$ left and right endpoints, respectively.
As all points of $X$ are $H_{+}$ equivalent all points are both
left and right endpoints.  That is, every point of $X$ is
isolated.  If, in addition, $X$ is complete, choose $f(0) \in X$
and inductively define $f(n+1) =  inf  (f(n),\infty)$ and $f(-n-1)
=  sup  (-\infty,f(-n))$ for every nonnegative integer $n$.  This
defines an order injection $f : \Z \to X$.
Since $X$ is discrete, $f(\Z)$ has no $sup$ or $inf$
and so is $\pm$cofinal  in $X$.  For $x \in X$, $f(n) < x$
for some $n \in \Z$ but not for all. If $n$ is the
largest such, then $f(n+1) \leq x$  since the interval
$(f(n),f(n+1))$ is empty.  By maximality of $n$, $x = f(n+1)$.
Hence, $f$ is surjective and so is an order isomorphism of
$\Z$ with $X$.  On the other hand, if $X$ is order dense
and complete then it is connected.

$\Z$ is first countable and in the connected case there
exists a bounded increasing sequence which converges to some point
$a$ by completeness. Similarly, some bounded decreasing sequence
converges to a point $b$.  As every $x \in X$ is $H_{+}$
equivalent to both $a$ and $b$, every $x$ is the limit of
increasing and decreasing sequences.  Thus, $X$ is first
countable.

(d) (i)$\Rightarrow$(iii)\&(iv):  If $a < b < c$ in $X$ then the pair
$a,b$ is $H_{+}$ equivalent to $b,c$ by double transitivity.
Hence, $a$ is not the $min$ and $c$ is not the $max$.  Thus, $X$
is unbounded.  If $x_{1},x_{2} \in (a,\infty)$ and $f \in H_{+}$
maps the pair $a,x_{1}$ to the pair $a,x_{2}$, then $f$ restricts
to an automorphism of $(a,\infty)$ which maps $x_{1}$ to $x_{2}$,
proving (iii).  Similarly, for (iv).

(ii)$\Rightarrow$(iii)\&(iv):  By (b) $X$ is unbounded.  Hence,
$(a,\infty)$ and $(-\infty,a)$ are open, nonempty subintervals of
$X$ and they are transitive by assumption (ii).

(iii)$\Leftrightarrow$ (v) and (iv)$\Leftrightarrow$ (vi): If $X$ is transitive,
then it is unbounded and all $(a,\infty)$ intervals
are isomorphic.
Hence (v) $\Rightarrow$ (iii). On the other hand, if $x < y$ in $X$ and $X$ has
no $min$, then there exists $a < x$ and so (iii) implies there
exists an automorphism of $(a, \infty)$ which maps $x$ to $y$. Extend by the
identity on $(-\infty,a]$. Hence, $X$ is transitive. The
proof of (iv)$\Leftrightarrow$ (vi) is similar.

(iii)$\Rightarrow$(viii): Given $x_{1} < ...< x_{n}$ and $y_{1} < ...<
y_{n}$ choose $a$ smaller than all of them.  We construct, by
induction on $n$, $f \in H_{+}((a,\infty))$ such that $f(x_{i}) =
y_{i}$ for $i = 1,...,n$.  Assume that $f_{1} \in
H_{+}((a,\infty))$ satisfies $f(x_{i}) = y_{i}$ for $i = 1,...n-1$
and let $\tilde{x}_{n} = f_{1}(x_{n}) > f_{1}(x_{n-1}) = y_{n-1}$.
Choose $f_{2} \in H_{+}((y_{n-1},\infty))$ such that
$f_{2}(\tilde x_{n}) = y_{n}$ and extend $f_{2}$ by the identity
on $(a,y_{n-1}]$ to get $f_{2} \in H_{+}((a,\infty))$.  Let $f =
f_{2} \circ f_{1}$.  Having obtained $f$, extend by the identity
on $(-\infty,a]$ to get $f \in H_{+}(X)$ mapping $x_{1} < ...<
x_{n}$ to $y_{1} < ...< y_{n}$.

(iv)$\Rightarrow$(viii):  Use a similar proof or apply (iii)$\Rightarrow$
(viii) to $X^{*}$.

(viii)$\Rightarrow$(vii):  Since condition (viii) clearly implies double
transitivity, X is unbounded by (i)$\Rightarrow$(iii)\&(iv).  If
$x_{1} < x_{2}$ and $y_{1} < y_{2}$, then $f \in H_{+}(X)$ which
maps the pair $x_{1}, x_{2}$ to $y_{1}, y_{2}$ restricts to an
isomorphism from $[x_{1},x_{2}]$ to $[y_{1},y_{2}]$, proving (vii).

(vii)$\Rightarrow$ $X$ is order dense and (i)\&(ii): Given
$x_{1}, x_{2} \in X$ we can choose $a, b \in X$ with $a <  x_{1}, x_{2}  < b$
because $X$ is unbounded. Put together order isomorphisms
\begin{equation}\label{eq3.3}
[a,x_{1}] \cong [a,x_{2}]\qquad [x_{1},b]
\cong [x_{2},b]
\end{equation}
 and extend by the identity outside $(a,b)$ to get
$f \in H_{+}(X)$ which maps $x_{1}$ to $x_{2}$.  This proves that
$X$ is transitive.

If $X$ were discrete, then we could choose
points $a < b < c$ with $(a,b) = (b,c) = \emptyset$.  But then
$[a,b] = \{a, b\}$ is not be isomorphic to $[a,c] = \{a, b, c \}$.
So by (c) $X$ is order dense.

Now we show that if  $J$ is a nonempty,
open convex subset of $X$, then $J$ is doubly transitive.  and $x_{1},x_{2} \in J$, then because $X$
is unbounded and order dense, we can choose $a,b \in J$ such that
$a < x_{1},x_{2} < b$ and just as before get $f \in H_{+}((a,b))$
mapping $x_{1}$ to $x_{2}$.  Extend by the identity outside of
$(a,b)$ to get $f \in H_{+}(J)$.  Hence, $J$ is transitive.
Assume $x_{1} < x_{2}$ and $y_{1} < y_{2}$ in $J$, choose $a$
smaller and $b$ larger than all four of them.  Put together
isomorphisms
\begin{equation}\label{eq3.4}
[a,x_{1}] \cong [a,y_{1}] \qquad [x_{1},x_{2}] \cong [y_{1},y_{2}]
\qquad [x_{2},b] \cong [y_{2},b]
\end{equation}
and extend by the identity outside $(a,b)$ to get $f \in H_{+}(J)$ which
maps the pair $x_{1},x_{2}$ to $y_{1},y_{2}$.  Observe that all
these intervals are nonempty because $X$ is order dense.  Thus,
$J$ is doubly transitive.

In particulary, all this proves (e).

\end{proof} \vspace{.5cm}

 We call
$X$ \emph{weakly homogeneous}\index{LOTS!weakly homogeneous} if it has at least three points and is order isomorphic with every
nonempty, bounded, open subinterval of itself.  Clearly, the reverse LOTS $X^*$ is  weakly homogeneous when
 $X$ is.

\begin{prop}\label{prop3.2a}Let $X$ be a LOTS with at least three points.
\begin{enumerate}
\item[(a)]  $X$ is weakly homogeneous iff for every $a \in X$
\begin{equation}\label{eq3.2}
(a,\infty) \ \cong \ X \ \cong \ (-\infty,a).
\end{equation}
If $X$ is weakly homogeneous, then it is doubly transitive and it
is order isomorphic with every nonempty, open subinterval of
itself (whether bounded or not).
\item[(b)]  If $X$ is doubly transitive,  then every nonempty, bounded, open
subinterval of $X$ is weakly homogeneous.  If $X$ is doubly transitive and
first countable, then it is weakly homogeneous iff it is $\sigma$-bounded.
\item[(c)]  Assume $X$ is transitive.  If for $a \in X$, $(-\infty,a)$ is
symmetric, then $X$ is doubly transitive.  If for $a \in X$, $(-\infty,a)$
and $(a,\infty)$ are symmetric, then $X$ is weakly homogeneous.
\end{enumerate} \end{prop}

\begin{proof}
(a): Assume that for all $a \in X$, $(-\infty,a) \cong X \cong
(a,\infty)$.  If $a < b$ in $X$, let $f : (a,\infty) \rightarrow
X$ be an order isomorphism and let $\tilde{b} = f(b)$.  Let
$\tilde{f} : X \rightarrow (-\infty,\tilde{b})$ be an order
isomorphism. Then $f^{-1} \circ \tilde{f} : X \rightarrow (a,b) $
is an order isomorphism.  It follows that $X$ is weakly
homogeneous.

Conversely, assume that $X$ is weakly homogeneous.  We prove first
that if $a < x < b$  in $X$ then $x$ is not an isolated point.  If
it were then we could choose $a$ and $b$ so that $(a,x) = (x,b) =
\emptyset$ and so $(a,b) = \{x\}$ could not be isomorphic to $X$.
It follows that $X$ is infinite.  If $b$ were the $max$ of $X$
then we could choose $a$ and $x$ so that $(a,x)$ and $(x,b)$ are
infinite.  Since $(a,x) \cong X$ the interval $(a,x)$ has a $max$
which we call $y$. Hence, $(y,x) = \emptyset $ and $(a,y)$ is
infinite.  Since $(a,y) \cong X$ the interval $(a,y)$ has a $max$
and so the point $y$ is isolated.  This contradiction implies that
$X$ has no $max$.  Similarly, there is no $min$. Now Proposition \ref{prop3.2}(d)
(vii)$\Rightarrow$(i) shows that $X$ is doubly transitive and so it
is order dense by Proposition \ref{prop3.2}(e).

Now given $a \in X$, $(a,\infty)$ and $(-\infty,a)$ are nonempty
since $a$ is neither $max$ nor $min$.  Choose $f : X \rightarrow
(b,c)$ an isomorphism with $b < c$ in $X$ and let $\tilde{a} =
f(a)$.  $f$ induces isomorphisms $(-\infty,a) \cong (b,\tilde{a})$
and $(a,\infty) \cong (\tilde{a},c)$.  Since neither is empty,
each is isomorphic with $X$.  Thus, $X$ is isomorphic with every
nonempty, open subinterval of itself.

(b):  If $X$ is doubly transitive and $a < x_{1} < x_{2} < b$ then
there exists $f \in H_{+}(X)$ which maps the pair $a,b$ to
$x_{1},x_{2}$.  This restricts to an order isomorphism $(a,b)
\cong (x_{1},x_{2})$.  Hence, $(a,b)$ is weakly homogeneous.

If $X$ is first countable then we can use Proposition \ref{prop2.8}(d) to
construct an order injection $f : \Z \to (a,b)$
whose image is $\pm$cofinal.  Similarly, if X is
$\sigma$-bounded we can construct an order injection $\tilde{f}
: \Z \to X$ whose image is $\pm$cofinal. If, in addition, $X$ is doubly transitive, then we can choose for each
$i \in \Z$ an order isomorphism $[f(i),f(i+1)] \cong
[\tilde{f}(i),\tilde{f}(i+1)]$ and put them together to get an
order isomorphism $(a,b) \cong X$.  Hence, $X$ is weakly
homogeneous.

If $X$ is weakly homogeneous then $X \cong (a,b)$ for $a < b$ and
the latter is $\sigma$-bounded if $X$ is first countable.

(c):  Since $X$ is transitive it is unbounded.  In addition, given
$a_{1} < b_{1}$ and $a_{2} < b_{2}$ in $X$ we can apply an element
of $H_{+}(X)$ to move the pair $a_{1},b_{1}$ so that we can assume
$b_{1} = b_{2}$ (hereafter denoted $b$). Let $q$ be an orientation
reversing homeomorphism of $(-\infty,b)$.  Define $\tilde{a}_{i} =
q(a_{i}) $ for $ i = 1,2$.  Since $X$ is transitive there exists an
 isomorphism $f : (-\infty,\tilde{a}_{1}) \rightarrow
(-\infty,\tilde{a}_{2})$.  Then $q^{-1} \circ f \circ q$ restricts
to an isomorphism $(a_{1},b) \cong (a_{2},b)$.  So $X$ is
doubly transitive by Proposition \ref{prop3.2}(d) (vii)$\Rightarrow$(i).  If, instead, we use
an  isomorphism $f_{1} : (-\infty,\tilde{a}_{1}) \rightarrow
(-\infty,b)$ then $q^{-1} \circ f_{1} \circ q$ restricts to an
 isomorphism $(a_{1},b) \cong (-\infty,b)$.

Now if $(b,\infty)$ is symmetric as well then we can similarly
construct an  isomorphism $(b,a_{3}) \cong (b,\infty)$ with
$b < a_{3}$.  Putting this together with the isomorphism $(a_{1},b) \cong (-\infty,b)$ above we get an isomorphism
$(a_{1},a_{3}) \cong X$.  Thus, $X$ is isomorphic to nonempty,
bounded, open subintervals of itself and so is  weakly homogeneous.
Since $X \cong (a,b) \cong (-\infty,b) $ it is symmetric as well.

\end{proof} \vspace{.5cm}

Since a nontrivial transitive LOTS $X$ is unbounded, it cannot be compact. There is a
condition due to G. D. Birkhoff, \cite{Bi} page 47, which we will call the
\emph{closed interval condition}\index{closed interval condition}. A LOTS $X$ satisfies the closed interval condition
when any two nontrivial, closed bounded subintervals are isomorphic. Clearly, any nontrivial convex subset of such a LOTS satisfies the closed
interval condition as well.

By Proposition \ref{prop3.2}(d) an unbounded LOTS satisfies the closed interval
condition iff it is doubly transitive. If $X$ is bounded with $min = m$ and $max = M$ and $X$ satisfies the closed interval condition, then
$X \setminus \{m. M \}$ is weakly homogeneous. Conversely, if $X_0$ is weakly homogeneous and $X = \{ m \} + X_0 + \{ M \}$, then
$X$ satisfies the closed interval condition.

 \begin{ex}\label{example3.3a} Linearly ordered groups, fields and products \end{ex}

 An ordered group is a group with a linear order such that the translation maps are order isomorphisms. Since a group acts transitively
 on itself by translation, it follows that an ordered group, like $\Z$, is a transitive LOTS. An ordered field is a field
 with a linear order such that the additive translation maps, and multiplication by positive elements are order isomorphisms. An ordered
 field, like $\Q$ and $\R$, is doubly transitive. Observe that if $x_1 < x_2$ in the field, then $x \mapsto (x - x_1)\cdot(x_2 - x_1)^{-1}$ maps
 $(x_1,x_2)$ isomorphically onto $(0,1)$.  In fact, an ordered field is weakly homogeneous, because the map:
 \begin{equation}\label{eqfieldmap}
 f(x) \ = \ \begin{cases} x \cdot (1 + x)^{-1} \ \text{ for } \ -1 < x \le 0, \\ x \cdot (1 - x)^{-1} \ \text{ for } \ 0 \le x  < 1, \end{cases}
 \end{equation}
 is an isomorphism from $(-1,1)$ to the entire field.

 If $\{ X_i : i \in I \}$ is a family of  LOTS indexed by an ordinal $I$, and $h_i$ is an automorphism of $X_i$ for each $i$,
 then $\prod_i \ h_i$ is an automorphism of the order space product $\prod_i X_i$. It follows if each $X_i$ is transitive, then
 the product  is transitive as well. If each $X_i$ is an ordered group, then the product is an ordered group with coordinate-wise
 addition.

 On the other hand, since a discrete LOTS with at least three points is never doubly transitive, it follows that the discrete, linearly
  ordered groups $\Z$, $\Z \times \Z$  and $\R \times \Z$ are not doubly transitive.

  The groups $\Q \times \Q$ and  $\R \times \R$ are order dense. Because $\Q \times \Q$ is a countable order dense LOTS, it is order isomorphic
  to the ordered field $\Q$ by Proposition \ref{prop2.8} (a) and so it is doubly transitive. On the other hand, $\R \times \R$ is not doubly transitive.
  Let $x_1 = (0,0), x_2 = (0,1),  x_3 = (1,0)$. The map $t \mapsto (0,t)$ is an order isomorphism from the unit interval $(0,1)$ in $\R$ onto
  the interval $(x_1,x_2)$ and, in particular, the latter interval is separable. On the other hand, if $D$ is a countable subset of
  $(x_1,x_3)$, then we can choose $t$ such that $0 < t < 1$ and $t$ is not the first coordinate of a point of $D$. Hence, the open
  interval $((t,0),(t,1))$ is disjoint from $D$. It follows that $(x_1,x_3)$ is not separable and so is not isomorphic to $(x_1,x_2)$.

   With $-1,+1 \in \R$
 the order space product $\R \times (-1,+1)$ is isomorphic to $\R \times \R$ and so is
transitive, locally connected and of countable type.  Its completion
is the order product $\R \times [-1,+1]$ which is
connected and of countable type, but not transitive. Observe that the union of the separable open subsets of $\R \times [-1,+1]$ is
$\R \times (-1,+1)$. That is, a point of the form $(t,\pm 1)$ does not have any separable neighborhood.

The product $\Z \times \R$ is isomorphic to the complement of $\Z$ in $\R$ and $\Q \times \R$ is isomorphic to the complement of the Cantor Set
$C$ in $[0,1]$. Both of these are transitive, but not doubly transitive.  Any homeomorphism maps components to components and so cannot map
a pair of points contained in a component to a pair in different components. Both of these provide examples of dense, transitive subsets of $\R$ which are
not doubly transitive. However, neither of these has dense holes.

Finally, we will see below that because its complement $\Q$ is doubly transitive, the set $\I$ of irrationals is doubly transitive. Hence, the
product $\Q \times \I$ is transitive. It is isomorphic to $X = [0,1] \setminus (\Q \cup C)$. Between any two points of $X$ in the same
component of $[0,1] \setminus C$ there are only countably many holes, whereas between two points of $X$ in different components of
$[0,1] \setminus C$ there are uncountably many holes.   Hence, $X$ provides an example in $(0,1) \cong \R$ of a dense, transitive subset with dense holes,
but which is not doubly transitive.

 We do not know whether or not a
  dense additive subgroup of $\R$ is necessarily doubly transitive. \vspace{.5cm}

  In contrast with the above examples Treybig proved the following beautiful result which we state using our language.

  \begin{theo}\label{theotreybig} If $X$ is a connected, unbounded transitive LOTS, then $X$ is doubly transitive.\end{theo}

  \begin{proof} See \cite{Tr}. As we will not need the result below, we include the proof in an Appendix.

  \end{proof} \vspace{.5cm}

  The following uses the easy part of an argument due to Babcock \cite{B}.

  \begin{theo}\label{theobabcock} If $\a$ is a  tail-like ordinal and $X$ is doubly transitive, then $X^{\a}$ is doubly transitive. \end{theo}

  \begin{proof} Choose two points of $X$ which we label $-1 < +1$. Let $x^{\pm}$ be defined by $x^+_i = +1, x^-_i = -1$ for all $i < \a$.
  Using transitivity of $X$ we choose for any $a \in X$ automorphisms $g^{\pm}_a$  of $X$ with $g^+_a(+1) = a, g^-_a(-1) = a$. Given $y < z \in X^{\a}$ we construct
  an isomorphism $ f : [x^-,x^+] \to [y,z]$.

  Let $\b = min \{ i : y_i \not= z_i \} < \a$ so that $y_i = z_i$ for all $i < \b$ and $y_{\b} < z_{\b}$. Let $g^0 : [-1,+1] \to [y_{\b},z_{\b}]$
  be an isomorphism. Because $\a$ is tail-like there is an isomorphism $\g : \a \to \a \setminus \b$. For $x \in [x-, x+]$
 define $f(x)$ by:
  \begin{itemize}
  \item $f(x)_i = y_i = z_i$ for $i < \b$.

  \item $f(x)_{\b} = g^0(x_0) = g^0(x_{\g(\b)})$.

  \item $f(x)_i = g^-_{y_i}(x_{\g(i)})$ if $\b < i < \a$ and $x_{\g(j)} = -1$ for all $\b \le j < i$.

   \item $f(x)_i = g^+_{z_i}(x_{\g(i)})$ if $\b < i < \a$ and $x_{\g(j)} = +1$ for all $\b \le j < i$.

   \item $f(x)_i = x_{\g(i)}$ otherwise.
   \end{itemize}

   For $\b \le i < \a$, $f(x)_i$ depends only on the values of $x_k$ with $k \le \g(i)$ and so it easily follows that
   $f$ is an order injection of $[x-,x+]$
into    $[y,z]$. On the other hand, it is easy to reverse the procedure to see that
   $f$ is surjective and so is an isomorphism.

\end{proof} \vspace{.5cm}

\noindent {\bfseries Remark.} Since $\a$ need not be countable, beginning with $X = \R$ this constructs doubly transitive LOTS of arbitrary
cardinality. \vspace{.5cm}

Using the completion $\hat{X}$ for an order dense LOTS $X$, we now
 describe the strong homogeneity condition that we want to focus on.

\begin{prop}\label{prop3.3} Let $\hat{X}$ be the completion of $X$, a nontrivial, order dense LOTS.
\begin{enumerate}
\item[(a)] The following conditions are equivalent.
\begin{itemize}
\item[(i)]  $X$ and $\hat{X}$ are doubly transitive.
\item[(ii)]  $X$ is doubly transitive and $\hat{X}$ is first countable.
\item[(iii)]  Any two nonempty, open, bounded, convex subsets of $X$ are order isomorphic.
\end{itemize}
\item[(b)]  The following conditions are equivalent.
\begin{itemize}
\item[(i)]  $X$ and $\hat{X}$ are weakly homogeneous.
\item[(ii)]  $X$ is doubly transitive and of countable type.
\item[(iii)]  $X$ is order isomorphic with every nonempty, open, convex subset of $X$.
\end{itemize}
\end{enumerate}
\end{prop}

\begin{proof} Since $X$ is order dense and nontrivial, it is
infinite.

(a)(i)$\Rightarrow$(ii): Since $\hat{X}$ is doubly transitive, it is
transitive.  Completeness implies first countability by
Proposition \ref{prop3.2}(c).

(ii)$\Rightarrow$(iii):  Bounded open convex sets are of the form
$(z_{1},z_{2}) \cap X$ and  $(w_{1},w_{2}) \cap X$ with $z_{1} <
z_{2}$ and  $w_{1} < w_{2}$ in $\hat{X}$.  Because $\hat{X}$ is
first countable and connected there exist order injections
$\tilde{g},\tilde{f} : \Z \rightarrow \hat{X}$ whose
images are $\pm$cofinal in $(z_{1},z_{2})$ and
$(w_{1},w_{2})$, respectively.  We choose for each $i \in \Z$,
$f(i) \in (\tilde{f}(i),\tilde{f}(i+1)) \cap X$ and $g(i)$
similarly so as to get $ g, f : \Z \rightarrow X $. In
the now familiar way we put together isomorphisms $(f(i),f(i+1))
\cong (g(i),g(i+1))$ to get the required isomorphism between the
open convex sets.

(iii)$\Rightarrow$(i):  If $X$ had a $max = b$ and $ a < b$ then
$(a,b)$ and $(a,b]$ are both bounded, open convex subsets of $X$.
Since $X$ is order dense $(a,b)$ has no $max$ and so cannot be
isomorphic to $(a,b]$.  Thus, condition (iii) implies that $X$ has
no $max$ and similarly no $min$.  If $z_{1} < z_{2}$ and $w_{1} <
w_{2}$ in $\hat{X}$, then any isomorphism $(z_{1},z_{2}) \cap X
\cong (w_{1},w_{2}) \cap X$ extends as in (\ref{eq2.20}) to an isomorphism
$(z_{1},z_{2}) \cong (w_{1},w_{2})$ in $\hat{X}$. Hence, both $X$
and $\hat{X}$ satisfy condition (vii) of Proposition \ref{prop3.2}(d). Since
(vii)$\Rightarrow$(i) there, $X$ and $\hat{X}$ are doubly
transitive.

(b): By Proposition \ref{prop3.2a}(a) (i) here implies condition (i) of (a).  By
Proposition \ref{prop2.8}(d) (ii) here  implies condition (ii) of (a).  Clearly, (iii) here implies condition (iii)
of (a).  Hence, in proving the equivalences we can use all of the
conditions of (a).

(i)$\Rightarrow$(ii): By weak homogeneity of $\hat{X}$, $\hat{X}$ is
isomorphic to any bounded, nonempty, open interval and such an interval is
$\sigma$-bounded because $\hat{X}$ is first countable.  Hence,
$\hat{X}$ is first countable and $\sigma$-bounded.  Thus, $X$ is
of countable type by Proposition \ref{prop2.8}(d).

(ii)$\Rightarrow$(iii): Because $\hat{X}$ is first countable, $X$ is and $X$ is
$\sigma$-bounded as well as doubly transitive by assumption. Therefore, it is weakly
homogeneous by Proposition \ref{prop3.2a}(b).  Let $J$ be any nonempty, open,
convex subset of $X$ and $a < b$ in $X$.  By weak homogeneity
there is an isomorphism $f : X \rightarrow (a,b)$ and $f(J)$ is a
bounded, open, convex subset of $X$.  So by condition (iii) of (a)
$f(J) \cong (a,b)$.  So $J \cong f(J) \cong (a,b) \cong X$.

(iii)$\Rightarrow$(i):  Clearly, $X$ is weakly homogeneous.  If
$w_{1} < w_{2}$ in $\hat{X}$ then $(w_{1},w_{2}) \cap X$ is a
nonempty, open, convex subset which is isomorphic to $X$.  The
isomorphism extends to the completions so that the interval
$(w_{1},w_{2})$ in $\hat{X}$  is isomorphic to $\hat{X}$.  Thus,
$\hat{X}$ is weakly homogeneous as well.

\end{proof}\vspace{.5cm}

We call $X$ a \emph{homogeneous LOTS}\index{LOTS!homogeneous}\index{HLOTS}\index{LOTS!HLOTS}, written HLOTS, if $X$ contains at least three points
and it is order isomorphic with every nonempty, open,
convex subset of itself.  If a HLOTS is complete then we call it a
CHLOTS, a \emph{complete homogeneous LOTS}\index{LOTS!complete homogeneous}\index{CHLOTS}\index{LOTS!CHLOTS}.  Otherwise, we call it
an IHLOTS, an \emph{incomplete homogeneous LOTS}\index{LOTS!incomplete homogeneous}\index{IHLOTS}\index{LOTS!IHLOTS}.

\begin{prop}\label{prop3.4} Let $X$ be a nontrivial LOTS.
\begin{enumerate}
\item[(a)] $X$ is a HLOTS if and only if $X$ is doubly transitive, of countable type and contains at least three points.
\item[(b)]  If $X$ is a HLOTS, then it is order dense, unbounded, $\sigma$-bounded and first countable.
\item[(c)]  The completion $\hat{X}$ of a HLOTS $X$ is a CHLOTS and a HLOTS  $X$
is a CHLOTS iff $X = \hat{X}$.  If $X$ is a CHLOTS, then it is connected, locally compact and $\sigma$-compact.
\item[(d)]  A HLOTS $X$ is an IHLOTS iff it is a proper subset of $\hat{X}$.  If $X$
is an IHLOTS then it has dense holes. Furthermore, $\hat{X} \setminus X$ is an IHLOTS
which is dense in $\hat{X}$, and so the completion of $\hat{X} \setminus X$ is $\hat{X}$.
\end{enumerate}
\end{prop}

\begin{proof} (a), (b):  A HLOTS is clearly
weakly homogeneous.  A finite LOTS is not isomorphic to its
singleton subsets and so a HLOTS is infinite.  It follows from
Proposition \ref{prop3.2a}(a) that a HLOTS is doubly transitive and so by
Proposition \ref{prop3.2}(e) it is order dense and unbounded.  Now
Proposition \ref{prop3.3}(b) (ii)$\Leftrightarrow$(iii) implies the equivalence
in (a).  Since a HLOTS is of countable type, it is first countable
and $\sigma$-bounded by Proposition \ref{prop2.5}(d).

(c), (d): $X$ is complete iff $X = \hat{X}$.  Completeness implies local compactness. If, in addition, it is
order dense, then it is connected and then
$\sigma$-boundedness implies $\sigma$-compactness.  In general, if
$X$ is a HLOTS, then by (a) and (b) it is order dense, doubly
transitive and of countable type.  By Proposition \ref{prop2.8}(e) $\hat{X}$
is of countable type and by Proposition \ref{prop3.3}(a) it is doubly
transitive.  By (a) above, $\hat{X}$ is a HLOTS and so is a
CHLOTS.

Now assume that $X$ is a proper subset of $\hat{X}$ with $z \in
\hat{X} \setminus X$.  Choose $a,b \in X$ so that $a < z < b$.
For any $c < d$ in $X$ there exists $f \in H_{+}(X)$ mapping the
pair $a,b$ to $c,d$.  The completion $\hat{f} \in H_{+}(\hat{X})$
maps $z$ to a point of $\hat{X} \setminus X$ between $c$ and $d$.
Thus, $X$ has dense holes.  That is, $\hat{X} \setminus X$ is
dense in $\hat{X}$ and so it is order dense with completion
$\hat{X}$.  Now if $\tilde{J}$ is a nonempty, open, convex subset
of $\hat{X} \setminus X$, then there exists an open, convex subset
$J$ of $\hat{X}$ such that $\tilde{J} = J \cap (\hat{X} \setminus
X)$.  There exists an isomorphism $f$ of the IHLOTS $X$ with the
nonempty, open, convex subset $J \cap X$ of itself. The completion
$ \hat{f} : \hat{X} \cong J $ restricts to an isomorphism $\hat{X}
\setminus X \cong \tilde{J}$.  Thus, $\hat{X} \setminus X$ is a
HLOTS, indeed an IHLOTS since $X$ is nonempty.

\end{proof}\vspace{.5cm}

\noindent {\bfseries Remark:} From Treybig's Theorem \ref{theotreybig} it follows
that a nontrivial, transitive, connected, $\s$-compact LOTS is
a CHLOTS. \vspace{.5cm}

The motivating example of an IHLOTS is the set of rationals
$\Q$ with completion the CHLOTS $\R$ and with complementary IHLOTS the irrationals $\I$ in $\R$. \vspace{.5cm}

 \begin{ex}\label{example3.4a} HLOTS ultraproduct \end{ex}

On $\omega$, the first infinite ordinal, we choose an ultrafilter $\mathcal{U}$ which
is nonprincipal,\index{ultrafilter}\index{ultrafilter!nonprincipal}
i.e. $\mathcal{U}$ contains all cofinite sets.
For a CHLOTS $F$ we define the ultraproduct construction on $F$
associated with  $\mathcal{U}$.   

On the, unordered, set of
maps $F^{\omega}$ we define the equivalence relation
$\equiv_{\mathcal{U}}$
\begin{equation}\label{eq3.36}
a  \equiv_{\mathcal{U}} b \qquad \Longleftrightarrow \qquad \{i :
a_{i} = b_{i}\} \in \mathcal{U},
\end{equation}
which is an equivalence relation because $\mathcal{U}$ is a
filter.  We let $F^{\mathcal{U}}$ denote the LOTS of equivalence
classes with the order
\begin{equation}\label{eq3.36a}
a < b  \qquad \Longleftrightarrow \qquad \{i : a_{i} < b_{i}\} \in
\mathcal{U}.
\end{equation}
This definition is independent of the choice of $a$ and $b$ in the
equivalence classes, but we can then choose representatives such
that $a_{i} < b_{i}$ for all $i \in \omega$.  The relation on
$F^{\mathcal{U}}$ is a total order because $\mathcal{U}$ is an
ultrafilter, i.e. if $D_{1} \cup ...\cup D_{n} \in \mathcal{U}$,
then $D_{k} \in \mathcal{U}$ for some $k = 1,...,n$.  In
particular, if $a_{i} < b_{i}$ for all $i$, then we can choose for
each $i$ an order isomorphism $f_{i} : F \cong (a_{i},b_{i})$.
The product $\prod_{i} \ f_{i}$ is an order
isomorphism $F^{\mathcal{U}} \cong (a,b) \subset F^{\mathcal{U}}$.
Thus, $F^{\mathcal{U}}$ is weakly homogeneous.  In the case $F =
\R$ this is the ordered field of hyperreal numbers.

$F^{\mathcal{U}}$ is  order dense but not complete. Its completion is not
even transitive.

We show that if $A = \{a^{1},a^{2},..\}, B = \{b^{1},b^{2},..\}$ are nonempty
countable subsets of $F^{\mathcal{U}}$ with $a < b$ for all $a \in A, b \in B$, then
there exists $w \in F^{\mathcal{U}}$ with $a < w < b$ for all $a \in A, b \in B$.

By replacing $a^{k}$ by $max_{\ell \leq k } a^{\ell}$ and eliminating any repeats we may assume that $\{ a^{k} \}$ is
a finite or infinite increasing sequence. Similarly, we may assume that $\{ b^{k} \}$ is a finite or infinite decreasing
sequence. Assume that representatives of $a^k$ and $b^k$ have been chosen so that $a^k_i < b^k_i$ for all $i \in \om$.
Now choose a representative $a^{k+1}$ so that $a^k_i < a^{k+1}_i < b^k_i$ for all $i$.  Then choose a representative
$b^{k+1}$ so that $a^{k+1}_i < b^{k+1}_i < b^k_i$ for all $i$. We do the obvious adjustments if either sequence is finite,
ending at the $k$ level.

We have thus obtained representatives of the elements of $A$ and $B$ so that $a_i < b_i$ for all $a \in A, b \in B, i \in \om$.
Choose $w_{i}$ so that
\begin{equation}\label{eq3.38}
a^{k}_{i} \ <  \ w_{i} \  <  \ b^k_{i} \qquad \mbox{for} \  k = 1,...,i.
\end{equation}

For each $k$, $\{i: a^{k}_{i} < w_{i} \}$ includes all $ i
\geq k $ and so is in $\mathcal{U}$. Hence, each $a^{k} < w$ for each $k$.
Similarly, $b^{k} > w$ for each $k$.

An obvious adjustment covers the case when $A$ or $B$ is empty.

Using the case when $A$ or $B$ is a singleton, we see that no increasing or decreasing
sequence converges in $F^{\mathcal{U}}$.

If $x$ in the completion $\widehat{F^{\mathcal{U}}}$ is the limit of an increasing
sequence $\{ a^k \}$, then we can choose the sequence
to be in $F^{\mathcal{U}}$ because the latter meets $(a^k,a^{k+1})$. So if $\{ b^k \}$
is any decreasing sequence with $b^k > x$ for all $k$,
then there exists $w \in F^{\mathcal{U}}$ with $a^k < w < b^k$ for all $k$. Since
$x \not\in F^{\mathcal{U}}$, $w \not= x$.
Since $\{a^k \}$ converges to $x$, it cannot be that $w < x$. Hence, $w > x$ and
so $\{b^k \}$ does not converge to $x$.

Thus, $\widehat{F^{\mathcal{U}}}$ is partitioned by three dense sets, namely
$\widehat{F^{\mathcal{U}}}_-$ consisting of the limits of increasing
sequences, $\widehat{F^{\mathcal{U}}}_+$ the limits of decreasing sequences and
the remaining set $\widehat{F^{\mathcal{U}}}_0$
which contains $F^{\mathcal{U}}$. Each automorphism of
$\widehat{F^{\mathcal{U}}}$ preserves each of these sets.

In general, Hausdorff defines an \emph{$\a$-set}\index{$\a$-set} $X$ to be
a linearly ordered set of cardinality $\aleph_{\a}$ such that
$A, B \subset X$ with $A < B$ and the cardinality of $A \cup B$ less than
$\aleph_{\a}$ implies there exists $x \in X$ such that $A < x < B$.
See \cite{Haus} and also section 4 or \cite{H}.
The isomorphism between any two $\a$-sets is proved using the transfinite
analogue of the back and forth argument which gives the uniqueness
up to isomorphism of countable, unbounded, order dense sets, see, e.g. Proposition \ref{prop2.8}(a).

The cardinality of  $F^{\mathcal{U}}$ is ${\mathbf c}$ which equals $\aleph_1$
if one assumes the Continuum Hypothesis. It then
follows that $F^{\mathcal{U}}$ is a $1$-set (assuming CH). Any open convex
subset of $F^{\mathcal{U}}$ which has neither a
countable cofinal nor a countable coinitial subset is also a $1$-set and
so is isomorphic to $F^{\mathcal{U}}$ by uniqueness.
Of course, we already know that if $a < b \in F^{\mathcal{U}}$, then the
interval is a isomorphic to $F^{\mathcal{U}}$ and so is a $1$-set,
but for $a < b \in \widehat{F^{\mathcal{U}}}_0$ the interval $(a,b)$ is a $1$-set as well.

Now suppose that $a < b \in \widehat{F^{\mathcal{U}}}_-$. Let $\{ a_n \}$ be
an increasing sequence in $F^{\mathcal{U}}$  which converges to $a$
and $\{ b_n \}$ be an increasing sequence in $(a,b) \cap F^{\mathcal{U}}$
which converges to $b$. Each of the open intervals
\begin{equation}\label{eq3.38aa}
(-\infty,a_1), \dots, (a_k,a_{k+1}), \dots, (a,b_1), \dots, (b_k,b_{k+1}), \dots, (b,\infty)
\end{equation}
is a $1$-set.

It easily follows that for $a < b, c < d \in \widehat{F^{\mathcal{U}}}_-$ we can complete
an automorphism of $F^{\mathcal{U}}$ to
obtain an automorphism of  $\widehat{F^{\mathcal{U}}}$ which maps $[a,b]$ onto $[c,d]$.

With a similar argument for $\widehat{F^{\mathcal{U}}}_+$, it follows that three sets
$\widehat{F^{\mathcal{U}}}_{\pm}, \widehat{F^{\mathcal{U}}}_0$
are the orbits of the action of $H_+(\widehat{F^{\mathcal{U}}})$ on $\widehat{F^{\mathcal{U}}}$.
Separately, each is a doubly transitive LOTS.

Thus, assuming CH, we see that the LOTS  $\widehat{F^{\mathcal{U}}}_-$ is a doubly transitive LOTS
with every increasing sequence convergent, but with no
convergent decreasing sequences.

\vspace{1cm}

\subsection{The Double Arrow of a LOTS}

For any LOTS $X$ we define the \emph{Alexandrov-Sorgenfrey Double
Arrow} \index{Alexandrov-Sorgenfrey Double
Arrow} of $X$, hereafter the \emph{AS double}\index{AS double}\index{LOTS!AS double} of $X$, to be the
order space product
\begin{equation}\label{eq3.5}
X' \quad =  \quad X \times \{ -1, +1 \},
\end{equation}\index{$X'$}
regarding $\{ -1, +1 \}$ as a two point LOTS.  For $x$ in $X$ we
denote by $x^{\pm}$ the pair $(x,\pm 1)$ and we define the first
coordinate projection map
\begin{equation}\label{eq3.6}
\begin{split}
\pi' : X' \rightarrow X  \\  \pi'(x^{\pm}) = x.
\end{split}
\end{equation}
Clearly, we have for $a < b$ in $X$:
\begin{equation}\label{eq3.7}
\begin{split}
(\pi')^{-1}((a,b)) = (a^{+},b^{-})  \\(\pi')^{-1}([a,b]) =
[a^{-},b^{+}].
\end{split}
\end{equation}
Since $(\pi')^{-1}(a)$ is the compact set $\{ a^{-}, a^{+} \}$,
Proposition \ref{prop2.1}(a) implies that the surjective order map $\pi'$ is closed and continuous.

For each $x \in X$, $x^{-} < x^{+}$ is a gap pair in $X'$ so that
$x^{-}$ is a left endpoint and $x^{+}$ is a right endpoint.
$x^{-}$ is also a right endpoint, and so is an isolated point of
$X'$, iff $x$ is a right endpoint in $X$, i.e. $[x,\infty)$ is
open.  For if $[x.\infty)$ is open then by (\ref{eq3.7}), $\{x^{-}\} =
(-\infty,x^{+}) \cap [x^{-},\infty)$ is open, while if $A \subset
X \setminus \{x\}$ with $x =  sup  A$ then $x^{-} =  sup   (\pi')^{-1}(A)$. Similarly, $x^{+}$ is isolated iff $x$ is a left
endpoint in $X$.  In particular, if $X$ has a $max = M$ (or a $min
= m$) then $M^{+}$ (resp. $m^{-}$) is an isolated point in $X'$.

If $f : X_{1} \rightarrow X_{2}$ is an order injection then we
define
\begin{equation}\label{eq3.8}
\begin{split}
f' : X_{1}' \rightarrow X_{2}'  \\  f'(x^{\pm}) = f(x)^{\pm}.
\end{split}
\end{equation}
Clearly, $f'$ is the unique order injection such that the diagram

\[  \begin{CD}
X_{1}'    @>f'>>    X_{2}' \\
\pi' @VVV         @VVV\pi'\\
X_{1}     @>f>>     X_{2}
\end{CD}  \]
commutes.   If $f$ is not
injective then the map defined by (\ref{eq3.8}) does not preserve order.
If $f$ is a continuous, noninjective, order map we can obtain many continuous
 order maps $f'$ which make the above diagram commute. If $X_1$ is complete, then
for each $y \in f(X_{1})$ the set $f^{-1}(y)$ is a closed interval $[x_{1},x_{2}]$. Choose a point $\tilde{x} \in
[x_{1},x_{2}] = f^{-1}(y)$ and then map
$[x^{-}_{1},\tilde{x}^{-}]$ to $y^{-}$ and
$[\tilde{x}^{+},x^{+}_{2}]$ to $y^{+}$. In general,$f^{-1}(y)$ is a closed, convex set and with the choice of
$\tilde{x} \in f^{-1}(y)$ we map $(-\infty,\tilde{x}^{-}] \cap (f \circ \pi')^{-1}(y)$ to $y^{-}$ and
$[\tilde{x}^{+},\infty) \cap (f \circ \pi')^{-1}(y)$ to $y^{+}$.

To check that such maps $f'$ are continuous, observe that
\begin{align}\label{eq3.8a}
\begin{split}(f')^{-1}((-\infty,y^-)) &= (f \circ \pi')^{-1}((-\infty,y)) \\
(f')^{-1}((-\infty,y^+)) &= (f')^{-1}((-\infty,y^-)) \ \ \text{if} \ y \not\in f(X), \\
(f')^{-1}((-\infty,y^+)) &= (-\infty,\tilde x^+) \ \ \text{if} \ y \in f(X).
\end{split}
\end{align}
Use a similar argument for the open intervals unbounded above.

If $r : X_{1} \rightarrow X_{2}$ is an injective order* map then
we define the continuous injective order* map
\begin{equation}\label{eq3.9}
\begin{split}
r^* : X_{1}' \rightarrow X_{2}'  \\  r^*(x^{\pm}) =
r(x)^{\mp}.
\end{split}
\end{equation}

If $f$ (or $r$) is bijective, then $f'$ (resp. $r^*$) is.

\begin{lem}\label{lem3.5} Let $X$ be an unbounded, order dense LOTS.
\begin{enumerate}
\item[(a)] Let $X_{1}$ be a LOTS and $g: X' \rightarrow X_{1}'$ be a homeomorphism.
The LOTS $X_{1}$ is order dense and unbounded.  If $g$ is an order map,
then there exists $f : X \rightarrow X_{1}$ an order isomorphism such that $g = f'$.
If $g$ is an order$^*$ map, then there exists $r : X \rightarrow X_{1}$ an order$^*$
homeomorphism such that $g = r^{*}$.

In particular, we have
\begin{equation}\label{eq3.09a}
X' \ \cong \ X_1' \quad \Longleftrightarrow \quad X \ \cong \ X_1.
\end{equation}

\item[(b)]  Assume that $X$ is complete and so is connected.  The bounded clopen
subintervals of $X'$ are compact sets which form a base for $X'$, and so $X'$ is
zero-dimensional.  If $C$ is any nonempty, bounded, clopen subset of $X'$, then
there is a unique finite sequence $a_{1} < b_{1} < a_{2} ... < b_{n}$ in $X$ such that
\begin{equation}\label{eq3.10}
C = \bigcup _{i=1}^{n}[a_{i}^{+},b_{i}^{-}].
\end{equation}
\end{enumerate}
\end{lem}

\begin{proof} (a): $X'$ has no isolated points
and so $X_{1}'$ has none.  Hence, $X_{1}$ has no left or right
endpoints.  If $g$ is either order preserving or order reversing
then gap pairs are mapped to gap pairs.  Hence $g$ induces a
bijection from $X$ to $X_{1}$ which preserves or reverses order
according to which $g$ does.

(b):  Since $X$ is connected with no $max$ or $min$, any bounded
clopen interval of $X'$ is of the form $[a^{+},b^{-}]$ with $a <
b$ in $X$.  These form a base.

Let $C$ be a nonempty, bounded, clopen subset of $X'$.  Since $C$ is open it is
a union of such subintervals and since $C$ compact it is a finite
union of them. If two such intervals intersect, or if the $max$ of
one and the $min$ of another form a gap pair, then the union of
the two is a clopen interval.  Combining in this way we obtain $C$
as the finite, disjoint union of clopen intervals with points of
$X' \setminus C $ between any two successive intervals, i.e.
(\ref{eq3.10}) holds. Furthermore, the intervals $[a_{i},b_{i}]$ are the
components of the image $\pi'(C)$ in $X$.  Uniqueness follows.

\end{proof}\vspace{.5cm}

We now describe the gap between $\pm$transitivity and transitivity
in the complete case.

\begin{prop}\label{prop3.6} A  complete LOTS  $X$ is $\pm$transitive but not transitive
iff $X \cong X_{0}'$ with $X_{0}$ a connected, symmetric, transitive LOTS.
\end{prop}

\begin{proof} If $X$ were trivial it would be
transitive.  If $X$ is the two point LOTS, then $X \cong X_{0}'$
with $X_{0}$ trivial.  So we can assume that $X$ has at least
three points so that Proposition \ref{prop3.2} applies.

Assume that $X_{0}$ is a nontrivial, order dense, symmetric,
transitive LOTS.  If $f \in H_{+}(X_{0})$ maps $x$ to $y$, then
$f'$ in $H_{+}(X_{0}')$ maps $x^{+}$ to $y^{+}$ and  $x^{-}$ to
$y^{-}$.  Thus, the right endpoints $x^{+}$ and the left endpoints
$x^{-}$ in $X_{0}'$ form $H_{+}(X_{0}')$ classes which are
distinct because $X_{0}'$ has no isolated points.  Thus, $X_{0}'$
is not transitive.  If $r$ is an order reversing homeomorphism on
$X_{0}$, then $r^{*}$ maps $x^{+}$ to $r(x)^{-}$ and so $X_{0}'$
consists of a single $H_{\pm}(X_{0}')$ equivalence class.  If, in
addition, $X_{0}$ is connected, then $X_{0}'$ is complete.

Conversely, if $X$ is $\pm$transitive, complete, not transitive
and contains at least three points, then by Proposition \ref{prop3.2}(b) it is
unbounded.  If there were any isolated points, then by transitivity $X$ would be discrete and we could apply
the argument in the proof of Proposition \ref{prop3.2}(c) to show that $X
\cong \Z$ which is transitive.  If $X$ were connected,
then by Proposition \ref{prop3.2}(a) it would be transitive.  Hence, $X$
contains gap pairs but no isolated points.  It follows from
$\pm$transitivity that every point is either a left or right
endpoint and that no point is both.  Since $X$ is unbounded, all of
the endpoints occur in gap pairs.  Let $X_{0}$ be the LOTS of gap
pairs.  Every element of $H_{\pm}(X)$ induces an element of
$H_{\pm}(X_{0})$ and so $X_{0}$ is $\pm$transitive.  Completeness
of $X_{0}$ follows from completeness of $X$.  For each pair $z$ in
$X_{0}$ let $z^{+}$ be the right and $z^{-}$ be the left endpoint
of the pair.  This yields an isomorphism $X \cong X_{0}'$.  Since
$X$ has no isolated points, $X_{0}$ is order dense and so is
connected.  Hence, $X_{0}$ is transitive by Proposition \ref{prop3.2}(a).
Finally, since X is $\pm$transitive but not transitive, it admits
some order reversing homeomorphism $g$.  By Lemma \ref{lem3.5}(a) $g = r^{*}$
for $r $ an order reversing homeomorphism of $X_{0}$.  Hence,
$X_{0}$ is symmetric.

\end{proof}\vspace{.5cm}

\noindent {\bfseries Remark.} Notice that for the case of $\Z$ we can define the order isomorphism:
\begin{equation}\label{eq3.11}
\begin{split}
f : \Z'  \rightarrow \Z \hspace{1.2in}  \\
f(n^{-}) = 2n \qquad \mbox{and} \qquad f(n^{+}) = 2n+1.
\end{split}
\end{equation}
\vspace{.5cm}

For $X$ a complete LOTS with no $max$ or $min$ the
\emph{two-point compactification}\index{two-point compactification},
denoted $\bullet X \bullet$\index{$\bullet X \bullet$}, is
the LOTS obtained by attaching a $max$ and a $min$, i.e. the order
sum
\begin{equation}\label{eq3.12}
\bullet X \bullet \quad =  \quad \{m\} + X  + \{M\}.
\end{equation}
\vspace{.5cm}

\begin{prop}\label{prop3.7}  If $X$ is a CHLOTS and $C$ is a nonempty, bounded, clopen subset of $X'$ then
\begin{equation}\label{eq3.13}
C \quad \cong \quad \bullet X'\bullet \quad.
\end{equation}
\end{prop}

\begin{proof}  Let  $a_{1} < b_{1} < a_{2} ... <
b_{n}$ be the sequence in $X$ such that (\ref{eq3.10}) holds.  Let $f_{i}
: [a_{i},b_{i}] \cong [a_{i},a_{i+1}]$ for $i = 1,...,n-1$.  Each
$f_{i}'$ restricts to an isomorphism $[a^{+}_{i},b^{-}_{i}] \cong
[a^{+}_{i},a^{-}_{i+1}]$. Put them together to get an isomorphism:
\begin{equation}\label{eq3.14}
C \quad \cong \quad (\bigcup^{n-1}_{i}[a^{+}_{i},a^{-}_{i+1}])
\cup [a^{+}_{n},b^{-}_{n}] \quad = \quad [a^{+}_{1},b^{-}_{n}].
\end{equation}
If $g : X \rightarrow (a_{1},b_{n})$ is an order isomorphism, then
so is
\begin{equation}\label{eq3.15}
g' : X' \rightarrow (a_{1},b_{n})' = (a_{1}^{+},b_{n}^{-}).
\end{equation}
Attaching the endpoints we obtain an isomorphism $\bullet X'\bullet \cong [a_{1}^{+},b_{n}^{-}]$.

\end{proof}\vspace{.5cm}

If $X$ has no isolated points, then the isolated points of $X'$ are
of the form $x^{+}$ where $x$ is a left endpoint of $X$ and
$x^{-}$ where $x$ is a right endpoint.  We define for a LOTS $X''$
with no isolated points:
\begin{equation}\label{eq3.16}
X''  =   X' \setminus (\{ x^{+}: x\  \mbox{ a left endpoint}\}
\cup \{x^{-} : x \ \mbox{ a right endpoint}\}).
\end{equation}\index{$X''$}
That is, $X''$ is the subset of $X'$ obtained by removing the isolated points.
Note that if $X$ is unbounded and order dense then $X'' = X'$.

In general, if $x_{1} < x_{2}$ is a gap pair in $X$, then since
$x_{1}$ is not a right endpoint and $x_{2}$ is not a left endpoint
(neither is isolated) we see that $x_{1}^{-} < x_{2}^{+}$ is a gap
pair in $X''$. If $M =  max  X$, then $M^{-}=  max  X''$ and
similarly, $m^{+} =  min  X''$ if $m =  min  X$.  If $x$ is not an
endpoint in $X$, then $x^{-} < x^{+}$ is a gap pair in $X''$.  Thus,
every point of $X''$ is a right or a left endpoint. Let
\begin{equation}\label{eq3.17}
\pi'' : X'' \rightarrow X
\end{equation}
denote the restriction of the projection $\pi' : X' \rightarrow
X$.  For every $x \in X$ either $x^{+}$ or $x^{-}$ or both lies in
$X''$ and so $\pi''$ is an order surjection.  Since $(\pi'')^{-1}(x)$
is either a singleton or a pair for every $x \in X$ it follows
from Proposition \ref{prop2.1}(a) that $\pi''$ is continuous and closed.

If $x =  sup  A$ with $A \subset (-\infty,x)$ in $X$, then $x$ is
not a right endpoint and so $x^{-} \in X''$.  Furthermore,
\begin{equation}\label{eq3.18}
x^{-} \quad = \quad  sup  \pi''^{-1}(A).
\end{equation}
Similarly, $x =  inf  B$ with $B \subset (x,\infty) \subset X$
implies
\begin{equation}\label{eq3.19}
x^{+} \quad = \quad  inf  \pi''^{-1}(B).
\end{equation}

In particular, $X''$ has no isolated points.

Now for a useful construction:

Assume $X$ is a complete LOTS with no isolated points, $X_{1}$ is
an unbounded LOTS with no isolated points, $A$ is a dense subset
of $X_{1}$ and $g : A \rightarrow X$ is an order injection.

First, define
\begin{equation}\label{eq3.20}
\begin{split}
G : X'_{1} \rightarrow X, \qquad \mbox{by} \hspace{1.5cm} \\
G(x^{-}) =  sup  g((-\infty,x) \cap A)   \\
G(x^{+}) =  inf  g((x,\infty) \cap A).
\end{split}
\end{equation}
So that $ G(x^{-}) \le G(x^{+}) $.

Assume $x_1 < x_2$ in $X_1$.

If $x_1, x_2$ is not a gap pair, then the
interval $(x_{1},x_{2})$ in $X_{1}$ is nonempty and so is
infinite because $X_1$ has no finite nonempty open set. Because $A$ is dense, $(x_{1},x_{2}) \cap A$ is infinite.
Thus, there exist $a_{1}, a_{2} \in A$ such that $x_{1} < a_{1} < a_{2} < x_{2}$. It follows that
\begin{equation}\label{eq3.25}
G(x_1^-) \le G(x_1^+) \le g(a_1) < g(a_2) \le G(x_2^-) \le G(x_2^+).
\end{equation}

If $x_1, x_2$ is a gap pair, then at least
\begin{equation}\label{eq3.23}
G(x_1^-) \quad \le \quad G(x_2^+).  \hspace{2cm}
\end{equation}
It can happen that in the gap pair case, $G(x_2^-) < G(x_1^+)$. On the other hand,
 if $x_1, x_2$ is a gap pair, then  $x_1^+, x_2^- \not\in X_1''$.

 It follows that the restriction of $G$ to $X_1''$ is an order map.

Now define the lift $g''$.
\begin{equation}\label{eq3.21}
\begin{split}
g'' : X_{1}'' \rightarrow X'', \qquad \mbox{by} \hspace{1.5cm} \\
g''(x^{-}) = sup \ (\pi'')^{-1}(g((-\infty,x) \cap A))   \\
g''(x^{+}) = inf \ (\pi'')^{-1}(g((x,\infty) \cap A)).
\end{split}
\end{equation}
Notice that if $x^{-} \in X_{1}''$, then $x$ is not a right endpoint
and so $x = sup \ ((-\infty,x)\cap A) $ and similarly for $x^+$.

It follows from (\ref{eq3.18})
and (\ref{eq3.19}) that on $X_{1}''  \subset X'_{1} $ we have
\begin{equation}\label{eq3.22}
\begin{split}
g''(x^{-}) \quad = \quad G(x^{-})^{-} \hspace{2cm} \\
g''(x^{+}) \quad = \quad G(x^{+})^{+}. \hspace{2cm} \\
\mbox{and so} \qquad \pi'' \circ g'' \quad = \quad G. \hspace{1.5cm}
\end{split}
\end{equation}

Assume that $x_1 < x_2$ in $X_1$.

If $x_{1}^{-},x_{2}^{+} \in X_{1}''$, then by (\ref{eq3.22}) together with (\ref{eq3.25}) or (\ref{eq3.23})
\begin{equation}\label{eq3.24}
g''(x_{1}^{-}) \quad < \quad g''(x_{2}^{+}). \hspace{2cm}
\end{equation}

If $x_{1}^{+}$ or $x_{2}^{-} \in X_{1}''$, then $x_1, x_2$ is not a gap pair
So from (\ref{eq3.25}) and (\ref{eq3.22}) again, we have
\begin{equation}\label{eq3.26}
g''(x_1^-) \ < \ g''(x_{1}^{+}) \  < \  g''(x_{2}^{-}) \ < \  g''(x_{2}^{+})
\end{equation}
omitting whichever terms are undefined.

It follows that $g''$ is an order injection.

Now assume that $x \in X_1$ is not a right endpoint so that $x^- \in X_1''$. By definition, $G(x^-)$ is not a right endpoint
and so $g''(x^-) = G(x^-)^- = sup \{ y^- : y < G(x^-) \}$. For $y < G(x^-)$ the interval $(y,G(x^-))$ is infinite and we need only
consider such points $y$ which are also not right endpoints. There exists $a \in (-\infty,x) \cap A$ such that $y < g(a)$. There exists
$x_1 \in (a,x) \subset X_1$ not a right endpoint so that $x_1^- \in X_1''$ as well.
Since $x_1, x$ is not a gap pair, (\ref{eq3.25}) implies that
 \begin{equation}\label{eq3.27}
 \begin{split}
 y < g(a) \le G(x_1^{-}) < G(x^-), \quad \text{and so} \hspace{2cm}\\
 y^- < G(x_1^{-})^-  < G(x^{-})^- , \ \text{i.e.} \  y^- < g''(x_1^{-})  < g''(x^{-}),
 \end{split}
 \end{equation}
and so
\begin{equation}\label{eq3.28}
g''(x^{-}) \quad = \quad sup \ \{g''(x_1^-) : x_1^- \in X_1'', \ \text{and} \ x_1 < x \}.
\end{equation}

From this together with a similar argument for $g''(x^+)$ when $x$ is not a left endpoint, it follows that $g''$ is continuous and so is
an order embedding. Since $\pi''$ is continuous, it follows from (\ref{eq3.22}) that $G : X_1'' \to X$ is a continuous order map.

In particular, we obtain:

\begin{lem}\label{lem3.7a} If $X$ and $X_1$ are connected LOTS, $A \subset X$ is dense and $g : A \to X_1$ is an order injection,
then $g'' : X' \to X_1'$ is an order embedding. \end{lem}
\vspace{.5cm}

Now we apply this construction.

\begin{theo}\label{theo3.8} (a) If $X$ is a connected, nontrivial LOTS, then $X$ contains a closed subset
$A$ which is a compact, separable LOTS with no isolated points.
\begin{enumerate}
\item[(b)] If $C$ is a separable, compact LOTS with no isolated points, then $C''$ has
the order type of the \emph{Fat Cantor Set}\index{Fat Cantor Set}, i.e.
\begin{equation}\label{eq3.30}
C'' \quad \cong \quad \bullet \R' \bullet . \hspace{1.5cm}
\end{equation}
\end{enumerate}
\end{theo}

\begin{proof}  (a): Let $a < b$ in $X$.  By
Proposition \ref{prop2.8}(a) there exists an order injection $ g: {\Q}
\rightarrow (a,b)$ where $\Q$ is the LOTS of rationals in
$\R$. By replacing $a,b$ if necessary we can assume that
$g(\Q)$ is $\pm$cofinal in $(a,b)$.  Since
$\R$ and $(a,b)$ are order dense and unbounded, $\R' = \R''$ and $(a,b)' = (a,b)''$. The above construction
yields an order embedding
\begin{equation}\label{eq3.31}
g'' : \R' \rightarrow (a,b)'.
\end{equation}
Extend to $\bullet\R'\bullet$ by mapping $m$ to $a^{+}$
and $M$ to $b^{-}$.  Because the image of $g$ is $\pm$cofinal
 we obtain an order embedding of $\bullet \R'
\bullet$ into $[a^{+},b^{-}] \subset X'$.  $A =
(\pi'\circ g'')(\bullet \R' \bullet)$ is compact and
separable by continuity.  Since the Fat Cantor Set, $ \bullet
\R' \bullet$, has no isolated points and $\pi'$ has
finite point inverses, $A$ has no isolated points.

(b): If $C$ is compact with no isolated points, let
\begin{equation}\label{eq3.32}
 L \quad = \quad C \setminus (\{ \mbox{right endpoints of} \   C \} \cup \{ max \ C\}).
\end{equation}
Regarded as a LOTS $L$ is unbounded and it is order dense because
if $x_{1} < x_{2}$ in $L$, then $x_{2}$ is not a right endpoint and
so $(x_{1},x_{2})$ is infinite in $C$.  If $y_{1} < y_{2}$ are
both right endpoints in this interval, then between them there is a
left endpoint.  Since $C$ has no isolated points it follows that
$(x_{1},x_{2})$ meets $L$.  Similarly, if $y$ is a right endpoint of $C$ and $y < x$ in $C$, then
the infinite interval $(y,x)$ meets $L$.  Finally, if $x \in C$ is not the maximum, then the infinite
interval $(x, max)$ meets $L$. It follows that $L$ is dense in $C$.

If $C$ is separable, then by Proposition \ref{prop2.5}(f) $L$ is separable and
so contains a countable dense set $D$.  So $D$ is countable, order
dense and unbounded.  By Proposition \ref{prop2.8}(a) there
exists an order isomorphism $g : \Q \rightarrow D \subset
C$. Regarded as a not necessarily continuous order injection into
$C$, $g$ induces, as above, $g'' : \R' \rightarrow C''$ an order
embedding.  Extend by mapping $m$ to $(min\ C)^{+}$ and $M$ to
$(max\ C)^{-}$.  We thus have an order embedding $g'' : \bullet
\R'\bullet \rightarrow C''$.  By continuity the image is
compact and hence closed in $C''$.  If $x^{-} \in C''$, then $x$ is
the limit of an increasing sequence in $D$.  If $x \not= max \ C$,
then the sequence is bounded in $D$.  Let $t \in \R$ be
the limit of the corresponding increasing sequence in $\Q$.
Clearly, $x^{-} = g''(t^{-})$.  With a similar argument for
$x^{+}$ in $C'' $ with $x \not= min \ A$, we see that $g'' $ is
onto.  It is thus an order isomorphism.

\end{proof} \vspace{.5cm}

\begin{ex}\label{ex3.8a} Embedding the Fat Cantor Set. \end{ex}

Suppose that $g : \Q \to X$ is an order injection with image a bounded, countable discrete subset. For example, if $X = \R \times [-1,1]$ we
can let $g(a) = (a,0)$ for $a \in \Q \cap (0,1) \subset \R$. Then for each $t \in (0,1)$, $G(t^-) < G(t^+)$ and so
$G : (0,1)' = (0,1)'' \to X$ is injective by (\ref{eq3.25}).
Since $g'' : \R' \to X'' = X'$ is an order embedding, $G = \pi'' \circ g''$ is continuous and so is an order embedding of $(0,1)'$. Extending to
the two-point compactification we obtain an embedding of the Fat Cantor Set into $X$. In the case of $X = \R \times [-1,1]$ with
 $g(a) = (a,0)$ the embedding  is given by $G(t^{\pm}) = (t, \pm 1)$ for $t \in [0,1]$ (omitting $0^-$ and $1^+$).
\vspace{1cm}

\subsection{The CHLOTS Long Line}

In a CHLOTS $F$ we label two points $-1 < +1$.  With $\Omega$ the
first uncountable ordinal we use the order product to define
\begin{equation}\label{eq3.33}
F^{\Omega \infty} \quad = \quad (\Omega \times [-1,+1)) \setminus
\{(0,-1)\}
\end{equation}\index{$F^{\Omega \infty}$}
In the case when $F = \R$ this is the \emph{Long Line}.

If $\alpha < \Omega$, then $\alpha $ is a countable ordinal and so
by Proposition \ref{prop2.8}(d) there exists an order embedding $f_{\alpha} :
\alpha + 1 \rightarrow [-1,+1]$ with $f_{\alpha}(0) = -1$ and
$f_{\alpha}(\alpha) = +1$.  Choose for each $i \in \alpha$ an
order isomorphism $\{i\} \times [-1,+1) \cong
[f_{\alpha}(i),f_{\alpha}(i+1))$ and put them together to get an
order isomorphism between the intervals $(-\infty,(\alpha,-1))
\subset F^{\Omega \infty}$ and $(-1.+1) \subset F$ which is
isomorphic to $F$ itself.  It follows that every nonempty, open
interval which is bounded above in $F^{\Omega \infty}$ is
isomorphic to $F$.  By Proposition \ref{prop3.2}(d)(vii)$\Rightarrow$(i)
$F^{\Omega \infty}$ is doubly transitive.  It is connected and
first countable, but not $\sigma$-compact.  It has countable
coinitial subsets but every countable subset is bounded above.

If we apply the above construction to $F^{*}$ and then reverse it,
then we obtain $^{-\infty \Omega}F$ \index{$^{-\infty \Omega}F$}which has countable cofinal
sets but no countable coinitial sets.  Using the order sum we
define
\begin{equation}\label{eq3.34}
^{-\infty \Omega}F^{\Omega \infty } \quad = \quad ^{-\infty
\Omega}F + \{0\} + F^{\Omega \infty }
\end{equation}\index{$^{-\infty \Omega}F^{\Omega \infty}$}
with every countable subset bounded.  In each case all of the
nonempty, bounded, open subintervals are isomorphic to $F$ and the
entire space is doubly transitive but not weakly homogeneous.

We can extend Proposition \ref{prop3.2a} to completely describe the gap
between double transitivity and homogeneity in the complete case.

\begin{theo}\label{theo3.9}  Let $X$ be a complete, doubly transitive LOTS and let $F$ be a nonempty,
bounded, open subinterval of $X$.  $F$ is a CHLOTS and $X$ is order isomorphic to exactly
one of the four spaces: $ F, F^{\Omega \infty }$, $^{ -\infty \Omega}F,$ or $^{-\infty \Omega}F^{\Omega \infty }$.
\end{theo}

\begin{proof} For a complete LOTS homogeneity is
equivalent to weak homogeneity.  Since $X$ is complete, $F$ is
complete and so by Proposition \ref{prop3.2a}(b) $F$ is a CHLOTS.  Also if $X$
has countable coinitial and cofinal subsets, then it is
$\sigma$-bounded and so is homogeneous by Proposition \ref{prop3.2a}(b) again.
In that case, $X \cong F$. Similarly, if $X$ has countable
coinitial subsets, then $X \cong (-\infty,a)$ for every $a \in X$.

Now assume that $X$ has countable coinitial subsets but that every
countable subset is bounded above.  Inductively, we can construct
an order embedding $f : \Omega \rightarrow X$.  Since $X$ is first
countable the image of $\Omega$ cannot be bounded in $X$ and so it
is cofinal in $X$. By Proposition \ref{prop2.8}(c) continuity of $f$ implies
that
\begin{equation}\label{eq3.35}
[f(1),\infty) = \bigcup \{ [f(i),f(i+1)) : 0 < i < \Omega \}.
\end{equation}
With points $-1 < +1$ fixed in $F$, choose for each $i$ an
isomorphism $\{i\} \times [-1,+1) \cong [f(i),f(i+1))$ and put
them together to get an isomorphism from $[(1,-1),\infty) \subset
F^{\Omega \infty}$ onto $[f(1),\infty) \subset X$. Put this
together with an isomorphism $\{0\} \times (-1,+1) \cong
(-\infty,f(1))$ to get an isomorphism $F^{\Omega \infty} \cong X$.

Applying this result to $ X^{*}$ we see that if $X$ has countable
cofinal subsets but no countable coinitial subsets, then $^{-\infty
\Omega}F \cong X$.  Finally, if every countable subset is bounded
in $X$, then we can pick $a \in X$ and apply the previous results
to show $(a,\infty) \cong F^{\Omega \infty}$ and $(-\infty,a)
\cong$ $  ^{-\infty \Omega}F$.  So directly from (\ref{eq3.34}) we see
that $X \cong $ $ ^{-\infty \Omega}F^{\Omega \infty}$ in this
case.

\end{proof}\vspace{.5cm}

\begin{ex} A first countable, doubly transitive LOTS whose completion is not first countable and so is not transitive.\end{ex}

With $X$ one of the long line examples above and $\a$ a countable tail-like ordinal with $\a > 1$,
 $Z = X^{\a}$ is first countable by
Proposition \ref{prop2.7a} and
 is doubly transitive by Theorem \ref{theobabcock}.

Consider $X = F^{\Omega \infty }$. There is an embedding $g : \Omega \to X$. Now choose $x \in X$ and let $\tilde g : \Omega \to Z$
be defined for $\b < \Omega$ by $\tilde g(\b)_1 = g(\b)$ and $ \tilde g(\b)_i = x$ for $i \not= 1$.
In particular,  $ \tilde g(\b)_0 = x$ for all $\b < \Omega$
and so the image of $\tilde g$ is bounded in $Z$. Let $\hat x$ be the supremum of the image in the completion $\hat Z$. Clearly,
$\hat x$ is not the limit of any countable increasing sequence.
\vspace{1cm}

\section{\textbf{Towers of CHLOTS}}
\vspace{.5cm}

\subsection{The LOTS $X_{\a}$}

In an unbounded LOTS $X$ we  pick out a distinguished
closed, bounded subinterval $J$\index{$J$} containing at least three points.  We label the endpoints $\pm
1$.  Thus, $J = [-1,+1]$ and its interior $J^{\circ} = (-1,+1)$ is nonempty.

For every positive ordinal $\alpha$ we define the subset of the
order space product $X^{\alpha}$ \index{$X^{\alpha}$}
\begin{equation}\label{eq4.1}
X_{\alpha} = \{ x \in X^{\alpha} : x_{i} \in J \ \mbox{for all}\
0 < i < \alpha \}.
\end{equation}
Thus, $X_{\alpha}$\index{$ X_{\a}$} is the order space product indexed by $\alpha$
with the first factor $X$ and the remaining factors copies of $J$.
In particular,
\begin{equation}\label{eq4.2}
\begin{split}
X_{1} \qquad = \qquad X \qquad \\
X_{2} \qquad = \qquad X \times J \quad \\
X_{\omega} \quad = \qquad X \times J \times J \times ...
\end{split}
\end{equation}

When $X$ is weakly homogeneous, e.g. if $X$ is a HLOTS, then it is order isomorphic to $J^{\circ}$ and so
the space $X_{\a}$ is independent of the choice of interval $J$.

For $0 < \beta \leq \alpha$ we have projections
$\pi_{\beta}^{\alpha} : X_{\alpha} \rightarrow X_{\beta}$.\index{$\pi_{\beta}^{\alpha} $}
Identifying $X_{1}$ with $X$ we have the special case
$\pi^{\alpha}: X_{\alpha} \rightarrow X$,\index{$\pi^{\alpha}$} the projection to the
first coordinate.

Following (\ref{eq2.14}) we define for $z \in X_{\beta}$ with $0 < \beta <
\alpha$ the points $z+$ and $z-$\index{$z-$}\index{$z+$} in $X_{\alpha}$ by
\begin{equation}\label{eq4.3}
(z\pm)_{i}  = \begin{cases} z_{i} \qquad i < \beta \\  \pm 1
\qquad \beta \leq i < \alpha.  \end{cases}
\end{equation}

As in (\ref{eq2.15}) we have for $a < b$ in $X_{\beta}$
\begin{equation}\label{eq4.4}
\begin{split}
 (\pi_{\beta}^{\alpha})^{-1}((a,b)) = (a+,b-) \\
 (\pi_{\beta}^{\alpha})^{-1}([a,b]) = [a-,b+].
\end{split}
\end{equation}

Using (\ref{eq4.3}) we define for $0 < \beta < \alpha$
\begin{equation}\label{eq4.5}
\begin{split}
j_{\alpha}^{\beta} : (X_{\beta})' \rightarrow X_{\alpha} \\
j_{\alpha}^{\beta}(z^{\pm}) = z\pm.
\end{split}
\end{equation}\index{$j^{\alpha}_{\beta}$ }

\begin{prop}\label{prop4.1}  Let $X$ be an unbounded, $\s$-bounded, order dense  LOTS with distinguished subinterval
$J = [-1,+1]$ and let $\alpha$ be a positive ordinal.
\begin{enumerate}
\item[(a)]  $X_{\alpha}$ is unbounded, order dense and $\sigma$-bounded.
\item[(b)]  If $X$ is an IHLOTS, then $X_{\alpha}$ has dense holes.
\item[(c)]  If $X$ is connected, then $X_{\alpha}$ is connected and $\sigma$-compact.
\item[(d)]  If $X$ has countable type and $\alpha$ is countable, then  $X_{\alpha}$ is
of countable type. If $\alpha$ is uncountable then $X_{\alpha}$ is not even first countable.
\item[(e)]  If $ 0 < \beta < \alpha$,  then $\pi_{\b}^{\a}$ is a continuous order
surjection, closed when $X$ is connected, and $j_{\a}^{\b}$ is an order embedding onto a closed subset of $X_{\alpha}$.
\item[(f)]  If $X$ is a HLOTS and $a < b$ in $X$,  then $(a+,b-)  \subset X_{\alpha}$ is
order isomorphic with $X_{\alpha}$ itself.
\end{enumerate}
\end{prop}

\begin{proof}  (a):  $X_{\alpha}$ is order dense
by Proposition \ref{prop2.3}(a) which also shows that each
$\pi_{\beta}^{\alpha}$ is a continuous order surjection.  Since
$X$ is unbounded and $\sigma$-bounded, $X_{\alpha}$ is unbounded
and $\sigma$-bounded by Proposition \ref{prop2.1}(a) applied to the order
surjection $\pi^{\alpha} : X_{\alpha} \rightarrow X$.

(b):  If $z < w$ in $X_{\alpha}$,  let $\beta = min \{j : z_{j}
\not= w_{j} \} $ so that $z_{\beta} < w_{\beta} $ and $ z_{i} =
w_{i} $ for $ i < \beta$.  By Proposition \ref{prop3.4}(d) $X$ has dense holes and so
 there exists a clopen subset $A$ of $X$ such that
$z_{\beta}  \in A, w_{\beta}  \not\in A $, and $x \in A \
\Rightarrow \ (-\infty,x) \subset A$.  Define
\begin{equation}\label{eq4.6}
\tilde{A} = \{ x \in X_{\alpha} : x < z \} \cup \{ x \in
X_{\alpha} : x_{i} = z_{i} \ \mbox{for} \ i < \beta \ \mbox{and}\
x_{\beta} \in A \}.
\end{equation}
It is clear that $\tilde{A}$ defines a hole in $X_{\alpha}$
between $z$ and $w$.

(c):  If $X$ is complete, then $J$ is compact and so $X_{\alpha}$ is
complete by Proposition \ref{prop2.3}(b).  As $X_{\alpha}$ is order dense and
$\sigma$-bounded, it is connected and $\sigma$-compact.

(d):  Since $X$ and $J$ are of countable type, the first result
follows from Proposition \ref{prop2.7}(b). If $\alpha$ is uncountable,  then
define $f : \alpha + 1 \rightarrow X_{\alpha} $ as in (\ref{eq2.24}) with
$a = (-1)-$ and $b = (+1)+$.  Then $f$ is an order embedding and
since $\Omega \leq \alpha$ the point $f(\Omega)$ is defined in
$X_{\alpha}$.  It is not the limit of any increasing sequence.

(e):  For $z \in X_{\b}$ the pre-image $(\pi^{\alpha}_{\beta})^{-1}(z) = [z-,z+]$ which is closed and is compact
when $X$, and therefore $X_{\a}$, are connected. So $\pi^{\alpha}_{\beta}$ is continuous, and is closed in the connected case, by
Proposition \ref{prop2.1}(a).

The complement of the image of $j^{\b}_{\a}$ is the
union of the collection of open intervals $ \{ (z-,z+) : z \in
X_{\beta} \} $.  In the case when $X$ is complete, continuity
follows from Proposition \ref{prop2.1}(b).  In the general case, one checks
directly that condition (\ref{eq2.2}) holds for the image, i.e. if $ a \in
j^{\b}_{\a}(X_{\beta})$ and $x \in X_{\alpha}$ with $a < x
$ and $ (a,x] \cap j^{\b}_{\a}(X_{\beta} )= \emptyset$ then
$a = z- $ and $z_{i} = x_{i} $ for all $i \in \beta $.  The
required $b$ is $z+$.

(f):  If $ \tilde{f} : X \rightarrow (a,b) $ is an order
isomorphism, then
\begin{equation}\label{eq4.7}
f(x)_{i} = \begin{cases} \tilde{f}(x_{0})  \qquad  i = 0 \\  x_{i}
\qquad  0 < i < \alpha   \end{cases}
\end{equation}
defines the required isomorphism  $f : X_{\alpha} \rightarrow
(a+,b-) $.

\end{proof} \vspace{.5cm}

\begin{theo}\label{theo4.2}  Let $\alpha$ be a countable, tail-like ordinal and $X$ be a LOTS.
\begin{itemize}
\item[(i)]  If $X$ is doubly transitive and first countable, then  $X_{\alpha}$ is
doubly transitive and first countable.

\item[(ii)]  If $X$ is an IHLOTS, then
$X_{\alpha}$ is an IHLOTS.

\item[(iii)]  If $X$ is a CHLOTS, then $X_{\alpha}$ is a CHLOTS.
\end{itemize}
\end{theo}

\begin{proof} In any case $X_{\a}$ is first countable by Proposition \ref{prop2.7a}. By Proposition \ref{prop4.1} $X_{\alpha}$
has dense holes if $X$ is an IHLOTS and it is complete if $X$ is a
CHLOTS.  Furthermore, in the HLOTS cases $X_{\alpha}$
has countable type.  We will show that if $\alpha$ is countable
and tail-like and $X$ is doubly transitive, then  $X_{\alpha}$ is
doubly transitive.  The HLOTS results then follow from Proposition \ref{prop3.4}.

If $\alpha = 1 $ then $X_{\alpha} = X$ which is doubly transitive.
Thus, we can assume that $\alpha > 1$ and so that $\a$ is a limit
ordinal.

Choose $ 0 \in J^{\circ}$ so that $ -1 < 0 < +1$.  Define $\tilde{a}
< \tilde{c} < \tilde{b}$ in $ X_{\alpha}$  by
\begin{equation}\label{eq4.8}
\begin{split}
\tilde{a}_{0} = -1 \qquad \tilde{c}_{0} = 0 \qquad \tilde{b}_{0} = +1  \\
\tilde{a}_{i} = +1 \qquad \tilde{c}_{i} = 0 \qquad \tilde{b}_{i} =
-1 \quad 0 < i < \alpha.
\end{split}
\end{equation}
Given an arbitrary pair $ a < b $ in  $X_{\alpha}$,  it suffices to
prove $ (a,b) \cong (\tilde{a},\tilde{b})$.

Let $\beta = min \ \{ j : a_{j} \not= b_{j} \}$, so that $a_{\beta}
< b_{\beta}$ and $a_{j} = b_{j} $ for $ j < \beta$.  Because $X$
is order dense, we can choose $c \in X_{\alpha}$ such that
\begin{equation}\label{eq4.9}
\begin{split}
a_{i} = c_{i} = b_{i} \qquad \mbox{for} \ i < \beta  \\
a_{\beta} < c_{\beta} < b_{\beta} \qquad \mbox{for} \ i = \beta  \\
\qquad c_{i} = 0  \qquad \mbox{for} \ \beta < i < \alpha.
\end{split}
\end{equation}

We will construct an order isomorphism $ [c,b) \cong
[\tilde{c},\tilde{b}) $.  Applying a similar argument (or the same
argument to $X^{*}$) we obtain an isomorphism $(a,c] \cong
(\tilde{a},\tilde{c}] $.  Putting them together we get $ (a,b)
\cong (\tilde{a},\tilde{b})$ as required.

Now define for the purposes of this proof (i.e. ignore (\ref{eq3.5}) for
the duration)
\begin{equation}\label{eq4.10}
\begin{split}
K = \{ \beta \} \cup \{ k : \beta < k < \alpha \ \mbox{and} \
b_{k} > -1 \ \mbox{in} \ J \}  \\  \qquad  K' \ = \ K \cup \{ \alpha
\}. \hspace{2cm}
\end{split}
\end{equation}
Thus, $K'$ is a subset of $\alpha + 1 $ and so is a countable
well-ordered set.  We define $f : K' \rightarrow  X_{\alpha} $ by
\begin{align}\label{eq4.11}
\begin{split}
f(\beta) \ &= \ c  \\
f(k)_{i} \ &= \  \begin{cases}  b_{i} \quad i < k \\
-1 \quad k \leq i < \alpha \end{cases}
\quad \mbox{ for } \ \beta < k \le \alpha.
\end{split}
\end{align}
In particular, $f(\a) \ = \ b.$

For $ k \in K' $ let $k'$
denote its successor in the well-ordered set $K'$.  For $ k \in K'
\setminus \{\beta \} $  we have
\begin{align}\label{eq4.12}
\begin{split}
f(k)_{i} \ = \ f(k')_{i} \  = \  b_{i} \qquad &\mbox{for} \quad i < k  \\
-1 \  = \ f(k)_{k} \  < \  f(k')_{k} \  =  \  b_{k}&   \\
f(k)_{i} \  = -1  \  = \ b_{i} \ = \ f(k')_{i} \qquad &\mbox{for}
\quad k < i < k' \\  f(k)_{i} \  = -1  \   = \ f(k')_{i} \qquad
&\mbox{for} \quad k' \leq i < \alpha .
\end{split}
\end{align}
Thus, $f$ is an order injection.

Notice that if $x \in X_{\a}$ and $k \in K$ with $k > \b$, then
$x_j \geq f(k)_j = -1$ for all $j \geq k$ and
$x_j \geq -1$ for all $j > \b$ with $j \not\in K$.

Thus, if  $k$ is a limit element of $K$ and $x \in X_{\alpha} $ with $ x
< f(k) $  and $j$ the minimum at which $x_j \not= f(k)_j$, then  $x_{j} < f(k)_{j}$.
So either  $j \leq \beta < k $ or else $j \in K $ with
$j < k$.  In either case, there exists some $\tilde{k} \in K$ with
$ max (j, \b) < \tilde{k} < k $  so that $x_j < f(\tilde{k})_j = f(k)_j$ and so $x < f(\tilde{k})$.  It follows
that $f$ is continuous and so Proposition \ref{prop2.8}(c) implies that
\begin{equation}\label{eq4.13}
[c,b) \ = \ \bigcup _{k\in K} [f(k),f(k')).
\end{equation}

By Proposition \ref{prop2.8}(d) there exists an order embedding
$\tilde{f}_{0} : K' \rightarrow [0,+1] $ with $\tilde{f}_{0}(\beta)
= 0  $ and  $\tilde{f}_{0}(\alpha) = +1  $.  Now define $\tilde{f}
: K' \rightarrow X_{\alpha} $ by
\begin{equation}\label{eq4.14}
\begin{split}
\tilde{f}(k)_{0} \ = \ \tilde{f}_{0}(k) \qquad \mbox{for} \quad  k \in K'  \hspace{2cm}\\
\tilde{f}(\beta)_{i} \ = \ 0 \qquad \mbox{for}\quad    0 < i < \alpha   \hspace{2cm} \\
\tilde{f}(k)_{i} \ = \ -1 \qquad \mbox{for} \quad   0 < i < \alpha
\quad \mbox{and} \quad k \in K' \setminus \{\beta \}.
\end{split}
\end{equation}
Clearly, $\tilde{f}$ is an order embedding and again continuity
implies
\begin{equation}\label{eq4.15}
[\tilde{c},\tilde{b}) \ = \ \bigcup _{k\in K}
[\tilde{f}(k),\tilde{f}(k')).
\end{equation}
Notice that $ \tilde{f}(\beta) = \tilde{c} $  and
$\tilde{f}(\alpha) = \tilde{b}$.

Because $X$ is doubly transitive, we can choose for each $k \in K$
an order isomorphism between intervals in $J$ :
\begin{equation}\label{eq4.16}
q_{k} : [f(k)_{k},f(k')_{k}) = [-1,b_{k}) \rightarrow
[\tilde{f}(k)_{0},\tilde{f}(k')_{0}) =
[\tilde{f}_{0}(k),\tilde{f}_{0}(k')).
\end{equation}

Because $\alpha$ is tail-like, there exists for each $k \in K$ a
unique order isomorphism $\tau^{k} : \alpha \rightarrow \alpha
\setminus k = \{ \epsilon : k \leq \epsilon < \alpha \} $.  Define
for each  $k \in K$ the map between intervals of $X_{\alpha}$
\begin{equation}\label{eq4.17}
\begin{split}
Q_{k} : [f(k),f(k')) \rightarrow [\tilde{f}(k),\tilde{f}(k')) \\
Q_{k}(x)_{i} = \begin{cases} q_{k}(x_{k}) \qquad \mbox{for} \ i = 0  \\
x_{\tau^{k}(i)} \qquad \mbox{for} \ 0 < i < \alpha . \end{cases}
\end{split}
\end{equation}
Notice that by (\ref{eq4.12}) $x_{i} = b_{i} $  for all  $i < k$  when $ x
\in [f(k),f(k')) $.  It follows that each $Q_{k}$  is an order
isomorphism.

Putting together these isomorphisms we obtain the required
isomorphism $[c,b) \cong [\tilde{c},\tilde{b}) $.

\end{proof}\vspace{1cm}

\subsection{Size Comparisons}

In distinguishing between CHLOTS we define a rough order of size.

\begin{df}\label{df4.3} For LOTS $X$ and $X_{1}$ we say that $X$  \emph{ injects into } $X_{1}$ if
there exists an order injection $ g : X \rightarrow X_{1} $.  We say that $X_{1}$ is
\emph{ bigger than}\index{size!bigger} $X$ if $X$ injects into $X_{1}$ but not the reverse. When neither injects into the other
we say their sizes are \emph{not comparable}. \index{size!not comparable}  On the other
hand, we say that $X$ \emph{ has the same size as}\index{size!same} $X_{1}$ when each injects into the other.
Finally, we say that the size of $X$ \emph{lies between}\index{size!between} $X_{1}$ and $X_{2}$ when $X_{1}$
injects into $X$ and $X$ injects into $X_{2}$.
\end{df}
\vspace{.5cm}

This is the usual crude partial ordering used to compare order
types.  We will now see that for CHLOTS $X$ and $X_{1}$, $X$
injects into $X_{1}$ iff there exists an order surjection $f :
X_{1} \rightarrow X$.  By Proposition \ref{prop2.1}(a) such a surjection is
always continuous.  On the other hand, there exists a continuous
order injection, i.e. an order embedding, $ g: X \rightarrow
X_{1}$ iff $X \cong X_{1}$.

\begin{prop}\label{prop4.4}
\begin{enumerate}
\item[(a)] If $f: X_{1} \rightarrow X$  is an order surjection of LOTS then there exists
a map $g : X \rightarrow X_{1} $ such that $f \circ g = 1_{X}$.  Any such map $g$ is an order injection.
\item[(b)] Assume $X$ is connected.  If $g : X \rightarrow X_{1} $ is an order injection
with image $g(X)$ $\pm$cofinal in $X_{1}$, then there exists a continuous order
surjection $f: X_{1} \rightarrow X$  such that  $f \circ g = 1_{X}$.
\item[(c)] Assume that $X_{1}$ is a HLOTS and that $X$ is unbounded.  If there exists an
order injection of $X$ into $X_{1}$, then there exists an order injection with image $\pm$cofinal  in $X_{1}$.
\item[(d)]  Assume that $X_{1}$  and  $X$ are CHLOTS.  If there exists a nonconstant,
continuous order map from $X_{1}$ to $X$, then $X$ injects into $X_{1}$.  If there exists
an order embedding of $X$  into $X_{1}$, then $X_{1} \cong X$.
\end{enumerate}
\end{prop}

\begin{proof}(a):  This is a repeat of Proposition \ref{prop2.1}(d).

(b):  Define for $x \in X$ the closed convex set $J_{x} \subset X_1$ by
\begin{equation}\label{eq4.18}
X_{1} \setminus J_{x} \ =\ (\bigcup \{(-\infty,g(a)) : a < x
\})\cup(\bigcup\{(g(b),\infty) : b > x \}).
\end{equation}

If $x_{1} < x_{2}$ in  $X$ then because  $X$ is order dense we can
choose $a, b $ such that $ x_{1} < b < a < x_{2}$.  Since $g(b) <
g(a) $
\begin{equation}\label{eq4.19}
(-\infty,g(a)) \cup (g(b),\infty) \quad = \quad X_{1}
\end{equation}
and so $ J_{x_{1}} \cap J_{x_{2}} = \emptyset$.

If $y \in X_{1}$ equals $g(x)$ then $ y \in J_{x} $.  In particular, $J_x$ is nonempty. If $y
\not\in g(X)$ then because the image is cofinal and coinitial, the
pair $g^{-1}((-\infty,y)),g^{-1}((y,\infty))$  is a partition of
$X$ by nonempty convex sets. By completeness of $X$ we can define $ x =
inf \ g^{-1}((y,\infty))$.  If $g(b) < y $, then $b$ is a lower
bound for $g^{-1}((y,\infty))$ and so $b \leq x $.
Contrapositively, $b > x$ implies $g(b) > y $.  If $g(a) > y $,
then $a \in g^{-1}((y,\infty))$ and so $x \leq a$.
Contrapositively, $a < x $ implies $g(a) \leq y$.  Thus, $y \in
J_{x}$.

It follows that $\{ J_{x} : x \in X \} $ is an $X$ indexed family
of nonempty, closed convex sets with union $X_{1}$.  So mapping
$J_{x}$ to $x$ defines an order surjection $f : X_{1} \rightarrow
X $ which is continuous by Proposition \ref{prop2.1}(a) because each
$f^{-1}(x) = J_{x}$ is closed.  If $x \in X$, then $g(x) \in J_{x}$
implies $f(g(x)) = x $.

(c):  Define $J = \{ y \in X_{1} : $ for some $x_{1},x_{2} \in X
\quad g(x_{1}) < y < g(x_{2}) \} $.  $J$ is an open, convex subset
of $X_{1}$.  Because $X$ is unbounded, $g(X)$ is a subset of $
J$ which is $\pm$cofinal  in $J$.  Because $X_{1}$ is
a HLOTS there exists an order isomorphism $q : J \rightarrow
X_{1}$.  Replace $g$ by $q \circ g$.

(d):  If $f : X_{1} \rightarrow X$ is a continuous order map and $c
< d$ in the image $f(X_{1})$, then let $a = sup f^{-1}(c)$ and $b =
inf f^{-1}(d)$.  The interval $[a,b]$ is compact and connected and
so by continuity of $f$ the image is as well.  Since $f$ is an
order map $f([a,b])$ is a compact, connected subset of $[c,d]$
which contains $c$ and $d$.  Hence, $f([a,b]) = [c,d]$ and so, by
definition of $a$ and $b$, $f((a,b)) = (c,d)$.  Let $h_{1} : X_{1}
\rightarrow (a,b)$ and $h : (c,d) \rightarrow X$ be order
isomorphisms. Then $ f_{1} = h \circ f \circ h_{1}$ is an order
surjection of $X_{1}$ onto $X$ and so $X$ injects into $X_{1}$ by
(a).

If, in addition, $f$ is injective then the restriction to $(a,b)$
is an order isomorphism to $(c,d)$ and so $f_{1}$ is an
isomorphism.

\end{proof}\vspace{.5cm}

As an immediate consequence we obtain the following.

\begin{cor}\label{cor4.5} If $X$ and $X_{1}$ are CHLOTS then $X$ injects into $X_{1}$ iff there exists an
order surjection $g : X_{1} \rightarrow X$.  Such an order surjection is necessarily continuous.
\end{cor}
 \vspace{.5cm}

 We can strengthen this result.

 \begin{theo}\label{theo4.5a} If $X$ and $X_1$ are connected LOTS and $f : X_1 \to X$ is a non-constant
 continuous map, then either $f$ is an order$^*$ map onto a non-trivial interval in $X$, or there exists
 an order  surjection from $X_1$ onto a  non-trivial interval in $X$. If $X$ is a CHLOTS,  then there exists
 an injection from $X$ to $X_1$ which is either order-preserving or order-reversing. \end{theo}

 \begin{proof} Because the LOTS are connected and $f$ is continuous and non-constant the image of $f$ is an
 interval in $X$ in any case. If $f$ is order-reversing then by applying Proposition \ref{prop4.4} to the reverse order, we
 obtain an order$^*$ injection  to $X_1$ from an interval in $X$.

  Assume $f$ is not order-reversing so that there exist $a <  b$ in $X_1$ such that
 with $c = f(a), d = f(b)$ we have $c < d$. We will construct an order preserving injection $g : [c,d] \to [a,b]$.
Let $g(c) = a, g(d) = b$. For $x \in (c,d)$, let $g(x) =  sup f^{-1}(x) \cap [a,b]$ which is not empty because $f([a,b])$ is
connected and contains $[c,d]$. Because $f$ is continuous, $f^{-1}(x)$ is closed and so $f(g(x)) = x$.
If $x < x_1 < d$ then $f([g(x),b])$ contains $[x,d]$ and so there exist points of $f^{-1}(x_1)$
between $g(x)$ and $b$.  Hence, $g(x) < g(x_1) < b$. That is, $g$ is an order injection.

 It now follows from Proposition \ref{prop4.4}(b) that with $a_1 =  inf g((c,d))$ and $b_1 =  sup g((c,d))$ there exists a - necessarily
 continuous - order surjection from $[a_1,b_1]$ to $[c,d]$. Extending by constants below $a_1$ and above $b_1$ we obtain an
 order surjection from $X_1$ onto $[c,d]$.

 Thus, we obtain an  injection to $X_1$ from a non-trivial open interval in $X$.  If $X$ is a CHLOTS we can precede by an isomorphism from
 $X$ to the open interval and get an injection from $X$ to $X_1$.

 \end{proof} \vspace{.5cm}

\begin{prop}\label{prop4.6} Let $X$ and $X_{1}$ be LOTS.
\begin{enumerate}
\item[(a)] Let $f : X \rightarrow X_{1} $ be an order map.  Assume that $X$ is order
 dense and that $D $ is a dense subset of $X$. If the restriction $f|D$ is injective,
 then $f$ is an order injection, i.e. it is injective on all of $X$.
\item[(b)]  Assume that $X$ is order dense and unbounded and that $X_{1}$ is complete.
 If $X$ injects into $X_{1}$, then the completion $\hat{X}$ injects into $X_{1}$.
\item[(c)]  Let $\alpha \leq \beta $ be positive ordinals.  If $X$ injects into $X_{1}$
then $X^{\alpha}$ injects into $X_{1}^{\beta}$.  If, in addition, $X$ and $X_{1}$ are HLOTS
then $X_{\alpha}$ injects into $(X_{1})_{\beta}$.
\item[(d)]  If $\alpha$ and $\beta$ are positive ordinals then
\begin{equation}\label{eq4.20}
(X^{\alpha})^{\beta} \quad \cong X^{\alpha \cdot \beta }.
\hspace{2cm}
\end{equation}
If, in addition, $X$ is a HLOTS then
\begin{equation}\label{eq4.21}
(X_{\alpha})_{\beta} \quad \cong X_{\alpha \cdot \beta }.
\hspace{2cm}
\end{equation}\end{enumerate}
\end{prop}

\begin{proof}(a): This is a repeat of Proposition \ref{prop2.1} (e).

(b): If $f : X \rightarrow X_{1} $ is an order injection, then
$\hat f$ is an order injection by Proposition \ref{prop2.4}.

(c): If $f : X \rightarrow X_{1} $ is an order injection and $z \in
X_{1}$ then we can define the injection $\tilde{f} : X^{\beta}
\rightarrow X_{1}^{\alpha}$ by
\begin{equation}\label{eq4.23}
\tilde{f}(x)_{i} = \begin{cases}  f(x_{i}) \qquad \mbox{for}\ i < \beta \\
z \qquad \mbox{for} \ \beta \leq i < \alpha.  \end{cases}
\end{equation}
In the HLOTS case, we can assume that the distinguished intervals
have been chosen so that $f(J) \subset J_{1}$ and that $z \in
J_{1}$.  Then $\tilde{f}$ restricts to an injection of $X_{\beta}$
into $(X_{1})_{\alpha}$.

(d):  We see first that the natural bijection $q :(X^{\alpha})^{\beta}
\rightarrow X^{\beta \times \alpha}$ given by $q(x)(i,j) =
x(i)(j)$ preserves the orders.

For $q(x_1) < q(x_2)$ means that for some $(i,j) \in \beta \times \a$
$q(x_1)(i,j) < q(x_2)(i,j)$ and $q(x_1)(k,\ell) = q(x_2)(k,\ell)$ for all
$(k,\ell) < (i,j)$. So for $k < i, x_1(k) = x_2(k) $. Furthermore,
$x_1(i)(\ell) = x_2(i)(\ell)$ for $\ell < j$, while $x_1(i)(j) < x_2(i)(j)$.
That is, $x_1(i) < x_2(i)$. Consequently, $x_1 < x_2$.  Since $q$ is a bijection
it is an order isomorphism.

The result follows because $\beta \times \alpha =
\alpha \cdot \beta$.

In the HLOTS case with $J$ the distinguished
interval in $X$ we choose $J^{\a} = (\pi^{\alpha})^{-1}(J)$ as the distinguished
interval in $X^{\alpha}$. It is an interval by (\ref{eq4.4}). Then $q$ maps $(X_{\alpha})_{\beta}$
onto $ X_{\alpha \cdot \beta } $.

\end{proof} \vspace{.5cm}

\noindent{\bfseries Remark.}  Notice that for an unbounded,
connected LOTS $X$ the projection $\pi' : X' \rightarrow X$ is
injective on the dense subset $D = X \times \{-1\}$ but not itself
injective. Thus, the order dense hypothesis in (a) is required.
\vspace{.5cm}

We will now show that for a CHLOTS $X$ and positive ordinals
$\alpha > \beta $ it is always true that $X_{\alpha}$ is bigger
than $X_{\beta}$.

\begin{df}\label{df4.7}  A LOTS $X$ is called \emph{ order simple}\index{order simple}\index{LOTS!order simple}
if $X'$ is bigger than $X$ where $X'$ is the AS double of $X$.
\end{df}
\vspace{.5cm}

It is always true that  $X$ injects into $X'$, e.g. use $x \mapsto x^-$.  So $X$ is
order simple when $X'$ does not inject into $X$.

\begin{prop}\label{prop4.8}
\begin{enumerate}
\item[(a)] If $X$ is an uncountable, order dense LOTS which satisfies the countable
chain condition, then $X$ is order simple.
\item[(b)] If $X_{1}$ and $X_{2}$ are LOTS of the same size, then $X_{1}$ is order simple iff $X_{2}$ is.
\item[(c)] If $X$ is order simple, then the reverse $X^{*}$ is order simple.
\item[(d)] If $X$ is order simple, then there does not exist an injective order* map from $X'$ into $X$.
\end{enumerate}
\end{prop}

\begin{proof} (a) If $ f : X' \rightarrow X$
is an order injection for any LOTS $X$, then  $\{ (f(z^-),f(z^+))
: z \in X \} $ is a family of open intervals in $X$ and each
is nonempty if $X$ is order dense. If $z_{1} < z_{2}$ in $X$,
then $f(z_{1}^+) < f(z_{2}^-) $ so that the intervals are pairwise
disjoint.  If $X$ satisfies c.c.c., then $X$ must be countable.

(b)   If $f : X_{2} \rightarrow X_{1}$ is an order injection, then
from (\ref{eq3.8}) we obtain the order injection $f' : X_{2}' \rightarrow
X_{1}' $.  So if $X_{1}$ is at least as big as $X_{2}$ then
$X_{1}'$ is at least as big as $X_{2}'$.  Thus, if $X_{1}$ and
$X_{2}$ have the same size and $X_{1}'$ has the same size as
$X_{1}$, then $X_{2}'$ has the same size as well.

(c) It is clear that
\begin{equation}\label{eq4.24}
(X^{*})' \qquad = \qquad (X')^{*}.
\end{equation}
So any order injection of $X'$ into $X$ is an order injection of
$(X^{*})'$ into $X^{*} $.

(d) Define the map
\begin{equation}\label{eq4.25}
\begin{split}
q : X' \rightarrow (X')'   \\
q(x^{\pm}) \quad = \quad (x^{\pm})^{\pm}.
\end{split}
\end{equation}
This is an order embedding of $X'$ onto the closed set whose
complement is the open set of isolated points $\{
(x^{+})^{-},(x^{-})^{+} : x \in X \} $. Now if $g : X' \rightarrow
X$ is an order reversing injection, then we use (\ref{eq3.9}) to define an
order injection from $X'$ to $X$ as the composition:
\[ \begin{CD}
X' @>q>> (X')' @>g^{*}>> X' @>g>> X.
\end{CD}  \]

\end{proof} \vspace{.5cm}

\noindent{\bfseries Remark.}  If $X$ is unbounded and order dense so that
$X'$ has no isolated points, then the map $q$ of (\ref{eq4.25}) is actually
an order isomorphism of $X'$ onto $(X')''$ defined via (\ref{eq3.16}).
\vspace{.5cm}

In the Remark after Proposition \ref{prop3.6} we observed that $\Z'
\cong \Z$ and so $\Z$ is not order simple.  By
Proposition \ref{prop2.8}(a) $\Q$ is not order simple.  There exist
connected LOTS which are not order simple as well.  With $J =
[-1,+1] \subset \R $ define
\begin{equation}\label{eq4.26}
\begin{split}
X_{n} \quad = \quad [n,n+1) \times J^{n}  \\
X \quad = \quad \Sigma_{n \in \omega} X_{n}.\hspace{1cm}
\end{split}
\end{equation}
It is easy to see that $X$ is connected with $min = 0$ in $X_{0} =
[0,1) $.  The order isomorphisms
\begin{equation}\label{eq4.27}
\begin{split}
f_{n} : X_{n} \times J \rightarrow X_{n+1} \qquad \mbox{for} \quad n \in \omega  \\
f_{n}(x,t)_{i} =
\begin{cases}
 x_{0} + 1 \qquad i = 0  \\
x_{i} \qquad 0 < i \leq n
\\  t \qquad  i = n+1
\end{cases}
\end{split}
\end{equation}
can be put together to get an order embedding of $X \times J$ into
$X$.  Since $X' \subset X \times J $, it follows that $X$ is not
order simple.  Notice that $ X $ has no $max$.

\begin{lem}\label{lem4.9}{\bfseries (The Shift Lemma)} \index{Shift Lemma}Let $X$ be a complete LOTS with $min = m$.
If $ f : X' \rightarrow X$ is an order injection, then for all $x \in X$
\begin{equation}\label{eq4.28}
f(x^{+}) \quad > \quad x.
\end{equation}
In particular, $X$ has no $max$.
\end{lem}

\begin{proof} Define $S = \{ a \in X : f(x^{+})
> x  \ \mbox{for all} \  x \leq a \}$. Hence, $a \in S$ implies
$(-\infty,a] \subset S$.  Since $m = min \ X$
\begin{equation}\label{eq4.29}
m \quad \leq  \quad f(m^{-}) \quad < \quad f(m^{+})
\end{equation}
and so $m \in S$.  It suffices to show that $S$ is unbounded as
this implies $S = X$.

Assume $S$ is bounded and let $z = sup \ S$.

First we show that $z \in S$.  If not, then $a < z $ for all $a \in
S$ and so $a^{+} < z^{-}$.  Hence
\begin{equation}\label{eq4.30}
a \quad < \quad f(a^{+}) \quad < \quad f(z^{-}) \quad < \quad
f(z^{+}).
\end{equation}
Thus, $f(z^{-})$ is an upper bound for $S$.  Since $z = sup \ S$
\begin{equation}\label{eq4.31}
z \quad \leq \quad f(z^{-}) \quad <  \quad f(z^{+}).
\end{equation}
As all $x < z$ are not upper bounds for $S$ they are elements of
$S$. It follows that $z \in S$ after all.

Now let $y = f(z^{+}) > z$.  If $z < x \leq y$ then $z^{+} <
x^{+}$ and so
\begin{equation}\label{eq4.32}
x \quad \leq \quad y \quad = \quad f(z^{+}) \quad < \quad
f(x^{+}).
\end{equation}
Thus, $y \in S$.  Since $y > z $, this contradicts the assumption
that $z = sup \ S$.

$X$ has no $max$ because $S$ is unbounded.  Also,  $x = max$ could
not satisfy (\ref{eq4.28}).

\end{proof} \vspace{.5cm}

\begin{cor}\label{cor4.10} A LOTS $X$ is order simple if it satisfies one of the following conditions:
\begin{enumerate}
\item[(i)]  $X$ is compact.
\item[(ii)]  $X$ is doubly transitive and complete.
\item[(iii)]  $X \cong (X_1)_{\alpha} $ where $X_1$ is any CHLOTS and $\alpha$ is any positive ordinal.
\end{enumerate}
\end{cor}

\begin{proof} (i) The Shift Lemma implies that
if a complete LOTS $X$ is not order simple and has a $min$, then it
has no $max$.  In particular, a compact LOTS is order simple.

(ii)  Assume that $f : X' \rightarrow X$ is an order map with $X$
doubly transitive and complete.  For any pair $a < b$ in $X$
define $\tilde{a} = f(a^{-})$ and $\tilde{b} = f(b^{+})$.  If $\tilde{a} = \tilde{b}$, then $f$ is not injective.
If $\tilde{a} \not= \tilde{b}$, then
$f$ maps $[a,b]'$ into $[\tilde{a},\tilde{b}]$.  Because $X$ is
doubly transitive we can choose an order isomorphism $g :
[\tilde{a},\tilde{b}] \rightarrow [a,b]$.  The composition $g
\circ f : [a,b]' \rightarrow [a,b]$ is an order map.  Since
$[a,b]$ is compact, (i) implies that this map is not
injective and so $f$ is not injective.  It follows that $X$ is
order simple.

(iii)  For the closed interval $J$ in $X_1$,the subset  $J^{\alpha}$
of $(X_1)_{\alpha}$ is compact and so is order simple. If $a < b$ in
$J$ then Proposition \ref{prop4.1}(f) implies that  $ (X_1)_{\alpha}$ is order isomorphic to the subset $(a+,b-) $ of $J^{\alpha}$.  Hence,
$(X_1)_{\alpha}$ has the same size as $J^{\alpha}$ and so is order simple by Proposition \ref{prop4.8}(b).

\end{proof} \vspace{.5cm}

\begin{theo}\label{theo4.11} Assume that $X$ is a CHLOTS.  If $\alpha > \beta$ are positive ordinals
then $X_{\alpha}$ is bigger than $X_{\beta}$.  Furthermore,   $X_{\alpha}$ is not homeomorphic to  $X_{\beta}$.
\end{theo}

\begin{proof} Recall that $j^{\b}_{\a}$
injects $X_{\beta}'$ into $X_{\a}$, see (\ref{eq4.5}). If $f : X_{\alpha} \rightarrow
X_{\beta}$ is an injective map then the composite $f \circ j^{\b}_{\a} : X_{\beta}' \rightarrow X_{\beta}$ is an
injection which is order preserving or order reversing if $f$ is.
Because $X_{\beta}$ is order simple by  Corollary \ref{cor4.10}, $f$ cannot
be order preserving and by Proposition \ref{prop4.8}(d) it cannot be order reversing
either.

If $f : X_{\alpha} \to X_{\beta}$ were a homeomorphism, then by Lemma \ref{lem3.1} it would be
either order preserving or reversing.

\end{proof}\vspace{.5cm}

\begin{cor}\label{cor4.12} If $X$ is a CHLOTS and $\alpha$ is a positive ordinal such that $X_{\alpha}$
is transitive, then $\alpha$ is countable and tail-like.
\end{cor}

\begin{proof} Since $X_{\alpha}$ is connected,
transitivity implies first countability by Proposition \ref{prop3.2}(c).  By
Proposition \ref{prop4.1}(d), $\alpha$ is countable.

Now assume that $\alpha$ is not tail-like.  This means that there
exists $\beta < \alpha$ and $\epsilon < \alpha$ such that
$\epsilon \cong \alpha \setminus \beta = \{ i : \beta \leq i <
\alpha \}$. Choose $b \in X_{\alpha}$ such that $b_{\beta} \in
J^{\circ}$, i.e. $-1 < b_{\beta} < +1 $ in $X$.  Define $a = (-1)-$
 in $X_{\alpha}$, i.e. $a_{i} = -1 $ for all $i < \alpha$. We
will obtain a contradiction from the assumption that there exists
$f \in H_{+}(X_{\alpha}) $ such that $f(a) = b $.

Define $c,d \in X_{\alpha}$ by
\begin{equation}\label{eq4.33}
\begin{split}
c_{i} \  = \  d_{i} \  = \  b_{i} \qquad \mbox{for} \quad i < \beta \hspace{1cm}\\
c_{\beta} \  =\  -1, \qquad d_{\beta} \  = \  +1   \hspace{1cm}\\
c_{i} \  = \  +1, \qquad d_{i}\  = \  -1 \qquad \mbox{for} \quad
\beta < i < \alpha.
\end{split}
\end{equation}
Clearly, we have $c < b < d$ and since
$\epsilon \cong \alpha \setminus \beta $  Proposition \ref{prop4.1} (f) implies
\begin{equation}\label{eq4.34}
(c,d) \quad \cong \quad X_{\epsilon}.
\end{equation}
By continuity of $f$, $f^{-1}((c,d))$ contains a neighborhood
of $a$ and so contains some interval $(x,a)$ with $x < a$ in
$X_{\alpha}$.  This implies $x_{0} < a_{0} = -1 $ in $X$ and so we
can choose $\tilde{a}, \tilde{b} \in X$ such that $x_{0} <
\tilde{a} < \tilde{b} < a_{0}$ in $X$.  Because $x < \tilde{a}+ <
\tilde{b}- < a$ in $X_{\alpha}$ we have
\begin{equation}\label{eq4.35}
(\tilde{a}+,\tilde{b}-) \quad \subset \quad f^{-1}((c,d)).
\end{equation}

By Proposition \ref{prop4.1}(f) again
\begin{equation}\label{eq4.36}
X_{\alpha} \quad \cong \quad (\tilde{a}+,\tilde{b}-).
\end{equation}

Composing the isomorphism of (\ref{eq4.36}) with the restriction of $f$
and the isomorphism of (\ref{eq4.34}) we get an order injection of
$X_{\alpha}$ into $X_{\epsilon}$ contradicting Theorem \ref{theo4.11}.

\end{proof}\vspace{.5cm}

For any connected unbounded LOTS $X$, e.g. any CHLOTS,  the associated \emph{ Cantor Space for} $X$, \index{Cantor Space for $X$}
\index{CHLOTS!Cantor Space}
denoted $C(X)$, \index{$C(X)$} is the two point compactification
\begin{equation}\label{eq4.37}
C(X) \quad = \quad \bullet X' \bullet
\end{equation}
where $X'$ is the AS double. In particular, the Cantor Space for $\R$ is the Fat Cantor Set.

Because an isomorphism maps $max$ to $max$ and $min$ to $min$,(\ref{eq3.09a}) implies for connected unbounded LOTS $X, X_1$
\begin{equation}\label{eq3.37a}
C(X) \ \cong \ C(X_1) \quad \Longleftrightarrow \quad X \ \cong \ X_1.
\end{equation}

For any positive ordinal $\alpha$ we
have
\begin{equation}\label{eq4.38}
C(X_{\a}) \quad = \quad \bullet (X_{\alpha})' \bullet.
\end{equation}

If $X$ is a CHLOTS, then for $a < b$ in $X$ the
isomorphism $f : X_{\alpha} \rightarrow (a+,b-)$ of Proposition
\ref{prop4.1}(f) induces the isomorphism $f' :(X_{\alpha})' \rightarrow
((a+)^{+},(b-)^{-})$ which extends to the two-point
compactification to show
\begin{equation}\label{eq4.39}
C(X_{\a}) \quad \cong \quad [(a+)^{+},(b-)^{-}] .
\end{equation}\vspace{.25cm}

\begin{theo}\label{theo4.13}  Assume that $X$ is a CHLOTS.  If $\alpha > \beta$ are positive ordinals
then $C(X_{\a})$ is bigger than $X_{\alpha}$ which is bigger than $C(X_{\b})$.
\end{theo}

\begin{proof} Because $X_{\alpha}$ is order
simple, there is no order injection from $(X_{\alpha})'$ into it.
$C(X_{\alpha})$ projects onto $X_{\alpha}$  and  contains
$(X_{\alpha})'$.  Hence, $C(X_{\alpha})$ is bigger than
$X_{\alpha}$.

By (\ref{eq4.39}) $C(X_{\beta})$ is the same size as $(X_{\beta})'$.  Now
choose $a < b$ in $J^{\circ} \subset X$. Clearly, we have
\begin{equation}\label{eq4.40}
\begin{split}
\{-1,+1\} \times \{-1,+1\} \quad \cong \quad \{-1,a,b,+1\} \qquad \mbox{and so} \\
((X_{\beta})')' \quad \cong \quad X_{\beta} \times \{-1,a,b,+1\}
\quad \subset \quad X_{\beta +1}.
\end{split}
\end{equation}
Hence, there is an order injection from  $C(X_{\beta})'$ into
$X_{\beta +1}$ which injects into $X_{\alpha}$.  If $X_{\alpha}$
were to inject into $C(X_{\beta})$, then $C(X_{\beta})$ would not be
order simple, contradicting Corollary \ref{cor4.10}.

\end{proof} \vspace{.5cm}

\noindent{\bfseries Remark.}
Because $C(X_{\alpha})$ and $C(X_{\beta})$ are not connected, we cannot use Lemma \ref{lem3.1} to show that
they are topologically
distinct. As far as we know it may happen that for some $\a \not= \b$ $C(X_{\alpha})$ and $C(X_{\beta})$ are homeomorphic.
On the other hand, it is clear that $C(X)$ is separable or satisfies
c.c.c. iff $X$ satisfies the corresponding property. So
 the original Fat Cantor Set $C(\R)$ is the
only  Cantor Space of a CHLOTS which is separable.

In summary, we have the following.

\begin{theo}\label{theo4.14}  If $F$ is a CHLOTS, then $F_{\om^{\g}}$ is a tower of CHLOTS strictly
increasing in size and no two distinct members of which
are homeomorphic. The tower $C(F_{\om^{\g})}$  of CHLOTS Cantor Spaces is also strictly increasing in size. \end{theo}
\vspace{1cm}

\section{\textbf{Trees}}
\vspace{.5cm}

\subsection{Trees and Bi-Ordered Trees}

For the theory of trees we follow \cite{J} Section 22.

If $(T,\succ)$ is a partially ordered set then for $p \in T$ we
define the \emph{tail set}, the \emph{predecessor set} and the
\emph{successor set}\index{tail set}\index{prededessor set}\index{successor set} of $p$:
\begin{equation}\label{eq5.1}
\begin{split}
T_{p} \quad = \quad \{ q \in T : q \succeq p \} \hspace{4cm}  \\
A_{p} \quad = \quad \{ q \in T : q \prec p \} \hspace{4cm}  \\
S_{p} \quad  = \quad  \{ q \in T : q \succ p \  \mbox{and}\
\not\exists r \in T\  \mbox{such that}\  q \succ r \succ p \}.
\end{split}
\end{equation}\index{$T_p$}\index{$A_p$}\index{$S_p$}

A non-empty, partially ordered set  $(T,\succ)$ is called a \emph{tree}\index{tree} when $A_{p}$ is well-ordered by $\prec$ for each $p \in T $.
The elements of $T$ are then called the \emph{vertices}\index{tree!vertex}\index{vertex} of $T$.
If $p \in T$ then the \emph{order}\index{vertex!order} of $p$, denoted $o(p)$ is the
ordinal whose order type is that of $A_{p}$, i.e. there is a unique order isomorphism from $o(p)$ onto $A_p$. The bijection can be extended
to $o(p) + 1$ by mapping $o(p)$ to $p$. For any ordinal
$\alpha$ the \emph{level} $\alpha$ set\index{level $\alpha$ set} is
\begin{equation}\label{eq5.2}
L_{\alpha} \quad = \quad \{ p \in T : o(p) = \alpha \}.
\end{equation}\index{$L_{\a}$}
The successors of $p$ are the points of $T_{p} $ at the next
level, i.e.
\begin{equation}\label{eq5.3}
o(p) = \alpha \quad \Longrightarrow \quad S_{p} \  = \  T_{p} \cap
L_{\alpha +1}.\hspace{1cm}
\end{equation}
If $A$ is a nonempty subset of $T$ then we define its \emph{
height}\index{height} by
\begin{equation}\label{eq5.4}
\begin{split}
h(A) \  = \   sup  \{ o(p)+1 : p \in A \} \hspace{1cm} \\  = \   min  \{\alpha :
o(p) < \alpha \ \mbox{for all}\  p \in A \}.
\end{split}
\end{equation}\index{$h(A)$}

Any subset of a tree is a tree in its own right, leading to
different notions of order and height for elements of the subset.
The two concepts agree when $R \subset T$ is a
\emph{subtree}\index{subtree}\index{tree!subtree} defined by the condition
\begin{equation}\label{eq5.5}
p \in R\quad \Longrightarrow \quad A_{p} \subset R.\hspace{1cm}
\end{equation}
For example, for each positive ordinal $ \alpha$ we define the $\a$ \emph{truncation}\index{truncation}\index{tree!truncation}
subtree
\begin{equation}\label{eq5.6}
T^{\alpha} \  = \  \{ p \in T : o(p) < \alpha \} \ = \  \bigcup
_{\beta < \alpha} \ L_{\beta}.
\end{equation}\index{$T^{\a}$}

A \emph{branch}\index{branch}\index{tree!branch} $x$ of a tree $T$ is a maximal, linearly ordered
subset of $T$.  Because the set of predecessors of any vertex $p$
of $T$ is totally ordered, any branch of $T$ is a subtree by maximality.  A
branch $x$ is a well-ordered set whose order type is that of the
ordinal $h(x)$.   By Zorn's
Lemma every subset of $T$  linearly ordered by $\succ$ is
contained in some branch.  In particular, each vertex lies in some
branch.

We denote by $X(T)$ the \emph{branch space}\index{branch space}\index{tree!branch space}\index{$X(T)$}
of the tree $T$, i.e. the set of branches of $T$.

 Since there is a branch through every vertex
\begin{equation}\label{eq5.25}
h(T) \quad = \quad  sup  \{ h(x) : x \in X(T) \}.
\end{equation}
We define for a branch $x$ and any ordinal
$\alpha < h(x)$ the vertex $x_{\alpha}$\index{$x_{\a}$} to be the - unique - level $\alpha$
element of $x$, i.e.
\begin{equation}\label{eq5.7}
\{x_{\alpha}\} \quad = \quad x \cap L_{\alpha}.
\end{equation}

Let $p$ be a vertex of $T$.  The tail set $T_{p}$ is a tree
but not a subtree of $T$.  For $q \in T_{p}$ let $o_{p}(q)$ denote
its order in the tree $T_{p}$.  Clearly,
\begin{equation}\label{eq5.8}
o(q) \quad = \quad o(p) \ + \ o_{p}(q),
\end{equation}
using ordinal addition.

If $x$ is a branch of $T$ then
\begin{equation}\label{eq5.9}
\begin{split}
x \cap T_{p} \not= \emptyset \quad \Longleftrightarrow \quad p \in x \hspace{1cm} \\
\mbox{in which case} \quad p = x_{o(p)}.\hspace{1cm}
\end{split}
\end{equation}
In that case, $x \cap T_{p}$ is a branch of $T_{p}$. On the other
hand, if $y$ is a branch of $T_{p} $ then
\begin{equation}\label{eq5.10}
j_{p}(y) \quad =  \quad y \cup A_{p}
\end{equation}
is a branch of $T$. Hence, (\ref{eq5.10}) defines a bijection
\begin{equation}\label{eq5.11}
j_{p} : X(T_{p}) \rightarrow \{ x \in X(T) : p \in x \},
\end{equation}
and we have for all $y \in X(T_{p})$
\begin{equation}\label{eq5.12}
h(j_{p}(y)) \quad = \quad o(p) + h_{p}(y),
\end{equation}
where $h_{p}$ denotes the height with respect to $T_{p}$.

We call a subset $A$ of a tree $T$ an \emph{antichain}\index{antichain} if no two
vertices in $A$ are comparable with respect to $\succ$.  For
example, $L_{\alpha}$ is an antichain for any ordinal $\alpha$.

Let $\#S$ denote the cardinality of a set $S$. Following
\cite{J} Chapter 4:

\begin{df}\label{df5.1} A tree $T$ is called a \emph{semi-normal tree}
\index{tree!semi-normal}\index{semi-normal tree} when it satisfies the following conditions
\begin{enumerate}
\item[(i)] $\#L_{0} = 1$, i.e. $T $ has a \emph{root}\index{root}\index{tree!root} which we denote $0 \in T$.
\item[(ii)] For all $p \in T$, $\#S_{p} \not= 1$.
\item[(iii)]  If $p,q \in T$ with $o(p) = o(q)$ a limit ordinal and with $A_{p} = A_{q}$, then $ p = q $.
\end{enumerate}

$T$ is called a \emph{normal}\index{tree!normal}\index{normal tree} tree when it is semi-normal and, in
addition, satisfies the condition
\begin{enumerate}
\item[(iv)] If $p \in T$ and $\alpha$ is an ordinal with $ o(p) \leq \alpha < h(T) $,
then there exists $q \in T_{p}$ with $ o(q) = \alpha $.
\end{enumerate}
This implies
\begin{enumerate}
\item[(v)] If $p \in T$ and $ o(p) + 1 < h(T) $, then $S_{p} \not= \emptyset$ and so  $\#S_{p} > 1$.
\end{enumerate}

We call a tree \emph{$\Omega$-bounded}\index{$\Omega$-bounded tree}\index{tree!$\Omega$-bounded} when it satisfies
\begin{enumerate}
\item[(vi)]   $h(x) < \Omega$ for all $x \in X(T)$ and so $h(T) \le \Omega$.
\end{enumerate}

An \emph{Aronszajn tree}\index{Aronszajn tree}\index{tree!Aronszajn} is a normal tree of height $\Omega$, the
first uncountable ordinal, which satisfies (vi) and
\begin{enumerate}
\item[(vii)] If $p \in T$ and $ o(p) + 1 < h(T) $, then $S_{p}$ is an infinite set.
\item[(viii)]  $L_{\alpha}$ is a countable set for each $\alpha < \Omega$.
\end{enumerate}

A \emph{Suslin tree}\index{Suslin tree}\index{tree!Suslin} is a normal tree of height $\Omega$ which
satisfies (vii) and
\begin{enumerate}
\item[(ix)]  Every antichain in $T$ is a countable set.
\end{enumerate}
\end{df}
\vspace{.5cm}

{\bfseries N. B. From now on we will assume that all trees are at
least semi-normal, unless otherwise mentioned.}

Notice that a branch $x$ is a subtree which is not semi-normal as each vertex $p \in x$ with $o(p) + 1 < h(x)$
has a single successor in the branch. \vspace{.5cm}

Condition (ii) says that if $p \in T$ has any successors then it
has at least two.  Thus, a vertex is either terminal, i.e. $S_{p}
= \emptyset$ or the tree branches in at least two directions after
$p$. For a normal tree, the latter always happens unless $o(p) + 1 = h(T)$.

If $T$ is a normal tree, $0 < \alpha < h(T)$  and $p \in
L_{\alpha}$ then $T^{\alpha}$ and $T_{p}$ are normal trees with
\begin{equation}\label{eq5.13}
\begin{split}
h(T^{\alpha}) \quad = \quad \alpha \hspace{2cm} \\
h(T_{p}) \quad \cong \quad h(T) \setminus \alpha.\hspace{1cm}
\end{split}
\end{equation}

Since $L_{\alpha}$ is an antichain, (ix) implies (viii) in
Definition \ref{df5.1}.  Furthermore, if $x \in X(T)$ then we can choose
for each $\alpha$ with $\alpha + 1 < h(x)$ a successor
$y_{\alpha}$ of $x_{\alpha}$ different from $x_{\alpha +1}$.  The
set $\{ y_{\alpha} : \alpha + 1 < h(x) \}$ is an antichain.  It
follows that (ix) implies (vi).  Thus, every Suslin tree is
Aronszajn.

If $T$ is normal and $x \in X(T)$, then by (v) either $h(x) =
h(T)$ or $h(x)$ is a limit ordinal less than $h(T)$.

If $x$ and $y$ are two distinct branches of a  tree $T$,
then by (i) the set
\begin{equation}\label{eq5.14}
Eq(x,y) \quad = \quad  \{ i : x_{i} = y_{i} \} \quad \not= \
\emptyset.
\end{equation}
Furthermore, $ i \in Eq(x,y) $ and  $ j < i $ imply $ j \in
Eq(x,y) $ and so $Eq(x,y)$ is an ordinal and by condition (iii) it
is not a limit ordinal.  So we can define
\begin{equation}\label{eq5.15}
\begin{split}
\epsilon (x,y) \quad =  \quad   max  Eq(x,y) \qquad \mbox{so that}  \\
Eq(x,y) \quad = \quad \epsilon(x,y) + 1 \quad \cong \quad  x \cap
y.
\end{split}
\end{equation}
We call $\epsilon(x,y)$ the \emph{equality level}\index{equality level}\index{$\ep(x,y)$} of the pair
$x,y$.  Clearly, the equality level $\epsilon$ is the unique
ordinal $\epsilon$ such that
\begin{equation}\label{eq5.16}
x_{\epsilon}\  =\  y_{\epsilon} \qquad \mbox{and} \qquad
x_{\epsilon +1}\ \not= \ y_{\epsilon +1}.
\end{equation}
Since both $x$ and $y$ extend to the $\epsilon + 1 $ level we have
\begin{equation}\label{eq5.17}
\epsilon(x,y) + 1 \quad < \quad h(x), h(y).
\end{equation}

\begin{df}\label{df5.2} A \emph{bi-ordered tree}\index{bi-ordered tree}\index{tree!bi-ordered} $T$ is a  tree of height greater
than $1$ with a linear order on each nonempty set of successors $S_{p}$.

We will say that a bi-ordered tree $T$ is \emph{of $Y$ type}\index{tree!of $Y$ type} for a
LOTS $Y$ if for $p \in T$
\begin{equation}\label{eq5.18}
S_{p} \ \not= \emptyset \qquad \Rightarrow \qquad S_{p} \quad
\cong \quad Y.
\end{equation}
We will say that a bi-ordered tree $T$ is \emph{of unbounded type},   \emph{of dense type}, \emph{of separable type} or \emph{of
countable type} if each nonempty successor set is a LOTS which is unbounded, order dense, separable
or of countable type, respectively.\index{tree!of unbounded type}
\index{tree!of dense type}\index{tree!of separable type}\index{tree!of countable type}
\end{df}
\vspace{.5cm}

For example, a normal bi-ordered tree is of $\Q$ type
if $p \in T$ and $o(p) + 1 < h(T)$ implies that $S_{p}$ is an unbounded,
countable, order dense LOTS.

For a bi-ordered tree $T$ the \emph{induced order}\index{induced order}\index{branch space!induced order} on the branch
space $X(T)$ is defined by
\begin{equation}\label{eq5.19}
x < y \quad \Longleftrightarrow \quad x_{\epsilon} = y_{\epsilon}
\  \mbox{and} \ x_{\epsilon + 1} < y_{\epsilon +1} \
\mbox{for some ordinal} \ \epsilon.
\end{equation}
Note that at the $\ep + 1$ level we are using the LOTS ordering.
By (\ref{eq5.16}) the ordinal $\epsilon$ is the equality level of $x,y$.
Conversely, if $x \not= y$ and $x,y$ have equality level $\epsilon$,
then both $x_{\epsilon +1}$ and $y_{\epsilon +1} $ are successors
of $x_{\epsilon}$ and so either $x_{\epsilon +1} < y_{\epsilon
+1}$ or the reverse.  Furthermore, it is easy to check that
\begin{equation}\label{eq5.20}
x < z < y \qquad \Rightarrow \qquad \epsilon(x,y) =
min (\epsilon(x,z),\epsilon(z,y))
\end{equation}
and from this that $x < y$.  Consequently, with the induced order
$X(T)$ is a LOTS.

If $T$ is a bi-ordered tree, then by retaining the LOTS structure
on each $S_{p}$ we give each tail tree $T_{p}$ and any subtree, e.g. any truncated tree
$T^{\alpha}$, the structure of a bi-ordered tree of the same type
(i.e. normal, of $Y$ type, or of unbounded, dense, separable or
countable type).

\begin{prop}\label{prop5.2a} Let $T$ be a bi-ordered tree with $p \in T$ and $R$ a subtree of $T$.
\begin{enumerate}

\item[(a)] The injection $j_{p} : X(T_{p}) \rightarrow
X(T)$ of (\ref{eq5.11}) is an order embedding with a convex image.

\item[(b)] The map $\pi : X(T) \to X(R)$ defined by $x \mapsto x \cap R$ is an
 order surjection.
\end{enumerate}\end{prop}

\begin{proof} (a): If $x < z < y$ in $T$ with $x,y \in T_p$, then $\epsilon(x,y) \ge o(p)$ and so by (\ref{eq5.20}),
$z \in T_p$ and so $x < z < y$ in $T_p$. Thus, $j_p$ is an order injection with a convex image and so an embedding by Proposition
\ref{prop2.1}(b).

(b): Because $R$ is a subtree it is clear that $\pi(x) < \pi(y)$ in $R$ implies $x < y$ in $T$. It follows that $\pi$ is an order map.
By Zorn's Lemma any branch of $R$ extends to a branch of $T$. Hence, $\pi$ is surjective.

\end{proof} \vspace{.5cm}

\noindent{\bfseries Remark.} It follows from  Proposition \ref{prop2.1}(d) that there exists an order injection from $X(R)$ into $X(T)$ for
any subtree $R$ of $T$. \vspace{.5cm}

In particular, for $0 < \beta \leq \alpha
\leq h(T) $ the projection map
\begin{equation}\label{eq5.21}
\begin{split}
\pi^{\a}_{\b} : X(T^{\alpha}) \rightarrow X(T^{\beta})  \\
 \pi^{\a}_{\b}(z) \quad = \quad z \cap T^{\beta}
\end{split}
\end{equation}
is an order surjection.  When $\alpha = h(T)$ we will omit the superscript,
writing $\pi_{\beta} : X(T) \rightarrow X(T^{\beta}) $.

If $T$ is a bi-ordered tree, then we denote by $T^*$\index{$T^*$} the same tree with the reverse LOTS order on each nonempty $S_p$. It is clear
the $T^*$ is a bi-ordered tree of $Y^*$ type if $T$ is of $Y$ type and is otherwise of the same type (normal, unbounded, dense, etc.).
The  branch space $X(T^*)$ is the branch space $X(T)$ with the reverse ordering, ie.  $X(T^*) = X(T)^*$.

\begin{prop}\label{prop5.3a}  If $T$ is a  bi-ordered  tree, then  $X(T)$ is an order dense LOTS when either of the
following two conditions hold.
\begin{itemize}
\item[(i)] The tree $T$ is of dense type, i.e. each nonempty $S_p$ is order dense.
\item[(ii)] $T$ is normal and there exists a limit ordinal $\a$
such that $h(T) = \a$ or $h(T) = \a + 1$ and  $S_p$ has no max for any $p \in T$ (e.g. if
$T$ is of unbounded type).
\end{itemize}\end{prop}

\begin{proof} For $x < y$ in $X(T)$ let
$\epsilon$ be the $x,y$ equality level. Because $x \not= y$, $\ep < h(x)$.

(i): If $S_{p}$ is order dense, then with $p = x_{\epsilon}
= y_{\epsilon}$ we can choose $q \in S_{p}$ so that
\begin{equation}\label{eq5.22}
x_{\epsilon +1} \quad < \quad q \quad < \quad y_{\epsilon +1}.
\end{equation}

(ii): If  $h(x) = \a + 1$, then condition (iii) of Definition \ref{df5.1} implies that $\ep < \a$. If $h(x) \le \a$, then
$h(x)$ is a limit ordinal because $T$ is normal. So in either case, $\ep + 2 < h(x)$. We have $x_{\epsilon +1} < y_{\epsilon +1}$
and with $p = x_{\epsilon +1}$, $o(p) + 1 < h(x)$ implies that $S_p$ is nonempty. If $S_p$ has no max, then we can choose $q \in S_p$
such that $x_{\epsilon +2} < q$.

In either case, if $z$ is a branch through $q$, then $x < z < y$.

\end{proof} \vspace{.5cm}

When $X(T)$ is an order dense LOTS, we will denote by $\widehat{X(T)}$ its completion.
In particular when $T$ is of dense type, we will write $\hat{\pi}^{\a}_{\b} : \widehat{X(T^{\alpha})}
\rightarrow  \widehat{X(T^{\b})} $ and  $\hat{j}_{p} :
\widehat{X(T_p)} \rightarrow \widehat{X(T)}$ for the extensions to the
completions of the maps defined above.

\begin{prop}\label{prop5.3}  Let $T$ be a  bi-ordered  tree.
\begin{enumerate}
\item[(a)]   For each $p \in T$ the inclusion $j_{p}$ is an order embedding onto a
convex subset of $X(T)$.

If $o(p) + 1 = h(T)$, then $T_{p} = \{ p \}$ and the image of $j_{p}$
is the unique branch through $p$.

Assume that for every $p \in T$ with $o(p) + 1 < h(T)$ , $\# S_p > 1$
(i.e. condition (v) of Definition \ref{df5.1} holds), e.g. $T$ is
 normal. If   $o(p) + 2 < h(T)$ or if $o(p) + 2 = h(T)$ and  $\# S_p > 2$,
 then  $j_p(X(T_p))$ has a nonempty interior in $X(T)$.

\item[(b)] The following equivalence holds:
\begin{align}\label{eq5.22a} \begin{split}
\text{For all} \ p \in T, \a \ \text{with}& \ o(p)+1 \le \a \le h(T) \\
 j_p(X(T^{\a}_p)) \ &\text{is a closed subset of} \ X(T^{\a}). \\
&\Longleftrightarrow   \\
\text{For all} \ \b \leq \a \leq h(T), \ \ \pi_{\b}^{\a} : \  &X(T^{\a}) \to X(T^{\b}) \ \text{ is  continuous }.
\end{split}\end{align}

These are both true when any of the following three conditions hold.
\begin{itemize}
\item[(i)] Every successor set $S_p$ with $o(p) > 1$ is either empty or bounded.
\item[(ii)] $T$ is normal and of dense type.
\item[(iii)] $T$ is normal and of unbounded type.
\end{itemize}

\item[(c)]  Assume $T$ is normal and of unbounded type,
and $p \in T$
with $o(p) + 1 < h(T)$.  The image $j_{p}(X(T_{p}))$ is a nonempty, infinite, clopen, convex set in
$X(T)$. If, in addition, $X(T)$ is order dense, then $X(T_p)$ is order dense
 and the image  $\hat{j}_{p}(\widehat{X(T_p)})$ is a nonempty, open interval in
$\widehat{X(T)}$.

\item[(d)]  If $S_{0}$, the successor set to the root, is unbounded, then $X(T)$ is unbounded.
If $S_{0}$ is unbounded and is $\sigma$-bounded, then $X(T)$ is
$\sigma$-bounded. If, in addition, $X(T)$ is order dense, then $\widehat{X(T)}$ is unbounded and $\sigma$-compact.

\item[(e)]  Assume there exists a limit ordinal $\a$ such that $h(T) = \a$ or $h(T) = \a + 1$.
If $T$ is normal and of unbounded type, then $X(T)$ is order dense and has dense holes.

\item[(f)]  Assume $T$ is normal and of dense type. For each $0 < \beta \leq \alpha \leq h(T) $ the projection $\pi^{\a}_{\b}$
and its extension $\hat{\pi}^{\a}_{\b}$ to the completions are continuous order surjections.
The extension $\hat{j}_{p}$ is an order embedding onto an interval in $\widehat{X(T)}$.

\item[(g)]  A subset $W$ of $X(T)$ is dense in $X(T)$ if for every $p \in T$ with
$o(p) + 1  < h(T)$ the set $\{q \in S_{p} : \exists x \in W $ with $q \in x\}$ is
dense in $S_{p}$.  In particular, if for every $p \in T$ there exists $x \in W$
such that $p \in x$, then $W$ is dense in $X(T)$.

Conversely, assume for every $p \in T$ with $o(p) + 1 < h(T)$, $\# S_p > 1$ and if $o(p) + 2 = h(T)$, then $\# S_p > 2$.
If $W \subset  X(T)$ is dense in $X(T)$, then for every $p \in T$ there exists $x \in W$
such that $p \in x$.
\end{enumerate}
\end{prop}

\begin{proof} (a): By Proposition \ref{prop5.2a}(a) $j_p$ is an order embedding with a convex image.

Now assume that (v) of Definition \ref{df5.1} holds. If $o(p) + 2 < h(T)$, then
we can choose $q_{1} < q_{4} \in S_{p}$ and $q_2 < q_3 \in S_{q_4}$. If $S_p$ contains at least three points, then there exist $q_1 < q_2 < q_3 $ in $S_p$.
In either case, let $z_i \in X(T)$ be a branch containing $q_i$ for $i = 1, 2, 3$. So $z_1 < z_2 < z_3$ and thus the interval
 $(z_1,z_3) \subset j_p(X(T_p))$ is nonempty. \vspace{.25cm}

(b): We are writing $T^{\a}_p$ for  $(T^{\a})_p$.

Assume that every $j_p$ has a closed image.

For $x \in T^{\b}$, $(\pi^{\a}_{\b})^{-1}(x) = \bigcap_{p \in x} j_p(X(T^{\a}_p))$. Hence, the order surjection $\pi^{\a}_{\b}$
has closed point-inverses and so is continuous by Proposition \ref{prop2.1}(a).

Assume that every $\pi^{\a}_{\b}$ is continuous. If $p = 0$, then $X(T_p) = X(T)$ is closed.  If $o(p) + 1 = h(T)$, then
$j_p(X(T_p))$ is a singleton and so is closed.

Assume that $1 \le o(p), o(p) +1 < h(T)$. Let $\b = o(p)+1$ and $x = x(p)$. $j_p(X(T^{\a}_p)) = (\pi^{\a}_{\b})^{-1}(x)$ and so is closed.\vspace{.25cm}

Assume (i):  Replacing $T$ by $T^{\a}$ we prove that $j_p(X(T_p))$ is closed. We may assume $1 \le o(p), o(p) +1 < h(T)$.

We  show that the convex set $j_p(X(T_p))$ has a max and min and so is closed.

Inductively, we define a collection of points $p_{\a}$ totally ordered by $\succ$ and use it to define the max $M$ of $j_p(X(T_p))$.

Begin with $\a = o(p)$ and $p_{\a} = p$.

 If $\a = \b+1$ and $S_{p_{\b}} = \emptyset$ then the process stops and we let $M = x(p_{\b})$. Otherwise, let $p_{\a}$ be the maximum
 element of $S_{p_{\b}}$.

 If $\a$ is a limit ordinal, then $\{ p_{\b} : o(p) \le \b < \a \}$ is contained in a branch $x \in X(T)$. Clearly, $x \in j_p(X(T_p))$.
 If $h(x) = \a$, then the process stops and we let $M = x$. Otherwise, let $p_{\a} = x_{\a}$.

 The process stops at or before $\a = h(T)$ and defines $M$. If $y \in X(T_p)$, then $\ep = \ep(y,M)$ has $o(p) \le \ep$.
 Hence, $y_{\ep} = p_{\ep}$ and so $y_{\ep + 1} < M_{\ep + 1}$ since the latter is the maximum element of $S_{p_{\ep}}$. Thus, $M$ is the
 maximum element of $j_p(X(T_p))$.\vspace{.25cm}

 Assume (ii): By Proposition \ref{prop5.3a}, $X(T^{\b})$ is order dense and so the surjection $\pi^{\a}_{\b}$ is continuous by
 Proposition \ref{prop2.1}(a).\vspace{.25cm}

 Assume (iii): Replacing $T$ by $T^{\a}$ we prove that $j_p(X(T_p))$ is closed. We may assume $1 \le o(p), o(p) +1 < h(T)$.

 If $ x < y, p \in x $ and $p \not\in y$,  then with $\ep = \ep(x,y)$, the equality level of $x,y$, $\ep + 1 \le o(p)$.
 By assumption $o(p) + 1 < h(T)$, Either $h(y) = h(T)$ or
$h(y)$ is a limit ordinal. In either case, with $p_1 = y_{\ep + 1}$, $y_{\ep + 2} \in S_{p_1}$ and, by assumption,
$S_{p_1}$ is unbounded. So  we can choose $q \in S_{p_1}$ with $q < y_{\ep + 2} $. A
branch $z$ through a point $q$ satisfies $x < z < y$ and $p \not\in z$, because $p_1 = y_{\ep + 1} \in z$.
It follows that  the open interval
$(z,\infty)$ in $X(T)$ contains $y$ and is disjoint from the image
of $j_{p}$.  Arguing similarly if $y < x $ we see that $j_{p}$ has
a closed image in X(T).\vspace{.25cm}

(c): Since $o(p) + 1 < h(T), p \in x $ and $S_{p}$ is unbounded, for $x \in X(T_p)$
we can choose $q_{1},q_{2} \in S_{p}$ so that $q_{1} < x_{\alpha
+1} < q_{2}$ where $\alpha = o(p)$.  Choosing branches $z_{i}$
through $q_{i}$ for $i = 1,2$ we see that the image of $j_{p}$
contains the open interval $(z_{1},z_{2})$.  Hence, the image of
$j_{p}$ is open. It is closed by (b).

If  $X(T)$ is order dense, then the convex subset $j_p(X(T_p))$ is order dense and it is isomorphic to $X(T_p)$ via $j_p$.
Hence, $X(T_p)$ is order dense.

Furthermore, the image of
$\hat{j}_{p}$ is the completion of the image of $j_{p}$ and so it
is an interval in $\widehat{X(T)}$.  The completion of an open convex set in $X(T)$ is an open
interval in $\widehat{X(T)}$.  Hence, the image of $\hat{j}_{p}$ is
open in $\widehat{X(T)}$.

Since $o(p) + 1 < h(T) $ $S_p$ is unbounded and so
$T_{p}$ is infinite. Hence, the image of $j_{p}$ is infinite.\vspace{.25cm}

(d): Since the root $0$ is the unique element of level $0$ in $T$,
$X(T^{2}) \ \cong \ S_{0}.$

Since $ \pi_{2} : X(T) \rightarrow X(T^{2}) $ is an order
surjection the results for $X(T)$ follow from Proposition \ref{prop2.1}(a).
If $X(T)$ is unbounded, then its completion is. If $X(T)$ is
$\sigma$-bounded then its completion is,too, and so is
$\sigma$-compact.\vspace{.25cm}

(e): If $T$ is of unbounded type,  $X(T)$ is order dense by Proposition \ref{prop5.3a}. It is unbounded by (d).

Let $x,y$ be elements of $X(T)$ with $x < y$,  $\ep = \ep(x,y)$, and $p = x_{\ep} = y_{\ep}$.
 Since $X(T)$ is order dense we choose $z \in X(T)$ such that $x < z < y$. By (\ref{eq5.20})
$\ep_1 = \ep(z,y) \geq \ep$. With $q = z_{\ep_1 + 1}$ we have $ q < y_{\ep_1 + 1}$.

Define the set
\begin{equation}\label{eq5.23a}
G \ = \ (-\infty,z) \cup j_{q}(T_{q}) \ = \ (-\infty,z] \cup j_{q}(T_{q}) \ \subset \ X(T).
\end{equation}
$G$ is a  convex set in $X(T)$ which contains $x$ but not the upper bound $y$.

If $h(y) = \a + 1$, then $\ep_1 + 1 \le h(y)$ implies
$\ep_1 \le \a$. But (iii) of Definition \ref{df5.1} then implies $\ep_1 < \a$. If $h(y) \le \a$, then $\ep_1 + 1 \le h(y)$ implies $\ep_1 < \a$.
Consequently, $\ep_1 + 2 < \a \le h(T)$ because $\a$ is a limit ordinal.
Because $o(q) + 1 = \ep_1 + 2 < h(T) $, $j_{q}(T_{q})$ is clopen in $X(T)$ by (c). Since $(-\infty,z)$ is open and
$(-\infty,z]$ is closed, it follows that
$G$ is clopen in $X(T)$.

We show that $G$ has no supremum and so reveals a hole between $x$ and $y$.

Since $G$ is closed, if $w$ were the supremum of $G$, then $w$ would be an element of $G$ and so would be the max of $G$.

On the other hand, $X(T)$ is unbounded. So $w \in G$ and $G$ open implies there exist $z_1, z_2 \in X(T)$ such that
$w \in (z_1,z_2) \subset G$. Since $X(T)$ is order dense, we can take $z_1, z_2$ to lie in $G$. Hence, $G$ has no max.\vspace{.25cm}

(f): Each
$\pi^{\a}_{\b}$ is continuous by (b).  The extensions to the completions are
therefore well-defined, surjective and continuous by Proposition \ref{prop2.4a}.\vspace{.25cm}

(g):  Assume that $\{ q \in S_{p} : \exists z \in
W$  with  $ q \in z \}$ is dense in $S_{p}$ for all $p$.

For $x < y$ in $X(T)$ with $(x,y)$ nonempty there exists $z \in X(T)$ such that $x < z < y$. Let
$\epsilon_1 = \epsilon(z,y) \ge \epsilon(x,y)$ and $p =
z_{\epsilon_1} = y_{\epsilon_1}$.

 If $\ep_1 > \ep$, then $z_{\epsilon_1 + 1} < y_{\epsilon_1 + 1}$ implies that
there exists $w \in W$ with $w_{\epsilon_1 + 1} < y_{\epsilon_1 + 1}$ in $S_p$. Because $\ep_1 > \ep$,
$x_{\epsilon + 1} < y_{\epsilon + 1} = w_{\epsilon + 1}$ and so $x < w < y$.

If $\ep_1 = \ep$, then $x_{\epsilon + 1} < z_{\epsilon + 1} < y_{\epsilon + 1}$ implies that
there exists $w \in W$ with $x_{\epsilon + 1} < w_{\epsilon + 1} < y_{\epsilon_1 + 1}$ in $S_p$ and
so again $x < w < y$.  So the condition is sufficient for density.

In particular, if $ p \in T$ implies there exists $x \in W$ with $p \in x$, then $W$ is dense.

Now assume that $o(p) + 1 < h(T)$ implies $\# S_p > 1$ and  $o(p) + 2 = h(T)$ implies $\# S_p > 2$. Further, assume that
 $W$ is dense.

 By (a) $j_p(X(T_p))$ has a nonempty interior and so meets $W$. It follows that there
 exists $w \in W$ with $p \in w$.

\end{proof}\vspace{.5cm}

If $\alpha + 1 \leq h(T)$ and $p \in L_{\alpha}$, then
\begin{equation}\label{eq5.26}
x(p) \quad =  \quad \{p\} \cup A_{p}
\end{equation}
is a branch in $X(T^{\alpha +1})$ of height $\alpha + 1$.
So $p \mapsto x(p)$ defines an injective
map from $L_{\alpha}$ into
$X(T^{\alpha +1})$ and so into $\widehat{X(T^{\alpha +1})}$ when $T$ is of dense type.  Regarding
this map as an inclusion,  we will regard $L_{\alpha}$ as
a subset of these LOTS
and so induce an order upon it.  Since $(\pi_{\alpha +1})^{-1}(x(p))$ in
$X(T)$ consists of the branches which contain $p$ we have
\begin{equation}\label{eq5.27}
j_{p}(X(T_{p})) \  = \  (\pi_{\alpha +1})^{-1}(x(p))\  \subset\
X(T) \qquad \mbox{with} \ \alpha = o(p).
\end{equation}

If $\alpha$ is a limit ordinal with $\alpha + 1 \leq h(T)$ then
condition (iii) of Definition \ref{df5.1} implies that
\begin{equation}\label{eq5.28}
\pi^{\alpha+1}_{\alpha} : X(T^{\alpha +1}) \ \cong \
X(T^{\alpha}).
\end{equation}
$\pi^{\alpha+1}_{\alpha}$ maps $x(p) = \{p\} \cup A_{p}$ to
$A_{p}$ which is a branch in $X(T^{\alpha})$ of height $\alpha$.

When $T$ is of dense type, this order isomorphism extends to $\hat{\pi}^{\alpha+1}_{\alpha}$,
an order isomorphism between the completions.

If $\alpha \leq h(T)$ is a limit ordinal then we define
\begin{equation}\label{eq5.29}
\tilde{L}_{\alpha} \quad =  \quad \{ x \in X(T^{\alpha}) :
h(x) = \alpha \}.
\end{equation}\index{$\tilde{L}_{\alpha}$}
If $\alpha < h(T)$, then
\begin{equation}\label{eq5.30}
\pi^{\alpha+1}_{\alpha}(L_{\alpha}) \ \subset \
\tilde{L}_{\alpha}.
\end{equation}

\begin{lem}\label{lem5.4}  Let $T$ be a  bi-ordered normal tree and let $\alpha$ be an ordinal.
\begin{enumerate}
\item[(a)] If $\alpha < h(T)$, then for all $p \in T$ there exists $x \in X(T)$ such
that $p \in x$ and $h(x) > \alpha$.  Furthermore,
the subset $(\pi_{\alpha +1})^{-1}(L_{\alpha})$
is dense in $X(T)$ (and hence in $\widehat{X(T)}$ when $X(T)$ is order dense).
\item[(b)]  If $\alpha = h(T)$ and $\alpha$ is a countable limit ordinal, then
for all $p \in T$ there exists $x \in X(T)$ such that $p \in x$ and $h(x) = \alpha$.
Furthermore,
the subset $\tilde{L}_{\alpha}$ is dense in $X(T)$  (and hence in $\widehat{X(T)}$ when $X(T)$ is order dense).
\end{enumerate}
\end{lem}

\begin{proof} (a) If $o(p) < \alpha$ then by
condition (iv) of Definition \ref{df5.1} there exists $q \in L_{\alpha}$
such that $p \prec q$. If $o(p) \geq \alpha$ then there exists a
unique $q \in L_{\alpha}$ such that $p \succeq q$.  Any branch $x$
containing $p$ and $q$ satisfies $h(x) > \alpha$.  Density follows
from Proposition \ref{prop5.3}(g).

(b) Let $\{\alpha_{n}\}$ be an increasing sequence of ordinals
with $o(p) < \alpha_{1}$ and $ sup  \{\alpha_{n}\} = \alpha $.
Apply condition (iv) inductively to choose a sequence of vertices
$\{p_{n}\}$ such that
\begin{equation}\label{eq5.31}
\begin{split}
p \ = \ p_{0} \ \prec \ p_{1} \ \prec ...     \\
o(p_{n}) \quad = \quad \alpha_{n}.
\end{split}
\end{equation}
The unique branch $x$ which contains $\{p_{0},p_{1},...\}$ has
height $\alpha$.  Density again follows from Proposition \ref{prop5.3}(g)

\end{proof}\vspace{.5cm}

\begin{prop}\label{prop5.5a} Let $T$ be a bi-ordered tree.

If for each $p \in T$ the successor set $S_{p}$ is a, possibly empty, complete LOTS
and is bounded and so compact if $o(p) > 0$, then $X(T)$ is complete.
If, in addition, $T$ is of dense type, then $X(T)$ is connected.
\end{prop}

\begin{proof}By Proposition \ref{prop5.3a} $X(T)$ is order dense when $T$ is of
dense type. So $X(T)$ is  connected when it is complete.

We prove by induction on $\alpha$ that $X(T^{\alpha})$ is
complete and connected when the successor sets are connected.

\textbf{Case 1:}  $\alpha = 2$: $X(T^{\alpha})
\cong S_{0}$  and so it is complete.  \vspace{.25cm}

\textbf{Case 2:}  $\alpha = \beta + 1$ with $\beta$ a limit ordinal:  then
$\pi^{\a}_{\b} :  X(T^{\alpha}) \rightarrow
X(T^{\beta})$ is an order isomorphism and so
$X(T^{\alpha})$ is complete by induction hypothesis.\vspace{.25cm}

\textbf{Case 3:} $\alpha = \beta + 1$ with $\beta = \epsilon + 1 $:
we define a family of compact  LOTS indexed by the LOTS
$X(T^{\beta})$.  For $x \in X(T^{\beta})$
with $h(x) < \beta$ let $X_{x} = \{x\}$ the trivial
LOTS.  If $x \in X(T^{\beta})$ with $h(x) = \beta$, then $x = x(p)$
with $p \in T$ and $o(p)= \epsilon$. Let $X_{x} = \{ x \}$ if $S_p = \emptyset$ and $X_{x} = S_{p}$ otherwise.
It is easy to see that
\begin{equation}\label{eq5.33}
X(T^{\alpha}) \quad \cong \quad \Sigma \{ X_{x} :x \in
X(T^{\beta}) \}.
\end{equation}
$X(T^{\alpha})$ is complete by Proposition \ref{prop2.2}.\vspace{.25cm}

\textbf{Case 4:}  $\alpha$ is a limit ordinal:
Let $A \subset X(T^{\alpha})$ which is bounded. For each $\b < \a$ $A_{\b} = \pi_{\b}^{\a}(A)$ is bounded and
so has a supremum $s_{\b} \in X(T^{\b})$ by inductive hypothesis. If $\b_1 < \b < \a$, then $ A_{\b_1} = \pi_{\b_1}^{\b}(A_{\b})$.
By Proposition \ref{prop5.3}(b) each $\pi_{\b_1}^{\b}$ is continuous and so by Proposition \ref{prop2.1} (c)
$ s_{\b_1} = \pi_{\b_1}^{\b}(s_{\b})$.  Hence, $s = \bigcup_{\b < \a} s_{\b}$ is a branch of $X(T^{\a})$. If $y < s \in X(T^{\a})$,
then with $\ep = \ep(y,s)$, $y_{\ep + 1} < s_{\ep + 1}$ and so $ \pi_{\ep + 1}^{\a}(y) <  \pi_{\ep + 1}^{\a}(s) = s_{\ep + 1}$.
So there exists $a \in A$ such that $\pi_{\ep + 1}^{\a}(y) < \pi_{\ep + 1}^{\a}(a) $ and so $y < a$. Similarly, one shows that $y \in A$ with
$y \not= s$ implies $y > s$. Hence, $s =  sup A$. With a similar argument for the infimum we see that  $X(T^{\alpha})$ is complete.

\end{proof}\vspace{.5cm}

It will be helpful to describe the completion of the branch space
of a  bi-ordered normal tree $T$ of dense type as a branch space itself. Recall that the completion
$\hat{X}$ of an order dense LOTS $X$ is a connected LOTS which has
a $max$ or $min$ iff $X$ does. We call a LOTS $Y$ the
\emph{completion with endpoints of} $X$ if $Y$ is the completion
of $X$ with $max$ (and $min$) attached if $X$ did not already have
a $max$ (resp. a $min$).  So $Y$ is compact as well as connected.
If $X$ had neither $max$ nor $min$ to begin with, then the
completion with endpoints is the two-point compactification
$\bullet \hat{X} \bullet$. In general, if $X$ is a dense subset of
a compact LOTS $Y$, then $Y$ is isomorphic to the completion with endpoints of $X$.

\begin{df}\label{df5.5} Let $T$ be a
bi-ordered normal tree of dense type.  We define
its \emph{completion}\index{tree!completion}\index{completion (tree)},
denoted $\hat{T}$, to be the tree which contains $T$ as follows:
\begin{enumerate}
\item[(a)] For the root $0$ in $T$ and $\hat{T}$, the successor set in $\hat{T}$,
denoted $\hat{S}_{0}$, is the completion of the order dense LOTS $S_{0}$.
\item[(b)] For $p \in T$ with $o(p) > 0$ the successor set in $\hat{T}$, denoted
$\hat{S}_{p}$, is the completion with endpoints of the order dense LOTS $S_{p}$.
\item[(c)]  If $q \in \hat{S}_{p} \setminus S_{p}$ for any $p \in T$, then $q$ is
called a \emph{new vertex} of $\hat{T}$.  If $q$ is a new vertex of $\hat{T}$, then
its successor set in $\hat{T}$, denoted $\hat{S}_{q}$, is empty.
\end{enumerate}
\end{df}
\vspace{.5cm}

\begin{prop}\label{prop5.6} Let $T$ be a
bi-ordered normal tree of dense type and let $\hat{T}$ be its completion.
\begin{enumerate}
\item[(a)]  $\hat{T}$ is a bi-ordered tree with $h(\hat{T}) = h(T)$.
\item[(b)]  For each $p \in T$ the successor set $\hat{S}_{p}$ is a connected LOTS which is compact if $o(p) > 0$.
\item[(c)]  If $q \in \hat{S}_{p}$ is a new vertex in $\hat{T}$, then $q$ is the end-point
of a unique branch of $\hat{T}$ namely $x(q) = \{q\} \cup A_{q}$ with
\begin{equation}\label{eq5.33aa}
h(x(q))\  =\  o(q) + 1 \ = \ o(p) + 2,
\end{equation}
and so $h(x(q))$ is a successor ordinal.  No two new vertices lie
on the same branch, i.e. the set of new vertices is an anti-chain
in $\hat{T}$.
\item[(d)]  Each branch of $T$ is a branch of $\hat{T}$ of the same height, i.e. $X(T) \subset X(\hat{T})$.
\item[(e)]  If $S_{0}$ is unbounded, then so are $\hat{S}_{0}, X(T)$ and $X(\hat{T})$.
\item[(f)] If $S_{0}$ is unbounded, then $X(\hat{T})$ is the completion $\widehat{X(T)}$. If $S_0$ has both $max$ and $min$,
then $X(\hat{T})$ is compact and is the completion with endpoints of $X(T)$. In either case, $\hat{X}(T)$ is connected.
\item[(g)] If for each $p \in T$ the successor set $S_{p}$ is a connected LOTS which is compact if $o(p) > 0$, then $\hat{T} = T$ and so
$X(T)$ is connected.
\end{enumerate}
\end{prop}

\begin{proof} It is clear that $\hat{T}$ is
semi-normal and so is a bi-ordered tree, although it is usually not normal.  If $q$ is a new vertex,
then $x(q)$ is a branch since the successor set for $q$ is empty.
The results of (a),(b),(c) and (d) follow easily.

(e):  If $S_{0}$ is unbounded, then its completion $\hat{S}_{0}$ is
and so $X(T)$ and $X(\hat{T})$ are  by Proposition \ref{prop5.3}(d).

(f): From Proposition \ref{prop5.5a} (b) it follows that $X(\hat{T})$ is connected.

If $q \in \hat{S}_{p}$ is a new vertex and $x(q) < x \in X(\hat{T})$, then let $\ep = \ep(x(q),x)$ and
$r = x_{\ep} \in \{ p \} \cup A_p$. If $r = p$, then because $S_p$ is order dense and dense in $\hat{S_p}$ we can choose
$q_1 \in S_p$ and $x_1 \in X(T)$ with $q_1 \in x_1$ so that $q < q_1 < x_{\ep+1}$ and so $x(q) < x_1 < x$.
On the other hand, if $\ep = o(r) < o(p)$, then because $S_r$ is order dense we can choose $q_1 \in S_r$
so that $x(q)_{\ep+1} = x(p)_{\ep+1} < q_1 < x_{\ep+1}$. Again if $x_1 \in X(T)$ with $q_1 \in x_1$, then $x(q) < x_1 < x$.
Similarly, if $x(q) > x \in X(\hat{T})$, then there exists $x_1 \in X(T)$ with $x(q) > x_1 > x$.
Hence, $X(T)$ is dense in
$X(\hat{T})$.

Thus, $X(\hat{T})$ is connected and
contains the order dense LOTS $X(T)$ as a dense subset. If $S_0$ is unbounded, it follows
that $X(\hat{T})$ is the completion of $X(T)$.

Now assume that $S_0$ has a maximum $M$. Define  $x_M \in X(\hat{T})$ inductively with $(x_M)_1 = M$.
If $(x_M)_i$ is defined and is a new vertex or $i+1 = h(T)$, then $x_M$ terminates with height
$i + 1$. If $p = (x_M)_i \in T$ with $i+1 < h(T)$,
then $(x_M)_{i+1}$ is chosen to be the maximum of $\hat{S}_p$. If for a limit ordinal $\a \le h(T)$ $ \ (x_M)_i$ is defined
for all $i < \a$, then
all such  $(x_M)_i \in T$. By (iii) of Definition \ref{df5.1} there exists at most one
vertex $p \in T$ with $o(p) = \a$ and $p \succ (x_M)_i$ for all $i < \a$. If no such  $p$ exists, then $\{ (x_M)_i \}$
 defines the branch $x_M$ of height $\a$. If such a $p$ exists, let $(x_M)_{\a} = p$.

It is easy to see that
$x_M =  max X(\hat{T})$ and so is $ max X(T)$ if it does not terminate at a new vertex. It it does terminate at a new vertex
$(x_M)_{i+1}$, then with $p = (x_M)_i \in T$, $S_p$ has no maximum and so $X(T)$ has no maximum.

With a similar construction for the minimum, we see that if $S_0$ has both a maximum and a minimum, then $X(\hat{T})$ does so as well
and so is compact as well as connected. It follows that $X(\hat{T})$ is the completion with endpoints of $X(T)$.

(g): Obvious.

\end{proof} \vspace{.5cm}

\noindent{\bfseries Remark.} If $T$ is a bi-ordered normal
tree of dense type and  $x \in X(\hat{T})$ with $h(x) < h(T)$, then $h(x)$ is a limit ordinal iff $x \in
X(T)$.
\vspace{.5cm}

Clearly, for any ordinal $\alpha$ the completion of $T^{\alpha}$
is $(\hat{T})^{\alpha}$ with the new vertices attached which have
order less than $\alpha$. The situation for the tail trees
requires a bit of quibbling.

\begin{lem}\label{lem5.7}  Let $T$ be a bi-ordered normal
tree of dense, unbounded type.
 If $p \in T$ with $o(p) > 0$ then $X((\hat{T})_{p})$ is a compact, connected LOTS.  It is
 the two-point compactification of the completion of $X(T_{p})$. That is,
\begin{equation}\label{eq5.34}
 X((\hat{T})_{p}) \quad \cong \quad \bullet X(\widehat{T_{p}})\bullet.
\end{equation}
\end{lem}

\begin{proof}  Since $o(p) > 0$ the successor
space $\hat{S}_{p}$ is the compact LOTS with the $max = M_{p}$ and
$min = m_{p}$ attached as two new vertices.  In $T_{p}$ the vertex
$p$ is the root and so the $max$ and $min$ are not included in the
completion.  That is,
\begin{equation}\label{eq5.35}
(\hat{T})_{p}  \quad = \quad \widehat{T_{p}}  \cup \{m_{p},M_{p}\}
\end{equation}
from which the result clearly follows.

\end{proof} \vspace{.5cm}

If $x < y$ in $X(T)$ and $q \in L_{\ep + 1}$ is between  $x_{\epsilon +1}$ and $y_{\epsilon +1}$ with $\ep = \ep(x,y)$, then
every branch through $q$ lies between $x$ and $y$.  Furthermore,
since the completion  $\hat{\pi}_{\epsilon +2}$ is order
preserving we have
\begin{equation}\label{eq5.36}
x_{\epsilon +1}\  < \ q \ < \ y_{\epsilon +1} \quad \Longrightarrow
\quad (\hat{\pi}_{\epsilon +2})^{-1}(x(q)) \ \subset \ (x,y) \
\subset \  \widehat{X(T)}.
\end{equation}

\begin{lem}\label{lem5.8}  Assume that $T$ is a  bi-ordered normal tree of dense type and that $S_{0}$ is unbounded.
 Let $\alpha$ be a positive ordinal with $\alpha + 1 < h(T)$. Let $\hat{\pi}_{\alpha +1} :
  \widehat{X(T)}\rightarrow \widehat{X(T^{\alpha +1})}$ be the canonical projection.

If $x \in  \widehat{X(T^{\alpha +1})}$ but $x \not= x(p)$ for any $p
\in L_{\alpha}$ then $x$ is the image under $\hat{\pi}_{\alpha
+1}$ of a unique point in $\hat{X}(T)$, i.e.
\begin{equation}\label{eq5.37}
\#(\hat{\pi}_{\alpha +1})^{-1}(x)\quad = \quad 1.
\end{equation}

The collection $\{(\hat{\pi}_{\alpha +1})^{-1}(x(p)) : p \in
L_{\alpha} \}$ is a pairwise disjoint family of closed, nontrivial
subintervals of $\widehat{X(T)}$, and the open set
\begin{equation}\label{eq5.38}
O_{\alpha +1} \quad =  \quad \bigcup \{ [(\hat{\pi}_{\alpha
+1})^{-1}(x(p))]^{\circ} : p \in L_{\alpha} \}
\end{equation}
is dense in $ \widehat{X(T)}$.

Furthermore, the open set
\begin{equation}\label{eq5.39}
O^{\alpha +2}_{\alpha +1} \quad =  \quad \bigcup \{
[(\hat{\pi}_{\alpha +1}^{\alpha +2})^{-1}(x(p))]^{\circ} : p \in
L_{\alpha} \}
\end{equation}
is dense in $\widehat{X(T^{\alpha +2})}$  and the restriction
\begin{equation}\label{eq5.40}
\hat{\pi}_{\alpha +2} :  \widehat{X(T)} \setminus O_{\alpha +1}
\rightarrow \widehat{X(T^{\alpha +2})} \setminus O^{\alpha +2}_{\alpha
+1}
\end{equation}
is a homeomorphism.
\end{lem}

\begin{proof}  By Proposition \ref{prop5.6} $ \widehat{X(T)} =
X(\hat{T})$.  Let $x$ be a branch of $\hat{T}$.  If $h(x) \leq
\alpha $ then $x$ regarded as a branch of $\hat{T}$ is the unique
branch which contains $x$ regarded as a branch of $\hat{T}^{\alpha
+1}$. If $h(x) > \alpha$ then $p = x_{\alpha} \in \hat{T}$ is
defined and $\hat{\pi}^{\alpha +1}$ maps $x$ to $x(p)$. If $p$ is
a new vertex of $\hat{T}$ then $x = x(p)$ is the unique branch
which contains $p$. Otherwise, $p \in L_{\alpha}$.

By (\ref{eq5.27}) and Proposition \ref{prop5.3} (a) the interval $(\hat{\pi}_{\alpha +1})^{-1}(x(p))$ is
nontrivial for $p \in L_{\alpha}$ because $\alpha + 1 < h(T)$ and
by Lemma \ref{lem5.4}(a) the union of these intervals is dense in
$\widehat{X(T)}$. If the union of a family of nontrivial intervals is
dense then the union of the interiors is dense, since each interval is contained in
the closure of its interior, provided that the LOTS is order dense.

We can apply the result to $T^{\alpha +2}$ whose height is $\alpha
+ 2 > \alpha + 1$. By (\ref{eq5.37}) the map in (\ref{eq5.40}) is a bijection.
The closed subsets of the completions are locally compact spaces
and the map is topologically proper.  Hence, it is a
homeomorphism.

\end{proof} \vspace{1cm}

\subsection{Countability Conditions}

\begin{prop}\label{prop5.9}  Let $T$ be a bi-ordered normal tree of dense type with $S_{0}$ unbounded.
 Assume that the height of $T$ is not the successor of a limit ordinal.
Define $\epsilon = h(T)$ if $h(T)$ is a limit ordinal and by $\epsilon + 2 = h(T)$ if the height is a successor.
\begin{enumerate}
\item[(a)] The following conditions are equivalent.
\begin{itemize}
\item[(i)]  $X(T)$ is separable.
\item[(ii)]  $\widehat{X(T)}$ is separable.
\item[(iii)]  $\widehat{X(T)} \quad \cong \quad \R$.
\item[(iv)]  $h(T)$ is a countable ordinal, for all $\alpha \leq \epsilon$ the level
set $L_{\alpha}$ is countable, and for each $p \in L_{\epsilon}$ the LOTS $S_{p}$ is separable.
\item[(v)]  $T^{\epsilon +1}$ is a countable tree and for each $p \in L_{\epsilon}$ the LOTS $S_{p}$ is separable.
\end{itemize}
\item[(b)]  The following conditions are equivalent.
\begin{itemize}
\item[(i)]  $X(T)$ satisfies the countable chain condition.
\item[(ii)]  $\widehat{X(T)}$ satisfies the countable chain condition.
\item[(iii)]  The tree  $T^{\epsilon +1}$ satisfies condition (ix) of Definition \ref{df5.1},
i.e. every anti-chain is countable, and for each $p \in L_{\epsilon}$ the LOTS $S_{p}$
satisfies the countable chain condition.
\end{itemize}
These conditions imply that  $T^{\epsilon +1}$ satisfies
conditions (vi) and (viii) of Definition \ref{df5.1} and, in particular, that $h(T) \leq
\Omega$. If $h(T) = \Omega$, then $T$ is a Suslin tree.
\item[(c)]  The following conditions are equivalent.
\begin{itemize}
\item[(i)]  $X(T)$ is of countable type.
\item[(ii)]  $\widehat{X(T)}$ is first countable and $\sigma$-compact.
\item[(iii)]  The tree $T$ is of countable type and is  $\Omega$ bounded, i.e.  $h(x) < \Omega$ for all $x \in X(T)$
(condition (vi) of Definition \ref{df5.1}).
\item[(iv)]  The tree $T$ is of countable type and  $h(x) < \Omega$ for all $x \in X(\hat{T})$.
\end{itemize}
\end{enumerate}
\end{prop}

\begin{proof} When $h(T) = \epsilon + 2$, let
$S_{x} = S_{p}$ for $x = x(p)$ with $o(p) = \epsilon$ and let
$S_{x} = \{x\} $ for $x \in X(T^{\epsilon +1})\setminus
L_{\epsilon}$.  Each branch in $X(T)$ extends its projection $ x
\in X(T^{\epsilon +1})$ by a point in $S_{x}$.  Thus, in this case
we have the order sum isomorphism
\begin{equation}\label{eq5.41}
X(T) \quad \cong \quad \Sigma \{S_{x} : x \in X(T^{\epsilon +1})
\}.
\end{equation}

(b)  (i)$\Leftrightarrow$(ii): By Proposition \ref{prop2.1}(b) the inclusion of
$X(T)$ into $\widehat{X(T)}$ is continuous and so the equivalence
follows from Proposition \ref{prop2.5}(f),(h).

(iii)$\Rightarrow$(i): Let $\{ J_{i} :i \in I \}$ be a pairwise
disjoint family of nonempty open intervals in $X(T)$.  By (\ref{eq5.36})
there exists a vertex $q_{i} \in T$ such that every branch through
$q_{i}$ is contained in $J_{i}$.  If $o(q_{i}) < \epsilon +1 $
then let $p(q_{i}) = q_{i}$.  If $o(q_{i}) = \epsilon +1 $ then
let $p(q_{i})$ be the immediate predecessor of $q_{i}$.  In that
case, $J_{i}$ meets $S_{p(q_{i})}$ which satisfies c.c.c.
Consequently, for each $p \in L_{\epsilon}$ the set $\{ i \in I :
p(q_{i}) = p \}$ is countable.  The set $\{p(q_{i}) : i \in I \}$
is an anti-chain in $T^{\epsilon +1}$ since the intervals $J_{i}$ are
disjoint. So by assumption this set is countable and consequently
$I$ itself is countable.

(i)$\Rightarrow$(iii):  By (\ref{eq5.41}) each $S_{p}$ for $p \in
L_{\epsilon}$ is isomorphic to a subinterval of $X(T)$.  So if
$X(T)$ satisfies c.c.c. then each such $S_{p}$ does. If $A$ is an
anti-chain in $T^{\epsilon +1}$ then $o(p) + 1 < h(T)$  implies
that $\{ [(\pi^{o(p)+1})^{-1}(x(p))]^{\circ} : p \in A \}$ is a family of
nonempty open intervals in $X(T)$ by Lemma \ref{lem5.8}.  The family is
pairwise disjoint since $A$ is an anti-chain.  Because $X(T)$
satisfies c.c.c. $A$ is countable.

We proved after Definition \ref{df5.1} that condition (ix) implies (vi)
and (viii). Thus $h(X(T^{\ep + 1}))\leq \Omega$. So if $h(T)$ is not a limit ordinal, then
$\ep + 1 < \Omega$ and so $h(T) = \ep + 2 < \Omega$.

(a) (i)$\Rightarrow$(ii):  $X(T)$ is dense in $\widehat{X(T)}$.

(ii)$\Rightarrow$(iii):  This follows from Proposition \ref{prop2.8}(b) since
$\widehat{X(T)}$ is separable, connected and is unbounded.

(iii)$\Rightarrow$(i): By Proposition \ref{prop2.1}(b) the topology on $X(T)$ is
inherited from $\widehat{X(T)}$.  Any subset of $\R$ is
second countable and so is separable.

(i)$\Rightarrow$(iv):  Since separability implies c.c.c. the results
of (b) can be applied and so $T^{\epsilon +1}$ satisfies (vi) and
(viii) of Definition \ref{df5.1}.  In particular, for each $\alpha \leq
\epsilon$ the set $L_{\alpha}$ is countable. Let $D$ be a
countable dense subset of $X(T)$.  If $J$ is any nontrivial
interval in $X(T)$, then $J^{\circ} \cap D$ is dense in $J^{\circ}$ and so
in $J$.  By (\ref{eq5.41}) $S_{p}$ is separable for each $p \in
L_{\epsilon}$.  Now let $\alpha =  sup \{ h(x) : x \in D \}$.  By
condition (vi) each $h(x)$ is countable and so $\alpha < \Omega$.
I claim that $h(T) \leq \alpha + 1$,  for if $\alpha + 1 < h(T)$
and $o(p) = \alpha$, then by Lemma \ref{lem5.8} $[\pi^{\alpha +1})^{-1}(x(p))]^{\circ}$
is a nonempty open interval and so contains some
point $y \in D$.  But then $p \in y$ implies $h(y) \geq \alpha +1$
which contradicts the definition of $\alpha$.

(iv)$\Rightarrow$(v): $ T^{\epsilon +1}$ has countably many levels
and each level is countable.

(v)$\Rightarrow$(i):  If $h(T) = \epsilon +2$, then we choose for
each $p$ in the countable set $ L_{\epsilon}$ a countable dense
subset of $S_{p}$. The union is a countable set which is dense in
$O_{\epsilon +1}$
(see Equation (\ref{eq5.38})) which is dense in $X(T)$ by Lemma \ref{lem5.8}. If $h(T)
= \epsilon$ is a limit ordinal, then $T = T^{\epsilon +1}$ is
countable.  Choose a branch through each vertex to get a countable
set which is dense in $X(T)$ by Proposition \ref{prop5.3}(g).

(c) (i)$\Leftrightarrow$(ii): By Proposition \ref{prop2.8}(d).

(i)$\Rightarrow$(iii): Assume that $f : \Omega \rightarrow S_{p}$ is
an order preserving or reversing injection for some $p \in T$.
For each $i \in \Omega$ choose $\tilde{f}(i)$ a branch through
$f(i)$.  Then $\tilde{f} : \Omega \rightarrow X(T)$ is an
injection which similarly preserves or reverses order.  Hence,
$X(T)$ is not of countable type.

On the other hand, suppose that $x \in X(T)$ with $h(x) \geq
\Omega$.  Define
\begin{equation}\label{eq5.42}
\begin{split}
K_{+} \quad = \quad \{ j \in \Omega : x _{j+1} \not=  min S_{x_{j}} \} \   \\
K_{-} \quad = \quad \{ j \in \Omega : x _{j+1} \not=  max S_{x_{j}} \}.
\end{split}
\end{equation}
Either $K_{+}$ or $K_{-}$ is uncountable since $max \not= min $
for any $S_{p}$.  Suppose $K_{+}$ is uncountable.  For each $i \in
K_{+}$ choose $f(i) \in X(T)$ such that
\begin{equation}\label{eq5.43}
f(i)_{i+1} \in S_{x_{i}} \qquad \mbox{with} \qquad f(i)_{i+1} <
x_{i+1}.
\end{equation}
If $i < j$ in $K_{+}$, then the equality level for $f(i),f(j)$ is
$i$ and
\begin{equation}\label{eq5.44}
f(i)_{i+1} \quad < \quad x_{i+1} \quad = \quad f(j)_{i+1}.
\end{equation}
Thus, $f : K_{+} \rightarrow X(T)$ is an order injection.  Since
$K_{+} \cong \Omega$, $X(T)$ is not of countable type.

(iii)$\Leftrightarrow$(iv):  It is clear from Proposition \ref{prop5.6} that
$h(x) < \Omega$ for all $x \in X(T)$ implies $h(x) < \Omega$ for
all $x \in X(\hat{T})$ and the converse is obvious.

(iii)$\Rightarrow$(i):  If $X(T)$ is not of countable type, then there
exists an injective map $f : \Omega \rightarrow X(T)$ which we can
assume without loss of generality to be order preserving.

We now prove by induction on $\alpha \in \Omega $ that there exist
$\epsilon(\alpha) \in \Omega$ and $p(\alpha) \in L_{\alpha}$ such
that
\begin{equation}\label{eq5.45}
\pi_{\alpha +1} \circ f (\beta) \quad = \quad x(p(\alpha)) \qquad
\mbox{for all}\  \beta \in \Omega \setminus \epsilon(\alpha).
\end{equation}

If $\alpha = 0$, then $\pi_{1} \circ f$ is constantly the root $0$.
Let $\epsilon(0) = 0$ and $p(0) = 0$.

Now assume that for all $i < \alpha$, $\epsilon(i) \in \Omega$ and
$p(i) \in L_{i}$ have been defined so that (\ref{eq5.45}) holds with
$\alpha$ replaced by $i$.\vspace{.25cm}

\textbf{Case 1:} If $\alpha = \tilde{\alpha} +1$, then $\pi_{\tilde{\alpha}
+1} \circ f $ is constant on the tail $\Omega \setminus
\epsilon(\tilde{\alpha})$ with value $x(p(\tilde{\alpha}))$.
Hence, if $j < k $ in the tail, then the equality level of $f(j)$
and $f(k)$ is at least $\tilde{\alpha}$.  Hence, $ f(j) < f(k) $
in $X(T)$ implies  $f(j)_{\alpha} \leq f(k)_{\alpha}$.  Thus, $j
\mapsto f(j)_{\alpha}$ defines an order map from $\Omega \setminus
\epsilon(\tilde{\alpha}) \cong \Omega$ to the LOTS of countable
type $S_{p}$ with $p = p(\tilde{\alpha})$.  By Corollary \ref{lem2.6} this map
is eventually constant.  Hence there exists $\epsilon(\alpha) \geq
\epsilon(\tilde{\alpha})$ such that for all $j \in \Omega
\setminus \epsilon(\alpha)$  $f(j)_{\alpha}$ is a common vertex
$p(\alpha) \in S_{p(\tilde{\alpha})}$ and so (\ref{eq5.45}) holds.\vspace{.25cm}

\textbf{Case 2:} If $\alpha$ is a limit ordinal, then define
$\epsilon(\alpha) = sup \{ \epsilon(i) : i < \alpha \}$.  For
$\beta > \epsilon(\alpha)$, $f(\beta)$ and $f(\epsilon(\alpha))$
are distinct points with $f(\beta)_{i} = f(\epsilon(\alpha))_{i}$
for all $i < \alpha$. Hence, the equality level is at least
$\alpha$.  In particular, $f(\beta)_{\alpha} =
f(\epsilon(\alpha))_{\alpha}$.  Let $p(\alpha)$ be this common
vertex.  Again (\ref{eq5.45}) holds.\vspace{.25cm}

Having proved (\ref{eq5.45}) for all $\alpha$ we can restate it as
\begin{equation}\label{eq5.46}
f(\beta)_{\alpha} \quad = \quad p(\alpha) \qquad \mbox{for all} \
\beta \in \Omega \setminus \epsilon(\alpha).
\end{equation}
It follows that if $\tilde{\alpha} < \alpha$ and $\beta \geq
\epsilon(\alpha)$, then $p(\tilde{\alpha})$ and $p(\alpha)$ both
lie on the branch $f(\beta)$.  Hence, $p(\tilde{\alpha}) \prec
p(\alpha)$.  Thus, $\{ p(\alpha) : \alpha \in \Omega \}$ is a
linearly ordered collection of vertices .  Hence, it is contained
in some branch $x$ and any such branch satisfies $h(x) \geq
\Omega$. This completes the proof that (iii) $\Rightarrow$ (i).

\end{proof}\vspace{.5cm}

\noindent{\bfseries Remark.}  In the case excluded by the hypothesis, $h(T)
= \epsilon +1 $ with $\epsilon$ a limit ordinal, $\pi_{\epsilon} :
X(T) \rightarrow X(T^{\epsilon})$ is an order isomorphism and
$h(T^{\epsilon}) = \epsilon$.  We obtain results for $T$ in this
case by applying the proposition to $T^{\epsilon}$. \vspace{.5cm}

\begin{cor}\label{cor5.10}  If $T$ is a bi-ordered Aronszajn tree of $\Q$ type, then
$X(T)$  is an order dense LOTS of countable type with dense holes and
$\widehat{X(T)}$  is a $\sigma$-compact, connected, first countable LOTS. Neither is
separable.  For each $1 < \alpha < \Omega $, $\widehat{X(T^{\alpha})} \cong \R$.
If $T$ is a Suslin tree, then $\widehat{X(T)}$ satisfies the countable chain condition.
\end{cor}

\begin{proof}  $X(T)$ is order dense, unbounded and with dense holes by Proposition \ref{prop5.3}(d) and (e).
It is of countable type by
Proposition \ref{prop5.9}(c). Hence, its completion $\widehat{X(T)}$ is $\sigma$-compact,
connected  and  first countable.  Because $h(T) = \Omega$, it is not separable by
Proposition \ref{prop5.9}(a) and the same result implies that for countable
$\alpha$, $\widehat{X(T^{\alpha})} \cong \R $.  By
Proposition \ref{prop5.9}(b), $\widehat{X(T)}$ for a Suslin tree satisfies c.c.c.

\end{proof} \vspace{1cm}

\subsection{Homogeneous and Reproductive Trees}

If $T_{1}$ and $T_{2}$ are bi-ordered trees, then an
\emph{isomorphism}\index{isomorphism}\index{tree!isomorphism} $f : T_{1} \rightarrow T_{2} $ is a bijection
which preserves both orders, i.e.
\begin{equation}\label{eq5.47}
\begin{split}
p \prec _{1} q \qquad \Longleftrightarrow \qquad f(p) \prec _{2} f(q) \hspace{1.5cm}  \\
q <_{1} r \quad \mbox{in} \ S_{p} \qquad \Longleftrightarrow
\qquad f(q) <_{2} f(r) \quad \mbox{in} \ S_{f(p)}.
\end{split}
\end{equation}
The first condition says that $f$ relates the tree structures. Hence, it
 maps the vertices of level $\alpha$ in $T_{1}$ to those of
level $\alpha$ in $T_{2}$ and also $f$ induces a bijection of branch
spaces, denoted $f_{*} : X(T_{1}) \rightarrow X(T_{2}) $. The
second condition implies that $f_{*}$ is an order isomorphism with
respect to the induced orders and so, in the order dense case, extends to an order
isomorphism on the completions, which we will also denote by
$f_{*}$.  When $T_{1} = T_{2}$ such an isomorphism is called an
\emph{automorphism}\index{automorphism}\index{tree!automorphism}.

\begin{df}\label{df5.11}  A bi-ordered tree $T$ is called \emph{homogeneous}\index{homogeneous tree}\index{tree!homogeneous} if for all
$p,q \in T$ such that $o(p) = o(q)$ there exists an automorphism $f$ of $T$ such that $f(p) = q$.

A bi-ordered tree $T$ is called \emph{reproductive}\index{reproductive tree}\index{tree!reproductive} if for all $p
\in T$ the tail tree $T_{p}$ is isomorphic to $T$.
\end{df}
\vspace{.5cm}

\begin{theo}\label{theo5.12}  If $T$ is a homogeneous, bi-ordered Aronszajn tree of dense type, then
 $\widehat{X(T)}$ is a nonseparable CHLOTS.
\end{theo}

\begin{proof} Homogeneity clearly implies that
no $S_{p}$ has a $max$ or $min$ and so $T$ is of $\Q$
type.  By Corollary \ref{cor5.10} $\widehat{X(T)}$ is first countable,
$\sigma$-compact and nonseparable.  By Proposition \ref{prop3.4}(a) it
suffices to prove that $\widehat{X(T)}$ is doubly transitive.  Given $
x < y$ and $z < w$  in $\widehat{X(T)}$  we construct an order
isomorphism on $\widehat{X(T)}$ which maps the pair $x,y $ to $z,w$.

Because $T$ is an Aronszajn tree, it satisfies condition (vi) of
Definition \ref{df5.1} and so every height $h < \Omega$ on $\widehat{X(T)} = X(\hat{T}) $.
Choose a countable ordinal $\alpha  > h(x),h(y),h(z),h(w) $ so
that
\begin{equation}\label{eq5.48}
\hat{\pi}_{\alpha +1}(r) \  \in \ \widehat{X(T^{\alpha +1})} \setminus L_{\alpha}
\qquad \mbox{for} \ r = x,y,z,w.
\end{equation}
By Lemma \ref{lem5.8} $\hat{\pi}_{\alpha +1}$ is injective on the
complement of  $(\hat{\pi}_{\alpha +1})^{-1}(L_{\alpha})$ we have
\begin{equation}\label{eq5.49}
\hat{\pi}_{\alpha +1}(x) \  < \  \hat{\pi}_{\alpha +1}(y) \quad
\mbox{and} \quad \hat{\pi}_{\alpha +1}(z) \ < \ \hat{\pi}_{\alpha
+1}(w).
\end{equation}

By Corollary \ref{cor5.10}  $\widehat{X(T^{\alpha})} \cong \R$ and by
condition (vii) or Definition \ref{df5.1} and Lemma \ref{lem5.4}(a), $L_{\alpha}$ is
a countable dense subset of $\widehat{X(T^{\alpha})}$. We can choose
an order isomorphism of $L_{\alpha}$ which maps the convex set
$L_{\alpha} \cap (\hat{\pi}_{\alpha +1}(x),\hat{\pi}_{\alpha
+1}(y))$ to $L_{\alpha} \cap (\hat{\pi}_{\alpha
+1}(z),\hat{\pi}_{\alpha +1}(w))$ because $L_{\alpha}$ is order
isomorphic to the IHLOTS $\Q$. Extending  to the
completion we obtain an order isomorphism $\tilde{f} :
\widehat{X(T^{\alpha})} \rightarrow \widehat{X(T^{\alpha})}$ such that
\begin{equation}\label{eq5.50}
\tilde{f}(L_{\alpha}) = L_{\alpha}, \quad
\tilde{f}(\hat{\pi}_{\alpha +1}(x)) = \hat{\pi}_{\alpha +1}(y),
\quad  \tilde{f}(\hat{\pi}_{\alpha +1}(z)) = \hat{\pi}_{\alpha
+1}(w).
\end{equation}

We can regard the LOTS $\widehat{X(T)}$ as the order space sum
\begin{equation}\label{eq5.51}
\widehat{X(T)} \  = \  \Sigma \{ (\hat{\pi}_{\alpha +1})^{-1}(p) : p
\in  \widehat{X(T^{\alpha +1})} \},
\end{equation}
where for $ p \not\in L_{\alpha} $, $(\hat{\pi}_{\alpha
+1})^{-1}(p)$ is a single point and for $ p \in L_{\alpha} $ it is
a closed interval with nonempty interior, consisting of the
branches of $\hat{T}$ which contain $p$.

For each $p \in L_{\alpha}$ choose an automorphism $\tilde{g}_{p}$
of $T$ which maps the vertex $p$ to the vertex $\tilde{f}(p)$.  We
can restrict the order isomorphism $(\tilde{g}_{p})_{*}$ on
$\hat{X}(T)$ to define for $p \in L_{\alpha}$ an order isomorphism
\begin{equation}\label{eq5.52}
g_{p} : (\hat{\pi}_{\alpha +1})^{-1}(p) \rightarrow
(\hat{\pi}_{\alpha +1})^{-1}(\tilde{f}(p))
\end{equation}
and for $p \not\in L_{\alpha}$ let $g_{p}$ denote the unique map
between the singletons.

The required isomorphism on $\widehat{X(T)}$ is
\begin{equation}\label{eq5.53}
g \  = \  \Sigma \{ g_{p} : p \in \widehat{X(T^{\alpha})} \}.
\end{equation}

\end{proof} \vspace{.5cm}

\begin{prop}\label{prop5.13}  Let $T$ be a bi-ordered Aronszajn tree of $\Q$ type.
\begin{enumerate}
\item[(a)]  Let $g : Y \rightarrow \widehat{X(T)} $  be a continuous, injective map
with $Y$ an arbitrary separable topological space.  There exists a countable ordinal $\beta$
such that $\hat{\pi}_{\beta +1} \circ g : Y \rightarrow \widehat{X(T^{\beta +1})} $ is
injective.  In particular, any separable space which can be continuously injected
into $\widehat{X(T)}$ can be continuously injected as well into $\R$.
\item[(b)]  Let the LOTS $X_{1}$ be an uncountable, dense subset of $\R$ .
The AS double $X_{1}'$ does not inject into $\widehat{X(T)}$.
\end{enumerate}
\end{prop}

\begin{proof}   (a): Let $D$ be a countable dense
subset of $Y$.  Since $h < \Omega $ on $\widehat{X(T)} = X(\hat{T})$,
\begin{equation}\label{eq5.54}
\alpha \quad =  \quad  sup  \{ h(g(y)) : y \in D \}
\end{equation}
is a countable ordinal and we have
\begin{equation}\label{eq5.55}
\hat{\pi}_{\alpha +1}(g(y)) \ \in \ \widehat{X(T^{\alpha +1})}
\setminus L_{\alpha} \qquad \mbox{for} \ y \in D.
\end{equation}
In the notation of Lemma \ref{lem5.7}, this implies that $g(D)$ is disjoint
from the open set $O_{\alpha + 1}$.  Because $D$ is dense in $Y$
and $g$ is continuous it follows that $g(Y)$ is disjoint from
$O_{\alpha + 1}$.  By Lemma \ref{lem5.7} the restriction of
$\hat{\pi}_{\alpha +2}$ to $\widehat{X(T)} \setminus O_{\alpha + 1}$
is injective and so the result follows with $\beta = \alpha +2$.

(b):  We will assume that $G_{1} : X_{1}' \rightarrow \widehat{X(T)}$
is an order injection and use it to construct a separable, compact
nonmetrizable subset $C$ of $\widehat{X(T)}$.  Since a compact space
which continuously injects into $\R$ is second countable
this will contradict part (a).

We will use the fact that a second countable, compact space has only countably many clopen sets.  If $\B$ is a countable basis, then any
clopen set $U$ is a union of members of $\B$. Since $U$ is closed, and hence compact, it is a union of finitely many members of $\B$. So the cardinality
of the set of clopens is bounded by that of the collection of finite subsets of $\B$ which is countable.

The LOTS $X_{1}$ is order dense.  Let $D$ be a countable, dense
subset of $X_{1}$ which we can identify with the countable dense
subset $D \times \{-1\}$ of $X_{1}'$ so define the order injection
$g : D \rightarrow \widehat{X(T)}$ by $g(t) = G_{1}(t^{-})$. Use $g$
to define $G : X_{1}' \rightarrow \widehat{X(T)}$ by using (\ref{eq3.20}).
Because $G_{1}$ is injective, we have for any $x \in X_{1}$
\begin{equation}\label{eq5.56}
G(x^{-})\quad  \leq \quad  G_{1}(x^{-})\quad  < \quad
G_{1}(x^{+})\quad  \leq \quad G(x^{+}).
\end{equation}
Combined with (\ref{eq3.25}) we see that the continuous order map $G$ is
injective.  Choose $a < b \in X_{1}$ and let $C$ be the closure in
$\widehat{X(T)}$ of the image $G([a^{+},b^{-}])$.  For any $x \in
(a,b)$ (\ref{eq5.56}) implies that  $G(x^{-}) < G(x^{+})$ is a gap pair in
the image and so in $C$.  Since $X_{1}'$ is separable and $G$ is
continuous, $C$ is separable.  Since $\widehat{X(T)}$ is complete and
$C$ is bounded, it is compact.  Since it has uncountably many gap
pairs, it has uncountably many clopen subsets and so is not
metrizable.

\end{proof}\vspace{.5cm}

For the following recall from (\ref{eq4.1}) the definition of $Y_{\a}$ for a nontrivial connected LOTS $Y$.

\begin{cor}\label{cor5.14}  Let $T$ be a bi-ordered Aronszajn tree of $\Q$ type,  and $\alpha > 1$ be an ordinal.
\begin{enumerate}
\item[(a)] If $Y$ is a nontrivial, connected LOTS, then $Y_{\alpha}$ does not inject into $\widehat{X(T)}$.
\item[(b)] If a LOTS $Y$ is an uncountable, dense subset of $\R$, then
$\widehat{Y_{\alpha}}$ does not inject into $\widehat{X(T)}$.
\end{enumerate}
\end{cor}

\begin{proof} (a): By Theorem \ref{theo3.8}, $Y$ contains a
compact subset $A$ such that $A'' \subset A'$ is order isomorphic
to the Fat Cantor Set $\bullet \R' \bullet$.  Since
$\alpha > 1$, $\hat{j}_{\alpha}^{1} : Y' \rightarrow Y_{\alpha} $
restricts to an order embedding of $A''$ into $Y_{\alpha}$.  So an
order injection of $ Y_{\alpha}$  into $ \widehat{X(T)}$ would
restrict to an order injection of $A''$ into $ \widehat{X(T)}$. But by
Proposition \ref{prop5.13}(b) the Fat Cantor Set cannot be order injected into
$\widehat{X(T)}$.

(b):  Similarly, an order injection of $\widehat{Y_{\alpha}}$ into
$\widehat{X(T)}$ would restrict to an order injection of $Y'$ into
$\widehat{X(T)}$ which again contradicts Proposition \ref{prop5.13}(b).

\end{proof} \vspace{.5cm}

\begin{lem}\label{lem5.15}  Let $T$ be a bi-ordered normal tree, $p \in T$ and $\alpha$ a positive ordinal with $\alpha < h(T)$.
\begin{enumerate}
\item[(a)]  Assume $T$ is homogeneous.  If $o(p) + 1 < h(T)$, then $S_{p}$ is a transitive LOTS and so is unbounded,
i.e. $T$ is of unbounded type, and the tail tree $T_p$ is homogeneous.
  If $q \in T$ with $o(q) = o(p)$, then as LOTS
$S_{p} \cong S_{q}$.  The subtree $T^{\alpha}$ is homogeneous.  Regarded as a subset
of $X(T^{\alpha +1})$, $L_{\alpha}$ is a transitive LOTS.
\item[(b)]  Assume $T$ is reproductive.  The height $h(T)$ is a tail-like ordinal.
As LOTS $S_{p} \cong S_{0}$ where $0$ is the root of $T$.  Hence, with $Y \cong S_{0}$
 $T$ is of $Y$ type.   The tail tree $T_{p}$ is reproductive and if $\alpha$ is
tail-like, then $T^{\alpha}$ is reproductive.
\item[(c)]  Assume $T$ is homogeneous and reproductive.  The tail tree $T_p$ is
homogeneous and reproductive and if $\alpha$ is tail-like, then $T^{\alpha}$ is homogeneous and reproductive.
\end{enumerate}
\end{lem}

\begin{proof}(a): If $q_{1}, q_{2} \in S_{p}$,
then $o(q_{1}) = o(q_{2}) = o(p) +1$ and so there exists an
automorphism $f$ of $T$ such that $f(q_{1}) = q_{2}$. Since $f$
preserves $\prec$, $f(S_{p}) = S_{p}$ and so $f$ restricts to an
order automorphism on $S_{p}$.  Similarly, if $f(p) = q$ then $f$
restricts to an order isomorphism $S_{p} \cong S_{q}$.

In general, if $q_1, q_2 \in T_p$ with $o_p(q_1) = o_p(q_2)$, then
$o(q_1) = o(q_2)$ and so there exists an
automorphism $f$ of $T$ such that $f(q_{1}) = q_{2}$. Since $f$
preserves $\prec$, $f$ fixes $p$ and so maps $T_p$ to itself. Hence, $T_p$ is homogeneous.

 Since any
automorphism of $T$ preserves $T^{\alpha}$, the latter is
homogeneous. Using the induced isomorphisms on $X(T^{\alpha +1})$,
which preserve $L_{\alpha}$, we see that $L_{\alpha}$ is
transitive as well.

(b): An isomorphism of $T_{p}$ with $T$ restricts to an isomorphism
of $S_{p}$ with $S_{0}$.  If $h(T) = 1$, then $T = \{0\}$ and
$S_{0} = \emptyset $.  Otherwise, all successor sets are of
$S_{0}$ type and $h(T)$ is an infinite limit ordinal.
Furthermore, by (\ref{eq5.10})
\begin{equation}\label{eq5.57}
h(T) = h(T_{p}) \cong h(T) \setminus o(p).
\end{equation}
Thus, $h(T)$ is tail-like.

Now assume that $\alpha$ is an infinite tail-like ordinal (the
case $\alpha = 1$ is trivial). If $o(p) < \alpha  \leq \beta$, then
(\ref{eq5.8}) implies that $o(q) < \beta$ iff $o_{p}(q) < \beta$.  This
shows that for $\alpha$ tail-like:
\begin{equation}\label{eq5.58}
o(p)\ < \ \alpha \ \leq \ \beta \quad \Rightarrow \quad
(T^{\beta})_{p} \ = \ (T_{p})^{\beta}.
\end{equation}
In particular, if $p \in T^{\alpha}$ then an isomorphism from $T$
to $T_{p}$ restricts to an isomorphism from $T^{\alpha}$ to
$(T^{\alpha})_{p}$ which shows that $T^{\alpha}$ is reproductive.

Since $T_{p}$ is isomorphic to $T$ it is reproductive.

(c):  If $\alpha$ is tail-like, then $T^{\alpha}$ is homogeneous and
reproductive by (a) and (b).  Since  $T_{p}$ is isomorphic to $T$
it is homogeneous and reproductive.

\end{proof} \vspace{.5cm}

\begin{prop}\label{prop5.16} If $T$ is a bi-ordered tree which is reproductive and with $S_{0}$ a
transitive LOTS, then $T$ is a normal, homogeneous tree of unbounded type with height an infinite, tail-like ordinal.
Furthermore, if $x,y \in X(T)$ with $h(x) = h(y)$, then there exists an automorphism $f$ of $T$ such that $f_{*}(x) = y$.
\end{prop}

\begin{proof} Since $S_{0}$ is a transitive LOTS
it is nonempty.  By condition (ii) of Definition \ref{df5.1} it contains
at least two points.  As the two point LOTS is not transitive,
Proposition \ref{prop3.2}(b) implies that $S_{0}$ is infinite with no $max$ or
$min$.  For every $ p \in T$,  there exists an isomorphism $ f: T
\rightarrow T_{p}$ and so $S_{p} = f(S_{0})$ is infinite.
Furthermore, for every ordinal $\alpha < h(T)$, there exists $q
\in T$ with $o(q) = \alpha$ and so $f(q) \in T_{p}$ with $o(f(q))
= o(p) + o(q)$.  Thus, condition (iv) of Definition \ref{df5.1}
holds and so $T$ is a normal tree.  Since $T$ is reproductive,
$h(T)$ is tail-like and since $S_{0} \not= \emptyset $, $h(T) >
1$. Hence, $h(T)$ is infinite.

We now have to show that if $p,q \in T$ with $o(p) = o(q)$ then
there exists an automorphism $f$ of $T$ such that $f(p) = q$. We
will call this the \emph{vertex case} and use $x,y$ to stand for
the linearly ordered sets $x =  \{p\} \cup A_{p}$ and $y =  \{q\}
\cup A_{q}$.  Let $\alpha = h(x) = h(y) = o(p) + 1$ so that in the
vertex case $\alpha$ is a successor.  When $x$ and $y$ are
branches, that is the \emph{branch case},  then $\alpha = h(x) =
h(y)$ is a limit ordinal.  In the vertex case, let $T_{x} = T_{p}$
and $T_{y} = T_{q}$.  In the branch case, let $T_{x} = T_{y} =
\emptyset $.  In either case,

\begin{equation}\label{eq5.59}
\begin{split}
T_{x} \quad = \quad \{ r \in T : x \subset \{r\} \cup A_{r} \}   \\
T_{y} \quad = \quad \{ r \in T : y \subset \{r\} \cup A_{r} \}.
\end{split}
\end{equation}

For any $r \in T$ we can define by analogy with (\ref{eq5.15}) the $r,y$
equality level to be $\epsilon = \epsilon(r,y)$ where
\begin{equation}\label{eq5.60}
y   \cap  (\{r\} \cup A_{r}) \quad \cong \quad \epsilon + 1.
\end{equation}
So $\ep = \alpha$ if $r \in T_{y}$ and $\ep < \alpha$
otherwise.  Similarly, define $ \epsilon(r,x)$.

Since the results are obvious for the trivial tree we will assume
that $h(T)$ is an infinite, tail-like ordinal and $S_{r} \cong
S_{0}$ is infinite for each $r \in T$.

By induction on $\b
\leq \alpha$ we will construct automorphisms $f_{\b}$ of $T$
which satisfy:
\begin{equation}\label{eq5.61}
\delta \leq \b \quad \Rightarrow \quad
f_{\b}(y_{\delta}) = x_{\delta},
\end{equation}
and for all $r \in T $
\begin{equation}\label{eq5.62}
\epsilon(r,y)   \leq   \delta  <  \b \quad  \Rightarrow
\quad   f_{\delta}(r) = f_{\b}(r)
\end{equation}
In the final step, when $\b = \alpha$, the inequality in
(\ref{eq5.61}) is replaced by the strict inequality $\delta < \alpha$
since $y_{\alpha} $ and $x_{\alpha}$ are not defined.  The final
automorphism $f_{\alpha}$ is the one which maps $x$ to $y$ (and
$p$ to $q$ in the vertex case).

Begin with $f_{0} = 1_{T}$, the identity. Condition (\ref{eq5.61}) holds
since $y_{0} = x_{0} = 0 $ and condition (\ref{eq5.62}) holds vacuously.
Now assume that $f_{\delta}$ has been constructed for all $\delta
< \b$.\vspace{.25cm}

\textbf{Case 1: } If $\b = \delta + 1$, then the automorphism
$f_{\delta}$ maps $y_{\delta}$ to $x_{\delta}$ and so
$x_{\b}$ and $f_{\delta}(y_{\b})$ both lie in the
transitive LOTS $S_{p}$ with $p = x_{\delta}$.  Choose $g$ a LOTS
automorphism of $S_{p}$ which maps $f_{\delta}(y_{\b})$ to
$x_{\b}$.

Because $T$ is reproductive we can choose for each $r \in S_{p}$ a
tree isomorphism $g_{r} : T_{r} \rightarrow T_{g(r)}$.  We define
the tree automorphism $\tilde{g}$ to be $g_{r}$ on each such
$T_{r}$ and to be the identity on the rest of $T$. Then define
$f_{\b} \ = \ \tilde{g} \circ f_{\delta}$.  By construction
(\ref{eq5.61}) holds and if $\epsilon(r,y) < \b$ then
$f_{\delta}(r)$ does not lie in $T_{p} \setminus \{p\}$ and so
$\tilde{g}$ is the identity on $f_{\delta}(r)$ from which (\ref{eq5.62})
follows.\vspace{.25cm}

\textbf{Case 2:}  If $\b$ is a limit ordinal, then define
$f_{\b}(r) = f_{\delta}(r)$ whenever $\epsilon(r,y) \leq
\delta < \b $.  By (\ref{eq5.62}) this definition is independent of
the choice of $\delta$ and defines $f_{\b}$ whenever
$\epsilon(r,y) < \b$.  Furthermore, by (\ref{eq5.62})
$f_{\b}(y_{\delta}) = x_{\delta} $ for all $\delta <
\b$.

If $\b = \alpha$, then, since $\b$ is a limit ordinal,
we are in the branch case with $T_{y} = T_{x} = \emptyset$.  So
$f_{\b}$ is defined on all of $T$ and satisfies (\ref{eq5.62}) and
the adjusted version of (\ref{eq5.61}).

If $\b < \alpha$, then $f_{\b}$ maps the predecessors
of $y_{\b}$ to the corresponding predecessors of
$x_{\b}$. Because $T$ is reproductive we can finish our
definition of $f_{\b}$ by choosing an isomorphism
$f_{\b} : T_{y_{\b}} \rightarrow T_{x_{\b}}$.
This tree isomorphism maps the root $y_{\b}$ to the root
$x_{\b}$ and so (\ref{eq5.61}) holds.\vspace{.25cm}

This finishes the inductive construction and so completes the
proof.

\end{proof}\vspace{.5cm}

\begin{lem}\label{lem5.17}  Let $T$ be a reproductive, bi-ordered tree with $S_{0}$ a doubly
transitive LOTS containing at least three points.
\begin{enumerate}
\item[(a)]  Let $\alpha = h(T)$.  If $x,y,z,w \in \tilde{L}_{\alpha}$ with $x < y $
and $ z < w $, then there exists an order isomorphism $k : (x,y) \rightarrow (z,w) $ (intervals in $X(T)$) such that
\begin{equation}\label{eq5.63}
k((x,y) \cap \tilde{L}_{\alpha}) \quad = \quad (z,w) \cap
\tilde{L}_{\alpha}.
\end{equation}
\item[(b)]  Let $\alpha$  be an infinite, tail-like ordinal with $\alpha < h(T)$.
 If $x,y,z,w \in L_{\alpha}$ with $x < y $ and $ z < w$, then there exists an order
 isomorphism $k : (x,y) \rightarrow (z,w) $ (intervals in $X(T^{\alpha +1})$) such that
\begin{equation}\label{eq5.64}
k((x,y) \cap L_{\alpha}) \quad = \quad (z,w) \cap L_{\alpha}.
\end{equation}
\end{enumerate}
\end{lem}

\begin{proof} (a): With $\beta$ and $\epsilon$
the $x,y$ and $z,w$ equality levels, respectively, we let
\begin{equation}\label{eq5.65}
p  \  = \   x_{\beta}\  = \  y_{\beta} \qquad \mbox{and} \qquad q
\  = \ z_{\epsilon}\  = \  w_{\epsilon}.
\end{equation}
Note that $x < y$ and $z < w$ imply
\begin{equation}\label{eq5.66}
\begin{split}
\beta, \epsilon \ < \ \alpha \hspace{2cm}  \\
x_{\beta +1} \ < \ y_{\beta +1} \qquad \mbox{in} \ S_{p}  \\
z_{\epsilon +1} \ < \ w_{\epsilon +1} \qquad \mbox{in} \  S_{q}.
\end{split}
\end{equation}

Because $T$ is reproductive there exist isomorphisms $g_{1} :
T_{p} \rightarrow T$ and $g_{2} : T_{q} \rightarrow T$. Define
\begin{equation}\label{eq5.67}
\begin{split}
\tilde{x}\ = \ g_{1}(x \cap T_{p}) \qquad  \tilde{y}\ = \ g_{1}(y \cap T_{p})\\
\tilde{z}\ = \ g_{2}(z \cap T_{q}) \qquad  \tilde{w}\ = \ g_{2}(w
\cap T_{q}).
\end{split}
\end{equation}
Because $\alpha$ is tail-like, $\tilde{x}, \tilde{y}, \tilde{z},
\tilde{w}$ are branches in $X(T)$ of height $\alpha$. The pairs
$\tilde{x},\tilde{y}$ and $\tilde{z},\tilde{w}$ have equality
level $0$ and by (\ref{eq5.66})
\begin{equation}\label{eq5.68}
\tilde{x}_{1} \ < \ \tilde{y}_{1} \quad \mbox{and} \quad
\tilde{z}_{1} \ < \ \tilde{w}_{1} \qquad \mbox{in}\  S_{0}.
\end{equation}

Because $S_{0}$ is doubly transitive there exists a LOTS
automorphism $g $ of $S_{0}$ such that
\begin{equation}\label{eq5.69}
g(\tilde{x}_{1}) \ = \ \tilde{z}_{1} \quad \mbox{and} \quad
g(\tilde{y}_{1}) \ = \ \tilde{w}_{1} .
\end{equation}

Now choose for each $r \in S_{0}$ a tree isomorphism
\begin{equation}\label{eq5.70}
g_{r} : T_{r} \rightarrow T_{g(r)}.
\end{equation}
Before putting these together in the now familiar way, we make one
last pair of adjustments.

With $r_{1} = \tilde{x}_{1}$ and $r_{2} = \tilde{z}_{1}$ we have
that $(g_{r_{1}})_{*}(\tilde{x} \cap T_{r_{1}})$ and $\tilde{z}
\cap T_{r_{2}}$ are both branches of height $\alpha$ in the tree
$T_{r_{2}} \cong T$.  By Proposition \ref{prop5.16} there is an automorphism
of $T_{r_{2}}$ which maps one branch to the other.  By composing
with such an automorphism we can adjust $g_{r_{1}}$ so that
\begin{equation}\label{eq5.71}
(g_{\tilde{x}_{1}})_{*}(\tilde{x}\cap T_{\tilde{x}_{1}}) \ = \
\tilde{z}\cap T_{\tilde{z}_{1}}.
\end{equation}
Similarly, since $\tilde{w}_{1} \not= \tilde{z}_{1}$ we can make
an independent adjustment on $T_{\tilde{w}_{1}}$ to get
\begin{equation}\label{eq5.72}
(g_{\tilde{y}_{1}})_{*}(\tilde{y}\cap T_{\tilde{y}_{1}}) \ = \
\tilde{w}\cap T_{\tilde{w}_{1}}.
\end{equation}

Now put together the tree isomorphisms $g_{r}$ to get a tree
automorphism $g$ of $T$ such that
\begin{equation}\label{eq5.73}
g_{*}(\tilde{x}) \ = \ \tilde{z}\qquad \mbox{and} \qquad
g_{*}(\tilde{y}) \ = \ \tilde{w}.
\end{equation}

Define $f = (g_{2})^{-1}\circ g \circ g_{1} : T_{p} \rightarrow
T_{q} $ a tree isomorphism such that
\begin{equation}\label{eq5.74}
f_{*}(x \cap T_{p}) \ = \ z \cap T_{q}\qquad \mbox{and} \qquad
f_{*}(y \cap T_{p}) \ = \ w \cap T_{q}.
\end{equation}

Let $k = j_{q} \circ f_{*} \circ (j_{p})^{-1} $ so that $k$ is an
order isomorphism from $j_{p}(X(T_{p})) \subset X(T)$ to
$j_{q}(X(T_{q})) \subset X(T)$.  From (\ref{eq5.74}) and the definition
(\ref{eq5.7}) and (\ref{eq5.8}) of $j_{p}$ and $j_{q}$ we have
\begin{equation}\label{eq5.75}
k(x) \ = \  z \qquad \mbox{and} \qquad k(y) \ = \ w.
\end{equation}
By (\ref{eq5.27}) the open interval $(x,y)$ is contained in
$j_{p}(X(T_{p}))$ and $k $ restricts to an isomorphism of $(x,y)$
to $(z,w)$.  Since $j_{p}$,  $f_{*}$ and $j_{q}$ map branches of
height $\alpha$ to branches of height $\alpha$ we have that
\begin{equation}\label{eq5.76}
k(j_{p}(X(T_{p}) \cap \tilde{L}_{\alpha}))\ = \ X(T_{q}) \cap
\tilde{L}_{\alpha}
\end{equation}
which implies (\ref{eq5.63}).

(b):  In this case the branches $x,y,z,w$ in $X(T^{\alpha +1})$
correspond to vertices in $T$ of level $\alpha$. We mimic the
proof of part (a) defining $\beta$, $\epsilon$ and vertices $p,q$
as before.  $\tilde{x},\tilde{y},\tilde{z},\tilde{w}$ are branches
of height $\alpha$ in $X(T^{\alpha +1})$. However, when we define
the tree isomorphisms $g_{r}$ as before they are isomorphisms of
the entire tail tree not just of the level $\alpha$ truncation.
In making the adjustment to $g_{r_{1}}$ with $r_{1} =
\tilde{x}_{1}$, we think of $\tilde{x}$ and $\tilde{z}$ not as
branches but as vertices of the tree $T$ at level $\alpha$.  Then
$g_{r_{1}}(\tilde{x})$ and $\tilde{z}$ are both level $\alpha$
vertices of $T_{r_{2}} \cong T$.  By Proposition \ref{prop5.16} we may use
homogeneity of $T$ to adjust $g_{r_{1}}$ and similarly adjust
$g_{r_{2}}$ so that
\begin{equation}\label{eq5.77}
 g_{\tilde{x}_{1}}(\tilde{x}) \ = \ \tilde{z} \qquad \mbox{and} \qquad
   g_{\tilde{y}_{1}}(\tilde{y}) \ = \ \tilde{w}.
\end{equation}
Assemble the maps $g_{r}$ to form the automorphism $g$ of $T$ and
define $f = (g_{2})^{-1}\circ g \circ g_{1} : T_{p} \rightarrow
T_{q} $ as before. Regarding $x,y,z,w$ as vertices, $f$ satisfies
\begin{equation}\label{eq5.78}
f(x) \ = \ z \qquad \mbox{and} \qquad f(y) \ = \ w .
\end{equation}

Because $\alpha$ is tail-like and $o(p) < \alpha$ (\ref{eq5.58}) implies
that $(T^{\alpha +1})_{p} = (T_{p})^{\alpha +1}$ and similarly for
$q$ since $o(q) < \alpha$ as well.

Thus, we can define $f_{*} : X((T^{\alpha+ 1})_{p}) \rightarrow
X((T^{\alpha+ 1})_{q})$ and let $ k = j_{q} \circ f_{*} \circ
(j_{p})^{-1} $ an order isomorphism between the intervals
$j_{p}(X((T^{\alpha +1})_{p}))$ and $j_{q}(X((T^{\alpha
+1})_{q}))$  in  $X(T^{\alpha +1})$ so that
\begin{equation}\label{eq5.79}
k(j_{p}(X((T^{\alpha +1})_{p})) \cap L_{\alpha}) \ = \
j_{q}(X((T^{\alpha +1})_{q})) \cap L_{\alpha}.
\end{equation}

As before, (\ref{eq5.78}) implies that $k$ maps the branches $x$ and $y$
to $z$ and $w$, respectively. Finally, (\ref{eq5.27}) again implies that
the interval $(x,y)$ is contained in $j_{p}(X((T^{\alpha
+1})_{p}))$.

\end{proof} \vspace{.5cm}

We illustrate the use of this result with the following.

\begin{prop}\label{prop5.18}  Let $T$ be a reproductive, bi-ordered tree with $S_{0}$ a doubly transitive
LOTS with at least three elements.  Let $\alpha = h(T)$.
If $\tilde{L}_{\alpha} = \{ x \in X(T) : h(x) = \alpha \} $ is nonempty, then
it is a doubly transitive LOTS which is dense in $X(T)$.
\end{prop}

\begin{proof} If $x \in \tilde{L}_{\alpha}$ and
$p \in T$, then there exists $q \in x $ such that $o(q) = o(p)$. By
Proposition \ref{prop5.16}  there exists an automorphism $f$ of $T$ such
that $f(q) = p$.  Then $f_{*}(x) \in \tilde{L}_{\alpha}$ with $p
\in f_{*}(x) $.  It follows from Proposition \ref{prop5.3}(g) that
$\tilde{L}_{\alpha}$ is dense in $X(T)$.  Since $S_{0}$ is doubly
transitive and infinite it is unbounded.  Hence, by Proposition
\ref{prop5.3}(d), $X(T)$ is unbounded. Since $\tilde{L}_{\alpha}$ is dense in
$X(T)$ it is unbounded.  That $\tilde{L}_{\alpha}$ is doubly
transitive then follows from Lemma \ref{lem5.17}(a) and Proposition \ref{prop3.2}(d)
(vii)$\Rightarrow$(i).

\end{proof}\vspace{.5cm}

Our main application of Lemma \ref{lem5.17} requires the following
technical result.

\begin{lem}\label{lem5.19}  Let $X$ be an unbounded, order dense LOTS of countable type
and let $W$ be a dense subset of $X$.  If for every $x < y$, $z < w$ in $W$ there exists
an order isomorphism $k : (x,y) \rightarrow (z,w) $ (intervals in X) such that
\begin{equation}\label{eq5.80}
k((x,y) \cap W) \  = \  (z,w) \cap W,
\end{equation}
then the same is true for any $x < y$, $z < w$ in $X$.  In
particular, $W$ and $X$ are HLOTS.
\end{lem}

\begin{proof} $W$ is unbounded and so is doubly
transitive by Proposition \ref{prop3.2}(d)(v)$\Rightarrow$(i).  By Proposition \ref{prop2.5}(f),
$W$ has countable type and so is a HLOTS by Proposition \ref{prop3.4}(a).

If $x < y$, $z < w$ in $X$, then since $X$ is order dense, $W$ is
dense in $X$ and $X$ has countable type, there exist $f,g :
\Z \rightarrow W$ with $f(\Z)$ $\pm$cofinal in $(x,y)$ and with $g(\Z)$ $\pm$cofinal
in $(z,w)$.  Put together isomorphisms $(f(i),f(i+1)) \cong
(g(i),g(i+1))$ to get the required $k$. $X$ is doubly transitive
and so is a HLOTS by Proposition \ref{prop3.2}(d) and Proposition \ref{prop3.4}(a) again.

\end{proof}\vspace{.5cm}

\begin{theo}\label{theo5.20}  If $T$ is a reproductive, $\Omega$ bounded, bi-ordered  tree with $S_{0}$ a HLOTS,
then $X(T)$ is an IHLOTS and $\widehat{X(T)}$ is a CHLOTS.
\end{theo}

\begin{proof} By Lemma \ref{lem5.15}(b) $T$ is of $S_{0}$
type.  Proposition \ref{prop5.9}(c) and Proposition \ref{prop2.8}(c) then imply that $X(T)$ and $\widehat{X(T)}$
are of countable type. By Proposition \ref{prop5.3} $X(T)$ is unbounded,
order dense and has dense holes. So the completion $\widehat{X(T)}$ is unbounded and connected.
If $W$ is any dense subset of $X(T)$, then it
is of countable type by Proposition \ref{prop2.5}(f).  Since such a $W$ is a
proper subset of its completion $\widehat{X(T)}$ it is an IHLOTS if it
is doubly transitive by Proposition \ref{prop3.4}.  Then the completion
$\widehat{X(T)}$ is a CHLOTS by the same proposition.

Since $T$ is $\Omega$ bounded, we have $h(T) \leq
\Omega$ with $h(T)$ an infinite tail-like ordinal.  There are two
cases.\vspace{.25cm}

\textbf{Case 1:}  If $\alpha = h(T)$ is countable, then by Lemma \ref{lem5.4}(b)
$\tilde{L}_{\alpha}$ is dense in $X(T)$. As in Proposition \ref{prop5.18},
Lemma \ref{lem5.17}(a) allows us to conclude that $\tilde{L}_{\alpha}$ is
doubly transitive, using Proposition \ref{prop3.2}(d).  Lemma \ref{lem5.17}(a) together with
Lemma \ref{lem5.19} implies that $\tilde{L}_{\alpha}$ and $X(T)$ are HLOTS
in this case.\vspace{.25cm}

\textbf{Case 2:} If $h(T) = \Omega $, then given $x < y$ and $z < w$ in
$X(T)$ we can choose $\alpha$ a countable tail-like ordinal such
that
\begin{equation}\label{eq5.81}
\alpha \quad > \quad  h(x), h(y), h(z), h(w)
\end{equation}
because $h$ takes only countable values on $X(T)$ and $h(T) =
\Omega$ is the limit of countable, tail-like ordinals.

Now we apply Lemma \ref{lem5.17}(b) to see that $W = L_{\alpha} \subset
X(T^{\alpha +1})$ satisfies the hypotheses of Lemma \ref{lem5.19}.  From
(\ref{eq5.81}) we clearly have
\begin{equation}\label{eq5.82}
\pi^{\alpha +1}(x) \ < \ \pi^{\alpha +1}(y) \quad \mbox{and} \quad
\pi^{\alpha +1}(z) \ < \ \pi^{\alpha +1}(w)\quad \mbox{in} \
X(T^{\alpha +1}).
\end{equation}
It follows that there exists an order isomorphism $k :
(\pi^{\alpha +1}(x),\pi^{\alpha +1}(y)) \rightarrow (\pi^{\alpha
+1}(z),\pi^{\alpha +1}(w)) $, intervals in $X(T^{\alpha +1})$,
which preserves $L_{\alpha}$.\vspace{.25cm}

We finish up as in the proof of Theorem \ref{theo5.12}. As we did there,
write
\begin{equation}\label{eq5.83}
X(T) \  = \  \Sigma \{ (\pi^{\alpha +1})^{-1}(r) : r \in
X(T^{\alpha +1}) \}.
\end{equation}
Recall that for $r = x(p)$ with $o(p) = \alpha$, $(\pi^{\alpha
+1})^{-1}(r) = j_{p}(T_{p})$ while for $r \not\in L_{\alpha}$
$(\pi^{\alpha +1})^{-1}(r)$ is a singleton.

For $x(p) \in  (\pi^{\alpha +1}(x),\pi^{\alpha +1}(y)) \cap
L_{\alpha} $ let $k(p)$ denote the vertex of level $\alpha$ such
that $k(x(p)) = x(k(p))$. For each such $p$ choose a tree
isomorphism $g_{p} : T_{p} \rightarrow T_{k(p)}$.

For $r = x(p)$ let
\begin{equation}\label{eq5.84}
k_{r} =  j_{k(p)} \circ (g_{p})_{*} \circ
(j_{p})^{-1}:(\pi^{\alpha +1})^{-1}(r) \rightarrow (\pi^{\alpha
+1})^{-1}(k(r)).
\end{equation}
and for $r \not\in L_{\alpha}$ let $k_{r}$ be the map of
singletons given by $(\pi^{\alpha +1})^{-1}(r) \mapsto$\\
$ (\pi^{\alpha +1})^{-1}(k(\pi^{\alpha +1}(r)))$.  Putting these together using (\ref{eq5.83}) we obtain an
isomorphism $\tilde{k} : (x,y) \rightarrow (z,w)$, intervals in
$X(T)$.

Thus, $X(T)$ is doubly transitive in this case as well by Proposition \ref{prop3.2}(d).

\end{proof} \vspace{.5cm}

Adapting the proof, we obtain the first part of the following.

\begin{theo}\label{theo5.21}  Assume that $T$ is a reproductive,  bi-ordered
tree of height $\Omega$ with $S_{0}$ a HLOTS.
\begin{itemize}
\item[(a)] If $\tilde L_{\Omega}$ is nonempty, i.e. there exist branches of height $\Omega$,
then  $\tilde L_{\Omega}$ is dense in $X(T)$ and is
 doubly transitive. If, in addition, $X(T) \setminus \tilde L_{\Omega}$ is nonempty,
 then $X(T) \setminus \tilde L_{\Omega}$ is dense in $X(T)$ and is
first countable and weakly homogeneous.

\item[(b)] Assume that  $\tilde L_{\Omega}$ and $X(T) \setminus \tilde L_{\Omega}$ are nonempty and that $X =
\widehat{X(T)}  \setminus \tilde L_{\Omega}$ with $\widehat{X(T)}$ the completion of $X(T)$. The LOTS $X$ is first countable and weakly homogeneous.
In addition, $X$ is countably complete. That is, if $A \subset X$ is a bounded, countable set, then $sup \ A, \ inf \ A \ \in X$. On the other hand,
$X$ is not complete.
\end{itemize}
\end{theo}

\begin{proof}  (a): If $\tilde L_{\Omega}$ is nonempty, then it is dense in $X(T)$ and is
doubly transitive by Proposition \ref{prop5.18}.

Since $S_0$ is a HLOTS, $X(T)$ is unbounded and is $\sigma$-bounded.
If $x \in X(T)$ with $h(x) < \Omega$, then because $X(T)$ is unbounded and $\tilde L_{\Omega}$ is dense
there exist $z_1, z_2 \in \tilde L_{\Omega}$ such that $z_1 < x < z_2$.  Now let
$w_1 < w_2 \in \tilde L_{\Omega}$ be arbitrary. By Lemma \ref{lem5.17}
there exists an order isomorphism $k : (z_1,z_2) \to (w_1,w_2)$ which maps
$(z_1,z_2) \cap \tilde L_{\Omega}$ to $(w_1,w_2) \cap \tilde L_{\Omega}$.
Hence, $k(x) \not\in \tilde L_{\Omega}$, i.e. $w_1 < k(x) < w_2$ with
$h(k(x)) < \Omega$. It follows that $X(T) \setminus \tilde L_{\Omega}$ is dense
in $X(T)$.

The proof of Case 2 in the proof of Theorem \ref{theo5.20} directly shows that
$X(T) \setminus \tilde L_{\Omega}$ is doubly transitive.

Let $\a$ be a countable tail-like ordinal with $\a > h(x)$. By Lemma \ref{lem5.15} (b),
$T^{\a}$ is reproductive and so by Theorem \ref{theo5.20}
$X(T^{\a})$ is a HLOTS and so is first countable. Hence, $x$ is the limit of
increasing and decreasing sequences in $X(T^{\a})$.
So there exist  sequences in $X(T)$ whose projections via $\pi^{\a}$ are
strictly monotone sequences $X(T^{\a})$ converging to $x$.
Furthermore, $(\hat \pi^{\a})^{-1}(x)$ is a singleton in the complete LOTS
$\widehat{X(T)}$. It follows that the sequences in $X(T)$ converge
to $x$.  Thus,  $X(T) \setminus \tilde L_{\Omega}$ is doubly transitive,
first countable and $\sigma$-bounded. It is weakly homogeneous by Proposition
\ref{prop3.2a}(b).

(b): Applying Proposition \ref{prop5.6} we identify the completion $\widehat{X(T)}$ with the branch space $X(\hat{T})$ where $\hat{T}$ is
the tree completion of Definition \ref{df5.5}. Every branch of $X(\hat{T}) \setminus X(T)$ is of the form $x(q)$ with $q$ a new vertex.
These branches have countable height and each is clearly the limit of increasing and decreasing sequences because $S_0$ is of countable type.
From (a) it follows that  $X$ is first countable.

On the other hand, let $x \in \tilde L_{\Omega}$ and $A$ be a countable subset of $X(\hat{T})$ with $a < x$ for all $a \in A$.
The equality level of $a$ and $x$ is a countable ordinal for each $a$. We can choose $\a < \Omega$ which is greater than all of these
equality levels. The successor set in $T$ for $x_{\a}$ is isomorphic to $S_0$ and so is unbounded. Choose $y$ in the successor set
with $y < x_{\a + 1}$. Thus, a branch $b$ of $X(T)$ which contains $y$ satisfies $b < x$ and $a < b$ for all $a \in A$. Thus,
$sup \ A$ in the complete space $X(\hat{T})$ is not equal to $x$. Hence, $sup \ A \in X$.
Similarly, if $A$ is bounded below then $inf \ A$ lies in $X$. That is, $X$ is
countably compact. In particular,
$x$ is not a limit of an increasing or decreasing sequence in $X(\hat{T})$.  It follows that
any automorphism of $X(\hat{T})$ maps $\tilde L_{\Omega}$ to itself and so restricts to an automorphism of $X$.

 If $a < b, c < d$ in $X(T) \setminus \tilde L_{\Omega}$, then from (a) there exists an automorphism of
$X(T) \setminus \tilde L_{\Omega}$ which maps $[a,b]$ to $[c,d]$. Because $X(T) \setminus \tilde L_{\Omega}$ is dense in $X(T)$ which is
dense in $X(\hat{T})$ the automorphism extends to an automorphism of $X(\hat{T})$ and then restricts to an automorphism of $X$.

Now let $a < b$ in $X$. Because $X(T)$ is $\s$-bounded and $X(T) \setminus \tilde L_{\Omega}$
is dense in $X(T)$ we can choose an embedding $q : \Z \to X(T) \setminus \tilde L_{\Omega}$ which is
$\pm$ cofinal in $X$. Similarly, because $X$ is first countable, we can choose an
embedding $r : \Z \to (a,b) \cap (X(T) \setminus \tilde L_{\Omega})$ which is $\pm$ cofinal in $(a,b)$. Choose an isomorphism from
$[q(i),q(i+1)]$ to $[r(i),r(i+1)]$ for each $i \in \Z$ and concatenate to obtain an isomorphism from $X$ to $(a,b) \cap X$. Thus,
$X$ is weakly homogeneous.

\end{proof} \vspace{.5cm}

\noindent {\bfseries Remark.} Since $X$ is invariant for every automorphism of $\widehat{X(T)}$ it follows that
$X(T)$ and its completion are not transitive.  If $x < z \in \tilde L_{\Omega}$, then the convex set
$X_1 = (x,z) \cap  X$ is not $\sigma$-bounded. So while $X_1$ is unbounded, every countable subset is bounded and
so has a sup and inf in $X_1$.  $X_1$ is doubly transitive, but not weakly homogeneous.
 \vspace{1cm}

\section{\textbf{Tree Constructions}}
\vspace{.5cm}

\subsection{The Simple Trees on a LOTS}

Throughout this section $X$ will be a nontrivial
LOTS, i.e. with at least two points.  For any ordinal $\alpha$, the set of maps $X^{\alpha}$ is a
LOTS with the order space product structure, see (\ref{eq2.12}).  If $p \in
X^{\alpha}$ and $\beta \leq \alpha$, then the restriction $p|\beta
\in X^{\beta}$.  For $p \in X^{\alpha}, q \in X^{\beta}$ we write
\begin{equation}\label{eq6.1}
q \ \prec p \quad \Longleftrightarrow \quad \beta < \alpha \
\mbox{and} \ p|\beta = q.
\end{equation}

Accordingly, we define the \emph{simple tree on $X,\alpha$}\index{simple tree on $X,\alpha$}
\index{tree!simple on $X,\alpha$} by
letting $X^{i}$ be the set of vertices of level $i$ for $i <
\alpha$ and using the restriction partial order\index{restriction partial order} of (\ref{eq6.1}).  If $p
\in X^{i}$ and $i+1 < \alpha$ then we can identify the successor
set $S_{p} \subset X^{i+1}$ with $X$ by

\begin{equation}\label{eq6.2}
\begin{split}
S_{p} \  \cong \  X   \\
q \ \mapsto \ q(i).
\end{split}
\end{equation}

We use the LOTS ordering on $X$ and the identification of (\ref{eq6.2}) to
give the simple tree on $X,\alpha$ the structure of a normal tree,
bi-ordered and  of $X$ type.   Conditions (i)-(iv) of Definition \ref{df5.1}
are easy to check.  Every branch has height $\alpha$.  Such a
branch is a coherent collection $\{ p_{\epsilon} : \epsilon <
\alpha \}$ of functions whose union is a unique
element $x \in X^{\alpha}$.  Thus, we can identify $X^{\alpha}$
with the branch space of the simple tree on $X,\alpha$.

For $p \in X^{\alpha}$ and $q \in X^{\beta}$ we define the
\emph{sum} $p + q \in X^{\alpha + \beta}$\index{sum}\index{$p + q$} by
\begin{equation}\label{eq6.3}
(p + q)(i) \ = \ \begin{cases}   p(i) \qquad i < \alpha  \\
q(i \setminus \alpha) \qquad \alpha \leq i < \alpha + \beta
\end{cases}
\end{equation}
where, as usual, we identify the tail set $i \setminus \alpha$
with the ordinal which has its order type.

For ordinals $ \alpha < \epsilon $ we define the \emph{translation
map}\index{translation map}
\begin{equation}\label{eq6.4}
\begin{split}
\tau_{\alpha} : X^{\epsilon} \rightarrow X^{\epsilon \setminus \alpha} \hspace{2cm}  \\
\tau_{\alpha}(r)(i) \ = \ r(\alpha + i) \qquad \mbox{for} \ i <
\epsilon \setminus \alpha.
\end{split}
\end{equation}

If $p \in X^{\alpha}$, $q \in X^{\beta}$ and $p + q \in
X^{\epsilon}$ so that $\epsilon = \alpha + \beta $, then $\beta =
\epsilon \setminus \alpha$ and
\begin{equation}\label{eq6.5}
(p + q)|\alpha \ = \ p \qquad \mbox{and} \qquad \tau_{\alpha}(p +
q) \ = \ q.
\end{equation}.
Conversely, if $r \in X^{\epsilon}$ and $\epsilon = \alpha +
\beta$, then $\beta = \epsilon \setminus \alpha$ and
\begin{equation}\label{eq6.6}
r \ = \ (r|\alpha) + \tau_{\alpha}(r).
\end{equation}

Clearly, if $o(p) = \alpha$ and $o(r) = \ep$, then
\begin{equation}\label{eq6.7}
p \  \prec \ r \qquad \Longleftrightarrow \qquad r \ = \ p + q
\end{equation}
for some $q \in X^{\epsilon \setminus \alpha}$ in which case, by
(\ref{eq6.5}) $\tau_{\alpha}(r) = q $.

The map $\tau_{\a}$ is clearly surjective, but is not order preserving and need not be continuous.

\begin{prop}\label{prop6.1} Let  $\alpha$ be a positive ordinal.
\begin{enumerate}
\item[(a)]  If $X$ is transitive, then the simple tree on $X,\alpha$ is homogeneous.
\item[(b)]  If $\alpha$ is tail-like, then the simple tree on $X,\alpha$ is reproductive.
\end{enumerate}
\end{prop}

\begin{proof} (a) If $p,q \in X^{i}$, then there
exists for each $j < i$, $g_{j} \in H_{+}(X)$ such that
$g_{j}(p(j)) = q(j)$.  For $i \leq j < \alpha$ let $g_{j}$ be the
identity. On $X^{k}$ for $k < \alpha$ the product map $\Pi_{j < k}
g_{j} $ defines a tree automorphism which maps $p$ to $q$.

(b) For $p$ in the simple tree $o(p) < \alpha$ and $\alpha$
tail-like implies that $o(p) + \beta < \alpha$ whenever $\beta <
\alpha$.  So we can define the \emph{canonical tree isomorphism at
$p$}\index{canonical tree isomorphism at $p$} by
\begin{equation}\label{eq6.8}
\begin{split}
a_{p} :  T \rightarrow T_{p}   \\
a_{p}(q) \ = \ p + q.
\end{split}
\end{equation}

\end{proof}\vspace{.5cm}

As a corollary we obtain a tree proof of Theorem \ref{theobabcock}.

\begin{cor}\label{cor6.2}  Let  $\alpha$ be an infinite, tail-like ordinal.
\begin{enumerate}
\item[(a)]  If $X$ is doubly transitive, then $X^{\alpha}$ is doubly transitive.
\item[(b)]  If $\alpha$ is countable and $X$ is a HLOTS, then $X^{\alpha}$ is a HLOTS.
\end{enumerate}
\end{cor}

\begin{proof} $X^{\alpha}$ is the branch space
of the simple tree on $X,\alpha$. By Proposition \ref{prop6.1}(b) the simple
tree is reproductive.

In this case every branch has height $\alpha$, i.e. $X^{\alpha} =
\tilde{L}_{\alpha}$. If $X$ is doubly transitive, then $X^{\alpha}$
is doubly transitive by Proposition \ref{prop5.18}.  If, in addition,
$\alpha$ is countable and $X$ has countable type, then by
Proposition \ref{prop2.7}(b), $X^{\alpha}$ has countable type.  So it is a
HLOTS by Proposition \ref{prop3.4}(a).

\end{proof} \vspace{.5cm}

Recall from Section 4 that if $X$ is a HLOTS, then we choose a
nontrivial, bounded subinterval $ J = [-1,+1]$ in $X$ to define
$X_{\alpha}$ for every positive ordinal $\alpha$.  If $\alpha$ is
an infinite limit ordinal, then
we define
\begin{equation}\label{eq6.9}
X_{<\alpha} \quad =  \quad \bigcup _{\beta < \alpha}
j_{\alpha}^{\beta}((X_{\beta})')
\end{equation}\index{$X_{<\alpha}$}
using the inclusion maps of (\ref{eq4.5}).  So $x \in X_{<\alpha}$  iff
there exists $\beta < \alpha$ such that
\begin{equation}\label{eq6.10}
x_{i}  =  -1 \ \mbox{for all} \ \beta \leq i < \alpha \qquad
\mbox{or}\qquad x_{i}  =   +1 \ \mbox{for all} \ \beta \leq i <
\alpha.
\end{equation}
Clearly, $X_{<\alpha}$ and its complement are dense in
$X_{\alpha}$.

\begin{prop}\label{prop6.3}  If $X$ is a HLOTS and $\alpha$ is an infinite, tail-like ordinal, then
\begin{equation}\label{eq6.11}
X^{\alpha} \quad \cong \quad X_{\alpha} \setminus X_{<\alpha}.
\end{equation}
In particular, $X^{\alpha}, X_{<\alpha}, X_{\alpha} \setminus
X_{<\alpha}$ and $X_{\alpha}$ all have order isomorphic
completions.
\end{prop}

\begin{proof} For $i \in \alpha \setminus 1 = \{
j : 0 < j < \alpha \}$ we denote by $i^{*}$ the corresponding
element of the reverse $(\alpha \setminus 1)^{*}$ and we write
\begin{equation}\label{eq6.12}
|i| \ = \ i \ = \ |i^{*}| \qquad \mbox{for} \ i \in \alpha.
\end{equation}
Apply Proposition \ref{prop2.8}(d)
to $\bullet X \bullet$ to get an order
embedding
\begin{equation}\label{eq6.13}
\begin{split}
\tilde{g} : (\alpha \setminus 1)^{*} + (\alpha \setminus 1) \rightarrow  X \quad \text{with}\hspace{2cm} \\
\tilde{g}(1) \ = \  +1, \qquad \  \tilde{g}(\alpha \setminus 1) \ \mbox{cofinal in} \ X, \hspace{1cm}\\
\tilde{g}(1^{*}) \ = \  -1, \quad \mbox{and}\quad  \tilde{g}((\alpha
\setminus 1)^{*})  \ \mbox{coinitial in} \ X.
\end{split}
\end{equation}
We define the intervals $\{ J_{i^{*}} : i^{*} \in (\alpha
\setminus 1)^{*} \} \cup \{J_{0}\} \cup \{ J_{i} : i \in (\alpha
\setminus 1)\} $ and choose isomorphisms:
\begin{equation}\label{eq6.14}
\begin{split}
g_{0} \  = \ \mbox{identity \quad  on} \ J_{0} \ = \ (-1,+1). \hspace{.5cm} \\
g_{i} : [-1,+1) \rightarrow J_{i} \ = \  [\tilde{g}(i),\tilde{g}(i+1))  \\
g_{i^{*}} : (-1,+1] \rightarrow J_{i^{*}} \ = \
(\tilde{g}((i+1)^{*}),\tilde{g}(i^{*})]
\end{split}
\end{equation}
This defines an $(\alpha \setminus 1)^{*} + \{0\} + (\alpha
\setminus 1)$ indexed family of pairwise disjoint, nontrivial
subintervals of $X$ with union $X$.

Now for each $x \in X_{\alpha} \setminus X_{<\alpha}$ we define an
order embedding from $\alpha$ into itself:
\begin{equation}\label{eq6.15}
\begin{split}
\beta(x,0) = 0  \hspace{4cm}  \\
\beta(x,i) = sup \{ \beta(x,j) : j < i \} \quad \mbox{when} \  i \
\mbox{is a limit ordinal}.
\end{split}
\end{equation}

The remaining values are defined inductively.  Assuming
$\beta(x,i)$ is defined we construct $\beta(x,i+1)$.

First, we set  $m(x,0) = 0$ and for $i \in \alpha \setminus
1$:
\begin{equation}\label{eq6.16}
\begin{split}
m(x,i) = 0 \quad  \mbox{if} \quad -1 < x_{\beta(x,i)} < +1, \hspace{1cm}\\
m(x,i) = j \quad \mbox{if} \quad x_{\beta(x,i)}= +1 \ \mbox{and} \ j = min \ \{k : x_{\beta(x,i)+k} < +1 \},  \\
m(x,i) = j^{*} \quad \mbox{if} \quad x_{\beta(x,i)}= -1 \
\mbox{and} \ j = min \ \{k : x_{\beta(x,i)+k} > -1 \}.
\end{split}
\end{equation}
Observe that $x \not\in X_{<\alpha}$ implies that the sets used to
define $j$ in the latter cases are nonempty.  Complete the
inductive definition of $\beta$ by
\begin{equation}\label{eq6.17}
\beta(x,i+1) \ = \ \beta(x,i) + |m(x,i)| + 1,
\end{equation}
so that, for example, $\beta(x,1) = 1$.

Since $x \in X_{\alpha}$, $x_{0} \in X$ and $x_{i} \in J =
[-1,+1]$ for all $i \in \alpha \setminus 1$.  We define $f(x) \in
X^{\alpha}$ by
\begin{equation}\label{eq6.18}
f(x)_{0} = x_{0} \quad \mbox{and} \quad f(x)_{i} =
g_{m(x,i)}(x_{\beta(x,i)+|m(x,i)|}) \quad \mbox{for} \ i \in
\alpha \setminus 1.
\end{equation}

If $x_{\beta(x,i)} \in (-1,+1) $, then $m(x,i) = 0$ and
$f(x)_{i} = x_{\beta(x,i)}$.

If $x_{\beta(x,i)} = +1$, then $j =
m(x,i) \in \alpha \setminus 1$  and $x_{\beta(x,i)+k} = +1$ for
all $k < j$ while $x_{\beta(x,i)+j} \in [-1,+1)$. The map $g_{j}$
moves $x_{\beta(x,i)+j}$ to the interval $J_{j}$ in $X$.

Similarly, if $x_{\beta(x,i)} = -1$ then $j^{*} = m(x,i) \in
(\alpha \setminus 1)^{*}$ and $g_{j^{*}}$ moves $x_{\beta(x,i)+j}$
to $J_{j^{*}}$.

Thus, we have defined a map $f : X_{\alpha} \setminus X_{<\alpha}
\rightarrow X^{\alpha}$.

Suppose that $x < y$.  If $x_{0} < y_{0}$, then $f(x) < f(y)$.
Otherwise, let $\epsilon = min \ \{j : x_{j} \not= y_{j} \}$ so
that $\epsilon \in \alpha \setminus 1$ and $x_{\epsilon} <
y_{\epsilon}$.  By Proposition \ref{prop2.8}(c) applied to the order embedding
$i \mapsto \beta(x,i)$ there exists $i \in \alpha$ such that
\begin{equation}\label{eq6.19}
\beta(x,i)\   \leq \  \epsilon \  < \  \beta(x,i+1).
\end{equation}
Since $\epsilon \in \alpha \setminus 1$ and $\beta(x,1) = 1$, we
have $i \in \alpha \setminus 1$.  Notice that the inductive
definitions imply that
\begin{equation}\label{eq6.20}
m(x,j) = m(y,j) \quad j < i \qquad \mbox{and} \qquad \beta(x,j) =
\beta(y,j) \quad j \leq i
\end{equation}
and so we have $f(x)_{j} = f(y)_{j}$ for all $j < i$ and by
considering the various cases we check that $f(x)_{i} <
f(y)_{i}$. \vspace{.25cm}

\textbf{Case 1:} If $\b(x,i) = \b(y,i) = \ep $, then $x_{\b(x,i)} < y_{\b(y,i)}$ and so $m(y,i) \in \a$ and $m(x,i) \in \a^*$.
Hence, $m(x,i) \leq m(y,i)$. \vspace{.25cm}

\textbf{Case 2:} If $\b(x,i) = \b(y,i) < \ep $, then $x_{\b(x,i)} = y_{\b(y,i)}$ and the common value cannot lie in $(-1,+1)$ since
$\ep < \b(x,i+1)$. If $x_{\b(x,i)} = y_{\b(y,i)} = +1$, then $m(x,i) < m(y,i) \in \a \setminus 1$.
If $x_{\b(x,i)} = y_{\b(y,i)} = -1$, then $m(x,i) < m(y,i) \in (\a \setminus 1)^*$.  \vspace{.25cm}

In either case, $f(x)_{i} <
f(y)_{i}$ and so
$f(x) < f(y)$.  Thus, f is an order injection.

To reverse the procedure, start with $z \in X^{\alpha}$ and let
$x_{0} = z_{0}$.  Then for any $i \in \alpha \setminus 1$,
$\{z_{j} : 0 \leq j < i \}$ determines, inductively, the sets $\{
m(x,j) : j < i \}, \{ \beta(x,j) : j \leq i \}$ and $\{ x_{k} : k
< \beta(x,i) \}$.  Now if $z_{i} \in J_{0}$, then $m(x,i) = 0$,
$x_{\beta(x,i)} = z_{i}$ and $\beta(x,i+1) = \beta(x,i)+1$.  If
$z_{i} \in J_{j}$ for some $j \in \alpha \setminus 1$ then
$m(x,i) = j$, $x_{\beta(x,i)+k} = +1$ for $0 \leq k < j$,
$x_{\beta(x,i)+j} = (g_{j})^{-1}(z_{i})$ and $\beta(x,i+1) =
\beta(x,i)+j+1$.  Similarly, if $z_{i} \in J_{j^{*}}$.

Thus, $f$ is surjective and so is an order isomorphism.

\end{proof} \vspace{.5cm}

We can use this result to get an alternative proof of the CHLOTS
portion of Theorem \ref{theo4.2}.

\begin{cor}\label{cor6.4}  Let $X$ be a CHLOTS and $\alpha$ be a countably infinite, tail-like ordinal.
$X^{\alpha}$ and $X_{<\alpha}$ are IHLOTS with completion isomorphic to the CHLOTS $X_{\alpha}$.
\end{cor}

\begin{proof} $X^{\alpha}$ is a HLOTS by
Corollary \ref{cor6.2}(b) and it is isomorphic to a dense proper subset of
$X_{\alpha}$ by Proposition \ref{prop6.3}.  Since $X$ is complete,
$X_{\alpha}$ is complete by Proposition \ref{prop2.3}(b) and so $X_{\alpha}$
is the completion of the image of $X^{\alpha}$. By Proposition \ref{prop3.4}
$X_{\alpha}$ is a CHLOTS and $X^{\alpha}$ and its complement in
$X_{\alpha}$, which is $X_{<\alpha}$, are IHLOTS.

\end{proof}\vspace{.5cm}

\begin{prop}\label{prop6.4a} Assume $X$ is a LOTS.  Let $T$ be a bi-ordered tree of $Y$ type.  If
$Y$ injects into $X$ and the height of the tree is at most
$\alpha$, then the branch space $X(T)$ injects into $X^{\a}$.
If, in addition,
$X$ is a CHLOTS, $Y$ is unbounded and $h(T)$ is a limit ordinal, then
$X(T)$  and its completion inject
into $X_{\alpha}$. \end{prop}

\begin{proof} We identify each nonempty successor set $S_p$ with $Y$.
Let $j: Y \to X$ be an order injection and let $z \in X$.

Define $j^{\a} : X(T) \to X^{\a}$ by
\begin{equation}\label{eq6.20a}
j^{\a}(x)(i) = \begin{cases} j(x(i)) \ \ \text{for} \ i < h(x), \\ \ z \ \ \text{for} \ h(x) \le i < \a. \end{cases}
\end{equation}
This is clearly an order injection.

If $X$ is a CHLOTS, then $X \cong J^{\circ}$ implies that $X^{\a}$ injects into $X_{\a}$.  Hence, $X(T)$ injects into $X_{\a}$.

If $Y$ is unbounded and $h(T)$ is a limit ordinal, then $X(T)$ is order dense by Proposition \ref{prop5.3a}. Since
$X_{\a}$ is complete, the order injection $\widehat{j^{\a}} : \widehat{X(T)} \to X_{\a}$ is defined by Proposition \ref{prop2.4a}.

\end{proof}
\vspace{.5cm}

Contrast these results with those for  $X_{\Omega}$.  We can split
$X_{<\Omega}$ into two disjoint pieces $X^{\pm}_{<\Omega}$ with
\begin{equation}\label{eq6.21}
x \in X^{+}_{<\Omega} \qquad \Longleftrightarrow \qquad x_{i} = +1
\quad \mbox{for} \ i \in \Omega \ \mbox{sufficiently large},
\end{equation}
and similarly, for $X^{-}_{<\Omega}$.  Each of these is a dense
subset of $X_{\Omega}$ and by Proposition \ref{prop6.3} we can identify
$X^{\Omega}$  with the complement of their union.  It is easy to
check that $z \in X_{\Omega}$ is the limit of some increasing (or
of some decreasing) sequence iff $z \in X^{-}_{<\Omega}$ (resp. $z
\in X^{+}_{<\Omega}$ ).  It follows that if $f \in
H_{+}(X_{\Omega})$, then
\begin{equation}\label{eq6.22}
f(X^{-}_{<\Omega}) \ = \ X^{-}_{<\Omega}, \quad f(X^{\Omega}) =
X^{\Omega}, \quad  f(X^{+}_{<\Omega}) \ = \ X^{+}_{<\Omega}.
\end{equation}
That is, the decomposition of $X_{\Omega}$ into three pairwise
disjoint, dense subsets is invariant with respect to the action of
the automorphism group $H_{+}(X_{\Omega})$.

By Corollary \ref{cor6.2}(a), $X^{\Omega}$ is doubly transitive. Using the
proof of Proposition \ref{prop3.4}(d) and $H_{+}$ invariance it is easy to
show that $X^{+}_{<\Omega}$ and $X^{-}_{<\Omega}$ are doubly
transitive as well.

Moreover, if $X$ is symmetric, e.g. $X = \R$, then
$X_{\Omega}$ is symmetric and any orientation reversing
isomorphism interchanges $X^{+}_{<\Omega}$ and $X^{-}_{<\Omega}$
while preserving $X^{\Omega}$.  In that case, $X_{<\Omega}$ is
$\pm$transitive but not transitive.  Any dense subset of
$X_{\Omega}$ is order dense and so has no gap pairs.

Contrast
this with the complete case in Proposition \ref{prop3.6}. \vspace{1cm}

\subsection{Additive Trees}

Now we introduce an important class of subtrees of the simple
tree.

\begin{df}\label{df6.5}  For  a positive ordinal $\alpha$ let $T$ be
a subset of the simple tree on $X,\alpha$. $T$ is called an \emph{additive $X$ tree}\index{tree!additive} when $T$
contains the root $0$, the level $1$ vertices $X^{1}$ of the simple tree and for all vertices $p, q$ of the simple tree
\begin{equation}\label{eq6.23}
p + q \in T  \qquad \Longleftrightarrow \qquad p, q \in T.
\end{equation}

$T$ is called a \emph{partially additive $X$ tree}\index{tree!partially additive}  when $T$
contains $\{0\} \cup X^{1}$ and for all vertices $p, q$ of the
simple tree
\begin{equation}\label{eq6.24}
p + q \in T \qquad \Longleftrightarrow \qquad p, q \in T \quad
\mbox{and} \quad o(p) + o(q) < h(T) .
\end{equation}
\end{df}
\vspace{.5cm}

Clearly, $p + q \in T$ always implies
\begin{equation}\label{eq6.25}
o(p) + o(q) \ = \ o(p + q) \ < \ h(T),
\end{equation}
and so, as the name suggests, condition (\ref{eq6.24}) is a weakening of
condition (\ref{eq6.23}).

\begin{prop}\label{prop6.6}   Let $\alpha, \beta  > 1$ be  ordinals.
\begin{enumerate}
\item[(a)]  The simple tree on $X,\alpha$ is a partially additive $X$ tree.  It is additive iff $\alpha$ is tail-like.
\item[(b)]  If $T$ is contained in the simple tree on $X,\alpha$, then $h(T) \leq \alpha$.
Conversely, if $T$ is a partially additive $X$ tree with $h(T) \leq \alpha$, then $T$ is a normal
subtree of the simple tree on $X,\alpha$.
\item[(c)]  If $T$ is a partially additive $X$ tree, then $T^{\beta}$ is a partially additive $X$ tree.
\item[(d)]  Assume $T$ is a partially additive  $X$ tree and $ p\in T$.  If $\beta < o(p)$, then
\begin{equation}\label{eq6.26}
p|\beta, \tau_{\beta}(p) \ \in \ T. \hspace{2cm}
\end{equation}
If $o(p) + 1 < h(T)$, then the successor set $S_{p}$ in $T$
consists of all successors of $p$ in the simple tree.
\item[(e)]  A partially additive tree $T$ is additive iff $h(T)$ is a tail-like ordinal.
\item[(f)]  If $T$ is an additive $X$ tree, then $T$ is a reproductive tree and for every
finite $n$ the set $X^{n}$ of level $n$ vertices in the simple tree is the set of level $n$ vertices in $T$.
\end{enumerate}
\end{prop}

\begin{proof}  (a),(b),(d): That the simple tree
is partially additive is obvious.  Condition (\ref{eq6.24}) and (\ref{eq6.6}) imply
(\ref{eq6.26}) which implies that a partially additive $T$ is a subtree of
the simple tree on $X,\alpha$ when $\alpha \geq h(T)$.  Clearly,
for $p$ in the simple tree the successor set is given by
\begin{equation}\label{eq6.27}
S_{p} \ = \ \{ p + q : q \in X^{1} \}.
\end{equation}
So  $p \in T$ and $o(p) + 1 < h(T)$  implies this is a subset of
$T$ by (\ref{eq6.24}). Furthermore, if $p \in T$ and $o(p) < \alpha <
h(T)$, then there exists $q \in T$ with $o(q) = \alpha \setminus
o(p)$ and so $p + q \in T$ by (\ref{eq6.24}).  Hence, $T$ is a normal tree.

(c):  Since $\beta > 1$, $X^{1} \subset T^{\beta}$.  Condition
(\ref{eq6.24}) for $T^{\beta}$ follows from the condition on $T$.

(e),(f):  If $T$ is additive, then the canonical isomorphism (\ref{eq6.8})
restricts to an isomorphism of $T$ with $T_{p}$ for any $p \in T$.
Hence, $T$ is reproductive and so $h(T)$ is tail-like by Lemma \ref{lem5.15}(b).
Since $h(T) > 1$ it is infinite and so by induction on
$n$ using (d), $X^{n} \subset T$ for every finite $n$.

On the other hand, if $h(T)$ is tail-like, then $p,q \in T$ implies
$o(p) + o(q) < h(T)$ and so $p+q \in T$ by (\ref{eq6.24}) and so (\ref{eq6.23})
holds.

In particular, the simple tree $T$ on $X,\alpha$ is additive iff
$\alpha = h(T)$ is tail-like.

\end{proof}\vspace{.5cm}

\begin{cor}\label{cor6.7}  If $X$ is a HLOTS and $T$ is an additive $\Omega$ bounded $X$ tree,
then $X(T)$ is an IHLOTS with completion $\widehat{X(T)}$ a CHLOTS.
\end{cor}

\begin{proof} $T$ is reproductive by Proposition
\ref{prop6.6}(f) and so the result follows from Theorem \ref{theo5.20}.

\end{proof}\vspace{.5cm}

\begin{lem}\label{lem6.8}  Let $T$ be an additive $X$ tree.  We can identify the branch space $X(T)$ with the set
\begin{equation}\label{eq6.28}
\begin{split}
\{x \in X^{\beta} : \beta \  \mbox{is an infinite limit ordinal},
x \not\in T, \\ \mbox{and} \ x|\epsilon \in T \ \mbox{for all} \
\epsilon < \beta \}. \hspace{2cm}
\end{split}\end{equation}
With this identification we have, for $p \in T$ and $x \in
X^{\beta}$
\begin{equation}\label{eq6.29}
x \in X(T) \qquad \Longleftrightarrow \qquad p + x \in X(T).
\end{equation}
For $x \in X^{\beta} $ and $\epsilon < \beta$
\begin{equation}\label{eq6.30}
x \in X(T) \qquad \Longleftrightarrow \qquad x|\epsilon \in T
\quad \mbox{and} \quad \tau_{\epsilon}(x) \in X(T).
\end{equation}
\end{lem}

\begin{proof} Any branch of height less than
$h(T)$ has height an infinite limit ordinal.  Since $T$ is
reproductive, $h(T)$ is also an infinite limit ordinal.  A branch
of height $\beta$ is a coherent collection $\{ x_{\epsilon} \in
X^{\epsilon} : \epsilon < \beta \}$ which fits together to define
an element $x \in X^{\beta}$ such that $x|\epsilon = x_{\epsilon}$
for all $\epsilon < \beta$. Hence the restrictions are all in $T$.
If $x$ itself were in $T$, then we could extend the branch by
adjoining $x$ which violates the maximality condition of the
branch.  Hence, $x\not\in T$.  Conversely, if $x$ lies in the set
described by (\ref{eq6.28}), then $\{ x|\epsilon : \epsilon < \beta \}$ is
a maximal totally ordered set of vertices of $T$ and so is the
corresponding branch.

The characterization of (\ref{eq6.28}) together with (\ref{eq6.23}) easily implies
(\ref{eq6.29}).  Then (\ref{eq6.30}) follows from (\ref{eq6.6}).

\end{proof}\vspace{.5cm}

We now present the inductive construction which shows how all
additive trees are built. We use the translation maps defined by
(\ref{eq6.4}).  Recall that if $\epsilon$ is tail-like and $\alpha <
\epsilon$,  then  $\epsilon = \epsilon \setminus \alpha$ and so
$\tau_{\alpha}$ maps $X^{\epsilon}$ to itself.  If $\epsilon$ is
tail-like and $W \subset X^{\epsilon}$, then we call $W$
\emph{translation invariant}\index{subset!translation invariant} if
\begin{equation}\label{eq6.31}
\tau_{\alpha}(W) \  \subset \ W \qquad \mbox{for all} \ \alpha <
\epsilon.
\end{equation}

A collection of trees $\{ T_{\delta} : \beta \leq \delta <
\epsilon \}$ is called a \emph{coherent collection of trees
indexed by} $[\beta,\epsilon)$  \index{coherent collection} if
\begin{equation}\label{eq6.32}
\begin{split}
h(T_{\delta}) \ = \ \delta \qquad \mbox{for} \quad \beta \leq \delta < \epsilon, \qquad \mbox{and} \\
T_{\delta} = (T_{\r})^{\delta} \qquad \mbox{for} \quad \beta
\leq \delta \leq \r < \epsilon.
\end{split}
\end{equation}

\begin{theo}\label{theo6.9}  Let  $\alpha$ be a tail-like ordinal
and $T$ be a partially additive $X$ tree with $h(T) = \beta < \alpha$.
\begin{enumerate}
\item[(a)]  Let $\epsilon$ be a limit ordinal with $\beta < \epsilon \leq \alpha$.
If  $\{ T_{\delta}  \}$ is a coherent collection of partially
additive $X$ trees indexed by  $[\beta,\epsilon)$ then
\begin{equation}\label{eq6.33}
T_{\epsilon} \ =  \ \bigcup  \{ T_{\delta} : \beta \leq
\delta < \epsilon \}
\end{equation}
is a partially additive $X$ tree.  It defines the unique tree such
that $\{ T_{\delta} : \beta \leq \delta \leq \epsilon \}$ is a
coherent collection of trees indexed by $[\beta,\epsilon] =
[\beta,\epsilon +1)$.
\item[(b)]  Let $\epsilon(\beta) = min \ \{ i : \beta \leq i \leq \alpha $\  and \ $i$ \
 is tail-like $\}$.  For each $\delta$ such that $\beta \leq \delta \leq \epsilon(\beta)$
 there is a unique partially additive $X$ tree $T_{\delta}$ of height $\delta$ and such that
\begin{equation}\label{eq6.34}
(T_{\delta})^{\beta} \ = \ T \qquad \mbox{for} \quad \beta \leq
\delta \leq \epsilon(\beta).
\end{equation}
The collection $\{ T_{\delta} : \beta \leq \delta \leq
\epsilon(\beta) \}$ is a coherent collection of trees indexed by
$[\beta,\epsilon(\beta)]$.  The tree $T_{\epsilon(\beta)}$ is an
additive $X$ tree
with
\begin{equation}\label{eq6.34a}
T_{\ep(\b)} \ = \ \{ p_1 + \dots + p_n : p_i \in T_{\b} \ \ \text{for} \ \ i = 1, \dots, n < \om \}.
\end{equation}
\item[(c)]  Assume that $T$ is an additive $X$ tree so that the height $\beta$ is tail-like.

The set $\tilde{L}_{\beta} = \{ x \in X(T) : h(x) = \beta \}$ is a
translation invariant subset of $X^{\beta}$.

If $\tilde{L}_{\beta} = \emptyset $, then the only subtree
$\tilde{T}$ of the simple tree on $X,\alpha$ which satisfies
$(\tilde{T})^{\beta} = T$ is $\tilde{T} = T$ itself.

If $\tilde{L}_{\beta} \not= \emptyset $  and $W$ is any nonempty,
translation invariant subset of $\tilde{L}_{\beta}$, then
\begin{equation}\label{eq6.35}
T_{\beta +1} \  = \  T \cup \{ p + x : p \in T \quad \mbox{and}
\quad x \in W  \}
\end{equation}
is a partially additive tree of height $\beta +1$ such that
\begin{equation}\label{eq6.36}
(T_{\beta +1})^{\beta} \ = \ T. \hspace{2cm}
\end{equation}

Conversely, if $T_{\beta +1} $ is a partially additive tree of height $\b + 1$ which satisfies
(\ref{eq6.36}), then
 \begin{equation}\label{eq6.36a}
 W \ = \ \{ q|\b : q \in T_{\beta +1} \ \text{with} \ o(q) = \b \}
 \end{equation}
 is a nonempty, translation invariant subset of $\tilde{L}_{\beta}$ and $p + x \in T_{\beta +1}$ for all $p \in T$ and
 $x \in W$.
 \end{enumerate}
\end{theo}

\begin{proof}  (a):  Since each $T_{\delta} $ has
height $\delta$, $T_{\epsilon}$ has height $\epsilon = sup \
\{\delta : \delta < \epsilon \}$ and $(T_{\epsilon})^{\delta} =
T_{\delta}$ is clear for $\beta \leq \delta \leq \epsilon$.

If $p + q \in T_{\epsilon}$, then for some $\delta < \epsilon$ $p +
q \in T_{\delta}$ and so $p, q \in T_{\delta} \subset T_{\epsilon}
$ by partial additivity of $T_{\delta}$. On the other hand, if $p,
q \in T_{\epsilon}$ and $ o(p) + o(q) = o(p + q) < \epsilon$, then
for some $\delta < \epsilon$ $o(p) + o(q) < \delta$.  Hence, $p +
q \in T_{\delta} \subset T_{\epsilon}$.  Thus, $T_{\epsilon}$ is
partially additive.

(b):  By induction on $\r$ with $\beta \leq \r \leq
\epsilon(\beta)$ we construct the trees $T_{\r}$ verifying
(\ref{eq6.34}) and uniqueness and coherence at each stage.

Begin with $T_{\beta} = T$ and assume that $T_{\delta}$ has been
constructed for  $\beta \leq \delta < \r$.\vspace{.25cm}

\textbf{Case 1:} If $\r$ is a limit ordinal, define $T_{\r}$
using (\ref{eq6.33}) with $\r = \epsilon$.  By (a) it is a partially
additive X tree and the collection is coherent, indexed by
$[\beta,\r]$.

If $\tilde{T}$ is a partially additive X tree of height $\epsilon$
with $\tilde{T}^{\beta} = T$, then for $\delta \in [\beta,\epsilon)$
$\tilde{T}^{\delta}$ is partially additive and satisfies
$(\tilde{T}^{\delta})^{\beta} = \tilde{T}^{\beta} = T$.  So by
uniqueness at the $\delta$ level $\tilde{T}^{\delta} = T_{\delta}$.
Since $\r$ is a limit ordinal, $h(\tilde{T}) = \r$
implies
\begin{equation}\label{eq6.37}
\tilde{T} \ = \ \bigcup \{ \tilde{T}^{\delta} \} \ = \  \bigcup \{
T_{\delta} \} \ = \ T_{\epsilon}.
\end{equation}
Uniqueness follows.\vspace{.25cm}

\textbf{Case 2:} If $\r = \delta + 1$, define
\begin{equation}\label{eq6.38}
T_{\r} \ = \ T_{\delta} \cup \{ p + q : p,q \in T_{\delta} \
\mbox{and} \  o(p) + o(q) = \delta \}.
\end{equation}

Since $\beta \leq \delta < \ep(\beta)$,  $\delta$  is not a
tail-like ordinal and so there exist ordinals $i, j < \delta$ with
$i + j = \delta$.  Because $h(T_{\delta}) = \delta $  there exist
$p, q \in T_{\delta}$ with $o(p) = i$ and $o(q) = j$.  Hence, $p +
q \in T_{\r}$ with $o(p + q) = \delta$.  Hence,
$T_{\r}$ has height $\r$ and coherence is clear. In
particular, $(T_{\r})^{\beta} = T$.

If $o(p) + o(q) < \delta$, then $p, q \in T_{\r}$ iff $p, q
\in T_{\delta}$ iff $p + q \in T_{\delta}$ (by partial additivity
of $T_{\delta}$) iff $p + q \in T_{\r}$.  Thus in checking
(\ref{eq6.24}) for $T_{\r}$ we can restrict attention to the case
$o(p + q) = o(p) + o(q) = \delta$ and $o(q) > 0$ so that $o(p) <
\delta$.

If $p, q \in T_{\r}$, then $o(p) < \delta$ implies that $p
\in T_{\delta}$.  If $o(q) < \delta$, then $q \in T_{\delta}$ and
so $p + q \in T_{\r}$ by (\ref{eq6.38}).  If $o(q) = \delta$, then by
(\ref{eq6.38}) $q = \tilde{p} + \tilde{q}$ with $  \tilde{p}, \tilde{q}
\in T_{\delta}$.  Since $o(q) = \delta$, $o( \tilde{q}) > 0$.
\begin{equation}\label{eq6.39}
p + q \ = \ (p +  \tilde{p}) +  \tilde{q}.  \hspace{2cm}
\end{equation}
Since $o( \tilde{q}) > 0$, $o(p +  \tilde{p}) = o(p) + o(
\tilde{p}) < \delta$.  By partial additivity of $T_{\delta}$, $p +
\tilde{p} \in T_{\delta}$.  By (\ref{eq6.39}) and (\ref{eq6.38}) $p + q \in
T_{\r}$.

On the other hand, if $p + q \in T_{\r}$ with $o(p) + o(q) =
\delta$, then by (\ref{eq6.38}) there exist $ \tilde{p},  \tilde{q}  \in
T_{\delta} $ such that
\begin{equation}\label{eq6.40}
p + q \ = \  \tilde{p} + \tilde{q}.  \hspace{2cm}
\end{equation}
Now if $o(p) = i \geq o( \tilde{p})$ and $ \tilde{i} = i \setminus
o( \tilde{p})$, then
\begin{equation}\label{eq6.41}
p \ = \  \tilde{p} + ( \tilde{q}| \tilde{i}) \quad \mbox{and}
\quad q \ = \ \tau_{ \tilde{i}}( \tilde{q}).
\end{equation}
Since $\tilde{p}, \tilde{q}|\tilde{i} \in T_{\delta}$ and $o(p) <
\delta$, $p \in T_{\delta}$ because $T_{\delta}$ is partially
additive. By  (\ref{eq6.26}) $q \in T_{\delta}$ as well.

If, instead, $o(p) = i < o( \tilde{p})$,  then
\begin{equation}\label{eq6.42}
p \ = \ \tilde{p}|i \qquad \mbox{and} \qquad  q \ = \
\tau_{i}(\tilde{p}) + \tilde{q}.
\end{equation}
By (\ref{eq6.26}) $p, \tau_{i}(\tilde{p}) \in T_{\delta}$.  Since
$\tilde{q} \in T_{\delta} $, $o(q) < \delta $ implies $q \in
T_{\delta}$ because $T_{\delta}$ is partially additive, while
$o(q) = \delta$ implies $q \in T_{\r} $ by (\ref{eq6.42}).

Thus, $T_{\r}$ is partially additive.

On the other hand, if $\tilde{T}$ has height $\r$ and
$\tilde{T}^{\beta} = T$,  then as in the limit case $\tilde{T}^{\delta}
= T_{\delta}$  by uniqueness at the $\delta$ height.  If $p, q \in
T_{\delta}$ with $o(p) + o(q) = \delta < \r $, then $p + q
\in \tilde{T} $ because $\tilde{T} $ is partially additive.  Hence,
$T_{\r} \subset \tilde{T}$.
 On the other hand, if $r \in
\tilde{T}$ with $o(r) = \delta$, then we can choose $ i < \delta $
such that $\delta \setminus i < \delta $, because $\delta$ is not
tail-like.  By (\ref{eq6.6})  $r = r|i + \tau_{i}(r)$.  Because $\tilde{T}$
is partially additive, (\ref{eq6.26}) implies that $r|i, \tau_{i}(r) \in
\tilde{T}^{\delta} = T_{\delta}$.  Hence, $r \in T_{\epsilon}$ by
(\ref{eq6.38}).  Hence, $\tilde{T} \subset T_{\r}$ which proves
uniqueness.\vspace{.25cm}

This completes the induction.  At the final stage,
$T_{\ep(\beta)}$ is partially additive with height
$\ep(\beta)$ tail-like so that it is additive by Proposition
\ref{prop6.6}(e).

Since $T_{\ep}$ is additive and contains $T_{\b}$ it clearly contains any finite sum
$p_1 + \dots + p_n$ with $p_i \in T_{\b}$.

Conversely, if $p \in T_{\ep}$, then $o(p) < \ep$ and we can use Cantor Normal Form, Proposition \ref{propCNF}, to write
$o(p) = \a_1 + \a_2 + \dots + \a_n$ with $\b \ge \a_1 \ge \a_2 \dots \ge \a_n$. Let $\s_0 = 0$ and $\s_i = \a_1 + \dots \a_i$
for $i = 1, \dots, n-1$. $p = p_1 + p_2 + \dots p_n$ where $p_i = (\t_{\s_{i-1}}(p))|\a_i$ for $i = 1, \dots, n$. By additivity, each $p_i \in T_{\ep}$.
Since $o(p_i) = \a_i \le \b$ and $(T_{\ep})^{\b} = T_{\b}$ it follows that each $p_i \in T_{\b}$.

(c):  If $x \in \tilde{L}_{\beta}$ and $i < \beta$, then $\beta =
\beta \setminus i$ and $\tau_{i}(x) \in X(T)$ by (\ref{eq6.30}).  Hence,
$\tau_{i}(x) \in \tilde{L}_{\beta}$.  Thus, $\tilde{L}_{\beta}$ is
translation invariant.

If $\tilde{T}$ is a subtree of the simple tree with
$\tilde{T}^{\beta} = T$ and $p \in \tilde{T} \setminus T$, then
$o(p) \geq \beta$ and so $A_{p} \cap \tilde{T}^{\beta} = A_{p}
\cap T$  is a branch of $T$ with height $\beta$.  So if such a
$\tilde{T}$ exists, then  $\tilde{L}_{\beta} \not= \emptyset$.

Assume $W$ is a nonempty, translation invariant subset of
$\tilde{L}_{\beta}$ and that $T_{\beta +1}$ is defined by (\ref{eq6.35}).
Clearly, (\ref{eq6.36}) holds.  It remains to verify (\ref{eq6.24}).

Consider vertices $p, q$ of the simple tree.  As before it
suffices to consider the case $o(p) + o(q) = \beta$  and  $o(q) >
0$ so that $o(p) < \beta $.

If $p, q \in T_{\beta +1}$, then $o(p) < \beta$  implies $p \in T$.
Since $o(p) + o(q) = \beta$ and $\beta$ is tail-like, $o(q) =
\beta$ and so $q = \tilde{p} + x $ with $\tilde{p} \in T$  and  $x
\in W$.  Since $T$  is additive,  $ p + \tilde{p} \in T $ and so
\begin{equation}\label{eq6.43}
p + q \ = \ (p + \tilde{p}) + x \ \in T_{\beta +1}, \hspace{1.7cm}
\end{equation}
by definition (\ref{eq6.35}).

On the other hand, if $p + q \in T_{\beta +1}$, then $o(p + q) =
o(p) + o(q) = \beta$ implies $p + q = \tilde{p} + x $ for some
$\tilde{p} \in T$ and $x \in W$.

As with (\ref{eq6.41}), $o(p) \geq o(\tilde{p})$ and $\tilde{i} = i
\setminus o(p)$ implies
\begin{equation}\label{eq6.44}
p \ = \  \tilde{p} + ( x| \tilde{i}) \quad \mbox{and} \quad q \ =
\ \tau_{ \tilde{i}}( x).
\end{equation}
Since $\tilde{i} < \beta$, $x|\tilde{i} \in T$ by Lemma \ref{lem6.8} and so
$p \in T$ by additivity.  Because $W$ is translation invariant $q
\in W \subset T_{\beta +1}$.

As in (\ref{eq6.42}), $ o(p) = i < o(\tilde{p})$ implies
\begin{equation}\label{eq6.45}
p \ = \ \tilde{p}|i \qquad \mbox{and} \qquad  q \ = \
\tau_{i}(\tilde{p}) + x.
\end{equation}
Hence, $p, \tau_{i}(\tilde{p}) \in T$ and so $q \in T_{\beta +1}$.

This completes the proof of partial additivity for $T_{\beta +1}$

For the converse, it is clear that $W$ of (\ref{eq6.36a}) is a subset of $\tilde{L}_{\beta}$.
It is nonempty since $T_{\b +1}$ has height $\b$.
For $q \in T_{\b + 1}$ with
$o(q) = \b$ and $i < \b$, $\tau_i(q) \in T_{\b + 1}$ by partial additivity and $o(\tau_i(q)) = \b \setminus i = \b$.
So $\tau_i(q)|\b = \tau_i(q|\b) \in W$. Thus, $W$ is translation invariant.
If $p \in T$ and $x \in T_{\b +1}$ with $o(x) = \b$, then by partial additivity $o(p + x) = o(p) + \b = \b$ implies $p + x \in T_{\b +1}$.

\end{proof} \vspace{.5cm}

\noindent{\bfseries Remark.} Recall that if $T$ is a tree of height $\alpha$
a countable limit ordinal, then $\tilde{L}_{\alpha} \not=
\emptyset$ by Lemma \ref{lem5.4}(b).  Thus, the extension process of part (c)
can go beyond any countable level.  On the other hand, an
additive tree of height $\Omega$ is $\Omega$-bounded precisely when it
cannot be continued to level $\Omega + 1$. \vspace{.5cm}

\begin{cor}\label{cor6.9a} There exists a LOTS $X$ which is weakly homogeneous, first countable and countably complete,
but which is not complete. \end{cor}

\begin{proof} Because an additive tree is reproductive, it suffices by Theorem \ref{theo5.21}
 to construct an additive tree of $\Q$ type which has height $\Omega$ and which contains
branches of countable height and branches of height $\Omega$.

Let $\bar 0 \in \Q^{\Omega}$ with $\bar 0_i = 0 \in \Q$ for every $i < \Omega$. For $p \in \Q^{\a}$ with
$\a$ a limit ordinal, we say that $p$ \emph{eventually equals $0$ } \index{eventually equals $0$ }
if there exists $\b < \a$ such that $p_i = 0$ for
all $i$ with $\b \le i < \a$.

We use that inductive construction of Theorem \ref{theo6.9} to build a coherent collection $T_{\a}$ of partially
additive trees such that $\bar 0|\b \in T_{\a}$ for all $\b < \a$. This condition is clearly preserved in steps (a) and (b).
For the choice step (c) with $\b$ tail-like, we have that $\bar 0|\b \in \tilde L_{\b}$. The set $W = \{ \bar 0|\b \}$ is
translation invariant and we let $T_{\b + 1} = T \cup \{ p + (\bar 0|\b ) : p \in T \}$. Thus, $L_{\b}$ consists of the
elements of $\tilde L_{\b}$ which eventually equal $0$.

With $ T = \bigcup T_{\a}$ we identify $X(T)$ as in Lemma \ref{lem6.8}. So $\bar 0 \in X(T)$ with $h(\bar 0) = \Omega$. On the
other hand suppose $x \in \Q^{\a}$ with $\a < \Omega$ a tail-like ordinal and $x$ is not eventually $0$ and $x|\b \in T$ for all
$\b < \a$, then $x$ is a branch with height $\a$, i.e. it is an element of $\tilde L_{\a}$ which is not in $L_{\a}$. In particular,
if $x \in \Q^{\om}$ is not eventually $0$, then $x \in X(T)$ with $h(x) = \om$.

\end{proof} \vspace{.5cm}

This result answers a question raised by Babcock \cite{B} Section 2. A closed, bounded interval in a LOTS $X$ given by Corollary \ref{cor6.9a}
satisfies his Linear Homogeneity Condition 2, but not his Linear Homogeneity Condition 3.

There is a different way of looking at additive trees.

Let $\tilde \Omega$ be the set of infinite tail-like ordinals in $\Omega$ so that
\begin{equation}\label{eq6.new01}
\tilde \Omega \ = \ \{ \om^{\g} : 0 < \g < \Omega \}.
\end{equation}
\index{$\tilde \Omega$}

We write $\ep(s,t)$ for the \emph{equality level}\index{equality level} for a pair $s, t \in X^{\Omega}$ so that with $\ep = \ep(s,t)$
\begin{equation}\label{eq6.new02}
s_i = t_i \ \ \text{for} \ \ i < \ep, \quad \text{and} \quad s_{\ep} \not= t_{\ep}.
\end{equation}
Thus, $s|\ep = t|\ep$ and $s|(\ep + 1) \not= t|(\ep +1)$.
\vspace{.25cm}

\begin{df}\label{dfnew01} Let $H$ be a function from $X^{\Omega}$ to $\Omega \setminus 1$. We call $H$ a \emph{height function} \index{height function}
when it satisfies:
\begin{itemize}
\item[(i)] For all $s, t \in X^{\Omega}, \ep(s,t) \ge H(s)$ implies $H(t) = H(s)$.
\end{itemize}
We call it an $\tilde \Omega$ height function when it takes values in $\tilde \Omega$.

The function $H$ is called an \emph{additive height function} \index{height function!additive} when it is a $\tilde \Omega$ height function which
also satisfies:
\begin{itemize}
\item[(ii)] For all $s \in X^{\Omega}, \a < H(s)$ implies $H(\t_{\a}(s)) = H(s)$.
\end{itemize}
\end{df}\vspace{.5cm}

\begin{lem}\label{lemnew02} Assume that $H$ is an additive height function. Let $s, t \in X^{\Omega}, p \in X^{\a}$.
\begin{itemize}
\item[(a)]If $p = t|\a$ with $\a < H(t)$, then $\a < H(p + s)$.

\item[(b)] If $\a < H(s)$, then
$\a < H(p + s)$ implies $H(p + s) = H(s)$.
\end{itemize} \end{lem}

\begin{proof} (a): $\ep(p + s,t) \ge \a$, so $H(p + s) \le \a$ would imply, by (i),  $H(p + s) = H(t)$ which is larger than $ \a$.

(b): By (ii), $\a < H(p + s)$ implies that $H(s) = H(\t_{\a}(p + s)) = H(p + s)$.

\end{proof} \vspace{.5cm}

\begin{theo}\label{theonew03} If $T$ is an $\Omega$-bounded subtree of the simple tree on $X, \Omega$, then the associated height
function is given by
 \begin{equation}\label{eq6.new03}
 H(s) \ = \  min \{ \a : s|\a \not\in T \}.
 \end{equation}

 Conversely, if $H$ is a height function, then $T = \{ s|\a :\a < H(s) \}$ is an $\Omega$-bounded (not necessarily semi-normal)
 subtree of the simple tree on $X, \Omega$. The branch space is given by
 \begin{equation}\label{eq6.new04}
 X(T) = \{ s|H(s) : s \in X^{\Omega} \}.
 \end{equation}

 The tree $T$ is an additive, $\Omega$-bounded subtree of  the simple tree on $X, \Omega$ iff the associated height function $H$ is
 an additive height function.\end{theo}

 \begin{proof} If $T$ is an $\Omega$-bounded subtree, and $\a < \b < \Omega$, then  $s|\a \not\in T $ implies $s|\b \not\in T$.
 If $T$ is  $\Omega$-bounded then for every $s \in X^{\Omega}$ $s|\a \not\in T$ for some $\a < \Omega$.  It follows that $H$ defined by
 (\ref{eq6.new03}) is a height function.

 Conversely, if $H$ is a height function, then $s|\a = t|\a$ implies $\ep(s,t) \ge \a$.  Hence, if $\a \ge H(s)$, then $H(s) = H(t)$ and so
 $\a \ge H(t)$. It follows that $T$ is a subtree with associated height function $H$.  The description of the branch space follows as in
 Lemma \ref{lem6.8}.

 If $T$ is additive, then the height of every branch is an infinite, tail-like ordinal. From the description (\ref{eq6.new04}) it follows
 that $H$ is an $\tilde \Omega$ height function. Furthermore, if $x = s|H(s)$ is a branch and $\a < H(s)$, then by Lemma \ref{lem6.8}
 $\t_{\a}(x) = (\t_{\a} (s))|H(s)$ is a branch and so $H(\t_{\a} (s)) = H(s)$.

 Now let $H$ be an additive height function and  that $p = t|o(p), q = s|o(q)$ with $s,t \in X^{\Omega}$ and $o(p), o(q) < \Omega$.

 First, assume  $p, q \in T$, so that $o(p) < H(t), o(q) < H(s)$.

 By Lemma \ref{lemnew02} (a) $H(p + s) > o(p)$ and so by Lemma \ref{lemnew02} (b) $H(p + s) = H(s)$.
 Similarly, $H(q + t) > o(q)$ and  $H(q + t) = H(t)$. Finally, $H(p + q + t) = H(q + t) = H(t)$.

 If $o(p) < H(s)$, then since $H(s)$ is tail-like, $o(p + q) = o(p) + o(q) < H(s) = H(p + s)$. Hence,
 $p + q = (p + s)|(o(p) + o(q)) \in T$.

 If, instead, $o(p) \ge H(s)$, then $H(t) > o(p) \ge H(s) > o(q)$ implies that $o(p) + o(q)$ is less than the
 tail-like ordinal $H(t) = H(p + q + t)$. Hence, $p + q = (p + q + t)|(o(p) + o(q)) \in T$.

Conversely, assume that $p + q \in T$ so that $p + q = r|(o(p) + o(q))$ with $o(p) + o(q) < H(r)$.
Since $o(p) \le o(p) + o(q) < H(r)$, $p = r|o(p) \in T$. By condition (ii) $H(\t_{o(p)}r) = H(r)$.
So $o(q) \le o(p) + o(q) < H(r)$ implies $q =  (\t_{o(p)}(r))|o(q) \in T$.

It follows that the tree associated to the additive height function $H$ is additive.

\end{proof} \vspace{.5cm}

 \begin{prop}\label{propnew04} Let $H_0 : X^{\Omega} \to \Omega \setminus 1$ be an arbitrary function. There exists $H$ a maximum
 height function with $H \le H_0$.  Furthermore, if $H_0$ takes values in $\tilde \Omega$, then $H$ is an $\tilde \Omega$ height function.
 \end{prop}

 \begin{proof} Let $\{ H^i \}$ be the set of all height functions with $H^i < H_0$. This set is nonempty since the  constant function with value $1$ is
 in it.  Define $H(s) = \sup_i H^i(s)$.

 If $\ep(s,t) \ge H(s)$ then $\ep(s,t) \ge H^i(s)$ for all $i$ and so $H^i(s) = H^i(t)$ since each $H^i$ is a height function. Hence,
 $H(t) = H(s)$. That is, $H$ is a height function.

 We can describe $H$ by a finitely inductive construction. For $n \ge 0$, define
 \begin{equation}\label{eq6.new05}
 H_{n+1}(s) =  min \{ H_n(t) : \ep(s,t) \ge min (H_n(s),H_n(t)) \}.
 \end{equation}
 In particular, $H_{n+1}(s) \le H_n(s)$.

 Define $H_{\infty}(s) = min_n H_n(s)$. By well-ordering, there exists for each $s$ a positive integer $N_s$ such that
 $H_{\infty}(s) = H_k(s)$ for all $k \ge N_s$.

 If $H_n \ge H^i$ for a height function $H^i$, then $\ep(s,t) \ge min (H_n(s),H_n(t)) \ge min (H^i(s),H^i(t))$ implies
 $H_n(t) \ge H^i(t) = H^i(s)$. Hence, $H_{n+1}(s) \ge H^i(s)$. It follows that $H_{\infty} \ge H$.

 Now suppose $\ep(s,t) \ge H_{\infty}(s)$ which equals $ H_k(s) \ge min (H_k(s),H_k(t))$ for $k \ge max (N_s,N_t)$.
 It follows that $H_{k+1}(s) \le H_k(t)$ and $H_{k+1}(t) \le H_k(s)$.
 Since $k \ge max (N_s,N_t)$, $H_k(s) = H_{k+1}(s) = H_{\infty}(s)$ and $H_k(t) = H_{k+1}(t) = H_{\infty}(t)$.
 Thus, $H_{\infty}(s) = H_{\infty}(t)$. This means that $H_{\infty}$ is a height function and so equals $H$.

 Finally, if $H_0$ takes values in $\tilde \Omega$, then so does each $H_n$ and so $H_{\infty} = H$ is an
 $\tilde \Omega$ height function.

 \end{proof} \vspace{.5cm}

 The difficulty with the inductive construction of Theorem \ref{theo6.9} is that we don't have a convenient method for  making the choices along
 the way so that the resulting tree is $\Omega$-bounded.

 There is an inconvenient method. If we begin with any function from $X^{\Omega}$ to
 $\tilde \Omega$, the construction of Proposition \ref{propnew04} yields an $\tilde \Omega$ height function with associated $\Omega$-bounded tree
 $R$ having all branch heights tail-like. Any $\Omega$-bounded additive tree is such a tree.  On the other hand, beginning with such a tree $R$ we
 can use the construction of Theorem \ref{theo6.9} to build an additive tree contained in $R$. The process will terminate at a countable
 tail-like ordinal $\b$ if either we cannot choose $W$ as a nonempty,
translation invariant subset of $\tilde{L}_{\beta}$, so that
$ \{ p + x : p \in T \  \mbox{and} \  x \in W  \}  \subset R$, or else, if such a choice is possible but the extension to the successor
 tail-like ordinal given by (\ref{eq6.34a}) is not contained in $R$. If the process does not terminate at a countable level, then we obtain an additive
 tree of height $\Omega$, but which is $\Omega$-bounded since it is contained in $R$.

\vspace{1cm}

\subsection{Special Trees for HLOTS}

For a HLOTS $X$ there is a special class of trees whose
construction parallels that of the additive trees.

For a HLOTS  $X$ the completion $\hat{X}$ is a CHLOTS and $\bullet
\hat{X} \bullet = \{m\} + \hat{X} + \{M\} $ is its two point
compactification.  If $ p \in X^{\alpha} $ for any positive
ordinal $\alpha$, define $\hat{p} \in ( \bullet \hat{X} \bullet
)^{\alpha +1}$ \index{$\hat{p}$} by
\begin{equation}\label{eq6.46}
\hat{p}(0) \  = \  m, \qquad \hat{p}(i) \ = \  sup \  \{ p(j) : j
< i \} \quad \mbox{for} \ 0 < i \leq \alpha.
\end{equation}

We say that $p \in X^{\alpha} $ is \emph{sharply increasing}\index{sharply increasing} when
\begin{equation}\label{eq6.47}
\hat{p}(i) \ < \ p(i) \quad \mbox{for} \ i < \alpha.  \hspace{1cm}
\end{equation}
Notice that if $p : \alpha \rightarrow X$ is an order map, then
\begin{equation}\label{eq6.48}
\hat{p}(i + 1) \ = \ p(i)  \quad \mbox{for} \ i < \alpha,
\hspace{1cm}
\end{equation}
and if $p$ is injective, then $p(i) < p(i + 1)$.  That is, if $p$
is any order injection then (\ref{eq6.47}) holds for any successor ordinal
$i < \alpha$.  On the other hand, if $p$ is an order embedding,
then $\hat{p}(i) = p(i)$ for any limit ordinal $i < \alpha$. Thus,
$p$ is sharply increasing exactly when it is an order injection
which is discontinuous at each limit ordinal.  In general, $\hat{p}$ is a continuous  order map and it is an order embedding
into $\bullet \hat{X} \bullet $  if $p$
is an injective order map.

We define the \emph{order tree on}\index{order tree}\index{tree!order} a HLOTS $X$, denoted $T(X)$,
whose vertices at level $\alpha$ are the bounded, sharply
increasing maps from $\alpha$ to $X$.  That is, we define
\begin{equation}\label{eq6.49}
L_{\alpha}(X) \ = \  \{ p \in X^{\alpha} : \hat{p}(i) < p(i)
\ \mbox{for} \  i < \alpha \ \mbox{and} \ \hat{p}(\alpha) < M \}.
\end{equation}

\begin{theo}\label{theo6.10}  For a HLOTS $X$ the order tree $T(X)$ is a subtree of the simple tree
on $X,\Omega$.  $T(X)$ is a reproductive tree of $X$ type.  The height of $T(X)$ is
$\Omega$ but every branch has countable height, i.e. $T(X)$ is $\Omega$-bounded.  We can identify the branch space with the set
\begin{equation}\label{eq6.50}
\begin{split}
\{ x \in X^{\beta} : \beta \ \mbox{is a countable, limit
ordinal},\hspace{.5cm} \\  \hat{x}(i) < x(i) \quad \mbox{for} \ i
< \beta \quad \mbox{and} \quad \hat{x}(\beta) = M \}.
\end{split}
\end{equation}
The branch space of $T(X) $ is an IHLOTS and its completion is a
CHLOTS.
\end{theo}

\begin{proof} Since $X$ is of countable type,
only countable ordinals admit order injections into $X$. Hence,
$L_{\alpha} = \emptyset $ if $\alpha$ is uncountable.  If $p \in
L_{\alpha}$ and $\beta < \alpha$, then $p|\beta \in L_{\beta}(X)$
because for any $p \in X^{\alpha}$
\begin{equation}\label{eq6.51}
\hat{p}|(\beta + 1) \  = \  (\widehat{p|\beta}). \hspace{2cm}
\end{equation}
Thus, $T(X)$ is a subtree of the simple tree on $X,\Omega$.

If $p \in L_{\alpha}$ and $s = \hat{p}(\alpha)$, then because
$X$ is a HLOTS there exists an order isomorphism $f_{s} : X
\rightarrow (s,M) \cap X $ and it extends to the isomorphism
$\hat{f}_{s} : \bullet \hat{X} \bullet  \rightarrow [s,M]$.  For
any $q \in X^{\beta}$
\begin{equation}\label{eq6.52}
\hat{f}_{s} \circ \hat{q} \ = \  (\widehat{f_{s} \circ q}) \quad \text{on} \ \b \setminus 1
\hspace{2cm}
\end{equation}
and so $q \in L_{\beta}$ iff $f_{s} \circ q \in L_{\beta}$.

Now define the analogue of the canonical inclusion of (\ref{eq6.8})
\begin{equation}\label{eq6.53}
\begin{split}
f_{p} : T(X) \rightarrow T(X)_{p}  \hspace{2cm} \\
f_{p}(q) \  = \  p + (f_{s} \circ q). \hspace{2cm}
\end{split}
\end{equation}
Notice that
\begin{equation}\label{eq6.54}
\hat{p}(\alpha) \  =  \  s  \  =  \hat{f}_{s}(m) \  < \
f_{s}(q(0)) \hspace{.5cm}
\end{equation}
and so $p + (f_{s} \circ q) \in L_{\alpha + \beta}$.

Conversely, $r \in T(X)_{p}$ implies
\begin{equation}\label{eq6.55}
\hat{r}(\alpha) \  =  \  \hat{p}(\alpha) \ = \ s  \hspace{.5cm}
\end{equation}
and so $f_{s}^{-1} \circ (\tau_{\alpha}(r)) \in T(X)$.  Thus,
$f_{p}$ is a tree isomorphism and so $T(X)$ is reproductive.

Since $L_{1} = X^{1} \cong X $ we see that the tree is of $X$
type because it is reproductive.

The identification of the branch space with the set described in
(\ref{eq6.50}) is now routine using an argument similar to that of Lemma \ref{lem6.8}.
Because $X$ is of countable type it follows that if $x \in
X^{\beta}$ is a branch, then $\beta$ is countable.  Hence, every
branch has countable height.

It follows from Theorem \ref{theo5.20} that the branch space is an IHLOTS
and its completion is a CHLOTS.

It remains to show that the height of $T(X)$ is $\Omega$, i.e.
$L_{\alpha} \not= \emptyset $ for every countable ordinal
$\alpha$.  We use a construction which we will apply again later.

By Proposition \ref{prop2.8}(d) there exists an order embedding $\tilde{p} :
\alpha + 1 \rightarrow X$.  Choose for each $i < \alpha$, an
isomorphism $g_{i} : X \rightarrow (\tilde{p}(i),\tilde{p}(i+1))$.
For any $z \in X^{\alpha} $ define
\begin{equation}\label{eq6.56}
p(z)(i) \ = \ g_{i}(z(i)). \hspace{2cm}
\end{equation}
It is easy to see that $p(z)$ is an order injection and for each
limit ordinal $i \leq \alpha$
\begin{equation}\label{eq6.57}
\widehat{p(z)}(i) \  = \  \tilde{p}(i)  \hspace{2cm}
\end{equation}
and so if $i < \alpha$, $\widehat{p(z)}(i) < p(z)(i) $.  Thus, $z
\mapsto p(z)$ is an injective map from $X^{\alpha}$ into
$L_{\alpha}$.

\end{proof} \vspace{.5cm}

\noindent{\bfseries Remark.} We can apply this last construction as well
when $\tilde{p} : \alpha + 1 \rightarrow \bullet \hat{X} \bullet $
is an order embedding with
\begin{equation}\label{eq6.58}
\tilde{p}(0) \  = \  m \qquad \mbox{and} \qquad \tilde{p}(\alpha)
\  =  \  M.
\end{equation}
If $\alpha$ is a limit ordinal, then $z \mapsto p(z)$ defines an
order injection from $X^{\alpha}$ into the branch space of $T(X)$
as identified in (\ref{eq6.50}). \vspace{.5cm}

In order to define the analogue for $T(X)$ of additive subtrees,
we need a piece of auxiliary equipment.

\begin{df}\label{df6.11}  For a HLOTS $X$ a set $\mathcal{S}$ of maps from $X$ to $X$ is called a
\emph{special semigroup}\index{special semigroup} when it satisfies the following conditions
\begin{enumerate}
\item[(i)] Each $f \in \mathcal{S}$ is either an order isomorphism  $ f : X \rightarrow X $
 or an order isomorphism $f : X \rightarrow (x,\infty)$ with $x \in X$.  The former are called
  the \emph{invertible} elements of $\mathcal{S}$ and included among them is the identity $1_{X}$.
\item[(ii)]  If $f,g \in \mathcal{S}$, then $f \circ g \in \mathcal{S}$.  If $f$ is
 an invertible element of $\mathcal{S}$, then $f^{-1} \in \mathcal{S}$.  Thus, the
 invertible elements of $\mathcal{S}$ form a group under composition.
\item[(iii)]  The group of invertible elements acts transitively on $X$.
\item[(iv)]  The set of noninvertible elements of $\mathcal{S}$ is nonempty.
\end{enumerate}
\end{df}
\vspace{.5cm}

Notice that each $f \in \mathcal{S}$ is an order embedding by
Proposition \ref{prop2.1}(b).

\begin{lem}\label{lem6.12}  If $\mathcal{S}$ is a special semigroup for a HLOTS $X$, $x < y $
in $X$ and $z \in X$, then there exists $ f \in \mathcal{S}$ such
that $f(X) \ = \ (x,\infty)$ and $f(z) = y$.
\end{lem}

\begin{proof} There exists $f_{1} :X \rightarrow
(\tilde{x},\infty)$ in $\mathcal{S}$ for some $\tilde{x} \in X$ by
condition (iv).  By (iii) there exists an invertible element
$f_{2} $ such that $f_{2}(\tilde{x}) = x$.  Let $\tilde{y} =
(f_{1}^{-1} \circ f_{2}^{-1}) (y) $.
By (iii) again there exists an
invertible element $f_{3}$ such that $f_{3}(z) = \tilde{y}$.  The
composite $ f = f_{2} \circ f_{1} \circ f_{3} $ is the required
element of $\mathcal{S}$.

\end{proof} \vspace{.5cm}

A subset $W $ of $X^{\beta}$ is called
$\mathcal{S}$ \emph{invariant}\index{$\mathcal{S}$ invariant}\index{subset!$\mathcal{S}$ invariant} if for all $f \in \mathcal{S}$
\begin{equation}\label{eq6.59}
x \in W \qquad \Longleftrightarrow \qquad f \circ x \in W.
\end{equation}

\begin{df}\label{df6.13}  Let $\mathcal{S}$ be a special semigroup for a HLOTS $X$ and let $T$ be a
subset of the order tree $T(X)$.  $T$ is called an $\mathcal{S}$ \emph{tree}\index{tree!special semigroup $\mathcal{S}$}
 when it satisfies the following conditions.
\begin{enumerate}
\item[(i)]  The root $0$ and the level $1$ vertices $X^{1}$ of $T(X)$ are contained in $T$.
\item[(ii)]  For every $p \in T$ with $o(p) = \alpha$
\begin{equation}\label{eq6.60}
\hat{p}(i) \in X \quad \mbox{for} \ 1 \leq i \leq \alpha.
\hspace{1cm}
\end{equation}\item[(iii)]  For all $f \in \mathcal{S}$ and $p \in T(X)$
\begin{equation}\label{eq6.61}
p \in T \qquad \Longleftrightarrow \qquad f \circ p \in T.
\end{equation}
\item[(iv)]  For all $p,q \in T(X)$ with $o(p) = \alpha$, $o(q) = \beta$
\begin{equation}\label{eq6.62}
p + q \in T   \qquad \Longleftrightarrow \qquad  p,q \in T, \
\hat{p}(\alpha) < q(0), \ \alpha + \beta < h(T) .
\end{equation}
\end{enumerate}
\end{df}
\vspace{.5cm}

\begin{prop}\label{prop6.14}  Let $T$ be an $\mathcal{S}$ tree with $\mathcal{S}$ a special semigroup for a HLOTS $X$.
\begin{enumerate}
\item[(a)]  $T$ is a normal subtree of $T(X)$.
\item[(b)]  If $p \in T$ and $\beta \leq o(p) $,  then
\begin{equation}\label{eq6.63}
p|\beta  \in T \qquad \mbox{and} \qquad  \tau_{\beta}(p) \in T.
\hspace{.5cm}
\end{equation}
If $o(p) + 1 < h(T)$,  then the successor set $S_{p}$ in $T$
consists of all of the successors of $p$ in $T(X)$.
\item[(c)]  $T$ is a tree of $X$ type and $T$ is reproductive iff $h(T)$ is tail-like.
\item[(d)]  If $ 1 < \beta \leq h(T)$, then $T^{\beta}$ is an $\mathcal{S}$ tree.
\item[(e)]  We can identify the branch space of $T$ with the set
\begin{equation}\label{eq6.64}
\begin{split}
\{ x \in X^{\beta} : \beta = h(T) \ \mbox{or an infinite limit
ordinal,} \hspace{.5cm} \\ x \not\in T \quad \mbox{and} \quad
x|\gamma \in T \ \mbox{for all} \  \gamma < \beta \}.
\end{split}
\end{equation}
\item[(f)]  Assume that $h(T)$ is tail-like.  For $p, q \in T(X)$  with $o(p) = \alpha$,
\begin{equation}\label{eq6.65}
p + q \in T \qquad \Longleftrightarrow \qquad p, q \in T \quad
\mbox{and} \quad \hat{p}(\alpha) < q(0).
\end{equation}
Furthermore, if $x \in X^{\beta}$ with $\beta > 0$, then
\begin{equation}\label{eq6.66}
p + x \in X(T)  \qquad \Longleftrightarrow \qquad p \in T, x \in
X(T) \quad \mbox{and} \quad \hat{p}(\alpha) < x(0),
\end{equation}
and so if $\gamma < \beta$, then
\begin{equation}\label{eq6.67}
x \in X(T) \qquad \Longrightarrow \qquad x|\gamma \ \in T, \  \tau_{\gamma}(x)  \in  X(T).
\end{equation}
\item[(g)]  If $h(T)$ is tail-like, then the branch space $X(T)$ is an IHLOTS and its
completion $\hat{X}(T)$ is a CHLOTS.
\end{enumerate}
\end{prop}

\begin{proof} (a), (b), (d)  and (\ref{eq6.63}) follow from
condition (iv) of Definition \ref{df6.13} and (\ref{eq6.6}). Furthermore, (\ref{eq6.63}) implies that
$T$ is a subtree of $T(X)$.  For $p \in L_{\alpha}$ the
successor set in $T(X)$ is given by
\begin{equation}\label{eq6.68}
S_{p} \  = \  \{ p + q : q \in X^{1} \quad \mbox{and} \quad
\hat{p}(\alpha) < q(0) \}.
\end{equation}
If $p \in T$ and $o(p) + 1 < h(T)$ , then conditions (i) and (iv)
imply that this set is contained in $T$.  Normality of the tree is
clear from this and condition (iv).

With $\beta > 1$ to retain condition (i), (d) is obvious.

(e) and (f) follow as in Lemma \ref{lem6.8}.  When $h(T)$ is tail-like $p, q
\in T$ implies $o(p) + o(q) < h(T)$ and $x \in X(T)$ implies $o(p)
+ h(x) \leq h(T)$.

(c),(g):  If $T$ is reproductive, then $h(T)$ is tail-like.  For the
converse let $p \in T$ with $o(p) = \alpha$.  By condition (ii) $s
= \hat{p}(\alpha) \in X$.  By Lemma \ref{lem6.12} there exists $f_{s} :X
\rightarrow (s,\infty)$ in $\mathcal{S}$.  By (\ref{eq6.65}) and condition
(iii) the map $f_{p}$ of (\ref{eq6.53}) restricts to a tree isomorphism of
$T$ onto $T_{p}$.  Thus, $T$ is reproductive.  Since $T \subset
T(X)$, $h(x) < \Omega$ for all $x \in X(T)$.  Hence, (g) follows
from Theorem \ref{theo5.20}.

\end{proof}\vspace{.5cm}

When $\hat{X} = X$, i.e. $X$ is a CHLOTS, then the order tree
$T(X)$ is an $\mathcal{S}$ tree for any special semigroup
$\mathcal{S}$.  However, if $X$ is an IHLOTS, then condition (ii)
fails for $T(X)$.  In that case, the \emph{reduced order tree} \index{reduced order tree}
\index{tree!order!reduced}
is
defined to be
\begin{equation}\label{eq6.69}
\{ p \in T(X) : \hat{p}(i) \in X \quad \mbox{for all} \ 1 \leq i
\leq o(p) \}
\end{equation}
is an $\mathcal{S}$ tree for any special semigroup $\mathcal{S}$
and it contains all other $\mathcal{S}$ trees.

In general,  $T(X)^{2} = \{0\} \cup X^{1} $ is the unique
$\mathcal{S}$ tree of height $2$ for any special semigroup and
from $T(X)^{2}$ we can build all $\mathcal{S}$ trees by using an
inductive construction completely analogous to that of Theorem \ref{theo6.9}.

\begin{theo}\label{theo6.15}  Let $\mathcal{S}$ be a special semigroup for a HLOTS $X$,
 $\alpha \leq \Omega$ be a tail-like ordinal and  $T$
be an $\mathcal{S}$ tree with $h(T) = \beta < \a$.
\begin{enumerate}
\item[(a)]  Let $\epsilon$ be a limit ordinal with $\beta < \epsilon \leq \alpha$.
If  $\{ T_{\delta}  \}$ is a coherent collection of $\mathcal{S}$
trees indexed by  $[\beta,\epsilon)$, then
\begin{equation}\label{eq6.70}
T_{\epsilon} \ =  \ \bigcup  \{ T_{\delta} : \beta \leq
\delta < \epsilon \}
\end{equation}
is an $\mathcal{S}$ tree.  It defines the unique tree such that
$\{ T_{\delta} : \beta \leq \delta \leq \epsilon \}$ is a coherent
collection of trees indexed by $[\beta,\epsilon] = [\beta,\epsilon
+1)$.
\item[(b)]  Let $\epsilon(\beta) = min \ \{ i : \beta \leq i \leq \alpha $\  and \ $i$ \
is tail-like $\}$.  For each $\delta$ such that $\beta \leq \delta \leq \epsilon(\beta)$
there is a unique $\mathcal{S}$ tree such that
\begin{equation}\label{eq6.71}
(T_{\delta})^{\beta} \ = \ T \qquad \mbox{for} \quad \beta \leq
\delta \leq \epsilon(\beta).
\end{equation}
The collection $\{ T_{\delta} : \beta \leq \delta \leq
\epsilon(\beta) \}$ is a coherent collection of trees indexed by
$[\beta,\epsilon(\beta)]$.
\item[(c)]  Assume that $h(T) = \beta$ is tail-like. Using the description (\ref{eq6.46})
 we define $\ddot{L}_{\beta} = \{ x \in X(T) : h(x) = \beta \ and \ \hat{x}(\beta) \in X \} $
 \index{$\ddot{L}_{\a}$},
 a subset of $\tilde{L}_{\beta} = \{ x \in X(T) : h(x) = \beta \}$.  Both $\ddot{L}_{\beta}$
 and $\tilde{L}_{\beta}$ are nonempty, $\mathcal{S}$ invariant, translation invariant subsets of $X^{\beta}$.

If  $W$ is any nonempty, $\mathcal{S}$ invariant, translation
invariant subset of $\ddot{L}_{\beta}$, then
\begin{equation}\label{eq6.72}
T_{\beta +1} \  = \  T \cup \{ p + x : p \in T, x \in W, o(p) =
\alpha  \  \mbox{and} \  \hat{p}(\alpha) < x(0) \}
\end{equation}
is an $\mathcal{S}$ tree of height $\beta +1$ such that
\begin{equation}\label{eq6.73}
(T_{\beta +1})^{\beta} \ = \ T. \hspace{2cm}
\end{equation}
\end{enumerate}
\end{theo}

\begin{proof}  The proof parallels that of
Theorem \ref{theo6.9} with a few adjustments which we will note.

(a):  Conditions (ii) and (iii) are obviously preserved in the
union.  Condition (iv) and uniqueness are proved as in Theorem
\ref{theo6.9}.

(b):  In the inductive construction the limit stage follows from
(a) as before and the uniqueness arguments in each case are
completely analogous.

If $\epsilon = \delta + 1 $, then we define
\begin{equation}\label{eq6.74}
\begin{split}
T_{\epsilon} \ = \ T_{\delta} \cup \{ p + q : p,q \in T_{\delta},
o(p) = \alpha, o(q) = \beta, \\  \hat{p}(\alpha) < q(0) \
\mbox{and} \ \alpha + \beta = \delta \}.\hspace{1cm}
\end{split}
\end{equation}
Condition (ii) is clear and condition (iv) is checked as was
partial additivity before.

If $f \in \mathcal{S}$ and $p, q \in T(X)$, then
\begin{equation}\label{eq6.75}
f \circ (p + q ) \quad = \quad (f \circ p) + (f \circ q),
\hspace{.5cm}
\end{equation}
and with $\hat{f} : \bullet \hat{X} \bullet \rightarrow  \bullet
\hat{X} \bullet $ the extension of $f$
\begin{equation}\label{eq6.76}
(\hat{f} \circ \hat{p}) \quad = \quad \widehat{f \circ p},
\hspace{1.5cm}
\end{equation}
except at $0$.

These let us verify condition (iii) for $T_{\epsilon}$ and so
complete the inductive construction.

(c): From (\ref{eq6.64}) we have for $f \in \mathcal{S}$ that
\begin{equation}\label{eq6.77}
x \in X(T) \qquad \Longleftrightarrow \qquad f \circ x \in X(T).
\hspace{.5cm}
\end{equation}
Notice that here we use the fact that
\begin{equation}\label{eq6.78}
\hat{f}(M) \quad = \quad M \qquad \mbox{for all} \ f \in
\mathcal{S}.
\end{equation}
It follows that $\tilde{L}_{\beta}$ and $\ddot{L}_{\beta}$ are
$\mathcal{S}$ invariant.  Translation invariance follows from
(\ref{eq6.67}).

From (\ref{eq6.72}), (\ref{eq6.73}) is clear.  We require $W \subset
\ddot{L}_{\beta}$ so that condition (ii) holds for $T_{\beta +
1}$.  Condition (iii) uses (\ref{eq6.75}), (\ref{eq6.76}) and (\ref{eq6.77}).  Condition
(iv) is checked just as partial additivity was in Theorem \ref{theo6.9}(c).

It remains to prove that $\ddot{L}_{\beta} \not= \emptyset $.

Choose a sequence  $ a_{0} < a_{1} < ... $ in $X$ with limit $a
\in X$.  Because $\beta $ is a countable limit ordinal we can
choose a sequence $ 1 =  \beta_{0} < \beta_{1} < ... $  in
$\Omega$ with limit $\beta$.  We inductively construct $p_{0},
p_{1},...$ in $T$ such that for $i = 0,1,...$
\begin{equation}\label{eq6.79}
o(p_{i}) \ = \ \beta_{i}, \quad p_{i+1}|\beta_{i} \ = \ p_{i},
\quad \hat{p}_{i}(\beta_{i}) \ = \ a_{i}.
\end{equation}
There is, then, a unique branch $x$ of height $\beta$ with
$x|\beta_{i} = p_{i}$.  Since $\hat{x}(\beta) = a$, we have $x \in
\ddot{L}_{\beta}$.

For the induction, begin with $p_{0} \in X^{1}$ with $p_{0}(0) =
a_{0} $. Assume that $p_{i}$ is defined and choose $q \in T$ such
that $o(q) = \beta_{i+1} \setminus \beta_{i}$ which exists because
$T$ has height $\beta$.  By Lemma 5.12 and condition (ii) applied
to $q$ there exists $f \in \mathcal{S}$ such that $f(X) =
(a_{i},\infty)$ and
\begin{equation}\label{eq6.80}
f(\hat{q}(\beta_{i+1} \setminus \beta_{i})) \quad = \quad a_{i+1}.
\hspace{1.5cm}
\end{equation}
Define
\begin{equation}\label{eq6.81}
p_{i+1} \quad = \quad p_{i} + f \circ q   \hspace{2cm}
\end{equation}
which lies in $T$ by conditions (iii) and (iv). By (\ref{eq6.76})
$p_{i+1}$ satisfies (\ref{eq6.79}) as required.

\end{proof}\vspace{.5cm}

\begin{cor}\label{cor6.16} There exist countable special semigroups for the HLOTS of rationals $\Q$.
 For any countable special semigroup $\mathcal{S}$ for $\Q$ there exist $\mathcal{S}$ trees
 which are Aronszajn.  In particular, reproductive Aronszajn trees exist.
\end{cor}

\begin{proof}  Choose $f_{0} : \Q
\rightarrow (0,\infty) \cap \Q$ an order isomorphism.
Let $\mathcal{H}$ be the countable group of translations on
$\Q$.  Let $\mathcal{S}$ be the smallest semigroup
containing $\{f_{0}\} \cup \mathcal{H}$.  $\mathcal{S}$ is
countable and is clearly a special semigroup for $\Q$.

Assume $\mathcal{S}$ is a countable special semigroup for
$\Q$.  We will use Theorem \ref{theo6.15} to construct a
coherent family $\{ T_{\delta} : \delta < \Omega \}$ of countable
$\mathcal{S}$ trees.

At the infinite, tail-like ordinal $\beta$ stage, it suffices to
construct $W \subset \ddot{L}_{\beta}$ which is countably
infinite, $\mathcal{S}$ invariant and translation invariant.

To do so, choose $x \in \ddot{L}_{\beta}$ and let $W_{0} =
\{\tau_{i}(x) : i < \beta \}$ so that $W_{0}$ is a countably
infinite, translation invariant subset of $\ddot{L}_{\beta}$.
Then define $W_{0} \subset W_{1} \subset ...$ by
\begin{equation}\label{eq6.82}
\begin{split}
W_{n+1} \  = \ \{ f \circ x : x \in W_{n} \ \mbox{and} \ f \in \mathcal{S} \} \cup  \\
\{ x : \mbox{there exists} \ f \in \mathcal{S} \ \mbox{such that}
\ f \circ x \in W_{n} \}.
\end{split}
\end{equation}
Each $W_{n+1}$ is translation invariant.  Since $1_{{\Q}}
\in \mathcal{S}$, $W_{n+1} \supset W_{n}$ and since
$\ddot{L}_{\beta}$ is $\mathcal{S}$ invariant, $W_{n+1} \subset
\ddot{L}_{\beta}$.  $W = \cup _{n} W_{n} $ which is $\mathcal{S}$
invariant as well as translation invariant.  Because $\mathcal{S}$
and $W_{0}$ are countable, each $W_{n}$ is. Thus, $W$ is the
required countable subset of $\ddot{L}_{\beta}$.

By Theorem \ref{theo6.15}(a) $T = \cup \{ T_{\delta} : \delta < \Omega \}
$ is an $\mathcal{S}$ tree of height $\Omega$. Since every branch
of an $\mathcal{S}$ tree has countable height it follows that $T$
is Aronszajn.
\end{proof}\vspace{1cm}

\subsection{The Omega Thinning Construction}

 Recall from Proposition \ref{propCNF} the Cantor Normal Form which implies that for an ordinal $\a \ge 1$:

  \begin{equation}\label{eq6.86} \begin{split}
\a = \om^{\g_1} + \dots + \om^{\g_{k + \ell}} \ \text{with} \ k, \ell \in \N, k+\ell \ge 1, \hspace{2.5cm}\\ \text{and} \
   \g_1 \ge \dots \g_k \ge \om > \g_{k+1} \ge \dots \g_{k+\ell} \ge 0, \hspace{2cm}\\
\Longrightarrow \quad \om \cdot  \a = \om^{\g_1} + \dots + \om^{\g_k} + \om^{\g_{k+1} + 1} + \dots + \om^{\g_{k + \ell} + 1}.
\end{split} \end{equation}

In particular, $\a$ is a limit ordinal iff $\a = \om \cdot \b$ for some $\b \ge 1$, see Corollary \ref{corCNF}.

\begin{df}\label{df6.17}  Assume that $T$ is a normal bi-ordered tree  with height $h(T)$  a tail-like ordinal satisfying
   $\om < h(T) \leq \Omega$ and so $\om^2 \le h(T)$.

   Let the \emph{Omega Thinning} be the tree  $\om T$ with
   \begin{equation}\label{eq6.83}
   \om T \ = \ \{ p \in T : o(p) \ \text{ is a    limit ordinal} \}.
   \end{equation}

   We will write $(\om S)_p$ for the successor set of $p$ in the tree $\om T$.

    So $(\om S)_0 =
L_{\om}$ and for $p \in \om T$, $(\om S)_p = T_p \cap L_{o(p) + \om}$. \end{df}

   If $p \in \om T$, then $o(p) = \om \cdot \b$ for some $\b \ge 1$ and it follows that its order in $\om T$, $o_{\om}(p) = \b$.

   So if $h(T) =  \om^{\g}$, then $h(\om T) = h(T)$ if $\g \ge \om$ and $h(\om T) = \om^{\g - 1}$ if $2 \le \gamma < \om$.

   For $p \in T$, $o(p) < h(T)$ and $\om < h(T)$.  Since $h(T)$ is tail-like, $o(p) + \om < h(T)$ and so normality implies there
   exists $q \in T_p$ with $o(q) = o(p) + \om$.
    Thus condition (iv) of Definition 5.1 follows for $\om T$ from the same condition for $T$. Thus, $\om T$ is
   a normal tree.

  \begin{prop} \label{prop6.18} If $S_p$ has no max for each $p \in T$, e.g. if $T$ is of unbounded type,
  then $\om T$ is of dense type, i.e.  $(\om S)_p$ is order dense for each $p \in \om T$.

  In particular, if $T$ is of $\Z$ type or $\N$ type, then $\om T$ is of dense type.\end{prop}

  \begin{proof} We proceed as in Proposition \ref{prop5.3a}. Let $q < r \in (\om S)_p$ for $p \in \om T$.
  So $q,r \in T_p$ and with $\ep = \ep(q,r)$,
  we have $ o(p) \le \ep < o(p) + \om$ and so $\ep + 2 < o(p) + \om$.

  $q_i = r_i$ for all $i \le \ep$ and
  $q_{\ep +1} < r_{\ep + 1}$. Because $S_{q_{\ep + 1}}$ has no max, there exists $a \in S_{ q_{\ep + 1}} $ with $q_{\ep + 2} < a$. By normality
  we can choose $s \in T$ with $o(s) = o(p) + \om$ and $ s_{\ep +2} = a$. It follows that $s \in (\om S)_p$ with $q < s < r$.

   \end{proof} \vspace{.5cm}

   Since $h(T)$ is assumed to be a limit ordinal, and $T$ is normal, it follows that for any branch $x \in X(T)$ $h(x)$ is a limit ordinal
    and so from (\ref{eq6.86}) $h(x) = \om \cdot \b$ with $\b \ge 1$.     A branch  of $\om T$ is of the form
   $ x \cap \om T$ for a unique branch $x$ of $T$ which has $h(x) = \om^2 \cdot \b$ with $\b \ge 1$ and so that $h(x)$ is a limit of limit ordinals.
   The order of $x \cap \om T$ in $\om T$ is then $\om \cdot \b$.

   The map defined by
$x \cap \om T \mapsto x$ is an order  injection from $X(\om T)$ into $X(T)$.
   We identify $X(\om T)$ with its image under this map and so regard $X(\om T) \subset X(T)$.

   \begin{cor}\label{cor6.19} If $S_p$ has no max for each $p \in T$, e.g. if $T$ is of unbounded type,
   then $X(T)$ and $X(\om T)$ are order dense LOTS and $X(\om T)$ is a dense subset of $X(T)$.
   In particular, they have a common connected completion. \end{cor}

   \begin{proof} Both $X(T)$ and $X(\om T)$ are order dense by Proposition \ref{prop5.3a} and
   Proposition \ref{prop6.18}.

   For any $p \in T$ we show there exists $x \in X(T)$ with $p \in x$ and with $h(x)$ a tail-like ordinal $\om^{\g}, \g \ge 2 $
   and so $x \cap \om T \in X(\om T)$. It then
   follows from Proposition  \ref{prop5.3} (g) that $X(\om T)$ is dense in $X(T)$.\vspace{.25cm}

   {\bfseries Case 1} $h(T)$ is a countable ordinal: By Lemma \ref{lem5.4}, there exists $x \in X(T)$ such that $h(x) = h(T)$ and $p \in x$. Since $h(T)$
   is a tail-like ordinal with $h(T) > \om$, $x \in X(\om T)$.\vspace{.25cm}

  {\bfseries Case 2}  $h(T) = \Omega = \om^{\Omega}$ : Choose $\gamma_1 > 2$ so that $o(p) < \om^{\gamma_1}$.
  Choose $p(\gamma_1)$ such that $p \prec p(\gamma_1)$ and
   $o(p(\gamma_1)) = \om^{\gamma_1}$. Assume that for all $\d$ with $\gamma_1 \le \d < \gamma$, $p(\d)$ has been chosen with
   $o(p(\d)) = \om^{\d}$ and if $\d_1 < \d_2$, then $p(\d_1) \prec p(\d_2)$.

   If $\gamma = \d +1$, then choose $p(\gamma)$ such that $p(\d) \prec p(\gamma)$ and $o(p(\gamma)) = \om^{\gamma}$.

   If $\gamma$ is a limit ordinal, then there exists a branch $z \in X(T)$ with $p(\d) \in z$ for all $\d < \gamma$.
   If $h(z) = \om^{\gamma}$, then the process stops and $z = x$ is the required element of $X(T)$.  Otherwise, $h(z) > \om^{\gamma}$ and we let
   $p(\gamma) = z_{\om^{\gamma}}$.

   If the process never stops, then we obtain $p(\gamma)$ for all $\gamma$ with $\gamma_1 \le \gamma < \Omega$ and these define a branch $z$ of $X(T)$
   with height $\Omega$ and  $z = x$ is the required element of $X(T)$.

   \end{proof} \vspace{.5cm}

   \begin{prop}\label{prop6.20} If $T$ is homogeneous or reproductive, then $\om T$ is homogeneous or reproductive, respectively. \end{prop}

   \begin{proof} Any automorphism $f$ of $T$ restricts to an automorphism of $\om T$ because $o(f(p)) = o(p)$. So if $T$ is homogeneous,
   then $\om T$ is. In particular, it then follows that each successor LOTS $(\om S)_p$ is transitive.

   Now assume that $T$ is reproductive and that $p \in T$. The isomorphism $j_p : T \to T_p$ induces an isomorphism of $\om T $ onto
   $\{ p \} \cup [(\om T) \cap T_p)]$.  Notice that if $o(q)$ is a limit ordinal if and only if  $o(j_p(q)) = o(p) + o(q)$ is a limit ordinal.
   If $p \in \om T$, then $(\om T)_p = \{ p \} \cup [(\om T) \cap T_p)]$.  It follows that $\om T$ is reproductive.

   \end{proof} \vspace{.5cm}

  We now assume  for a LOTS $X$ that $T$ is an additive subtree of the simple tree on $X, \Omega$.  We would like to say that $\om T$ is
  then an additive tree and indeed it is, but to make sense of this requires a bit of preliminary work.

  The set $L_{\om}$ of vertices of $T$ at level $\om$ is a subset of $T^{\om +1}$ and $\pi^{\om+1}_{\om} : T^{\om +1} \to T^{\om}$ is an isomorphism.
  Moreover, additivity implies that $T^{\om}$ is the entire simple tree on $X, \om$ by Proposition \ref{prop6.6} (f). So $p \mapsto A_p$ is an
  order injection of $L_{\om}$ into $X^{\om}$.

  Because it is additive, $T$ is reproductive by Proposition \ref{prop6.6} (f) and so $\om T$ is
  reproductive by Proposition \ref{prop6.20}. So it follows that $\om T$ is a tree of
  $(\om S)_0 = L_{\om}$ type. We want to identify $\om T$ with a subtree of the simple tree on $L_{\om}, \Omega$.

Using (\ref{eq4.20}) we can identify $(L_{\om})^{\a} \subset (X^{\om})^{\a} = X^{\om \cdot \a}$.

\begin{lem}\label{lem6.21} For all $\a$ with $1 \le \a < \Omega$, the set of level $\om \cdot \a$ vertices of $T$,
i.e. $L_{\om \cdot \a} \subset X^{\om \cdot \a}$, is
a subset of $ (L_{\om})^{\a}$. \end{lem}

\begin{proof} This is trivial for $\a = 1$.

If $\a$ is a limit ordinal, then $L_{\om \cdot \a} \subset \{ s \in X^{\om \cdot \a} : s|(\om \cdot \b) \in L_{\om \cdot \b} $
for all $\b < \a \}$. On the other hand, $(L_{\om})^{\a} = \{ s \in X^{\om \cdot \a} : s|(\om \cdot \b) \in (L_{\om})^{\b}$ for all $\b < \a \}$.

If $\a = \b + 1$, then $L_{\om \cdot \a} = \bigcup (T_p \cap L_{o(p) + \om})$ with $p$ varying over $L_{\om \cdot \b}$.  By additivity,
$ T_p \cap L_{o(p) + \om} = \{ p + q : q \in L_{\om} \}$. So $L_{\om \cdot \a} $ consists of the successors of $L_{\om \cdot \b}$
in the $L_{\om}, \Omega$ simple tree.

Thus, the result follows by induction.

  \end{proof} \vspace{.5cm}

  With these identifications we have:

  \begin{prop} \label{prop6.22} If $T$ is an additive subtree of the simple tree on $X, \Omega$, then $\om T$ is an additive subtree
  of the simple tree on  $L_{\om}, \Omega$. If, in addition,
   $X$ is transitive, then $\om T$ is homogeneous with $L_{\om}$ transitive and order dense. \end{prop}

  \begin{proof} Additivity is clear and only requires identifying addition in $X^{\om \cdot \a}$ with addition in $(X^{\om})^{\a}$.

  If $X$ is transitive, then $T$ is homogeneous by Proposition \ref{prop5.16} and so $\om T$ is homogeneous by Proposition \ref{prop6.20}.
  This implies that $(\om S)_0 = L_{\om}$ is transitive. Since $X$ is transitive it has no max and so $L_{\om}$ is order dense by
  Proposition \ref{prop6.18}.
 \end{proof} \vspace{1cm}

\section{\textbf{The Double Tower for a CHLOTS}}
\vspace{.5cm}

Given a CHLOTS $F$ we constructed in Section 4 a tower of CHLOTS,
defining for every positive ordinal $\alpha$
\begin{equation}\label{eq7.1}
F_{\alpha} \quad = \quad \{ x \in F^{\alpha} : x(i) \in J \quad
\mbox{for} \ 1 \leq i < \alpha  \},
\end{equation}
where $J = [-1,+1]$ is a distinguished, nontrivial, closed
interval in $F$.  By Theorem \ref{theo4.2} each $F_{\alpha}$ is a CHLOTS
when $\alpha$ is a countable, tail-like ordinal.  In particular,
with $\alpha = 1$, $F_{\alpha} = F$.  By Proposition \ref{prop6.3}
$F_{\alpha}$ is order isomorphic to the completion of $F^{\alpha}$.

In Section 6.3 we defined the order tree $T(F)$ associated with $F$.
  We define the \emph{arborization}\index{arborization}\index{CHLOTS!arborization} of a CHLOTS $F$ to be the
completion of the branch space of its order tree.  We use the
notation\index{$a(F)$}\index{$\hat{a}(F)$}
\begin{equation}\label{eq7.2}
a(F) \ = \ X(T(F)) \qquad \mbox{and} \qquad \widehat{a(F)} \
= \ \widehat{X(T(F))}.
\end{equation}

The order tree is of type $F$ with $X(T(F)^{2}) \cong S_{0} \cong
F $ in a natural way.  Thus, from the projection map $\pi_{2}$ of
(\ref{eq5.21}) we can define the continuous order surjections
\begin{equation}\label{eq7.3}
\pi_{F} : a(F) \rightarrow F \qquad \mbox{and} \qquad
\hat{\pi}_{F} : \widehat{a(F)} \rightarrow F.
\end{equation}

In Definition \ref{df4.3} we called a CHLOTS $F$ at least as big as a
CHLOTS $F_{1}$ if there exists an order injection of $F_{1}$ into
$F$.  By Corollary \ref{cor4.5} this is equivalent to the existence of a,
necessarily continuous, order surjection of $F$ onto $F_{1}$.  We
summarize and extend some earlier results using this comparison
concept.

\begin{prop}\label{prop7.1}  Let $F$ and $F_{1}$ be CHLOTS and $\R$ be the CHLOTS of real numbers.
\begin{enumerate}
\item[(a)]  If $F$ is not order isomorphic to $\R$, then $F$ is bigger than  $\R$.
\item[(b)]  If $F_{1}$ is at least as big as $F$ and $F_{1}$ satisfies the countable chain condition, then so does $F$.
\item[(c)]  If $F_{1}$ is the completion of the branch space of an Aronszajn tree and
$\alpha$ is an infinite, tail-like ordinal, then $F_{1}$ is not as big as $F_{\alpha} $.
\item[(d)]  If $\beta < \alpha$ are countable, tail-like ordinals, then $F_{\alpha}$ is
bigger than $F_{\beta}$ and $\widehat{a(F)}$ is bigger than
$F_{\a}$.
\end{enumerate}
\end{prop}

\begin{proof} (a), (b):  If $f : F_{2}
\rightarrow F_{1}$ is an order surjection and $F_{2}$ is separable
or satisfies c.c.c., then $F_{1}$ satisfies the corresponding
condition by Proposition \ref{prop2.5}(g).  If $F$ is any CHLOTS, then by
Proposition \ref{prop2.8} there exists an order surjection of $F$ onto
$\R$ and so $F$ is at least as big as $\R$.  If
$\R$ is as big as $F$, then $F$ is separable and so is
order isomorphic to $\R$ by Proposition \ref{prop2.8}.

(c): This is a restatement of Corollary \ref{cor5.14}.

(d):  Since $\pi^{\a}_{\b} : F_{\alpha} \rightarrow
F_{\beta}$ is an order surjection, $F_{\alpha}$  is at least as
big as $F_{\beta}$.  By Theorem \ref{theo4.11} $F_{\beta}$ is not as big as
$F_{\alpha}$.

Now choose an order embedding $\tilde{p} : \alpha + 1 \rightarrow
\bullet F \bullet $ with $\tilde{p}(0) = m$ and $\tilde{p}(\alpha)
= M$ as in (\ref{eq6.58}).  The map $z \mapsto p(z)$ of (\ref{eq6.56}) associates
to each $z \in F^{\alpha}$ a branch of height $\alpha$ for the
tree $T(F)$.  We thus have an order injection from $F^{\alpha}
\supset F_{\alpha}$ into $\tilde{L}_{\alpha} \subset a(F) \subset
\widehat{a(F)}$.  Thus, $\widehat{a(F)}$ is at least as big as
$F_{\alpha}$.  Furthermore, if $\gamma > \alpha$ is a countable,
tail-like ordinal, then $\widehat{a(F)}$ is at least as big as
$F_{\gamma}$ which is bigger than $F_{\alpha}$.  Hence,
$\widehat{a(F)}$ is bigger than $F_{\alpha}$.

\end{proof}\vspace{.5cm}

\begin{theo}\label{theo7.2}  Let $F$ be a CHLOTS. For every countable ordinal $\alpha$ there exists a
HLOTS $a_{\alpha}(F)$ with completion $\widehat{a_{\alpha}(F)}$ and for each pair of countable
ordinals $\beta < \alpha$ there exists an order surjection
$p^{\a}_{\b} : a_{\alpha}(F) \rightarrow \widehat{a_{\b}(F)} $
 with completion $\hat{p}^{\a}_{\b} : \widehat{a_{\alpha}(F)} \rightarrow \widehat{a_{\b}(F)} $
 so that the following conditions hold.
\begin{enumerate}
\item[(i)] For $\gamma < \beta < \alpha < \Omega$ we have
\begin{equation}\label{eq7.4}
p^{\a}_{\g} \quad = \quad  p^{\b}_{\g} \circ
p^{\a}_{\b }. \hspace{1cm}
\end{equation}
\item[(ii)]  For $\beta < \alpha < \Omega$ define $a^{\b}_{\a }(F) \subset a_{\a }(F)$ by
\begin{equation}\label{eq7.5}
a^{\b}_{\a }(F) \ = \ \begin{cases} a_{\alpha}(F) \qquad
\mbox{for} \ \beta + 1 = \alpha \\ \cap
\{(p^{\a}_{i})^{-1}(a_{i}(F)) : \beta < i < \alpha \} \quad
\mbox{for} \ \beta + 1 < \alpha.  \end{cases}
\end{equation}
The restriction of $p^{\a}_{\b}$ to $a^{\b}_{\a}(F)$
is surjective.  That is,
\begin{equation}\label{eq7.6}
p^{\a}_{\b}(a^{\b}_{\a}(F))
 \quad = \quad
\widehat{a_{\b}(F)}. \hspace{1cm}
\end{equation}
\item[(iii)] $a_{0}(F) = \widehat{a_{0}(F)} = F$.
\item[(iv)]  If $\alpha = \gamma +1$, then
\begin{equation}\label{eq7.7}
\begin{split}
a_{\alpha}(F) \quad = \quad a(\widehat{a_{\g}(F)}). \hspace{2cm} \\
\widehat{a_{\alpha}(F)} \quad = \quad \widehat{a}(\widehat{a_{\g}(F)}). \hspace{2cm} \\
p^{\a}_{\g} \quad = \quad
\pi_{\widehat{a_{\gamma}(F)}}.
\hspace{3cm}
\end{split}
\end{equation}
\item[(v)]  If $\alpha$ is a countable limit ordinal, then $(\{ \widehat{a_{i}(F)} : i < \alpha \},
 \{ \hat{p}^{i}_{j} : j < i < \alpha \}) $ is an unbounded, special inverse system and
\begin{equation}\label{eq7.8}
\begin{split}
\widehat{a_{\alpha}(F)} \quad = \quad \overleftarrow{Lim}(\{ \widehat{a_{i}(F)} : i < \alpha \},
\{ \hat{p}^{i}_{j} : j < i < \alpha \}).  \\
\hat{p}^{\a}_{\b} \quad = \quad \mbox{coordinate projection to} \ \widehat{a_{\b}(F)}. \hspace{.5cm} \\
a_{\alpha}(F) \quad = \quad \bigcup _{i < \alpha} a^{i}_{\alpha}(F). \hspace{2cm} \\
p^{\a}_{\b} \quad = \quad
\hat{p}^{\a}_{\b}|a_{\alpha}(F). \hspace{3cm}
\end{split}
\end{equation}
\end{enumerate}
\end{theo}

\begin{proof} Conditions (iii),(iv) and (v)
define an inductive construction and we show, inductively that at
each stage $a_{\alpha}$ the properties described above hold, i.e.
each $a_{\alpha}$ is a HLOTS and each $ p^{\a}_{\b}$ is an
order surjection which satisfies (\ref{eq7.4}) and (\ref{eq7.6}).

The construction begins with condition (iii).\vspace{.25cm}

\textbf{Case 1:} $\alpha = \gamma + 1$. $a_{\alpha}$ is an IHLOTS by
Theorem \ref{theo6.10}, and $ p^{\a}_{\g}$ is the surjective order
map to the level $1$ vertex set of the order tree.  By Proposition
\ref{prop2.1}(a) such a surjective order map is continuous and so its
completion is defined.  For $\beta < \gamma$ define
\begin{equation}\label{eq7.9}
p^{\a}_{\b} \quad = \quad  p^{\g}_{\b} \circ
p^{\a}_{\g }. \hspace{1cm}
\end{equation}
As the composition of order surjections, each $p^{\a}_{\b}$
is an order surjection and (\ref{eq7.4}) for $\alpha$ follows from the
corresponding condition for $\gamma$.

If $\beta = \gamma$, then $a^{\b}_{\a} = a_{\alpha}$ and (\ref{eq7.6})
is clear from what we have already shown.  If $\beta < \gamma$
then by (\ref{eq7.9}) and (\ref{eq7.5})
\begin{equation}\label{eq7.10}
a^{\b}_{\a}(F) \quad = \quad
(p^{\a}_{\g})^{-1}(a^{\b}_{\g}(F)). \hspace{2cm}
\end{equation}
Hence, (\ref{eq7.6}) for $\alpha$ follows from (\ref{eq7.6}) for $\gamma$ together
with (\ref{eq7.9}).\vspace{.25cm}

\textbf{Case 2:}  $\alpha$ is a limit ordinal.  Our inductive
hypothesis implies that  $(\{ \widehat{a_{i}} : i < \alpha \}, \{
\hat{p}^{i}_{j} : j < i < \alpha \}) $ is an unbounded, special inverse limit
system.  By Proposition \ref{prop2.4} the inverse limit $\widehat{a_{\a}}$
is an unbounded, connected LOTS and each
$\hat{p}^{\a}_{\b}$ is a continuous order surjection.  By
Proposition \ref{prop2.7}(c) $\alpha < \Omega$ implies that $\widehat{a_{\a}}$
is of countable type and hence, by Proposition \ref{prop2.5}(f), any subset
has countable type.

For any $\beta < \alpha$ we can choose an increasing sequence of
ordinals  $\beta = \beta_{0} < \beta_{1} < ...$ with supremum
$\alpha$.

If $z_{0} \in \widehat{a_{\b}} $, then because
$p^{\b_{n}}_{\b_{n-1}}$ satisfies (\ref{eq7.6}) we can choose a
sequence $z_{1}, z_{2},...$ such that
\begin{equation}\label{eq7.11}
z_{n} \in a^{\b_{n-1}}_{\b_{n}} \qquad \mbox{with} \qquad
p^{\b_{n}}_{\b_{n-1}}(z_{n}) \ = \ z_{n-1} \quad \mbox{for}
\quad i = 1,2,...
\end{equation}
Because the sequence $\{ \beta_{n} \}$  is cofinal in $\alpha$ the
sequence $\{ z_{n} \}$ defines a unique element $z$ of the inverse
limit. By (\ref{eq7.11})
\begin{equation}\label{eq7.12}
p^{\a}_{i}(z) \in a_{i} \qquad \mbox{for} \quad \beta_{n-1} <
i \leq \beta_{n}, \ n = 1,2,... \hspace{.5cm}
\end{equation}
So by definition (\ref{eq7.5}) $z \in a^{\b}_{\a}$ with
$\hat{p}^{\a}_{\b}(z) = z_{0}$.  Strictly speaking the maps
$\hat{p}^{\a}_{i}$ are used in (\ref{eq7.5}) to define
$a^{\a}_{\b}$ because $a_{\alpha}$ and $p^{\a}_{\b}$
are defined subsequently in (\ref{eq7.8}).  It follows that $a_{\alpha}$
defined in (\ref{eq7.8}) is dense in $\widehat{a_{\a}}$ and so has
completion $\widehat{a_{\a}}$.  Furthermore, (\ref{eq7.6}) holds at the
$\alpha$ level.  By definition of the inverse limit projections
(\ref{eq7.4}) holds for $\alpha$.\vspace{.25cm}

It remains to prove that $a_{\alpha}$ is doubly transitive which
will imply it is a HLOTS since it is of countable type.

Assume $x < y$ and $z < w$ in $a_{\alpha}$.  Because $\alpha$ is a
limit ordinal there exists $\beta < \alpha$ such that $x,y,z,w \in
a^{\beta}_{\alpha}$ with
\begin{equation}\label{eq7.13}
p^{\a}_{\b}(x) < p^{\a}_{\b}(y) \qquad \mbox{and}
\qquad p^{\a}_{\b}(z) < p^{\a}_{\b}(w).
\end{equation}
By definition of $a^{\b}_{\a}$ we have
\begin{equation}\label{eq7.14}
p^{\a}_{i}(x),p^{\a}_{i}(y),p^{\a}_{i}(z),p^{\a}_{i}(w)
\in a^{i}_{\a} \qquad \mbox{for} \quad \beta < i < \alpha.
\hspace{1cm}
\end{equation}
Since $\widehat{a_{\b}}$ is a CHLOTS there exists $\hat{f}_{\beta}
\in H_{+}(\widehat{a_{\b}})$ such that
\begin{equation}\label{eq7.15}
\hat{f}_{\beta}(p^{\a}_{\b}(x)) \  = \
p^{\a}_{\b}(z) \qquad \mbox{and} \qquad
\hat{f}_{\beta}(p^{\a}_{\b}(y)) \ = \ p^{\a}_{\b}(w).
\end{equation}

By induction we will construct for $\beta < i \leq \alpha$, $f_{i}
\in H_{+}(a_{i})$ so that
\begin{equation}\label{eq7.16}
p^{i}_{\b} \circ f_{i} =  \hat{f}_{\beta} \circ p^{i}_{\b}
\quad \mbox{and} \quad p^{i}_{j} \circ f_{i} =  f_{j} \circ
p^{i}_{j} \quad \mbox{for} \ \beta < j < i.
\end{equation}
and
\begin{equation}\label{eq7.17}
f_{i}(p^{\a}_{i}(x)) \ = \ p^{\a}_{i}(z) \qquad \mbox{and}
\qquad f_{i}(p^{\a}_{i}(y)) \ = \ p^{\a}_{i}(w).
\end{equation}

If $i$ is a limit ordinal, e.g. $i = \alpha$, then we define
$\hat{f}_{i} \in H_{+}(\widehat{a_{i}}) $ to be the inverse limit of
 $\{ \hat{f}_{j} : \beta < j < i |]$.  That is,
$\hat{f}_{i}(z)$ is the unique element of $\widehat{a_{i}}$ which
projects via $\hat{p}^{i}_{j}$ to
$\hat{f}_{j}(\hat{p}^{\a}_{j}(z))$.  Since each $\hat{f}_{j}$
is the completion of an isomorphism $f_{j}$, it follows that
$\hat{f}_{i}$ maps $a_{i}$ to $a_{i}$. The restriction of
$\hat{f}_{i}$ to $a_{i}$ defines $f_{i}$ so that (\ref{eq7.16}) holds and
the original $\hat{f}_{i}$ is the completion of $f_{i}$.  Because
$\widehat{a_{i}}$ is the inverse limit of its predecessors,
$p^{\a}_{i}(r)$ is determined by the maps $p^{\a}_{j}(r)$ for
$\beta < j < i$.  Hence, (\ref{eq7.17}) for $i$ follows from the
inductively assumed equations for $j < i$.

Finally, assume that $i = k + 1$ and that $f_{j}$  is defined for
all $\beta < j \leq k$ so that (\ref{eq7.16}) and (\ref{eq7.17}) hold.

Now we use the tree structure: $a_{i}$ is the branch space of the
order tree on $\hat{a}_k$.  We can regard $\hat{f}_{k}$ as an
order isomorphism on the level $1$ vertices of the tree. By (\ref{eq7.14})
$p^{\a}_{i}(x),...,p^{\a}_{i}(w) $ are branches in $a_{i}$
with vertices at the $1$ level
$p^{\a}_{k}(x),...,p^{\a}_{k}(w) $.

By Theorem \ref{theo6.10} the order tree $T(\widehat{a_{k}}) $ is reproductive and so
is homogeneous by Proposition \ref{prop5.16}.  We use a variation of the
proof of Lemma \ref{lem5.17}(a).

For each $r \in S_{0} = \widehat{a_{k}}$ we choose a tree isomorphism
$g_{r} :T_{r} \rightarrow T_{\hat{f}_{k}(r)}$.  With $r_{1} =
p^{\a}_{k}(x) $, $ \hat{f}_{k}(r_{1}) = r_{2} =
p^{\a}_{k}(z)$.  So the branches
$(g_{r_{1}})_{*}(p^{\a}_{i}(x) \cap T_{r_{1}})$ and
$p^{\a}_{i}(z)\cap T_{r_{2}}$ both lie in the branch space of
the homogeneous tree $T_{r_{2}}$.  We can adjust $g_{r_{1}}$ so
that these two are in fact equal. Similarly, for
$p^{\a}_{i}(y)$ and $p^{\a}_{i}(w)$.  Assemble the maps
$g_{r}$ to get a tree isomorphism $g$.  Then $g_{*} \in
H_{+}(a_{i})$ so that
\begin{equation}\label{eq7.18}
p^{i}_{k} \circ g_{*} \quad = \quad \hat{f}_{k} \circ p^{i}_{k}.
\hspace{2cm}
\end{equation}

We let $f_{i} = g_{*}$.
From the construction (\ref{eq7.17}) is clear and (\ref{eq7.16}) for $i$ follows
from (\ref{eq7.16}) for $k$ together with condition (i).

This completes the inductive construction.  The case $i = \alpha$
yields $f_{\alpha}$ which is the required order isomorphism of
$a_{\alpha}$ which maps the pair $x,y$ to the pair $z,w$.

\end{proof}\vspace{.5cm}

\begin{theo}\label{theo7.3}  For a CHLOTS $F$ and $(i,j) \in \Omega \times \Omega$ define the CHLOTS
\begin{equation}\label{eq7.19}
F_{(i,j)} \quad =  \quad (\widehat{a_{i}(F)})_{\omega^{j}}.
\hspace{2cm}
\end{equation}
$F_{(0,0)} \cong F$ and if $(i,j) < (\tilde{i},\tilde{j}) $ in
$\Omega \times \Omega$, then $F_{(\tilde{i},\tilde{j})}$ is bigger
than $F_{(i,j)}$. In particular, $F_{(\tilde{i},\tilde{j})}$ is not homeomorphic to
 $F_{(i,j)}$.
\end{theo}

\begin{proof}  Recall that as $j$ varies through
the countable ordinals, $\omega^{j}$ varies through the countable,
tail-like ordinals and $j < \tilde{j}$ iff $\omega^{j} < \omega
^{\tilde{j}} $.  Also, $\omega^{0} = 1$.  For any countable,
tail-like ordinal $\alpha$ $F_{\alpha}$ defined by (6.1) is a
CHLOTS and $F_{1} \cong F$.  So by Theorem \ref{theo7.2}, $F_{(i,j)}$ is a
CHLOTS for all $(i,j) \in \Omega \times \Omega$.

If $i = \tilde{i}$ and $\tilde{F} = \widehat{a_{i}(F)}$, then $j <
\tilde{j}$ implies $\tilde{F}_{\beta}$ is bigger than
$\tilde{F}_{\alpha}$ with $\beta = \omega^{\tilde{j}}$ and $\alpha
= \omega^{j}$ by Theorem \ref{theo4.11}.  Thus,
$F_{(\tilde{i},\tilde{j})}$ is bigger than $F_{(i,j)}$.

Now suppose that $i  < \tilde{i}$ so that  $ i + 1 \leq
\tilde{i}$.

First, $F_{(\tilde{i},\tilde{j})}$ is at least as big as
$F_{(\tilde{i},0)}$ and so by using the projections of Theorem \ref{theo7.2}
we see that $F_{(\tilde{i},0)}$ is as big as $F_{(i+1,0)} =
\widehat{a(F_{(i,0)})}$.  By Proposition \ref{prop7.1}(d), $F_{(i+1,0)}$ is bigger
than $F_{(i,j)}$.  Thus, $F_{(\tilde{i},\tilde{j})}$ is bigger
than $F_{(i,j)}$.

That $F_{(i,j)}$ and $F_{(\tilde{i},\tilde{j})}$ are not homeomorphic follows just as in Theorem \ref{theo4.11} because
$F_{(i,j)}'$ injects into $F_{(\tilde{i},\tilde{j})}$.

\end{proof}\vspace{.5cm}

For the associated Cantor spaces which are the compactifications of the AS doubles we have
the following extension of Theorem \ref{theo4.13}.

\begin{cor}\label{cor7.3a}  For a CHLOTS $F$ and $(i,j), (\tilde{i},\tilde{j}) \in \Omega \times \Omega$
if $(i,j) < (\tilde{i},\tilde{j}) $ in
$\Omega \times \Omega$, then $C(F_{(\tilde{i},\tilde{j})})$ is bigger
than $C(F_{(i,j)})$. In particular, $C(F_{(\tilde{i},\tilde{j})})$ is not isomorphic to
 $C(F_{(i,j)})$.\end{cor}

\begin{proof} If $i = \tilde i$ this follows directly from in Theorem \ref{theo4.13}.

If $i < \tilde i$, then by Theorem \ref{theo7.3} $F_{(\tilde{i},\tilde{j})}$ is bigger than $F_{(i,j + 1)}$ and so
$C(F_{(\tilde{i},\tilde{j})})$ is at least as big as $C(F_{(i,j + 1)})$. By Theorem \ref{theo4.13} again
$C(F_{(i,j + 1)})$ is bigger than $C(F_{(i,j)})$.

\end{proof} \vspace{.5cm}

We call a LOTS $X$ \emph{$\R$-bounded}\index{LOTS!$\R$-bounded}\index{$\R$-bounded} if it admits an order injection into $\R_{\d}$ for some countable
 ordinal $\d$.

\begin{cor}\label{cor7.4} For a CHLOTS $F$ and $(i,j) \in \Omega \times \Omega$  with $0 < i$,
the CHLOTS $F_{(i,j)}$ is not $\R$-bounded. \end{cor}

\begin{proof} $F_{(i,j)}$ is at least as large as $\widehat{a(F)} = F_{(1,0)}$ which is larger than
$F_{\a}$ for every countable $\a$ and $F_{\a}$ is at least as large as $\R_{\a}$.

\end{proof}\vspace{1cm}

\section{\textbf{The Tree Characterization of a CHLOTS}}

\subsection{A Tree  for a LOTS and the IHLOTS Tower}

We begin with a version of the partition tree for a LOTS, described in \cite{T} and in \cite{BLR}.

Throughout this section, all intervals in a LOTS $X$ will be assumed nonempty.
A singleton is an \emph{improper} closed interval\index{improper interval}\index{interval!improper} and so an interval with
at least two points is a \emph{proper} interval\index{proper interval}\index{interval!proper}. If $I$ is a proper interval with
endpoints $a < b$ then we let $I^{\circ} = (a,b)$.

For intervals $I_1, I_2$ in a LOTS
we will write
\begin{equation}\label{eqn.01}
I_1 < I_2 \quad \Longleftrightarrow \quad I_1 \not= I_2 \ \text{and} \ c_1 \le c_2 \ \text{for all} \ c_1 \in I_1, c_2 \in I_2,
\end{equation}
and so $I_1 \cap I_2 $ is either empty or consists of a single common endpoint.

Recall that we use $\# X$ for the cardinality of a set $X$.

\begin{theo}\label{theon.01} If $X$ is a connected, first countable, bounded LOTS $X $,  there is an $\Omega$ bounded subtree $T$ of
  the simple tree on $\Z, \Omega $ whose branch space is order isomorphic to a
 dense subset of $X$. In particular, $\# X \le 2^{\aleph_0}$. \end{theo}

 \begin{proof} We let $m, M$ be the minimum and maximum of $X$, respectively, so that $X = [m,M]$. If $X$ is a singleton, then the subtree consisting
 of the root of the simple tree has branch space isomorphic to $X$ and so we may assume that $m < M$.

 For $ \a \le \Omega, s \in \Z^{\a}$ we will associate a closed interval $I_s = [a^s,b^s]$ in $X$. We will call $s$ \emph{proper}
 \index{sequence!proper} when the
 interval $I_s$ is proper and so $a^s < b^s$.

 The construction will satisfy the following properties.

 \begin{enumerate}
 \item[(i)]  $I_{\emptyset} = [m,M]$.

 \item[(ii)] If $\b < \a \le \Omega$ and $s \in \Z^{\a}$, then $I_s \subset I_{s|\b}$ and
 if $s|\b$ is proper, then $I_s \subset I_{s|\b}^{\circ}$.

 \item[(iii)] If $\a \le \Omega$ is a limit ordinal and $s \in \Z^{\a}$, then $I_{s} = \bigcap_{\b < \a} I_{s|\b}$.

 \item[(iv)] If $s_1 < s_2$ in  $\Z^{\a}$, then $I_{s_1} < I_{s_2}$.

 \end{enumerate}\vspace{.5cm}

 If $s \in X^{\a}$ is proper with $I_s = [a^s,b^s]$, then in $I_s$ we choose a $\pm$cofinal embedding of $z :\Z \to (a^s,b^s)$.
 If $s' \in \Z ^{\a +1}$ with $s'|\a  = s$ and $s'(\a ) = n$, then
we let  $I_{s'} = [z_n,z_{n+1}]$. The set of successors $S_s \cong \Z$ and if $s'' \in S_s$ has $s''|\a  = s$ and $s''(\a ) = n+1$, then
 $b^{s'} = a^{s''}$. It follows that if $s$ is proper, then all of its successors are proper. Observe that $I_{s'} \subset I_s^{\circ}$.

 If $s \in X^{\a}$ is improper, then we let $I_{s'} = I_s$ for all $s' \in S_s$.

 The construction is completed by using (iii) for limit ordinals. \vspace{.25cm}

 Conditions (i) - (iv) are easy to check from the inductive construction.

 The subtree $T = \{ s : s $ is proper $\}$, i.e. $s \in T$ iff $I_s$ is a proper interval.
  It is clear from (i) and (ii) that $T$ is a nonempty subtree of the simple tree.

 We can identify the branch space $X(T)$ with
\begin{equation}\label{eqn.02}
\{ x \in \Z ^{\a } :   x|\b  \in T \ \text{for all} \ \b  < \a , \ \text{and} \ x \not\in T \}. \hspace{1cm}
\end{equation}
That is,
\begin{equation}\label{eqn.03}
x \in X(T) \  \Leftrightarrow \ I_{x|\b}   \ \text{is proper for all} \ \b  < \a \ \text{and} \ I_x \ \text{is not proper}.
\end{equation}
Since every successor of a proper element is proper, it follows that $h(x)$ is a limit ordinal for all $x \in X(T)$.

For $x \in X(T)$ we let $f(x) = a^x \in X$
so that $I_x$ is the singleton $ \{ a^x  \}$. Furthermore, $\b  \mapsto a^{x|\b }$ is an
embedding of the ordinal $h(x) $ into $X$ with limit (= sup) $a^x $ and  $\b  \mapsto b^{x|\b }$ is an
embedding of $h(x)^*$ into $X$ with limit (= inf) $a^x $.  Because $X$ is first countable, $\Omega$ does not inject into $X$.
It follows that $h(x)$ is a countable ordinal.  Hence, $T$ is $\Omega$ bounded.

Observe that if $s \in T^{\a }$ and $\b  < \a $, then $I_{s } \subset I_{s|\b }^{\circ}$
 and so $a^{s|\b } < a^s < b^s < b^{s|\b }$.
Also if $s_1 < s_2 \in T^{\a}$, then by (iv) $I_{s_1}^{\circ} < I_{s_2}^{\circ}$ in $X$.

It follows that $x_1 < x_2$ in $X(T)$ implies $a^{x_1} < a^{x_2}$ and so $f $ is an order injection from $X(T)$ into $X$.

For $\a $ a countable ordinal, we let $C_{\a } = \{ a^s : s \in T^{\a } \} \cup \{ b^s : s \in T^{\a } \}$ and
$C = \bigcup_{\a  < \Omega} \ C_{\a }$ so that $C$ is the set of endpoints of the proper intervals.
Because we are taking the union over the countable ordinals, we have
$\# C \le 2^{\aleph_0}$.

By induction on $\a $ we see that for each $\a  \le \Omega$ the set $X$ is partitioned into three subsets:
\begin{equation}\label{eqn.04}
 \bigcup_{\b  \leq \a } \ C_{\b }, \quad
 \{ a^x  : x \in X(T) \ \text{with height} \le \a  \}, \quad
 \bigcup_{s \in T^{\a }} \ I_s^{\circ}.
 \end{equation}

 So with $\a  = \Omega$ we obtain $X$ is the disjoint union of $C$ and $f(X(T))$.
 Since every $a^x $ is a limit of elements of $C$, it follows
 that $C$ is dense in $X$. Because every element of $X$ is thus a limit of a sequence in $C$ it follows that $\# X \le 2^{\aleph_0}$.

 If $s \in T^{\a }$ then $s$ is contained in some branch $x$ and so $a^s < a^x  < b^s$. If $s < s_1 \in T$ with
 $a^s < a^{s_1}$, then $a^{s} < b^{s} \le a^{s_1}$
  and so $a^s < a^x  < a^{s_1}$. It follows that $f(X(T))$ is dense in $X$ and so is order dense because $X$ is connected. It follows that
  $f : X(T) \to X$ is an order embedding onto a dense subset.

\end{proof} \vspace{.5cm}

The cardinality result is well-known.  See, e.g. \cite{B} Section 4.

  We will require a height estimate for a special case.

  \begin{theo}\label{theonew01} If $\a$ is a positive ordinal, then $\R_{\a} \cong \widehat{X(T)}$ with $T$ a tree of $\Z$ type with height
  $h(T) \le \om \cdot \a$. \end{theo}

  \begin{proof}  With $J = [-1,+1] \subset \R$ we can replace $\R_{\a}$ by the isomorph $\{ x \in J^{\a} : -1 < x_0 < +1 \}$, ie.
  the interval $(-1+,+1-) \subset J^{\a}$. We  let $X $ be the closed interval $ [-1+,+1-] \subset J^{\a}$, which is isomorphic to the two point
  compactification of $\R_{\a}$.

  For $a < b$ in $X$ let $\ep =  min \{ i : a_i \not= b_i \} $ so that $a_i = b_i$ for all $i < \ep$ and $a_{\ep} < b_{\ep}$. We define the
 \emph{ midpoint} of the interval $[a,b]$ to be the point $c$ with
 \begin{equation}\label{eqnew01}
   c_i \ = \ \begin{cases} \ a_i = b_i \ \ \text{for} \ i < \ep, \\
    \frac{1}{2}(a_{\ep} + b_{\ep})  \ \ \text{for} \ i = \ep, \\
    \ \ 0 \ \ \ \text{for} \  \ep < i < \a.\end{cases}
  \end{equation}
  Notice that on $J \subset \R$ algebraic notions and length are defined.

  Now we apply the above construction to $X$ with the proviso that if $I_s = [a^s, b^s]$ is a proper interval, then the $\pm$cofinal sequence
  $z : \Z \to (a^s,b^s)$ maps $0 \in \Z$ to the midpoint of $[a^s,b^s]$.  This then implies for all $n \in \Z$
\begin{equation}\label{eqnew02}
 (z_{n+1})_\ep - (z_{n})_\ep \ \le \ \frac{1}{2}(b^s_{\ep} - a^s_{\ep}).
 \end{equation}

 We now show, by induction, that for all $\b$ with $ 0 < \b \leq \a$
 \begin{equation}\label{eqnew03}
 s \in \Z^{\om \cdot \b} \quad \Longrightarrow \quad a^s_i = b^s_i \ \ \text{for all} \ \ i < \b.
 \end{equation}\vspace{.25cm}

 \textbf{Case 1} $\b = \gamma + 1$:

 Let $s \in \Z^{\om \cdot \b}$. The induction hypothesis applied to $s0= s|\om \cdot \g$  implies that
 $a^{s0}_i = b^{s0}_i \ \ \text{for all} \ \ i < \g$. Since, $[a^s, b^s] \subset [a^{s0},b^{s0}]$,
 $a^{s0}_{\g} = b^{s0}_{\g}$ implies $a^{s}_{\g} = b^{s}_{\g}$ and so $a^s_i = b^s_i$ for all $i < \b$.

 Assume, instead, that $a^{s0}_{\g} < b^{s0}_{\g}$. Let $sn = s|(\om \cdot \g + n)$ for $n < \om$. It follows
 from (\ref{eqnew02}) $b^{sn}_{\g} - a^{sn}_{\g} \le \frac{1}{2^n}(b^{s0}_{\g} - a^{s0}_{\g})$.

 Because $\om \cdot \b = \om \cdot \g + \om, \ I_s = \bigcap_n I_{sn}$ and it follows that $b^s_{\g} = a^s_{\g}$. Hence,
 again $a^s_i = b^s_i$ for all $i < \b$. \vspace{.25cm}

  \textbf{Case 2} $\b$ is a limit ordinal:

 If $\g < \b$, then with $s' = s|\om \cdot \g$ we have $I_s \subset I_{s'}$. So the induction hypothesis applied to $s'$
 implies that $b^s_i = a^s_i $ for all $i < \g$. Since $\b$ is a limit ordinal and $\g < \b$ is arbitrary, it follows that
 $a^s_i = b^s_i$ for all $i < \b$. \vspace{.25cm}

 From (\ref{eqnew03}) applied with $\b = \a$ we see that for all $s \in \Z^{\om \cdot \a}$, $a^s_i = b^s_i$ for all $i < \a$.
 Since $a^s, b^s \in J^{\a}$ this means $a^s = b^s$. Hence, the interval $I_s$ is improper.

Thus, the tree $T$ consisting of those $s$ with $I_s$ proper is a subtree of the simple tree on $\Z, \om \cdot \a$ and so has
 height at most $\om \cdot \a$.

 \end{proof} \vspace{.5cm}

 \begin{cor}\label{cornew02} If $\a$ is a positive ordinal and $X$ is an unbounded LOTS, then there exists an order injection
 from $R_{\a}$ into the completion $\widehat{X^{\om \cdot \a}}$.  \end{cor}

 \begin{proof}  We can identify $X^{\om \cdot \a}$ with the branch space on the simple
 tree on $X, \om\cdot \a$. Since $X$ is unbounded and $\om \cdot \a$ is a limit ordinal, $X^{\om \cdot \a}$ is order dense
 by Proposition \ref{prop5.3a}. In particular, $\Z^{\om \cdot \a}$ is order dense.

 By Theorem \ref{theonew01} $R_{\a}$ is the completion of the branch space of a tree $T$ which is a subtree of the simple tree on
 on $\Z, \om \cdot \a$. From the inclusion of $T$ into the simple tree we obtain an order injection from the branch space $X(T)$ into
 $\Z^{\om \cdot \a}$ from Proposition \ref{prop6.4a}. By Proposition \ref{prop2.4a} the extension
 to the completions is injective. That is, we obtain an order injection
 from $\R_{\a}$ to $\widehat{\Z^{\om \cdot \a}}$.

 Since $X$ is unbounded, there is an order injection from $\Z$ to $X$. From it we obtain an order injection from $\Z^{\om \cdot \a}$ to
$X^{\om \cdot \a}$ which extends to an injection between the completions.

Composing, we obtain the required order injection.

\end{proof} \vspace{.5cm}

 \begin{cor}\label{cornew03} Assume that $X, X_1$ are HLOTS and that $X_1$ admits an order injection into $\R_{\d}$ for some countable
 ordinal $\d$.
 If $\b$ is a positive ordinal and  $\a$ is a countable, tail-like ordinal such that  $\a > \om \cdot \d \cdot \b $,
 then there does not exist an order injection of the completions
 from $\widehat{X_{\a}}$ into  $\widehat{(X_1)_{\b}}$ nor
 an order injection of the Cantor Spaces
 from $C(\widehat{X_{\a}})$ into  $C(\widehat{(X_1)_{\b}})$. \end{cor}

\begin{proof} The injection of $X_1$ into $\R_{\d}$ induces an  injection from $\widehat{(X_1)_{\b}}$ to $(\R_{\d})_{\b}$ which is
isomorphic to $\R_{\d \cdot \b}$ by (\ref{eq4.21}). So there is an injection from $C(\widehat{(X_1)_{\b}})$ to $C(\R_{\d \cdot \b})$
which in turn injects into $\R_{(\d \cdot \b)+1}$.

Since $\a$ is tail-like and $\b, \d > 0$, $\a > \om \cdot \d \cdot \b $ implies $\a > \om \cdot \d \cdot \b + \om = \om \cdot (\d \cdot \b + 1)$.
Hence, there is an order injection from $\widehat{X^{\om \cdot (\d \cdot \b + 1)}}$ into $\widehat{X^{\a}}$.

By Proposition \ref{prop6.3},  $\widehat{X^{\a}} \cong \widehat{X_{\a}}$.

Hence, from an order injection from $\widehat{X_{\a}}$ to $\widehat{(X_1)_{\b}}$ we would obtain an order injection from
$\widehat{X^{\om \cdot (\d \cdot \b + 1)}}$ into $\R_{\d \cdot \b}$.

By Corollary \ref{cornew02} there exists an order injection from $\R_{\d \cdot \b + 1}$ into $\widehat{X^{\om \cdot (\d \cdot \b + 1)}}$.

The composition would contradict Theorem \ref{theo4.11}.

The order injection from $\R_{\d \cdot \b + 1}$ into $\widehat{X^{\om \cdot (\d \cdot \b + 1)}}$ induces an order injection from
$C(\R_{\d \cdot \b + 1})$ into $C(\widehat{X^{\om \cdot (\d \cdot \b + 1)}})$.

As above, from an order injection from $C(\widehat{X_{\a}})$ to $C(\widehat{(X_1)_{\b}})$ we would obtain an order injection from
$C(\widehat{X^{\om \cdot (\d \cdot \b + 1)}})$ into $\R_{(\d \cdot \b)+1}$.

Since $\R_{\d \cdot \b + 1}' \subset C(\R_{\d \cdot \b + 1})$, the composition would yield an injection from
$\R_{\d \cdot \b + 1}'$ to $\R_{\d \cdot \b + 1}$. This contradicts Corollary \ref{cor4.10}  which says that $\R_{\d \cdot \b + 1}$ is order simple.

\end{proof} \vspace{.5cm}

Recall that a LOTS $X$ is $\R$-bounded if it admits an order injection into $\R_{\d}$ for some countable
 ordinal $\d$. Of course, it then injects into any $\R_{\g}$ with $\d \le \g < \Omega$.

 \begin{prop}\label{propnew04a} Assume $X$ is an $\R$-bounded LOTS.

 \begin{enumerate}
 \item[(a)] For any countable ordinal $\a$, $X^{\a}$ and $X_{\a}$ are
 $\R$-bounded LOTS.

 \item[(b)] If $X$ is order dense, then its completion is $\R$-bounded.

 \item[(c)] If $T$ is a tree of $X$ type $h(T) < \Omega$, then $X(T)$ is  $\R$-bounded.
 \end{enumerate} \end{prop}

 \begin{proof} On the one hand, $\R_{\d} \subset \R^{\d}$.  On the other, $\R \cong J^{\circ}$
 implies that $\R^{\d}$ injects into $\R_{\d}$, i.e. they have the same size.  So  $X$ injects into $\R_{\d}$ iff it injects into
 $\R^{\d}$.

 (a): If $X$ injects into $\R^{\d}$, then $X^{\a}$ injects into $(\R^{\d})^{\a} \cong \R^{\d \cdot \a}$ (see Proposition \ref{prop4.6}(d)).
 So $X^{\a}$ is $\R$-bounded when $X$ is and $\a$ is countable.
 Since $X_{\a} \subset X^{\a}$ it is $\R$-bounded as well.

 (b): If $X$ is  an order dense LOTS, and $j : X \to \R_{\d}$ is an order injection, then $\hat j$ is an order injection on the completion by
 Proposition \ref{prop2.4a}.

 (c): By Proposition \ref{prop6.4a} $X(T)$ injects into $X^{\a}$ if $h(T) = \a$. So by (a), $X(T)$ is $\R$-bounded.

 \end{proof} \vspace{.5cm}

 In particular, Proposition \ref{propnew04a} implies that if $T$ is a tree of $\Z$ type
 with height a countable limit ordinal, then $X(T)$ and its completion
 are $\R$-bounded. By sharpening the proof of Theorem \ref{theonew01} we obtain the following converse.

  \begin{theo}\label{theonew01a} If $X$ is a connected LOTS which admits an order injection $f : X \to \R_{\a}$ with
  $\a$ a positive ordinal, then $X \cong \widehat{X(T)}$ with $T$ a tree of $\Z$ type with height
  $h(T) \le \om \cdot \a$. \end{theo}

  \begin{proof}  With $J = [-1,+1] \subset \R$ we can, as before, replace $\R_{\a}$ by the isomorph $\{ x \in J^{\a} : -1 < x_0 < +1 \}$, ie.
  the interval $(-1+,+1-) \subset J^{\a}$. We let $Z $ be the closed interval $ [-1+,+1-] \subset J^{\a}$, which is isomorphic to the two point
  compactification of $\R_{\a}$. Assume that $f : X \to Z$ is an order injection. Let $\pi_i : Z \to J$ be the projection to the $i$ coordinate with
  $i < \a$.

  For $a < b$ in $X$ let $\ep =  min \{ i : f(a)_i \not= f(b)_i \} $ so that $f(a)_i = f(b)_i$ for all $i < \ep$ and $f(a)_{\ep} < f(b)_{\ep}$.
 We again have that $\pi_{\ep} \circ f( [a,b]) \subset J \subset \R$ and
so again algebraic notions and length are defined there. Notice that the order preserving map $\pi_{\ep} \circ f$ need not be continuous and
so its image on $[a,b]$ need not be connected. So we have to consider some cases.

Let $t_1 = f(a)_{\ep} + \frac{1}{3}(f(b)_{\ep} - f(a)_{\ep})$ and $t_2 = f(a)_{\ep} + \frac{2}{3}(f(b)_{\ep} - f(a)_{\ep})$. \vspace{.25cm}

\textbf{Case 1} (midpoint case): If $\pi_{\ep} \circ f( [a,b]) \cap [t_1,t_2] \not= \emptyset$, then let $c \in (a,b)$ such that
$f(c)_{\ep} \in [t_1,t_2]$.  Because $\pi_{\ep} \circ f$ is order preserving we have that $[a_1,b_1] \subset [a,c]$ or
$[a_1,b_1] \subset [c,b]$ implies $f(b_1)_{\ep} - f(a_1)_{\ep} \leq \frac{2}{3}(f(b)_{\ep} - f(a)_{\ep})$.\vspace{.25cm}

\textbf{Case 2} (edge cases): The left edge case applies when $\pi_{\ep} \circ f( (a,b)) \cap (f(a)_{\ep},t_1) = \emptyset$ whereas the right edge case
applies when $\pi_{\ep} \circ f( (a,b)) \cap (t_2,f(b)_{\ep}) = \emptyset$.  For either edge case,
$[a_1,b_1] \subset (a,b)$ implies $f(b_1)_{\ep} - f(a_1)_{\ep} \leq \frac{2}{3}(f(b)_{\ep} - f(a)_{\ep})$.\vspace{.25cm}

\textbf{Case 3} (boundary cases): If neither Case 1 nor Case 2 applies, then $A_1 = [a,b] \cap (\pi_{\ep} \circ f)^{-1}[f(a)_{\ep},t_1)$ is a
proper convex set which contain $a$ and $A_2 = [a,b] \cap (\pi_{\ep} \circ f)^{-1}(t_2,f(b)_{\ep}]$ is a proper convex set which contains $b$
and their union is $[a,b]$ since Case 1 does not apply. Because $X$ is connected, there exists a unique $c \in (a,b)$ such that
$[a,c) \subset A_1$ and $(c,b] \subset A_2$. If $[a_1,b_1] \subset [a,c)$ or
$[a_1,b_1] \subset (c,b]$, then  $f(b_1)_{\ep} - f(a_1)_{\ep} \leq \frac{2}{3}(f(b)_{\ep} - f(a)_{\ep})$. The inequality also holds
if $b_1 = c$ with $c \in A_1$ and if $a_1 = c$ with $c \in A_2$.

Finally, if $c \in A_1$, then the inequality does not hold for $[c,b_1]$, but now the interval $[c,b_1]$ is itself a left edge case. Similarly,
if $c \in A_2$, then the inequality does not hold for $[a_1,c]$, but the interval $[a_1,c]$ is a right edge case. \vspace{.25cm}

  As we did for Theorem \ref{theonew01}, we apply the  construction from the proof of Theorem \ref{theon.01}to $X$
  with the proviso that if $I_s = [a^s, b^s]$ is a proper interval, then the $\pm$cofinal sequence
  $z : \Z \to (a^s,b^s)$ maps $0 \in \Z$ to  the choice $c \in (a^s,b^s)$ when either Case 1 or Case 3 applies.   This then implies for all $n \in \Z$
\begin{equation}\label{eqnew02a}
 (z_{n+1})_\ep - (z_{n})_\ep \ \le \ \frac{2}{3}(b^s_{\ep} - a^s_{\ep}),
 \end{equation}
 except that in the boundary Case 3 with $n = 0$, we have for $c = z_0 \in A_1$, the interval $[z_0,z_1]$ is a left edge case and for $c = z_0 \in A_2$,
 $[z_{-1},z_0]$ is a right edge case.

Now we proceed as before showing, by induction, that for all $\b$ with $ 0 < \b \leq \a$
 \begin{equation}\label{eqnew03a}
 s \in \Z^{\om \cdot \b} \quad \Longrightarrow \quad f(a^s)_i = f(b^s)_i \ \ \text{for all} \ \ i < \b.
 \end{equation}\vspace{.25cm}

 \textbf{Case 1} $\b = \gamma + 1$:

 Let $s \in \Z^{\om \cdot \b}$. The induction hypothesis applied to $s0= s|\om \cdot \g$  implies that
 $f(a^{s0})_i = f(b^{s0})_i \ \ \text{for all} \ \ i < \g$. Since, $[a^s, b^s] \subset [a^{s0},b^{s0}]$,
 $f(a^{s0})_{\g} = f(b^{s0})_{\g}$ implies $f(a^{s})_{\g} = f(b^{s})_{\g}$ and so $a^s_i = b^s_i$ for all $i < \b$.

 Assume, instead, that $f(a^{s0})_{\g} < f(b^{s0})_{\g}$ so that $\ep = \g$
 for the interval $[a^{s0},b^{s0}]$. Let $sn = s|(\om \cdot \g + n)$ for $n < \om$.
 At worst, every other step shrinks length by a factor of $2/3$.  This suffices to show, as before, that
 $\om \cdot \b = \om \cdot \g + \om, \ I_s = \bigcap_n I_{sn}$ implies $f(b^s)_{\g} = f(a^s)_{\g}$. Hence,
 again $f(a^s)_i = f(b^s)_i$ for all $i < \b$. \vspace{.25cm}

  \textbf{Case 2} $\b$ is a limit ordinal:

 If $\g < \b$, then with $s' = s|\om \cdot \g$ we have $I_s \subset I_{s'}$. So the induction hypothesis applied to $s'$
 implies that $f(b^s)_i = f(a^s)_i $ for all $i < \g$. Since $\b$ is a limit ordinal and $\g < \b$ is arbitrary, it follows that
 $f(a^s)_i = f(b^s)_i$ for all $i < \b$. \vspace{.25cm}

 From (\ref{eqnew03a}) applied with $\b = \a$ we see that for all $s \in \Z^{\om \cdot \a}$, $f(a^s)_i = f(b^s)_i$ for all $i < \a$.
 Since $f(a^s), f(b^s) \in J^{\a}$ this means $f(a^s) = f(b^s)$. Because $f$ is injective, $a^s = b^s$. Hence, the interval $I_s$ is improper.

Thus, the tree $T$ consisting of those $s$ with $I_s$ proper is a subtree of the simple tree on $\Z, \om \cdot \a$ and so has
 height at most $\om \cdot \a$.

 \end{proof} \vspace{.5cm}

 \begin{theo}\label{theonew04} If $X$ is an $\R$-bounded IHLOTS, then the tower of
 CHLOTS $\widehat{X_{\om^{\g}}}$ is nondecreasing in size and is strictly increasing in size for
 sufficiently large $\g$ and the tower of CHLOTS Cantor Spaces $C(\widehat{X_{\om^{\g}}})$ is nondecreasing
 in size and is strictly increasing in size for
 sufficiently large $\g$.

 To be precise, if $X$ injects into $\R_{\om^{\g_0}}$, then
 $\widehat{X_{\om^{\g_1}}}$ is strictly bigger than $\widehat{X_{\om^{\g_2}}}$
 and $C(\widehat{X_{\om^{\g_1}}})$ is strictly bigger than $C(\widehat{X_{\om^{\g_2}}})$
 when $\g_1 > 1 + \g_0 + \g_2$.\end{theo}

 \begin{proof}  The precise estimate is clear from Corollary \ref{cornew03}. In particular, if $\g_3$ is the smallest tail-like
 ordinal larger than $\g_0$, then $1 + \g_0 + \g_2 = \g_2$ when $\g_2 \ge \g_3$. That is, beyond $\g_3$ the sequences are
 strictly increasing in size.

 \end{proof}\vspace{.5cm}

 By Corollary \ref{cor7.4} there exist CHLOTS which are not $\R$-bounded, e.g. $\widehat{a(\R)}$.

  \vspace{.5cm}

   Now with $X$ a CHLOTS we want to sharpen the  tree construction from Theorem \ref{theon.01} so that we obtain $T$ as an additive tree.

   We begin by choosing two points
   labeled $-1 < 1 \in X$ and apply the construction to $I = [-1,1]$.
   The image $f(X(T))$ of the branch space of the tree constructed in Theorem \ref{theon.01} is dense in $(-1,1)$ which is isomorphic to $X$.

   In order to obtain an additive tree we will have to choose the $\pm$cofinal
   embeddings of $\Z $ in a coherent way.    We will inductively define $t_s : I \to I_s$ which is the constant map when $I_s$ is improper and
   is an isomorphism when $I_s$ is proper. \vspace{1cm}

   \subsection{The Alphabet Construction}

   Regard a set $A$ as an  \emph{alphabet}\index{alphabet} and for every $a \in A$ there is  defined  a map $t_a : I \to I_a$
 with $I_a $ an interval contained in $I^{\circ}$. If $I_a$ is a proper interval, then $t_a$ is an isomorphism. If $I_a$ is improper,
 then $t_a$ is the constant map. In either case, $t_a$ is continuous. In addition, we let $I_{\emptyset} = I$ with $t_{\emptyset}$ the
 identity $1_I$.  We call $A_+ = \{ a \in A: I_a$ is
 a proper interval $\} \subset A$ the \emph{proper alphabet }\index{proper alphabet}\index{alphabet!proper}.

 We will extend these to $A^k$ for $ k < \om$ and to $ A^{\om}$, i.e. the spaces of finite and infinite sequences.

   A \emph{ word} $w \in A^k$ is a finite sequence $ a_0 a_2 \dots a_{k-1}$ in $A$ with $\emptyset$ the empty word of length $0$.
   For  a word $w \in A^k$ with $k \ge 1$, we define the map
     $ t_w = t_{a_0} \circ t_{a_1} \circ \dots t_{a_{k-1}}$ from
   $ I$ onto its image denoted $I_w$. We call the word \emph{proper} when all of its letters lie in the proper alphabet. In that case, $t_w$ is an isomorphism.
   If the interval $I_{a_i}$ is improper, then $t_w$ is the constant map onto the singleton $I_w = (t_{a_0} \circ \dots \circ t_{a_{i-1}})(I_{a_i})$.

   If $w$ is the concatenation $ w_1 w_2$ and $w_1$ is proper, then $I_w \subset I_{w_1}^{\circ}$.
   Notice that for proper words $w_1, w_2, w$, the composition $t_{w_1w}\circ (t_{w_2w})^{-1}  : I_{w_2w} \to I_{w_1w}$ is an isomorphism
   which is the restriction of the isomorphism $t_{w_1}\circ (t_{w_2})^{-1}  : I_{w_2} \to I_{w_1}$.

   We now consider the space of infinite sequences $A^{\om}$. Let
   $\t$ denote the \emph{shift map} on $A^{\om}$  with $\t(s)_i = s_{i+1}$.

    For $s \in A^{\om}$ we define the associated interval $I_s = \bigcap_n \ I_{s|n}$ with $s|n$ equal to the word $s_0 \dots s_{n-1}$.
   Because a continuous map commutes with the decreasing intersection of  compacta,
   we have for a word $w$ and $z \in A^{\om}$
   \begin{equation}\label{compeq01}
   I_{wz} = t_w(I_z). \hspace{2cm}
   \end{equation}
   Hence, $I_{wz}$ is proper if and only if $I_{z}$ is proper and the word $w$ is proper.

    When $I_s$ is improper, we let $t_s : I \to I_s$ be the constant map. In particular, this applies to any sequence $s \in A^{\om} \setminus A_+^{\om}$.

Now we restrict attention to sequences in the proper alphabet $A_+$.

   We call two sequences $s_1, s_2 \in A_+^{\om}$ \emph{end-equivalent}\index{end-equivalence} if there exist
   $i_1, i_2 \in \om$ such that $\t^{i_1}(s_1) = \t^{i_2}(s_2)$,
   or, equivalently, when there exist  proper words $w_1, w_2$ and $z \in A_+^{\om}$ such that
   $s_1 = w_1 z, s_2 = w_2 z$. We call them \emph{level end-equivalent}\index{level end-equivalence} when, in addition, $i_1$ and $i_2$ can be chosen equal,
   or, equivalently, the words $w_1$ and $w_2$ have the same length. We call the end-equivalence class of $s$ the \emph{class}
  \index{class!end-equivalence} of $s$.
   The class of $s$ is subdivided into \emph{level classes}\index{class!level end-equivalence}.

    If  for a sequence $s$ the interval
   $I_s$  is proper, then by (\ref{compeq01}) $I_{s_1}$ is proper for every sequence $s_1$ end-equivalent to $s$.
   We will call the class \emph{proper}\index{class!proper} when it
   consists of sequences with proper associated intervals.

   A sequence $s$ is \emph{eventually periodic}\index{sequence!eventually periodic}
   if there exist $i, p \in \N$ with $p > 0$ such that $\t^i(s) = \t^{i+p}(s)$ and so
   $\t^i(s) = \t^{i+np}(s)$ for all $n \in \N$. The minimum such
   $p$ is called the \emph{period} of $s$. Clearly if $s_1$ is end-equivalent to $s$ and $s$ is eventually periodic then $s_1$ is eventually
   periodic with the same period. We call a class \emph{periodic}\index{class!periodic}
   if it consists of eventually periodic sequences.  Otherwise we call the class \emph{nonperiodic}\index{class!nonperiodic}.

   If $s$ is not eventually periodic then the sequences $\t^i(s)$ are all distinct. In fact, all lie in different level classes.
 So if $s_1, s_2$ are in a
   nonperiodic class and $k_1, k_2, j_1, j_2 \geq 0$, then
   \begin{equation}\label{compeq02}
   \t^{k_1}(s_1) = \t^{k_2}(s_2) \ \text{and} \ \t^{k_1+j_1}(s_1) = \t^{k_2+j_2}(s_2) \ \Longrightarrow \quad j_1 = j_2.
   \end{equation}

   When a class  is proper,  we will choose a representative element $s$ for the class
   and choose $t_s : I \to I_s$ an isomorphism for this representative. \vspace{.5cm}

   \textbf{Case 1} Nonperiodic Proper Class:  For a nonperiodic proper class let $t_s : I \to I_s$ be the chosen
   isomorphism for the chosen representative.

   If $s_1$ is end-equivalent to $s$ then we let $k$ be the minimum such that $\t^k(s) = \t^{k_1}(s_1)$ for some $k_1$ and then
   choose $k_1$ minimum. Thus, there are unique words $w, w_1$ of length $k,k_1$ and
    $z \in A_+^{\N}$  such that
   $s = wz$, $ s_1 = w_1z$. By (\ref{compeq01}) $t_w|I_z : I_z \to I_s$ and $t_{w_1}|I_z : I_z \to I_{s_1}$ are isomorphisms.
    We define $t_{s_1} =  t_{w_1} \circ (t_w)^{-1} \circ t_s$.

      If $s = w'z'$,  $s_1 = w'_1 z'$ for finite words $w', w_1'$ and $z' \in A_+^{\om}$ then by minimality there exists a word $u$
    such that $w' = wu$ and $z = uz'$.  So by (\ref{compeq02}) it follows that $w_1' = w_1u$.
    Again
    $t_w \circ t_u = t_{w'}$ restricts to an isomorphism from $I_{z'} $ to $I_s$ and
    $t_{w_1} \circ t_u = t_{w_1'}$ restricts to an isomorphism from $I_{z'} $ to $I_{s_1}$. So we have
     \begin{equation}\label{compeq03a}
     t_{w'_1} \circ (t_{w'})^{-1} \circ t_s = t_{w_1} \circ t_u \circ (t_w \circ t_u)^{-1} \circ t_s = t_{w_1} \circ (t_w)^{-1} \circ t_s = t_{s_1}.
     \end{equation}

  So if $s_2 = w_2$  $s_1 = w_2 w_1 z$,  then $s_2$ is end-equivalent to $s_1$ and
   \begin{equation}\label{compeq03}
   t_{s_2} = t_{w_2w_1}\circ (t_{w})^{-1}\circ t_s = t_{w_2} \circ t_{s_1}. \hspace{2cm}
   \end{equation}

   It follows that
   \begin{equation}\label{compeq04}
   s_1 = w_1 z,\;  s_2 = w_2 z \quad \Longrightarrow \quad   t_{s_1} = t_{w_1} \circ  t_z = t_{w_1} \circ (t_{w_2})^{-1} \circ t_{s_2}.
    \end{equation}\vspace{.5cm}

   \textbf{Case 2} Periodic Proper Class: If $s$ is eventually periodic with period $p$, and so for some $i$
   $\t^i(s) = \t^{i+p}(s)$ then there exists a finite word $ e = e_0 e_1 \dots e_{p-1}$ such that
   $\t^i(s) =  e\t^{i + p}(s)$. We let $\bar e \in A_+^{\om}$ be the periodic element
   in the class with $\bar e_i = e_j$ if $i$ is congruent to $j$ mod $p$. Since $ \bar e =  e \bar e$
   we have $t_{e}(I_{\bar e}) = I_{\bar e}$.
   That is, $t_{ e}$ restricts to an automorphism of $I_{\bar e}$ and so, in particular, fixes the endpoints of $I_{\bar e}$.

   While the period $p$ is uniquely associated with all members of the class, the minimum block $e$ is not unique when $p > 1$.
  For $i = 1, \dots, p-1$ cyclic permutations $  e_i \dots e_{p-1} e_0 \dots e_{i-1} $ are blocks of length $p$ with
  $\overline{e_i \dots e_{p-1} e_0 \dots e_{i-1}}$
   end-equivalent to $\bar e$. To be precise, $e_i \dots e_{p-1}\bar e = \overline{e_i \dots e_{p-1} e_0 \dots e_{i-1}}$.

   The periodic class
   of period $p$ consists of $p$ level classes, each containing one of the periodic elements $\overline{e_i \dots e_{p-1} e_0 \dots e_{i-1}}$.

   Assume that the periodic class is proper. We fix a minimum block $ e$ and use $\bar e$ as the representative of the
   class. Then we choose the isomorphism $t_{\bar e} : I \to I_{\bar e}$. We choose $\overline{e_i \dots e_{p-1} e_0 \dots e_{i-1}}$ as the representative of
   its level class and define
   $t_{\overline{e_i \dots e_{p-1} e_0 \dots e_{i-1}}} = t_{e_i \dots e_{p-1}} \circ t_{\bar e}$.

   Now we operate in each level class
   separately using its unique periodic element as representative. We look at the level class containing $\bar e$.

   A sequence $s_1$ is in the level class of $\bar e$  if and only if  there exists a word $w$ of length $|w| = np$ for some $n \in \om$ such that $s_1 = w \bar e$.
   We choose $w$ to be the unique such word with $n$ minimum and
   define $t_{s_1} = t_w \circ (t_{ e})^{-n} \circ t_{\bar e}$, with $(t_{e})^{-n}$ the $n$-fold iterate of $(t_{e})^{-1} $
   (= identity when $n = 0$).

   If $s_1 = w' \bar e$ then $|w'| = n'p$ for some $n' \in \om$. It follows that $n' = n + k$ for some $k \in \om$ and $w' = w (e)^k$.
  \begin{equation}\label{compeq05a}
  t_{w'}  \circ (t_{ e})^{-n'} \circ t_{\bar e} = t_w \circ t_{ e}^k \circ (t_{ e})^{-(n+k)} \circ t_{\bar e} = t_{s_1}.
  \end{equation}

   So we have
   \begin{equation}\label{compeq05}\begin{split}
   s_1 = w_1 \bar e, s_2 = w_2 \bar e, \ \ \text{and} \ \ |w_1| = np = |w_2|  \quad \Longrightarrow \\
      t_{s_1} = t_{w_1}\circ (t_{ e})^{-n} \circ t_{\bar e} = t_{w_1} \circ t_{w_2}^{-1} \circ t_{s_2}. \hspace{2cm}
    \end{split}\end{equation} \vspace{.25cm}

This completes the Alphabet Construction.

\begin{ex}\label{exn02} Proper intervals and fixed points.\end{ex}

 If, in the Alphabet Construction, $e$ is a finite word with $\bar e$ the associated periodic sequence and $I_{\bar e} = [a,b]$
    interval, then $t_e$ is an automorphism of $[a,b]$ and so the endpoints are  fixed points of $t_e$.  They are distinct if the interval is proper.

The set of fixed points $\{ x \in I : t_e(x) = x \}$ is closed with minimum $m$ and maximum $M$. If $m < M$, then
    $t_e(I)$ is an interval which contains $m$ and $M$ and so $[m,M] \subset t_e(I)$ and so by induction $[m,M] \subset t_e^k(I)$ for all
    $k \in \om$.  It follows that $[a,b] = [m.M]$. Thus, $\bar e$ is improper  if and only if  $t_e$ has a unique fixed point.
    Finally, if $[a,b] \subset [a_1,b_1]^{\circ}$ and $[a_1,b_1] \subset I^{\circ}$, then there exist isomorphisms
       $[-1,a] \to [a_1,a]$ and $[b,1] \to [b,b_1]$. Combining these with the identity on $[a,b]$ we can define $t_e : I \to [a_1,b_1]$ which
    fixes $[a,b]$. \vspace{1cm}

 \subsection{The Additive Tree for a CHLOTS}

   We now proceed with our inductive construction of the mappings $t_s$. As part of the construction  we will  prove the following: \vspace{.5cm}

   \textbf{Composition Property} Let $s_1, s_2 \in \Z^{\a}$. If $\b$ is an ordinal with $\b  < \a$ such that $\t_{\b}(s_1) = \t_{\b}(s_2)$, and
   $s_1|\b$ and $s_2|\b$ are proper, then $s_1$ is proper if and only if $s_2$ is proper and in that case
   then
   \begin{equation}\label{compeq06a}
   t_{s_1} \ = \ (t_{s_1|\b}) \circ (t_{s_2|\b})^{-1} \circ t_{s_2}. \hspace{1cm}
  \end{equation} \vspace{.5cm}

  If an ordinal $\a$ is a sum of ordinals $\a = \a_1 + \a_2 + \dots \a_k$, we will write $\s_i = \a_1 + \dots \a_i$ for $i = 1, \dots k$
   and for $s \in \Z^{\a }$ we will write
  \begin{equation}\label{compeq06b}\begin{split}
  s_1 = s|\a_1, s_2 = \t_{\a_1}(s)|\a_2, s_3 = \t_{\s_2}(s)|\a_3, \dots, s_k = \t_{\s_{k-1}}(s).
 \end{split}  \end{equation}
  So that $ s = s_1 + s_2 + \dots s_k $ in the simple tree. \vspace{.5cm}

   {\bfseries Step 1} $\a \le \om$ :

   We first let $t_{\emptyset}$ be the identity on $I$.

   We then choose a $\pm$cofinal embedding of $\Z $ into $(-1,1)$ and for each $n \in \Z $ we choose an isomorphism $t_n : I \to [z_n,z_{n+1}]$.

We apply the Alphabet Construction with $A = \Z$ and using $t_n$ for $n \in \Z$.
So in this case, the entire alphabet is
   proper and from the Alphabet Construction we obtain $t_s : I \to I_s$ for every $s \in \Z ^{\a}$ with $\a \le \om$. Because a class
   is either entirely proper or entirely improper, the Composition Property follows from
   (\ref{compeq04}) in the nonperiodic case and from (\ref{compeq05}) in the periodic case. \vspace{.25cm}

{\bfseries Step 2}  $\om^{\gamma} < \a < \om^{\gamma + 1}$ with $\gamma \ge 1$ :

  For $\a$ between $\om^{\gamma}$ and $ \om^{\gamma+1} $ Cantor Normal Form is
  $\a = \a_1 + \a_2 + \dots \a_k$ with  $\omega^{\gamma} = \a_1 \ge \a_2 \ge \dots \a_k$ tail-like.
   As in (\ref{compeq06b})  for
  $s \in \Z^{\a }$ we have the decomposition $ s = s_{1} + s_{2} + \dots s_{k}$.

   We use the induction hypothesis to define
  \begin{equation}\label{compeq00a}
  t_{s} = t_{s_{1}} \circ t_{s_{2}} \circ \dots t_{s_{k}} : I \to I_s
  \end{equation}
  with $I_s$ the image of $t_s$.

   Thus, $s$ is  proper  if and only if each $s_i$ is proper and the composition
   is then an isomorphism. If, instead, $i$ is the smallest index such that $s_i$ is improper, then $I_s$ is the singleton
   $t_{s_1} \circ t_{s_2} \circ \dots t_{s_{i-1}}(I_{s_i})$ and the constant map $t_s$ is the composition $t_{s_1} \circ t_{s_2} \circ \dots t_{s_{i}} $.

  For the Composition Property let $s_1, s_2 \in \Z^{\a}$ and for $\ep = 1,2$  we  decompose
  $ s_{\ep} = s_{\ep 1} + s_{\ep 2} + \dots s_{\ep k}$. Assume $\b < \a$ such that $\t_{\b}(s_1) = \t_{\b}(s_2)$ and $s_{\ep}|\b$ is proper
  for $\ep = 1,2$.

  Let $i \leq k$ be the minimum such that $ \b < \s_i$ and let
  $\b_1 = \b \setminus \s_{i-1}$. If $i = 1$, then $\s_0 = 0$ by convention and $\b_1 = \b$.

  Since $\t_{\b}(s_1) = \t_{\b}(s_2)$ we have $s_{1j} = s_{2j}$ for  $ j = i+1, \dots k$
  and for such $j$ we will write $s_j$ for $s_{1j} = s_{2j}$. In addition, $\t_{\b_1}(s_{1i}) = \t_{\b_1}(s_{2i})$.

  So
  from (\ref{compeq00a}) we obtain (since $\b_1 < \a_i$):
   \begin{align}\label{compeq06}
   \begin{split}
   t_{s_{\ep}|\s^i}\  &= \ t_{s_{\ep 1}} \circ \dots t_{s_{\ep i}}, \\
 t_{s_{\ep}|\b} \ &= \ t_{s_{\ep 1}} \circ \dots t_{s_{\ep (i-1)}} \circ t_{s_{\ep i}|\b_1}, \\
  t_{s_{\ep}} \ &= \ t_{s_{\ep}|\s^i} \circ t_{s_{i+1}} \circ \dots t_{s_k}.
  \end{split}\end{align}

  Because $s_1|\b$ and $s_2|\b$ are proper,
 the maps $t_{s_{\ep 1}},t_{s_{\ep (i-1)}},t_{s_{\ep i}|\b_1}$ are isomorphisms.

    By the Composition Property for $s_{1i}, s_{2 i} \in \Z^{\a_i}$ we have that $s_{1i}$ is improper if and only if $s_{2i}$ is improper
     in which case $s_1$ and $s_2$ are both improper. In addition if any of $s_{i+1}, \dots, s_k$ are improper, then both $s_1$ and $s_2$
     are improper.

    Assume, instead, that both $s_{1i}$ and $s_{2i}$ as well as $s_{i+1}, \dots, s_k$ are proper.
    The Composition Property for $s_{1i}, s_{2 i}$ implies
   $$t_{s_{1i}} = (t_{s_{1i}|\b_1}) \circ (t_{s_{2i}|\b_1})^{-1} \circ t_{s_{2i}}.$$
   From this and (\ref{compeq06}) it follows that
   $$t_{s_1|\s^i} = (t_{s_1|\b}) \circ (t_{s_2|\b})^{-1} \circ t_{s_2|\s^i}.$$
   and we compose with $t_{s_{i+1}} \circ \dots t_{s_k}$ to obtain the Composition Property for
   $s_1$ and $s_2$.

   If $\a' = \a  + 1$ is the successor  of  $\a $, then the Cantor Normal Form for $\a'$ is $\a' = \a_1 + \a_2 + \dots \a_k + \a_{k+1}$ with
   $\a_{k+1} = 1$. Recall that $1$ is the unique tail-like ordinal which is not a limit ordinal. If $s' \in \Z ^{\a'}$ with
   $s'|\a  = s$, then $t_{s'} = t_s \circ t_n$ where $n = s'(\a )$. Thus, in both Step 1 and Step 2, we are using the $\pm$cofinal
   map $t_s \circ z : \Z \to I_s^{\circ}$ for $I_s$, where $z : \Z \to I$ is the $\pm$cofinal map with which we began.
 so that, inductively, the intervals $I_s$ are those obtained in the construction of
    Theorem \ref{theon.01}.   As before, each successor $s'$ is proper when $s$ is.\vspace{.25cm}

  {\bfseries Step 3} $\a = \om^{\gamma+1} $ :

   We have $\Z^{\om^{\gamma+1}} = \Z^{\om^{\gamma}\cdot \om}$ which we identify with
  $(\Z^{\om^{\gamma}})^{\om}$ as in (\ref{eq4.20}). That is,
   we regard an element of $\Z^{\om^{\gamma+1}}$  as a sequence of elements of the alphabet $A = \Z^{\om^{\gamma}}$.
  The proper alphabet consists  of the proper  elements of $\Z^{\om^{\gamma}}$.

   We now apply the  Alphabet Construction  to
   define $t_s : I \to I_s$ for all $s \in \Z^{\om^{\gamma + 1}}$. In particular, from the Alphabet Construction we obtain (iii) of
  the construction for  Theorem \ref{theon.01}.

   For the Composition Property, let $s_1, s_2 \in \Z^{\om^{\gamma+1}}$ and $\b$ be an ordinal less than $\om^{\gamma+1}$. Assume that
   $\t_{\b}(s_1) = \t_{\b}(s_2)$ and that $s_1|\b$ and $s_2|\b$ are proper.

   Let $k$ be the minimum in $\om$ such that
   $\b \le  \om^{\gamma} \cdot k$.
  Regarding $s_1$ and $s_2$ as sequences,  then $\t_{\b}(s_1) = \t_{\b}(s_2)$ implies that
   $s_1 = w_1 z, s_2 = w_2 z$ with $z \in \Z^{\om^{\gamma+1}}$ and $w_1, w_2$ words of length $k$ in the alphabet, or equivalently
   elements of $\Z^{\om^{\gamma} \cdot k}$.
   Furthermore, $s_{\ep}|\b = w_{\ep}|\b$ for $\ep = 1,2$ and we are assuming that these are proper.

    From the Composition Property applied to $w_1, w_2 \in \Z^{\om^{\gamma}\cdot k}$ we have that $w_1$ is improper if and only if $w_2$ is improper
   in which case both $s_1$ and $s_2$ are improper. In addition if any of the $\Z^{\om^{\gamma}}$ terms of the sequence $z$ is improper, then
   both $s_1$ and $s_2$ are improper.

  Assume, instead, that all of the terms of the sequences $s_1$ and $s_2$ lie in the proper alphabet.
   The Composition Property applied to $w_1, w_2$ then implies
  \begin{equation}\label{compeq07a}
   t_{w_1} = t_{w_1|\b} \circ t_{w_2|\b}^{-1} \circ t_{w_2}.
\end{equation}

  Since $s_1$ and $s_2$ are end-equivalent sequences in the proper alphabet, we can apply the Alphabet Construction results.
    From (\ref{compeq04}), or in the eventually periodic case from (\ref{compeq05}), we obtain from (\ref{compeq07a})
    \begin{equation}\label{compeq07}
    \begin{split}
  t_{s_1} = t_{w_1} \circ t_{w_2}^{-1} \circ t_{s_2} = \hspace{1cm} \\
  t_{w_1|\b} \circ t_{w_2|\b}^{-1} \circ t_{w_2} \circ t_{w_2}^{-1} \circ t_{s_2} \\
  =  t_{w_1|\b} \circ t_{w_2|\b}^{-1}  \circ t_{s_2}, \hspace{1cm} \\
   =  t_{s_1|\b} \circ t_{s_2|\b}^{-1}  \circ t_{s_2}, \hspace{1cm}
   \end{split}\end{equation}
  proving the Composition Property for $s_1$ and $s_2$. \vspace{.25cm}

   {\bfseries Step 4}  $\a = \om^{\gamma} $ with $\gamma$ a countable limit ordinal :

   For $s \in \Z^{\om^{\gamma}}$ if $s|\a$ is improper for some $\a  < \om^{\gamma}$ then $s$ is improper with a constant $t_{s}$.
   We call $s$ \emph{limit proper}\index{sequence!limit proper} if, instead, $s|\a$ is proper for all $\a  < \om^{\gamma}$.

   Call two elements $s_1, s_2 \in \Z^{\om^{\gamma}}$ \emph{end-equivalent}\index{end-equivalence} if there exists $\b  < \om^{\gamma}$ such that
   $\t_{\b }(s_1) = \t_{\b }(s_2)$.

   We consider an end-equivalence class of limit proper elements. From the Composition Property, it follows that for every
   $\a$ with $\b < \a < \om^{\gamma}$, $t_{s_1|\a} = (t_{s_1|\b}) \circ (t_{s_2|\b})^{-1} \circ t_{s_2|\a}$
   and so $I_{s_1|\a} = (t_{s_1|\b}) \circ (t_{s_2|\b})^{-1}(I_{s_2|\a})$. Intersecting as $\a \to \om^{\gamma}$
   we obtain
     \begin{equation}\label{compeq08}
   I_{s_1} = (t_{s_1|\b}) \circ (t_{s_2|\b})^{-1}(I_{s_2})
   \end{equation}

   In particular, $s_1$ is improper if and only if $s_2$ is improper.

   Assume, instead, that the end-equivalence class  consists of proper elements. We choose
   a representative $s$ and an isomorphism $t_s : I \to I_s$. For $s_1$ end-equivalent to $s$
   let $\b_1 $ be the smallest  ordinal such that $\t_{\b_1 }(s) = \t_{\b_1 }(s_1)$. We define
   $t_{s_1} = (t_{s_1|\b_1}) \circ (t_{s|\b_1})^{-1} \circ t_{s}$.

   If $\b_2 > \b_1$ then $\t_{\b_2 }(s) = \t_{\b_2 }(s_1)$ and by the Composition Property for $s_1|\b_2$ and $s_2|\b_2$,
   $t_{s_1|\b_2} = (t_{s_1|\b_1}) \circ (t_{s|\b_1})^{-1} \circ t_{s|\b_2}$.
   That is, $(t_{s_1|\b_2}) \circ (t_{s|\b_2})^{-1}$ is a restriction of $(t_{s_1|\b_1}) \circ (t_{s|\b_1})^{-1}$
   and so
    \begin{equation}\label{compeq09}
   t_{s_1} = (t_{s_1|\b_2}) \circ (t_{s|\b_2})^{-1} \circ t_{s}.
   \end{equation}

   Finally, for the Composition Property let $s_1, s_2 \in \om^{\gamma}$  and $\b < \om^{\gamma}$ be an ordinal such that $\t_{\b }(s_1) = \t_{\b }(s_2)$
   with both $s_1|\b$ and $s_2|\b$  proper.

   If neither $s_1$ nor $s_2$ is limit proper then both $I_{s_1}$ and $I_{s_2}$ are improper. Now assume that $s_1$ is limit proper. If $\a \le \b$
   then $s_2|\a = (s_2|\b)|\a$ is proper.  If $\b < \a$, then by the Composition Property applied to $s_1|\a$ and $s_2|\a$, $s_2|\a$ is proper
   because $s_1|\a$ is proper.  So we may assume that both $s_1$ and $s_2$ are limit proper.
   They are clearly end-equivalent. If the end-equivalence class consists of improper elements, then both $s_1$ and $s_2$ are improper.

   Assume, instead, that the class consists of proper elements and let  $s$ be the representative of the class.

   We can choose $\b_1 > \b$ so that $\t_{\b_1 }(s_1) = \t_{\b_1 }(s_2) = \t_{\b_1 }(s)$.
   \begin{equation}\label{compeq10}
   \begin{split}
   t_{s_1} = (t_{s_1|\b_1}) \circ (t_{s|\b_1})^{-1} \circ t_{s} = \hspace{1cm}\\
   (t_{s_1|\b_1}) \circ (t_{s|\b_1})^{-1} \circ (t_{s|\b_1}) \circ (t_{s_2|\b_1})^{-1} \circ t_{s_2} = \\
   (t_{s_1|\b_1})  \circ (t_{s_2|\b_1})^{-1} \circ t_{s_2} = (t_{s_1|\b})  \circ (t_{s_2|\b})^{-1} \circ t_{s_2},
   \end{split}\end{equation}
   completing the proof of the Composition Property for $s_1$ and $s_2$.\vspace{.25cm}

   {\bfseries Step 5}  $\a = \Omega$ :

   Because the tree $T$ is $\Omega$ bounded, every $s \in \Z^{\Omega}$ is improper and so $t_s$ is the constant map onto $I_s$.\vspace{.25cm}

   This completes the inductive construction of the tree $T$. \vspace{.5cm}

     \begin{theo}\label{theon03}  The tree $T$ is an $\Omega$ bounded additive tree of $\Z$ type. \end{theo}

   \begin{proof} We saw in   Theorem \ref{theon.01} that the tree is $\Omega$ bounded.

   We must show that if $s \in \Z^{\a}$ and $\a = \a_1 + \a_2$, then $s$ is proper if and only if
   both $s|\a_1$ and $\t_{\a_1}(s)$ are proper. Since $s$ proper implies $s|\a_1$ is proper, we assume that $s_1 = s|\a_1$ is proper and
   prove that $s_2 = \t_{\a_1}(s)$  is proper if and only if $s = s_1 + s_2$ is proper. \vspace{.5cm}

   First we reduce to the case when $\a_1$ is tail-like.  That is, assume the result is true when $\a_1$ is tail-like and in general
   write $\a_1 = \a_{11} + \a_{12} + \dots \a_{1k}$ in Cantor Normal Form and so write $s_1 =  s_{11} + s_{12} + \dots s_{1k}$.

   Because $s_1$ is proper, the definition (\ref{compeq00a}) implies that $s_{11}, s_{12},\dots s_{1k}$ are all proper and so are the partial sums.
   We then complete the proof by induction on $k$ with $k=1$ the assumed, tail-like, case.
We have
     \begin{equation}\label{compeq10a}
     s = s_1 + s_2 = (s_{11} + \dots s_{1(k-1)}) + (s_{1k} + s_2).
     \end{equation}

   By the inductive hypothesis, $s$ is proper if and only if $s_{1k} + s_2$ is proper and so,
   by the initial case, if and only if $s_2$ is proper.\vspace{.5cm}

   Now we prove the result assuming that $\a_1$ is tail-like and so equals some $\om^{\gamma_1}$. Now we proceed by induction on $\a_2$. \vspace{.5cm}

   {\bfseries Step 1} $\a_2 < \om^{\gamma_1 + 1} = (\a_1)^{\om}$ :

   Write $\a_2 = \a_{21} + \a_{22} + \dots \a_{2 \ell}$ in Cantor Normal Form and
   so write $s_2 =  s_{21} + s_{22} + \dots s_{2 \ell}$.

  The assumption $\a_2 < (\a_1)^{\om}$ implies that $\a_1 \ge \a_{21}$ and so
   $ \a_1 + \a_{21} + \a_{22} + \dots \a_{2 \ell}$ is Cantor Normal Form for $\a = \a_1 + \a_2$. Hence,
   (\ref{compeq00a}) implies that $s$ is proper if and only if $s_1$ and $s_{21}, s_{22},\dots s_{2 \ell}$ are all proper. Furthermore,
   $s_2$ is proper if and only if  $s_{21}, s_{22},\dots s_{2 \ell}$ are all proper. Since $s_1$ is assumed proper, it follows that
   $s$ is proper if and only if $s_2$ is proper, as required. \vspace{.25cm}

   {\bfseries Step 2} $\a_2 = \om^{\gamma + 1}$ with $\gamma_1 \leq \gamma$ :

   As before we identify $\Z^{\om^{\gamma + 1}} = \Z^{\om^{\gamma} \cdot \om}$ with $(\Z^{\om^{\gamma} })^{\om}$
   and so regard $\Z^{\om^{\gamma + 1}}$ as sequences on the alphabet $\Z^{\om^{\gamma}}$ with
   the proper alphabet consisting of the proper elements of $\Z^{\om^{\gamma}}$.
   Now we split into two subcases. \vspace{.25cm}

   {\bfseries Step 2a} $\gamma_1 = \gamma$ :

   In this case,
   $s_1$  is  a letter in the proper alphabet which we label $w$,
   $s_2 = z$ is a sequence in the alphabet and  $s = wz$. If any of the terms in $z$ is not proper, then both $s$ and $s_2$ are
   improper.

   Assume instead that $z$ is a sequence in the proper alphabet, and so  $wz$ is a sequence in the proper alphabet as well.
   Because $wz$ and $z$ are end-equivalent
   sequences, it follows from the Alphabet Construction that $s_2$ is proper if and only if $s$ is. \vspace{.25cm}

   {\bfseries Step 2b} $\gamma_1 < \gamma$ :

   This time we write $s_2 = wz$ for the sequence in the alphabet where we use $w$ to label the
  first term. So $w \in \Z^{\om^{\gamma}}$ and since $\a_1 < \om^{\gamma}$, $w_1 = s_1 + w $ is an element of
   $ \Z^{\om^{\gamma}}$ as well, with
   $s = w_1z$.

   By the inductive hypothesis, $w$ is improper if and only $w_1$ is improper in which case both $s$ and $s_2$ are improper.
   In addition, if any term in $z$ is improper, then both $s$ and $s_2$ are improper.

   Assume, instead, that both $s$ and $s_2$ are
   sequences in the Proper Alphabet. Because $wz$ and $w_1z$ are end-equivalent
   sequences, it follows from the Alphabet Construction that $s_2$ is proper if and only if $s$ is. \vspace{.25cm}

   {\bfseries Step 3} $\om^{\gamma} < \a_2 < \om^{\gamma + 1}$ with $\gamma_1 < \gamma$ :

   We return to Cantor Normal Form
   $\a_2 =  \a_{21} + \a_{22} + \dots \a_{2 \ell}$ with $\a_{21} = \om^{\gamma}$ and we write $s_2 =  s_{21} + s_{22} + \dots s_{2 \ell}$.
   Again  $s_2$ is proper if and only if  $s_{21}, s_{22},\dots s_{2 \ell}$ are all proper.

   By inductive hypothesis, $s_{21}$ is proper if and
   only if $s_1 + s_{21}$ is proper. But $\a_1 + \a_{21} = \a_{21}$ because $\a_1 < \a_{21}$ and $\a_{21}$ is tail-like.
   Hence, $\a_1 + \a_2 = \a_{21} + \a_{22} + \dots \a_{2 \ell}$ in Cantor Normal Form with $s_1 + s_2 = (s_1 + s_{21}) + s_{22} + \dots s_{2 \ell}$.
   Thus, $s = s_1 + s_2$ if proper if and only if $(s_1 + s_{21}), s_{22},\dots s_{2 \ell}$ are all proper. It follows that
   $s_2$ is proper if and only if $s$ is proper. \vspace{.25cm}

  {\bfseries Step 4} $\a_2 = \om^{\gamma}$ with $\gamma$ a limit ordinal such that $\gamma_1 < \gamma$ :

  Because $\gamma$ is a limit ordinal larger than $\gamma_1$, it is larger than $\gamma_1 + 1$. So  $\b = \om^{\gamma_1 + 1}$ is
  a tail-like ordinal with  $\a_1 < \b < \a_2$. Hence $\a_1 + \b = \b$ and $\a_1 + \a_2 = \a_2$.

   Because $\a_1 + \a_2 = \a_2$, $s, s_2 \in \Z^{\om^{\gamma}}$.

   Because $\a_1 + \b = \b$,  $\t_{\b}(s) = \t_{\b}(s_2)$ and $s|\b = s_1 + (s_2|\b)$.

   By inductive hypothesis, $s|\b$ is proper if and only if $s_2|\b$ is proper. If these are improper
   then both $s$ and $s_2$ are improper.

   Assume, instead, that both $s|\b$  and $s_2|\b$ are proper.

   For any $\b_1$ with $\b < \b_1 < \a_2$, the Composition Property for $s|\b_1$ and $s_2|\b_1$ implies that
   $s|\b_1$ is proper if and only if $s_2|\b_1$ is proper. If these are improper
   then both $s$ and $s_2$ are improper.

    Assume, instead, that both $s|\b_1$  and $s_2|\b_1$ are proper for
   all $\b_1 < \a_2$. That is, $s$ and $s_2$ are end-equivalent limit proper elements of $\Z^{\om^{\gamma}}$. From Step 4 of the above
   construction of the tree,  it follows that $s$ is proper if and only if $s_2$ is proper.

   \end{proof} \vspace{.5cm}

   \begin{ex}\label{exn04} Trees for $X = \R$. \end{ex}

   If $F$ is a countable unbounded LOTS, like $\Z$ or $\Q$, then $F^{\om}$ can be regarded as the branch space of the simple tree on $F, \om$.
   By Proposition \ref{prop5.3a} it is order dense. We can also think of $F$ as an alphabet and $F^{\om}$ as the space of sequences on $F$.
   Any end equivalence class is a countable dense set and so $F^{\om}$ is separable. By Proposition \ref{prop2.8}(b)
   the compactification $\bullet \widehat{F^{\om}} \bullet$
   is isomorphic to the unit interval in $\R$ and so $\widehat{F^{\om}} \cong \R$. It follows that $F^{\om}$ is order isomorphic to a dense
   subset of $\R$.  \vspace{.25cm}

  In the $F = \Z$ case we let
   \begin{equation}\label{eqn.05}
    z_n = \begin{cases} \ -1 +  2^{n} \ \ \text{for } \ n < 0, \\ \ \ 1 - 2^{-n}\ \ \text{for } \ n \ge 0. \end{cases}.
    \end{equation}
   This defines a $\pm$cofinal embedding of $\Z$ into $I = (-1,1) \subset \R$.
   Let $t_n$ from $I$ onto $I_n = [z_n,z_{n+1}]$ be the restriction of the affine map on $\R$ given by
  \begin{equation}\label{eqn.06}
  t_n(z) = \frac{1}{2}[(z_{n+1} + z_n) + z \cdot (z_{n+1}-z_n)].
  \end{equation}
    So the derivative $t_n'$ (equals the slope) is bounded by $1/4$ for all $n$.
   If $w \in \Z^k$, then the length in $\R$ of the interval $I_w$ is bounded by $4^{-k}$. It follows that
    for every $s \in \Z^{\om}$, $I_s$ is improper. Hence, $T$
   is the simple tree on $\Z, \om$ and $X(T) = \Z^{\om}$.

   Because there are countably many finite words, $C$ is a countable dense subset of $[-1,1]$, including $-1$ and $1$ and
   $I$ is the disjoint union of $X(T)$ and $C$.  Since $C$ is order-isomorphic to $\Q \cap I$, it follows that $X(T) \cong \Z^{\om}$
   is order-isomorphic to the
   set of irrationals in $I$ and so to the set of irrationals in $\R$. \vspace{.25cm}

   If $F = \N$, the map $z : \N \to [0,1)$ defined by the restriction of the map in (\ref{eqn.05}) is a cofinal embedding into $\tilde I = [0,1)$.
Now we apply the Alphabet Construction with $A = \N$ and $\tilde t_n$ given by
  \begin{equation}\label{eqn.06a}
  \tilde t_n(z) =  z_n + z \cdot (z_{n+1}-z_n).
  \end{equation}
   For every finite word $w$ let
$\tilde I_w = \tilde t_w(\tilde I)$. For each $k < \om$,  $ [0,1)$ is the disjoint union of
$\{ \tilde I_w : w \in \N^k \}$.
As above for every $s \in \N^{\om}$, $\tilde I_s = \bigcup_{k < \om} \tilde I_{s|k}$ is a singleton
 and so the branch space of the tree can be identified with $\N^{\om}$.
This time the map $\tilde f : X(T) \to \tilde I$ given by $x \mapsto a^x$ is surjective as well as injective.  Thus, $\N^{\om}$ is order isomorphic with
$[0,1) \subset \R$.

       \vspace{1cm}

       \subsection{Subsets of $\Z^{\om}$}

       We regard $\Z^{\om}$ on the one hand as the branch space of the simple tree on $\Z, \om$ and on the other as the space of sequences on
       the alphabet $\Z$.

We first characterize the subsets $W$ of $\Z^{\om}$ which can occur as the set $L_{\om}$ of $\om$ level vertices of some additive $\Z$ tree. Notice
that $L_{\om} \not= \emptyset$  if and only if  $h(T) > \om$.

\begin{prop}\label{propauto0} A subset $W \subset \Z^{\om}$ is equal to the set  $L_{\om}$ of level $\om$ vertices in some additive  tree $T$ of  $\Z$ type
with $h(T) > \om$  if and only if  $W$
is saturated by the end equivalence relation.  That is, if $s_1, s_2 \in \Z^{\om}$ are end equivalent, then
$s_1 \in W$ if and only if $s_2 \in W$. \end{prop}

\begin{proof} If $T$ is an additive tree on $\Z$, then for $k < \om$, $L_k = \Z^k$. That is, every finite sequence on $\Z$ of length $k$
is a vertex of $T$ with order $k$. So if $p_1, p_2 \in \Z^{\om + 1}$ so that $o(p_1) = o(p_2) = \om$, then the associated sequences in
$\Z^{\om}$ are end equivalent if and only if there exist $q_1, q_2 \in T$ with $o(q_1), o(q_2) < \om$ and $  p \in \Z^{\om + 1}$  such that
$p_1 = q_1 + p, p_2 = q_2 + p$. Since $q_1, q_2 \in T$, additivity implies that $p_1 \in T \ \Leftrightarrow \ p \in T \ \Leftrightarrow \ p_2 \in T$.
Hence, $W = L_{\om}$ implies that the set   is saturated by the end equivalence relation.

For the converse, assume that $W$ is saturated by the end equivalence relation.  We apply the inductive construction of Theorem \ref{theo6.9}.
 Observe first
that the simple tree on $\Z, \om$ is an additive tree with height the limit ordinal $\om$. Because $W$ is saturated, it is translation invariant in
the sense of (\ref{eq6.31}). So the construction of the Theorem \ref{theo6.9} allows us to build an additive tree of $\Z$ type with height $\om^2$ with
$L_{\om} = \{ p + x : o(p) < \om, x \in W \}$. Because $W$ is saturated, this set is equal to $W$.

\end{proof} \vspace{.25cm}

\noindent{\bfseries Remark.} Recall the shift map $\t$ on $\Z^{\om}$. A set $W \subset \Z^{\om}$ is
 invariant in the sense of (\ref{eq6.31}) precisely when it is
$\t$ invariant, i.e. $s \in W \ \Longrightarrow \ \t(s) \in W$. The set is saturated by the end equivalence relation  if and only if
$s \in W \ \Longleftrightarrow \ \t(s) \in W$. \vspace{.5cm}

An additive tree on the transitive LOTS $\Z$ is homogeneous by Proposition \ref{prop5.16}.  It follows that  an end-equivalence saturated subset
 $W \subset \Z^{\om}$ is always transitive.
In fact we will prove below that it is doubly transitive.  We begin by describing
certain order preserving maps on $\Z^{\om}$.\vspace{.5cm}

\textbf{The Tree Automorphisms:}

Because the only order preserving bijections of $\Z$ are translations, and any tree automorphism  $f$ maps $S_x$ bijectively onto
$S_{f(x)}$, we can construct an arbitrary tree automorphism of $\Z^{\om}$ as follows.  Choose for each finite word $w$ (including the empty word)
 $m(w) \in \Z$.  Define the associated automorphism by
\begin{equation}\label{auto1}
f(x)_i = x_i + m(x_0 \dots x_{i-1}),
\end{equation}
 so that $f(x)_0 = x_0 + m(\emptyset)$.

 Notice that, using pointwise addition, $\Z^{\om}$ is a group. Translation by an element of the group is the special case when
 $|w_1| = |w_2| \ \Longrightarrow \ m(w_1) = m(w_2)$. Note that for the group $\Z^{\om}$
 the level end equivalence class of $\bar 0$ (equals the end equivalence class of $\bar 0$)
is a subgroup and the level end equivalence classes are exactly the cosets for this subgroup. Hence,  the set of
level end equivalence classes has a group structure
induced by pointwise addition. Also it follows that any two level end equivalence classes are isomorphic via a translation map.

\begin{prop}\label{propauto1} Let $p, q \in \Z^{\om}$.  There exists a tree automorphism $f$ of $\Z^{\om}$ such that $f(p) = q$ and such that
for every $x \not= p$ in $\Z^{\om}$,  $x$ and $f(x)$ are level end equivalent.
\end{prop}

\begin{proof} Given $p \in \Z^{\om}$ we choose $m(p_0 \dots p_{i-1}) = q_i - p_i$ for $i = 0,1, \dots $  and $m(w) = 0$ for all other finite words.
 Clearly, $f(p) = q$. If $x \not= p$, then let
 $k$ be the equality level for $x$ and $p$.  That is, $x_i = p_i$ for all $i < k$ and
 $x_k \not= p_k$. We then have
 \begin{equation}\label{auto2}
 f(x)_i \ = \ \begin{cases} \quad q_i \qquad \qquad \text{for} \ i < k, \\
 x_k - p_k + q_k \ \text{for} \ i = k, \\
 \quad x_i \qquad \qquad \text{for} \ i > k. \end{cases}
 \end{equation}
 Thus, $x$ and $f(x)$ are level end equivalent.

 \end{proof} \vspace{.5cm}

 \begin{cor}\label{corauto2} If $W \subset \Z^{\om}$ is saturated by level end equivalence (and so, a fortiori, if it is saturated by end equivalence), then
 the group of tree automorphisms on $\Z^{\om}$ acts transitively on $W$. \end{cor}

 \begin{proof} If $p, q \in W$ then by Proposition \ref{propauto1} there exists a tree automorphism $f$ such that $f(p) = q$
 and such that $x$ and $f(x)$ are
 level end equivalent for all $x \not= p$. It follows that for all $x$, $x \in W$ if and only if $f(x) \in W$.

 \end{proof} \vspace{.5cm}

Notice that if $x, y, z \in \Z^{\om}$ with $x_0 = 0, y_0 = 1$ and $z_0 = 2$, then there does not exist a tree automorphism $f$
such that $f(x) = x$ and $f(y) = z$.  So for double transitivity of a subset $W$ of $\Z^{\om}$ we cannot hope to use tree automorphisms.
We require more general order preserving bijections which do not preserve the tree structure.\vspace{.5cm}

\textbf{The Reproduction Isomorphisms:}\index{reproductive automorphisms}

For a finite word $w$ of length $k$,  $a_w(z) = wz$ defines the canonical order isomorphism from
$\Z^{\om}$ onto $\{x \in \Z^{\om}: x_i = w_i $ for $i = 0, \dots k-1 \}$. That is, it is the canonical isomorphism from
$T$ to $T_w$ where $T$ is the simple tree on $\Z,\om$.

Furthermore, $z$ and $a_w(z)$ are end equivalent in $\Z^{\om}$.
So if $W \subset \Z^{\om}$ is end equivalence saturated, then $a_w$ is an isomorphism from $W$ onto the subset we denote $wW$, the copy of $W$
with foot $w$. \vspace{.5cm}

\textbf{The Lift Maps:}\index{lift maps}

For $x \in \Z^{\om}$ with $x \not= \bar 0$ let $k^*(x) \in \om$ such that $x_{k^*} \not= 0$ and $x_i = 0$ for all $i < k^*$. That is,
$k^*$ is the equality level for the pair $x, \bar 0$.

The \emph{lift maps}\index{lift maps} are defined as follows:

 \begin{equation}\label{auto3}\begin{split}
 \ell_+(x) \ = \ \begin{cases} 0x \ \ \text{if} \ x_{k^*} < 0, \\ x \ \text{otherwise}.\end{cases} \\ \\
 \ell_-(x) \ = \ \begin{cases} 0x \ \ \text{if} \ x_{k^*} > 0, \\ x \ \text{otherwise}.\end{cases}.
 \end{split}\end{equation}

 Let $(\Z^{\om})_+ = \{ x \in \Z^{\om} : x_0 \geq 0 \}, (\Z^{\om})_- = \{ x \in \Z^{\om} : x_0 \leq 0 \},$ and for
 $W \subset \Z^{\om}$ we let $W_{\pm} = W \cap (\Z^{\om})_{\pm}$.

 \begin{prop}\label{propauto3} $\ell_+$ is an order  isomorphism from $\Z^{\om}$ onto $(\Z^{\om})_+$ and
 $\ell_-$ is an order  isomorphism from $\Z^{\om}$ onto $(\Z^{\om})_-$. Furthermore,
 $x$, $\ell_+(x)$ and $\ell_-(x)$ are end equivalent for all $x \in \Z^{\om}$.

 If $W \subset \Z^{\om}$ is end equivalence saturated, then $\ell_+$ and $\ell_-$ restrict to isomorphisms from
 $W$ to $W_+$ and $W_-$, respectively.\end{prop}

 \begin{proof} We will do the proof for $\ell_+$.

 It is clear that $\ell_+$ is a bijection from $\Z^{\om}$ to $(\Z^{\om})_+$ and that
  $x$ and $\ell_+(x)$ are end equivalent. We must prove that $\ell_+$ is order preserving.

 We say that $x$ is \emph{fixed} when $\ell_+(x) = x$, and that $x$ is \emph{lifted} when $\ell_+(x) = 0x$. So $x$ is fixed when $x = \bar 0$ or
 $x_{k^*} > 0$ and is lifted when $x_{k^*} < 0$.

 Assume $x < y$. We say that $x$ and $y$ are on the same side when they are both fixed or both lifted. In that case it is clear that
 $\ell_+(x) < \ell_+(y)$.

 Let $k$ be the equality level for $x$ and $y$.  That is,  $x_i = y_i$ for $i < k$ and $x_k < y_k$. If for some $i < k$,
 $x_i = y_i \not= 0$, then $k^*(x) = k^*(y)$ and $x_{k^*} = y_{k^*}$. So in that case, $x$ and $y$ are on the same side.

 Assume now that $x_i = y_i = 0$ for $i < k$. If $x_k < y_k <0$, then $k^*(x) = k^*(y) = k$ and $x$ and $y$ are on the same side.

 If $y_k > 0$, then $y$ is fixed. If $x$ and $y$ are not on the same side then $x$ is lifted.
 Hence, $\ell_+(x)_k = 0 < y_k = \ell_+(y)_k$ and so $\ell_+(x) < \ell_+(y)$.

 If $x_k < y_k = 0$, then $x$ is lifted. If $x$ and $y$ are not on the same side, then $y$ is fixed. So either $y = \bar 0$ or
 $k^*(y) > k$ and $y_{k^*} > 0$.
 In either case we have $\ell_+(x)_{k+1} = x_k < 0 \leq y_{k+1} = \ell(y)_{k+1}$.  So in this case as well, $\ell_+(x) < \ell_+(y)$.

 This completes the proof that $\ell_+$ is an order map.

 Since $x, \ell_+(x)$ and $\ell_-(x)$ are all end equivalent, it follows that $\ell_{\pm}$ restricts to an isomorphism of $W$ onto $W_{\pm}$.

 \end{proof} \vspace{.5cm}

 \begin{theo} \label{theoauto4} If $W$ is an end equivalence saturated  subset of $\Z^{\om}$, then $W$ is doubly transitive and so is an IHLOTS
 with completion $\R$. \end{theo}

   \begin{proof} Using the reproduction isomorphisms we express various subsets of $\Z^{\om}$ as \emph{patterns} of copies of $W$.
  The decomposition $W = \bigcup_{i \in \Z} iW$ expresses $W$ as a $\Z$ pattern of copies of $W$, i.e. shows that it is
  isomorphic to the lexicographic product $\Z \times W$. Similarly,  $W_+ = \bigcup_{i \in \N} iW$ and $W_- = \bigcup_{i \in -\N} iW$
  express $W_+$ as an $\N$ pattern and $W_-$ as an $\N^*$ pattern.

  The isomorphism $\ell_+ : W \to W_+$ shows that $W$ is isomorphic to an $\N$ pattern of copies of $W$. Similarly, $\ell_-$
shows that $W$ is isomorphic to an $\N^*$ pattern of copies of $W$.

By replacing $W$ by a translate given by $x \mapsto x - y + \overline{(-1)}$ with $y \in W$, we may assume that $\overline{(-1)} \in W$. Since
$W$ is transitive, double transitivity follows from transitivity of the interval $\{x \in W : \overline{(-1)} < x \}$, see Proposition \ref{prop3.2}
(d)(v).

If $(-1)^k$ is the word $w$ of length $k$ with $w_i = -1$ for $i = 0,\dots, k-1$ and so $(-1)^0$ is the empty word, then
we have the decomposition $\{x \in W : \overline{(-1)}< x \} = \bigcup_{k \in \N} (-1)^k W_+$ with $(-1)^{k+1} W_+$ preceding $(-1)^k W_+$.
Each $(-1)^k W_+$ is isomorphic to $W$ and so  $\{x \in W : \overline{(-1)} < x \}$ is isomorphic to an $\N^*$ pattern of copies of $W$ which is,
in turn, isomorphic to $W$. Hence, $\{x \in W : \overline{(-1)} < x \}$ is transitive and so $W$ is doubly transitive.

For an alternative direct proof,   we can instead begin with $x < y \in W$. If the equality level is $k$ then there is a word $w$ of length $k$ and $x', z' \in W$
such that $x = wx', y = wy'$ and $x'_0 < y'_0$.  The interval $(x,y)$ in $W$ is contained in $wW$ and is isomorphic via $a_w$ to the
interval $(x',y')$ in $W$. So we may assume that $x_0 < y_0$.

For $k = 1,2, \dots$ let $W_{k+} = \{ z \in W : x_i = z_i $ for $ i = 0, \dots, k-1$ and $ x_k < z_k \}$ so that $W_{k+} \cong W_+ \cong W$ and
let $W_{k-} = \{ z \in W : y_i = z_i $ for $ i = 0, \dots, k-1$ and $ y_k > z_k \}$ so that $W_{k-} \cong W_- \cong W$. For $x_0 < n < y_0$
let $W_{n0} = nW$.

This expresses the interval $(x,y)$ in $W$ as a $\Z$ pattern of copies of $W$ and so it is isomorphic to $W$.  This directly shows that $W$ is
weakly homogeneous. It follows from Proposition \ref{prop3.2a} that $W$ is doubly transitive.

By Proposition \ref{prop2.5}
$W$ is of countable type because $\Z^{\om}$ is. It then
follows from Proposition \ref{prop3.4}  that $W$ is a HLOTS. Since it is a dense subset of the IHLOTS  $\Z^{\om}$ it is an IHLOTS with completion $\R$.

   \end{proof} \vspace{.5cm}

      \begin{prop}\label{propauto6} If $X$ is a proper subset of $\R$ which is invariant under the group of affine maps $x \mapsto q + 2^n \cdot x$ with
   $q$ varying over $\Q$ and $n$ varying over $\Z$, then $X$ is an IHLOTS. If $\Q \cap X = \emptyset$, then $X$
   is isomorphic to an end equivalence saturated subset $W$ of $\Z^{\om}$.
   \end{prop}

   \begin{proof}  We apply the construction of Example \ref{exn04}.

   Observe that $t_n$ is the restriction of an element of the group for all $n \in \Z$. The map $x \mapsto t_n(x-n)$ for $n \leq x \leq n+1$ defines an
   isomorphism from $X$ onto $X \cap (-1,1)$. Invariance implies that $X$ is dense in $\R$ and so has completion $\R$.

   First assume that $\Q \cap X = \emptyset$. Let $W = \{ s \in \Z^{\om} : I_s \subset X \}$.
   That is, the point in the singleton set $I_s$ lies in $X$. If $w$ is a finite word, then $t_w(I_s) = I_{ws}$. So the invariance assumption implies
   that $s \in W \ \Leftrightarrow \ ws \in W$.  That is, $W$ is end equivalence saturated. From the definition of the maps $t_n$ it is clear that
   the set of endpoints $C$ is contained in $\Q$. Because $X$ is disjoint from $\Q$ it follows that for every $x \in X \cap (0,1)$ there is a unique
   $s \in W$ such that $I_s = \{ x \}$. Thus, $W \cong X \cap (0,1) \cong X$. Hence, $X$ is an IHLOTS by Theorem \ref{theoauto4}.

   If $\Q$ meets $X$ then by invariance, $\Q \subset X$ and so the above result applies to the complement $\R \setminus X$. By
   Proposition \ref{prop3.4} $X$ is an IHLOTS
   because it complement is an IHLOTS.

   \end{proof}

   \vspace{1cm}

  \subsection{The Tree Characterizations}

 \begin{theo} \label{theoauto5}  For a LOTS $X$, the following are equivalent.
 \begin{itemize}
\item[(i)] $X$ is a CHLOTS.
\item[(ii)] There exists an additive, $\Omega$ bounded tree $T$ with $S_0 = \Z$ and $h(T) \ge \om$, and
such that the completion of the branch space $X(T)$ is isomorphic to $X$.
 \item[(iii)] There exists an additive, $\Omega$ bounded tree $T$ with $S_0 = Y$, an IHLOTS with completion $\R$,
such that the completion of the branch space $X(T)$ is isomorphic to
$X$.
 \item[(iv)] There exists a reproductive, $\Omega$ bounded tree $T$ with $S_0$ an HLOTS, and such that
 the completion of the branch space $X(T)$ is isomorphic to $X$.
\end{itemize}\end{theo}

 \begin{proof} (i) $\Rightarrow$ (ii): This is Theorem \ref{theon03} applied to our construction above. Note that we require $h(T) \ge \om$ because
 with $h(T) = 1$ the simple tree on $\Z$ has $\Z$ as branch space.

 (ii) $\Rightarrow$ (iii): If $h(T) = \om$, then $T$ is the simple tree on $\Z, \om$ and so $X \cong \R$. As in Example \ref{exn04} we can
 get $\R$ as the completion of the simple tree on $\Q, \om$, proving (iii) in this case. We can also use the simple tree on $\Q, 1$.

 When $h(T) > \om$ and so $h(T) \ge \om^2$, we use the Omega Thinning Construction. By Proposition \ref{prop6.22} the result is an
 additive tree of type $L_{\om}$.
 By Proposition \ref{propauto0} $W = L_{\om}$ is end equivalence saturated and so
 by Theorem \ref{theoauto4} $W$ is an IHLOTS with completion $\R$. Finally, Corollary \ref{cor6.19} implies that the $X(\om T)$ has the
 same completion as $X(T)$ and the latter is assumed to be isomorphic to $X$.

  (iii) $\Rightarrow$ (iv): An additive tree is reproductive.

   (iv) $\Rightarrow$ (i): Theorem \ref{theo5.20}.
   \end{proof} \vspace{.5cm}

   \noindent{\bfseries Remark.} For each of the trees described in (ii)-(iv), if $\a$ is a countable, tail-like ordinal, then the truncation $T^{\a}$
   is a tree of the same sort, e.g. additive or reproductive, and, in addition, with $h(T^{\a}) = \a$.
   It follows from the theorem together with  Proposition
   \ref{propnew04a} that each $\widehat{X(T^{\a})}$ is an $\R$-bounded CHLOTS
    (In (iv) we must assume that $S_0$ is $\R$-bounded). \vspace{.5cm}

   We have seen that if $T$ is an $\Omega$ bounded, additive tree of $\Z$ type, then the branch space $X(T)$ is order dense and transitive and
   it is a dense subset of its completion which is the same as the completion of its Omega thinned tree $\om T$. Furthermore,
   $X(\om T)$ is an IHLOTS and so the completion is a CHLOTS. We do not know
   whether $X(T)$ itself, though transitive and order dense,  is necessarily doubly transitive and so is a HLOTS. \vspace{.5cm}

   From the theorem we see that via the inductive construction of Theorem \ref{theo6.9} we can obtain every CHLOTS. We review the construction
   using  alphabet language.

   For $\g$ a positive countable ordinal, let $T(\g)$ be an additive tree of height $\om^{\g}$ with $\tilde L(\g) = \tilde L_{\om^{\g}}$
   the set of branches of height $\om^{\g}$. By additivity, we have, for $T = T(\g), \a = \om^{\g}, \tilde L = \tilde L(\g)$:
   \begin{equation}\label{eqfinalnew01}
   T \ = \ \{ x|\b : x \in \tilde L, \b < \a \}. \hspace{2cm}
   \end{equation}

   We select $A(\g) \subset \tilde L(\g)$ which is \emph{$\pm$ invariant} in the sense that with $A = A(\g), \a = \om^{\g}, T = T(\g)$
    \begin{equation}\label{eqfinalnew02}\begin{split}
    x \in A\ \  \text{and} \ \  \b < \a \qquad \Longrightarrow \qquad \t_{\b}(x) \in A,\\
  x \in A \ \  \text{and} \ \  p \in T  \qquad \Longrightarrow \qquad p + x \in A.
   \end{split}\end{equation}
   If we say that $x_1, x_2 \in \tilde L(\g)$ are \emph{end-equivalent} when there exist $\b_1, \b_2 < \om^{\g}$ such that
 $\t_{\b_1}(x_1) = \t_{\b_2}(x_1)$, then $A(\g)$ is $\pm$ invariant exactly when it is saturated by the end-equivalence relation.

   On the one hand, we let $A(\g)$ be the set of vertices of $T(\g + 1)$ of order $\om^{\g}$. On the other hand, we
   regard $A(\g)$ as an alphabet and let $\tilde L(\g + 1)$ be all of the elements of $X^{\om^{\g + 1}} = X^{\om^{\g} \cdot \om}$
   which correspond to infinite sequences on the alphabet $A(\g)$. Then $T = T(\g + 1)$ is defined by (\ref{eqfinalnew01}) with
   $\a = \om^{\g + 1}, \tilde L = \tilde L(\g + 1)$. It follows from $\pm$ invariance of $A(\g)$ that $T(\g + 1)$ is
   an additive tree of height $\om^{\g + 1}$ and with $\tilde L(\g + 1)$ the set of branches of height $\om^{\g + 1}$.
   Furthermore, $T(\g) = T(\g + 1)^{\a}$ with $\a = \om^{\g}$.

   A subset $A(\g + 1) \subset \tilde L(\g + 1)$ is $\pm$ invariant exactly when it corresponds to a set
   of $A(\g)$ sequences which is saturated by the
   end-equivalence relation on sequences, i.e. the two versions of end-equivalence agree.

   Let $\g$ be a limit ordinal. Assume that the  additive trees $T(\d)$ have been defined for all $\d < \g$ with
   $T(\d_1) = T(\d)^{\a}$ with $\a = \om^{\d_1}$ whenever $\d_1 < \d < \g$. Define
    \begin{equation}\label{eqfinalnew03}
   T(\g) \ = \ \bigcup_{\d < \g} \ T(\d),
   \end{equation}
   with $\tilde L(\g)$
   the branches of height $\om^{\g}$.

   The process stops at a countable ordinal $\g$ if we choose $A(\g) = \emptyset$.
   Otherwise, the process continues to $\g = \Omega = \om^{\Omega}$ and $A(\Omega) = \emptyset$.

   Having defined $T(\g), A(\g)$ for all positive $\g$ with $\om^{\g} \le h(T)$, the branch space is given by
    \begin{equation}\label{eqfinalnew04}
    X(T) \ = \ \bigcup \{ \tilde L(\g) \setminus A(\g) : \om^{\g} \le h(T) \}.
    \end{equation}

   The resulting tree is $\Omega$ bounded when $\tilde L(\Omega) = \emptyset$. As we saw when we considered the height function in Section 6,
   we can assure an $\Omega$ bounded result, by beginning with $R$ an arbitrary $\Omega$ bounded subtree of the
   simple tree of height $\Omega$ and only continuing the construction
   as long as $T(\g)$ remains a subset of $R$.
         \vspace{1cm}

   \subsection{Trees of Convex Sets}

   A convex set $J$ in an order dense LOTS $X$ is \emph{proper} when it has more than one point and so is infinite. In particular, a nonempty open convex set   is proper.

   Two convex sets $J_1, J_2$ in $X$
   \emph{overlap}\index{overlapping convex sets} when $J_1 \cap J_2$ is proper. If two convex sets do not overlap, then either they are disjoint
   or their intersection is a singleton consisting of a common endpoint.

   We write $J_1 \prec J_2$ when there exist $a, b \in J_1 ^{\circ}$
   such that $a < c < b$ for all $c \in \overline{J_2}$ and so
   \begin{equation}\label{eqcon01}
   \overline{J_2} \ \subset \ (a,b) \ \subset \ [a,b] \ \subset \ J_1^{\circ}.
   \end{equation}
   Here $\overline{J}$ and $J^{\circ}$ are the closure and interior, respectively, in the LOTS $X$.

   For $A \subset X$ we let $[A]$ denote the convex closure of $A$, i.e. the smallest closed, convex set which contains $A$.

   \begin{lem}\label{lemcon01}  Let $J, J_1, J_2$ be proper convex subsets of $X$ an order dense LOTS.
   \begin{itemize}
   \item[(a)] The set $\overline{J} \setminus J^{\circ}$ contains at most two points. The interior $J^{\circ}$ is dense in $J$ and so is itself
   a proper open convex set.

   \item[(b)] $J_1$ overlaps $J_2$ if and only if $J_1 \cap J_2^{\circ} \not= \emptyset$.

   \item[(c)] $J_1 \prec J_2$ if and only if there exist $a_2 < a_1 < b_1 < b_2 \in J_1$ such that $a_1 < c < b_1$ for all $c \in J_2$.

   \item[(d)] If $X$ is connected, then $J_1 \prec J_2$ if and only if $\overline{J_2} \subset J_1^{\circ}$.
   \end{itemize}
   \end{lem}

 \begin{proof} (a): Let $a \in \overline{J}$.

 \textbf{Case 1:} ($(-\infty,a) \cap J \not= \emptyset$ and $(a,\infty) \cap J \not= \emptyset$.) There exist $b_1, c_1 \in J$ with $b_1 < a < c_1$.
  Choose $b \in (b_1,a), c \in (a,c_1)$ and we have $a, b, c \in J^{\circ}$ with $ a \in (b,c)$.\vspace{.25cm}

 \textbf{Case 2:} ($(-\infty,a) \cap J \not= \emptyset$ and $(a,\infty) \cap J = \emptyset$.) There exists $b_1\in J$ with $b_1 < a$. Choose $b \in (b_1,a)$
 and we have $[b,a) \subset J^{\circ}$ and $a =  sup J$.\vspace{.25cm}

 \textbf{Case 3:} ($(-\infty,a) \cap J = \emptyset$ and $(a,\infty) \cap J \not= \emptyset$.) There exists $c$ such that
 $(a,c] \subset J^{\circ}$ and $a =  inf J$.\vspace{.25cm}

 Since $J$ is proper, one of these cases applies for every $a \in \overline{J}$. So the only points of $\overline{J} \setminus J^{\circ}$ are the
 supremum and infimum of $J$ if either of these exists. So $J^{\circ}$ is infinite and dense in $J$.

 (b): If $J_1 \cap J_2$ is proper, then by (a)  $(J_1 \cap J_2)^{\circ}$ is infinite. If $J_1 \cap J_2^{\circ} \not= \emptyset$, then because
 $J_1^{\circ}$ is dense in $J_1$, $J_1^{\circ} \cap J_2^{\circ}$ is a nonempty open convex set and so it is infinite.

 (c):  If $a_2 < a_1 < b_1 < b_2 \in J_1$ with $a_1 < c < b_1$ for all $c \in J_2$, then $\overline{J_2} \subset [a_1,b_1]$ and
 $[a_1,b_1] \subset (a_2,b_2) \subset J_1^{\circ}$. If $a \in (a_2,a_1)$ and $b \in (b_1,b_2)$, then $a$ and $b$ satisfy (\ref{eqcon01}).

 Conversely, if $a$ and $b$ satisfy (\ref{eqcon01}), then $[b, \infty) \cap J_1^{\circ}$ is nonempty and so by (b) it contains a nonempty open interval.
 Similarly for $(-\infty,a] \cap J_1^{\circ}$. Since these sets are infinite, we can choose the required $a_1, a_2, b_1, b_2$.

 (d): If $X$ is connected, then with $c_1 =  inf J_2, c_2 =  sup J_2$, $\overline{J_2} = [c_1,c_2]$. Similarly, $J_1^{\circ} = (b_1,b_2)$.
 If $\overline{J_2} \subset J_1^{\circ}$, then we can choose $a \in(b_1,c_1), b \in (c_2,b_2)$, and then  $a$ and $b$ satisfy (\ref{eqcon01}).

 \end{proof} \vspace{.5cm}

 Let $i : X_1 \to X_2$ be an order injection with $X_1$ order dense and $X_2$ connected. For $J$ a bounded, proper convex set in $X_1$ let $[i](J) = [i(J)]$,
 the convex closure of the image of $J$. Thus, $[i](J)$ is the closed, bounded, proper interval in $X_2$ with endpoints the infimum and supremum in $X_2$
  of $i(J)$. In particular, if $a < b \in X_1$, then  $[i]([a,b]) = [i(a),i(b)]$.

  \begin{lem}\label{lemcon02}  Let $J_1, J_2$ be bounded, proper convex subsets of $X_1$ an order dense LOTS and $i : X_1 \to X_2$ be an order injection
  with $X_2$ connected.
   \begin{itemize}
   \item[(a)] $J_1$ and $J_2$ overlap in $X_1$ if and only if $[i](J_1)$ and $[i](J_2)$ overlap in $X_2$.

   \item[(b)] $J_1 \prec J_2$ in $X_1$ if and only if $[i](J_1) \prec [i](J_2)$  in $X_2$.
 \end{itemize}
 \end{lem}

 \begin{proof} Let $[i](J_1) = [a_1,b_1]$ and $[i](J_2)= [a_2,b_2]$.

 (a): $[i](J_1 \cap J_2) \subset [i](J_1) \cap [i](J_2) = [max (a_1,a_2), min (b_1,b_2)]$. So if $J_1$ and $J_2$ overlap, then
 $[i](J_1) \cap [i](J_2)$ is proper. On the other hand, if $[i](J_1)$ and $[i](J_2)$ overlap, then $max (a_1,a_2) < min (b_1,b_2)$. Choose
 $c \in X_2$ between them. For $\ep = 1,2$, there exist $i(u_{\ep}) \in (a_{\ep}, c) \cap i(J_{\ep}), i(v_{\ep}) \in (c,b_{\ep}) \cap i(J_{\ep})$ by definition of  the sup and inf. Let $u = max (u_1,u_2), v = min (v_1,v_2)$. So $u, v \in [u_{\ep},v_{\ep}]$ for $\ep = 1,2$, Thus, $[u,v] \subset J_1 \cap J_2$.

 (b): $[i](J_1) \prec [i](J_2)$ if and only if $a_1 < a_2 < b_2 < b_1$. We can choose $u_1, u_2, v_1, v_2 \in J_1$
 with $a_1 < i(u_2) < i(u_1) < a_2$ and $ b_2 < i(v_1) < i(v_2) < b_1$. From Lemma \ref{lemcon01} (c) it follows that $J_1 \prec J_2$.

 Conversely, if $u_2 < u_1 < v_1 < v_2 \in J_1$ such that $u_1 < c < v_1$ for all $c \in J_2$, then
 $[i](J_2) \subset [i(u_1),i(v_1)] \subset (i(u_2),v(v_2)) \subset [i](J_1)^{\circ}$ and so $[i](J_1) \prec [i](J_2)$.

 \end{proof} \vspace{.5cm}

 \noindent{\bfseries Remark.} If $i$ is the inclusion of $X_1$ into its completion $X_2 = \widehat{X_1}$, then $[i](J)$ is the closure $\overline{J}$ of
 $J$ in $X_2$. \vspace{.5cm}

 \begin{df}\label{dfcon03} A collection $\T$ of proper convex subsets of an order dense LOTS $X$ is a \emph{tree of convex sets}
 \index{tree of convex sets} in $X$ when it satisfies:
  \begin{itemize}
 \item[(i)] $X \in \T$ with $J \in \T$ bounded when $J \not= X$.
 \item[(ii)] If $J_1$ and $J_2$ are distinct elements of $\T$, then $J_1$ overlaps $J_2$ if and only if either $J_1 \prec J_2$ or $J_2 \prec J_1$.
 \item[(iii)] With respect to the ordering $\prec$, $\T$ is a (not necessarily semi-normal) tree.
 \end{itemize}
 \end{df} \vspace{.5cm}

 The tree $\T$ has $X$ as its root and every other vertex $J$ of $\T$ is bounded by (\ref{eqcon01}) since $X \prec J$.
 We do not assume that (ii) of Definition \ref{df5.1} holds, i.e. $J \in \T$ may have a single successor and we do not assume
 that condition (iii) of Definition \ref{df5.1} holds.

 For example, for the tree $T$ constructed in the proof of Theorem \ref{theon.01}, the collection $\{ I_{s} : s \in T \}$ is a tree of convex sets
 in the connected LOTS $[m,M]$.

 \begin{prop}\label{propcon04} If $\T_1$ is a tree of convex sets in an order dense LOTS $X_1$ and $i : X_1 \to X_2$ is an order injection with
 $X_2$ connected, then $\T_2$ is a tree of closed intervals in  $X_2$ and $f : \T_1 \to \T_2$ is a tree isomorphism with
 \begin{equation}\label{eqcon02}
f(J) \ = \  \begin{cases} X_2 \ \text{if} \ J = X_1, \\ [i](J) \ \text{if} \ o(J) > 0, \end{cases}
\end{equation}
and $\T_2$ the image of $f$. \end{prop}

 \begin{proof} This is clear from Lemma \ref{lemcon02}.

 \end{proof} \vspace{.5cm}

 \begin{theo}\label{theocon05} If $T$ is a normal tree of unbounded type with $h(T)$ a limit ordinal, then $\T = \{ j_p(X(T_p)) : p \in T \}$ is a tree
 of clopen convex sets in the branch space $X(T)$. The map $j$ defined by $j(p) =  j_p(X(T_p))$ is a tree isomorphism from $T$ to $\T$. In addition,
 $\hat{\T} = \{ \overline{j_p((T_p))} : p \in T \}$ is a tree of closed intervals in the completion $\widehat{X(T)}$, with the map $f$ given by
 $f(j_p(X(T_p))) = \overline{j_p(X(T_p))}$ a tree isomorphism from $\T$ to $\hat{\T}$. \end{theo}

 \begin{proof} By Proposition \ref{prop5.3a} the branch space $X(T)$ is order complete. By Proposition \ref{prop5.3} (c) each $j_p(X(T_p))$ is a
 clopen convex set in $X(T)$ with completion the open interval $\hat{j}_{p}(\widehat{X(T_p)})$ in $\widehat{X(T)}$. Since $T$ is of unbounded type,
 each $j_p(X(T_p))$ is bounded for $p \not= 0$.
 By the Remark following
 Lemma \ref{lemcon02}, the common closure of $j_p(X(T_p))$ and $\hat{j}_{p}(\widehat{X(T_p)})$ in $\widehat{X(T)}$ is $[i](j_p(X(T_p)))$ where
 $i$ is the inclusion of $X(T)$ into its completion.

 We check (i)-(iii) for $\T$.

 For the root $0$, $T_0 = T$ and $j_0(X(T_0)) = X(T)$, verifying (i).

 If $p_1$ and $ p_2$ are distinct vertices of $T$, then $p_1 \prec p_2$ implies $T_{p_2} \subset T_{p_1}$ and there exist $q_1 < q_2 < q_3 \in S_{p_1}$
 with either $q_2 = p_2$ or $q_2 \prec p_2$. Let $a, b$ be branches through $q_1$ and $q_3$. Then (\ref{eqcon01}) is satisfied showing that
 $j_{p_1}(X(T_{p_1})) \prec  j_{p_2}(X(T_{p_2})) $ since the convex sets are clopen.

On the other hand, if $j_{p_1}(X(T_{p_1})) \cap j_{p_2}(X(T_{p_2}))$ is nonempty and so there exists
a branch which contains both $p_1$ and $p_2$, then either $p_1 \prec p_2$ or $p_2 \prec p_1$.

This implies (ii) and shows that $j : T \to \T$ is an order isomorphism. Hence, $\T$ is a tree, verifying (iii).

Finally, the completion results follow from Proposition \ref{propcon04}.

\end{proof} \vspace{.5cm}

Now we apply these results to show that the branch space of an Aronszajn tree is not $\R$-bounded.

\begin{theo}\label{theocon06} Let $\a$ be a countable ordinal and let $\T$ be a tree of convex sets in the connected LOTS $\R_{\a}$.
If every level of $\T$ is countable, then $\T$ is countable. \end{theo}

\begin{proof} By applying Proposition \ref{propcon04} to $i$ equal the identity map on $\R_{\a}$ we can assume that $\T$ is a
tree of closed intervals in $\R_{\a}$. So for $J \in \T$ with $o(J) > 0$, we have $J = [a(J),b(J)]$ with $a(J) < b(J) \in \R_{\a}$.

Since every level is countable, it suffices to prove that the height $h(\T)$ is countable.

If $x$ is a branch of $\T$, then $J \mapsto  a(J)$ is an embedding of the ordinal $h(x)$ into $\R_{\a}$. Since $\a$ is countable,
$\R_{\a}$ is first countable and so $\Omega$ cannot be injected into $\R_{\a}$.  It follows that $h(x)$ is countable.

For $J \in \T$, let $\ep(J)$ be the equality level for $a(J)$ and $b(J)$ so that $a(J)_i = b(J)_i$ for $i < \ep$ and
$a(J)_{\ep} < b(J)_{\ep}$ with $\ep = \ep(J)$. Hence, $\ep(J) < \a$.

We let $r(J) \in \R_{\ep}$ be the common restriction of $a(J)$ and $b(J)$ to $\ep$
and we call $r(J)$ the \emph{stem} of  $J$. We write  $span(J) = [a(J)_{\ep},b(J)_{\ep}]$ so that $span(J)^{\circ} = (a(J)_{\ep},b(J)_{\ep})$.
These are proper intervals in $\R$. Notice that if $q \in span(J)^{\circ}$, then any $c \in \R_{\a}$ with $c_i = a(J)_i = b(J)_i$ for $i < \ep$
and $c_{\ep} = q$ satisfies $a(J) < c < b(J)$ and so $c \in J^{\circ}$.  For example, $q+ \in J^{\circ}$.

Now we prove by induction that for every $\b < \a$, there exists $\xi_{\b} < \Omega$ such that $o(J) > \xi_{\b}$ implies $\ep(J) > \b$.

This is trivial for $\b = 0$.

Define $\rho_{\b}$ to be the smallest ordinal greater than $\xi_{\g}$ for all $\g < \b$. Thus,
\begin{equation}\label{eqcon05}
\rho_{\b} \ = \  sup \ \{\xi_{\g} : \g < \b \} + 1
\end{equation}

Now suppose $r \in \R_{\b}$ and there exists $J \in \T$ such that $stem(J) = r$. Let $J$ be an element of $\T$ with minimum order
such that  $stem(J) = r$. If $J_1 \prec J$ in $\T$, then $a(J_1) < a(J) < b(J) < b(J_1)$ and so $\ep(J_1) \le \ep(J)$. Furthermore,
$\ep(J_1) = \ep(J)$ would imply $stem(J_1) = stem(J)$ violating the minimality condition on $J$.  Hence, $\g = \ep(J_1) < \ep(J) = \b$.
From the inductive hypothesis it follows that  $o(J_1) \le \xi_{\g}$. Therefore,
\begin{equation}\label{eqcon06}
o(J)\ \le \  sup \ \{ o(J_1) : J_1 \prec J \} + 1 \ \le \  sup \ \{ \xi_{\g} : \g < \b \} + 1 \ = \ \rho_{\b}.
\end{equation}

Because each level of $\T$ is countable, there are only countably many $J \in \T$ with $o(J) \le \rho_{\b}$ and so the set
 $A_{\b} = \{ r \in \R_{\b} : r = r(J)$ for some $J \in \T \}$ is countable. \vspace{.25cm}

\textbf{Claim:} For each $r \in A_{\b}$ there are only countably many $J \in \T$ such that $r(J) = r$.

Suppose instead that for some $r \in A_{\b}$ there are uncountably many such $J$.
It follows that for some $q \in \Q$ there is an uncountable $\X \subset \T$ such
that for $J \in \X$  $r(J) = r$ and $q \in span(J)^{\circ}$. So $q+ \in J^{\circ}$ for all $J \in \X$. Thus, any two members of $\X$ overlap and so
by condition (ii) of Definition \ref{dfcon03}, any two are $\prec$ comparable. It follows that $\X$ is contained in an uncountable branch of $\T$
which is, as we have seen above, impossible.  This proves the Claim.\vspace{.25cm}

Now let $\xi_{\b} \ge \rho_{\b}$ and $\xi_{\b} \ge o(J)$ for all $J$ such that $r(J) = r$ for some $r \in A_{\b}$.

If $\ep(J) = \g < \b$, then $o(J) \le \xi_{\g} < \rho_{\b} \le \xi_{\b}$.

If $\ep(J) = \b$, then $r(J) \in A_{\b}$ and so $o(J) \le \xi_{\b}$.

This completes the inductive step.

The height of $\T$ is bounded by $sup \{ \xi_{\b} : \b < \a \}$ and so is countable.

 \end{proof} \vspace{.5cm}

 \begin{cor}\label{corcon07} If $T$ is an Aronszajn tree of unbounded type, then $X(T)$ is not $\R$-bounded. \end{cor}

 \begin{proof} Assume that $T$ is a normal tree of unbounded type with height a limit ordinal and that every level of $T$ is countable.
 Suppose that $i : X(T) \to \R_{\d}$ is an order injection with $\d$ countable.

 By Theorem \ref{theocon05}, there is a tree $\T_1$ of convex sets in $X(T)$ which is isomorphic with $T$ itself.
By Proposition \ref{propcon04} the injection $i$ induces an isomorphism of $\T_1$ onto a tree $\T_2$ of closed intervals in
$\R_{\d}$. By Theorem \ref{theocon06} the tree $\T_2$ is countable since every level is countable.  This in turn implies that $T$
is countable and so with countable height.

However, an Aronszajn tree is uncountable with height $\Omega$.
 \end{proof} \vspace{1cm}

\section{\textbf{HLOTS in $\R$ }}
\vspace{.5cm}

\subsection{Comparisons Along the Tower}

We begin with a pair of useful isomorphisms.

Recall that if $X$
is an IHLOTS with completion $F = \hat{X}$, then by Proposition
\ref{prop3.4}(d), the complement $F \setminus X$ is an IHLOTS with the same
completion.  Also if $\alpha$ is an infinite, tail-like ordinal
then by Proposition \ref{prop6.3}, $X^{\alpha}$ is order isomorphic to a
dense subset of $X_{\alpha}$. If, in addition, $\a$ is countable, then
$X^{\alpha}$ and $X_{\alpha}$ are HLOTS by Corollary \ref{cor6.2} and Theorem \ref{theo4.2}.

\begin{lem}\label{lem8.1}  Let $X$ be a IHLOTS with completion the CHLOTS $F$ and let $\alpha $
be a countably infinite, tail-like ordinal. The complement $\widehat{X^{\alpha}} \setminus X^{\alpha}$
is an IHLOTS and the  IHLOTS $(\widehat{X^{\alpha}} \setminus X^{\alpha})^{\alpha}$ has
completion isomorphic to $F_{\alpha}$.
\end{lem}

\begin{proof} Let $Y = F \setminus X$ be the
complementary IHLOTS.  Assume that the interval $J = [-1,+1]
\subset F$ has been chosen so that $-1,+1 \in Y$. Let $ J(X) =  X
\cap J = X \cap J^{\circ} $ and $ J(Y) = Y \cap J $.  Since $X$ is
a HLOTS, $X \cong J(X)$ and so $X^{\alpha} \cong J(X)^{\alpha}$.
By applying Proposition \ref{prop5.6} to the simple tree, it is easy to see that the complement of
$J(X)^{\alpha}$ in its completion can be identified with
\begin{equation}\label{eq8.1}
\begin{split}
Z \quad = \quad \{ z \in J^{\beta + 1}:\mbox{for some} \
\beta < \alpha \hspace{1cm} \\  \mbox{with} \ z(j) \in J(X) \  \mbox{for all}\
j < \beta \ \mbox{and} \ z(\beta) \in J(Y) \}
\end{split}
\end{equation}
except that $Z$ includes the endpoints $m,M \in J^{1}$ with $m(0)
= -1$ and $M(0) = +1$. With $Z^{\circ} =  Z \setminus
\{m,M\}$ we apply Proposition \ref{prop6.3} to the HLOTS $Z^{\circ}$, to see
that the completion of $(Z^{\circ})^{\alpha}$ is the same as the
completion of $(Z^{\circ})_{\alpha}$.  Since the distinguished
closed bounded interval in $Z^{\circ}$ is isomorphic to $Z$ we can
identify $(Z^{\circ})_{\alpha}$ with
\begin{equation}\label{eq8.2}
\tilde{Z} \qquad = \qquad \{ x \in Z^{\alpha} : x(0) \in
Z^{\circ} \}.
\end{equation}

We will show that $\tilde{Z}$ is order isomorphic to
\begin{equation}\label{eq8.3}
\begin{split}
D \quad = \quad \{ x \in J^{\alpha} :x(0) \not= \pm 1 \
\mbox{and} \   \gamma(x) \cong \alpha  \} \\ \mbox{where} \qquad
\gamma(x) \quad = \quad  \{i \in \alpha : x(i) \in J(Y) \}.
\end{split}
\end{equation}
Since $\alpha$ is tail-like, $D$ includes every $x \in J^{\alpha}$
such that $x(0) \not= \pm 1$ and  $x(i) \in J(Y)$ for sufficiently
large $i \in \alpha$.  So $D$ is dense in $\{ x \in J^{\alpha} :
x(0) \not= \pm 1 \} \cong F_{\alpha}$.

Given $x \in D$ let $\tau_{x} : \alpha \rightarrow \gamma(x)$ be
the (unique) order isomorphism which exists by definition of $D$.

Now we use a variation of the construction of (\ref{eq6.47}).

If $\tau : \a \to \a$ is an order injection, then we define $\tilde{\tau}(0) = 0$ and for $0 < i < \a$\index{$\tilde{\tau}$}:
\begin{equation}\label{eq8.3a}
\begin{split}
\tilde{\tau}(i) =  sup \{ \tau(j) + 1 : j < i \} =  min \{ k : k > \tau(j) \ \text{for all} \ j < i \},\\
\text{ So that } \  \tilde{\tau}(i+1) = \tau(i) + 1, \ \text{and} \hspace{3cm}\\
\tilde{\tau}(i) =  sup \{ \tau(j) : j < i \} \ \text{for} \ i \ \text{a limit ordinal}.\hspace{2cm}
\end{split}
\end{equation}
By induction, $\tau(i) \geq i$ and so the image of $\tau$ is cofinal in $\a$. From Proposition \ref{prop2.8}(c) follows that
$\{ [\tilde{\tau}(i),\tau(i)] =  [\tilde{\tau}(i),\tilde{\tau}(i + 1)): i < \a \}$
is a $\a$ indexed partition of $\a$ by closed intervals. Following (\ref{eq2.7}) we can identify $\a$ with the associated sum:
\begin{equation}\label{eq8.3b}
\sum_{i < \a} \ [\tilde{\tau}(i),\tilde{\tau}(i+1)) \ \cong \ \a.
\end{equation}

For each $i \in \alpha$ identify the interval
$[\tilde{\tau_x}(i),\tilde{\tau_{x}}(i+1)) \subset \alpha $ with the ordinal
$\beta(x,i) +1 = \tilde{\tau}_x(i + 1) \setminus \tilde{\tau}_x(i)$ which has its order type.  The restriction of $x$
to this interval yields an element $z(i) \in Z$ and $i \mapsto
z(i)$ defines the element of $\tilde{Z}$ which we associate with
$x$. The procedure is order preserving and reversible and so
defines the required isomorphism.

\end{proof} \vspace{.5cm}

\begin{lem}\label{lem8.2}  Let $F$ be a CHLOTS and $\{ f_{n} : F \rightarrow  F \quad \mbox{with} \ n \in \omega \setminus 1 \}$
be a sequence of order surjections.  Let $Y$ be the connected LOTS which is the
inverse limit of the  special inverse system indexed by \ $\omega$, defined by the sequence $\{ f_{n} \}$. Thus, $x \in \prod_{n \in \om} \ F$
is in $Y$  if and only if  $f_{n+1}(x_{n+1}) = x_{n}$ for all $n \in \om$.

Let  $\alpha$ be a tail-like ordinal with $\alpha \geq \omega ^{\omega}$.  With $J$ the
distinguished subinterval in $F$ we let $J_{0} = J$ and inductively define
$J_{n+1} = (f_{n+1})^{-1}(J_{n})$.
The points of $Y$ whose $n^{th}$ coordinate lies in $J_{n}$ for all $n \in \om$ define a compact interval $J_{Y}$ of the
connected LOTS Y and we use it to define the product space $Y_{\alpha}$. Then:
\begin{equation}\label{eq8.4}
Y_{\alpha} \qquad \cong \qquad F_{\alpha}. \hspace{3cm}
\end{equation}
\end{lem}

\begin{proof}  A point $z \in  Y_{\alpha}$ is
indexed by the order space product $\alpha \times \omega $ which
is order isomorphic to the ordinal product $\omega \cdot \alpha$,
\begin{equation}\label{eq8.5}
\begin{split}
z(i,n) \in X(i,n)\quad = \quad  \begin{cases}  F \quad \mbox{if} \   i = 0 \\
J_{n} \quad \mbox{if} \ 0 < i < \alpha.  \end{cases} \\
\mbox{with} \quad f_{n+1}(z(i,n+1))
 \ = \ z(i,n) \qquad \mbox{for} \
(i,n) \in \a \times \omega.
\end{split}
\end{equation}
For $ (i,n) > (0,0) $ in $\a \times \omega$ we let $z(<(i,n))$
denote the projection of the point to the subproduct $\Pi \{ X(j,m) : (j,m) < (i,n)
\}$ and let
\begin{equation} \label{eq8.6}
\begin{split}
I(z( <(i,n))) \ = \ \{ w(i,n) : w \in Y_{\alpha} \mbox{with} \ w(<(i,n)) = z(<(i,n))\}. \\
 \mbox{So that} \quad I(z(<(i,n))) = \begin{cases}  J \quad \mbox{if} \ n = 0 \\
(f_{n})^{-1}(z(i,n-1)) \quad \mbox{if} \ n > 0. \end{cases} \hspace{1cm}
\end{split}
\end{equation}
Thus,  each $I(z( <(i,n)))$ is a closed bounded interval in $X(i,n)$.

Now let
\begin{equation}\label{eq8.7}
\gamma (z) \quad = \quad \{(0,0)\} \cup \{ (i,n) : I(z( <(i,n)))
\mbox{ is nontrivial} \}.
\end{equation}

By (\ref{eq8.6}) $(i,0) \in \gamma(z)$ for all $z \in Y_{\alpha}$.  It
follows that, identifying $\gamma(z)$ with the ordinal which is
its order type,
\begin{equation}\label{eq8.8}
\alpha \quad \leq \quad \gamma(z) \quad \leq \quad \omega \cdot
\alpha \quad = \quad \alpha,
\end{equation}
with the latter equation from $\alpha = \omega ^{\rho} $ with
$\rho \geq \omega$ so that $1 + \rho = \rho$.  Let
\begin{equation}\label{eq8.9}
\tau_{z} : \alpha \rightarrow  \alpha \times \omega   \hspace{2cm}
\end{equation}
denote the unique order isomorphism onto $\gamma (z)$.

Notice that $z(<(i,n)) = w(<(i,n)) $ implies $\gamma (z)$ and
$\gamma (w)$ agree through $(i,n)$ and
\begin{equation}\label{eq8.10}
\tau_{z}(\beta) = \tau_{w}(\beta) \qquad \mbox{when} \quad
\tau_{z}(\beta) \leq (i,n).
\end{equation}

For $(i,n) \in \gamma (z)$ choose $g_{z(<(i,n))} : I(z(<(i,n)))
\rightarrow J$ an order isomorphism which exists because $F$
is a CHLOTS.

We now define for $z \in Y_{\alpha}$ the associated point $Q(z)
\in F_{\alpha}$
\begin{equation}\label{eq8.11}
Q(z)_{\beta} \quad = \quad g_{z(<(i,n))}(z(i,n)) \qquad
\mbox{where} \ (i,n) \ = \ \tau_{z}(\beta).
\end{equation}

If $z < w$ in $Y_{\alpha}$ and $(i,n)$ is the smallest coordinate
where $z(i,n) \not= w(i,n)$, then $z(i,n) < w(i,n)$ and since both
of these points are in $I(z( <(i,n))) = I(w( <(i,n)))$ it follows that
$(i,n) \in \gamma (z) \cap \gamma (w) $.  With $(i,n) =
\tau_{z}(\beta) = \tau_{w}(\beta) $ we have $Q(z)_{\epsilon} =
Q(w)_{\epsilon}$ for all $ \epsilon < \beta $ and $Q(z)_{\beta} <
Q_{w}(\beta) $.  Thus, $Q$ is an order injection.

Conversely, given $x \in F_{\alpha}$ we inductively define the
associated point $z \in Y_{\alpha}$ and the order injection
$\tau_{z}$.  Begin with $z(0,0) = x_{0}$ and $\tau_{z}(0) =
(0,0)$.  Now for $0  < \beta < \alpha$ we will define
$\tau_{z}(\beta) = (i,n)$ and $z(k,m) \in X(k,m)$ for all $(k,m)
\leq (i,n)$  so that
\begin{equation}
\begin{split}\label{eq8.12}
(k,m) \ = \ \tau_{z}(\epsilon) \quad \mbox{for some} \ 0 <
\epsilon \leq \beta \\ \Longleftrightarrow \qquad I(z(<(k,m))) \
\mbox{is nontrivial}.
\end{split}
\end{equation}

Now assume that the definitions have been completed for all
$\epsilon < \beta$.\vspace{.25cm}

\textbf{Case 1:} If $\beta$ is a limit ordinal, then let
\begin{equation} \label{eq8.13}
\begin{split}
\tau_{z}(\beta) \quad = \quad  sup \ \{ \tau_{z}(\epsilon) : \epsilon < \beta \} = (i,0) \hspace{1cm} \\
z(i,0) \quad = \quad (g_{z(<(i,0))})^{-1}( x_{\beta} ).
\hspace{2cm}
\end{split}
\end{equation}
Notice that the only limit elements of $\alpha \times \omega$ are
of the form $(i,0)$.\vspace{.25cm}

\textbf{Case 2:}  If $\beta = \epsilon + 1$ and
$\tau_{z}(\epsilon) = (k,m)$,
then by finite induction we define for $n = 0,1,...$
\begin{equation}\label{eq8.14}
\{z(k,m+n+1)\} \ = \ (f_{m+n+1})^{-1}(z(k,m+n)) \ = \
I(z(<(k,m+n+1)))
\end{equation}
if the preimage is a singleton.  There are two possibilities.

If this procedure stops for some finite $n$, then we define
\begin{equation}\label{eq8.15}
\begin{split}
\tau_{z}(\beta) \quad = \quad (k,m+n+1) \hspace{2cm} \\
z(k,m+n+1) \quad = \quad  (g_{z(<(k,m+n+1))})^{-1}( x_{\beta} ).
\end{split}
\end{equation}

If this procedure continues for all finite $n$, then we define
\begin{equation}\label{eq8.16}
\begin{split}
\tau_{z}(\beta) \quad = \quad (k+1,0) \hspace{2cm} \\
z(k+1,0) \quad = \quad  (g_{z(<(k+1,0))})^{-1}( x_{\beta} ).
\end{split}
\end{equation}\vspace{.25cm}

Notice that in both (\ref{eq8.15}) and (\ref{eq8.16})  $0 < \beta$ implies $x_{\beta} \in J$.

The order injection $\beta \mapsto \tau_{z}(\beta)$ from $\alpha$
into $\alpha \times \omega \cong \alpha$  has a cofinal image and
so we obtain a point $z \in Y_{\alpha}$.  Clearly, $Q(z) = x$ and
so $Q$ is the required order isomorphism.

\end{proof}\vspace{.5cm}

For the special case where all of the maps $f_n$ are equal to a fixed map $f$ we can define on
 $Y$  the \emph{shift automorphism}\index{shift automorphism} and its inverse $f_*$ by
\begin{equation}\label{eq8.16a}
\t(x)_n \ = \ x_{n+1}, \qquad f_*(x)_n = f(x_n)
\end{equation}
From (\ref{eq2.19x}) it follows that these are order isomorphisms.

Now we apply these preliminary results.  Recall that the size of
$X$ lies between $X_{1}$ and $X_{2}$ if $X_{1}$ injects into $X$
and $X$ injects into $X_{2}$.

\begin{theo}\label{theo8.3}  Assume that $F_{1}$ and $F_{2}$ are CHLOTS
\begin{enumerate}
\item[(a)]  If $F_{1}$ and $F_{2}$ have the same size and if $\alpha $  is a
tail-like ordinal with $ \alpha \geq \omega^{\omega}$ , then
\begin{equation}\label{eq8.17}
(F_{1})_{\alpha} \quad  \cong  \quad (F_{2})_{\alpha}.
\hspace{2cm}
\end{equation}
\item[(b)]  If for some countable ordinal $\beta$, the size of $F_{1}$ lies
between $F_{2}$ and $(F_{2})_{\beta}$ and if $\alpha $  is a sufficiently large countable tail-like ordinal,
then the isomorphism of (\ref{eq8.17}) holds.  Specifically, for ordinals $i, j$
\begin{equation}\label{eq8.18}
\begin{split}
\mbox{if} \quad \beta \cdot \omega^{i} \ = \ \omega^{i} \qquad
\mbox{and} \qquad  j \geq (i + \omega) \\ \mbox{then} \quad
\alpha \ = \ \omega^{j} \ \mbox{implies} \quad (F_{1})_{\alpha}
\cong  (F_{2})_{\alpha}.
\end{split}
\end{equation}
\end{enumerate}
\end{theo}

\begin{proof}  (a):  By Corollary \ref{cor4.5} there exist
continuous order surjections $ g_{21} : F_{1} \rightarrow F_{2} $
and  $ g_{12} : F_{2} \rightarrow F_{1} $.  Let $ g_{1} = g_{12}
\circ g_{21} : F_{1} \rightarrow F_{1} $ and  $ g_{2} = g_{21}
\circ g_{12} : F_{2} \rightarrow F_{2} $.  Let $Y_{1}$ be the
inverse limit of the special inverse system indexed by $\omega$
with $f_{n} = g_{1} $ for all $n \in \omega$ and similarly for
$Y_{2}$.  By Lemma \ref{lem8.2},
\begin{equation}\label{eq8.19}
(F_{1})_{\alpha} \  \cong \  (Y_{1})_{\alpha} \qquad \mbox{and}
\qquad (F_{2})_{\alpha} \  \cong \  (Y_{2})_{\alpha}.
\end{equation}

Now define the order surjection $ \tilde{g}_{21} : Y_{1}
\rightarrow Y_{2} $  to be a copy of $g_{21}$ on each coordinate,
and similarly define $ \tilde{g}_{12} : Y_{2} \rightarrow Y_{1} $.
The compositions each way are the shift automorphisms on the
inverse limits. So we have $Y_{1} \cong Y_{2}$.  Hence,
$(Y_{1})_{\alpha} \cong (Y_{2})_{\alpha} $ which completes the
required chain of isomorphisms.

(b):  As in (\ref{eq2.11}) $\beta =  \omega^{k_{1}} + ... +\omega^{k_{N}} $
with  $k_{1} \geq ...\geq k_{N}$. If $i$ is a tail-like ordinal
with $i > k_{1} $, then
\begin{equation}\label{eq8.20}
\omega^{i} \ \leq \ \beta \cdot \omega^{i} \ \leq  \
\omega^{k_{1}+1} \cdot \omega^{i} \ = \ \omega^{k_{1}+1+i} \ = \
\omega^{i}.
\end{equation}

Now from Proposition \ref{prop4.6}(c),(d) we have, with $\gamma = \omega^{i}$ that
\begin{equation}\label{eq8.21}
\begin{split}
((F_{2})_{\beta})_{\gamma} \ \mbox{is  at least as big as} \
(F_{1})_{\gamma}, \hspace{2cm} \\
(F_{1})_{\gamma} \ \mbox{is at least as big as} \  (F_{2})_{\gamma}, \hspace{2cm} \\
\mbox{and} \qquad  ((F_{2})_{\beta})_{\gamma}  \cong  \
(F_{2})_{\gamma} \hspace{2.5cm}
\end{split}
\end{equation}
since $\beta \cdot \gamma = \gamma$. Thus, $(F_{1})_{\gamma}$ and $(F_{2})_{\gamma}$ have the same size.

Now let $\tilde{j} = j \setminus i$ so that $i + \tilde{j} = j$.
By assumption $\tilde{j} \geq \omega$ and so $\tilde{\alpha} =
\omega^{\tilde{j}} \geq \omega^{\omega}$. So by part (a) we have
\begin{equation}\label{eq8.22}
((F_{2})_{\gamma})_{\tilde{\alpha}} \quad \cong \quad
((F_{1})_{\gamma})_{\tilde{\alpha}}.  \hspace{1.5cm}
\end{equation}
Since $\a = \gamma \cdot \tilde{\alpha}$ the result  follows
from Proposition \ref{prop4.6}(d) again.

\end{proof} \vspace{.5cm}

\noindent{\bfseries Remark.}  The interest in part (a) of the theorem comes
from the fact that there exists CHLOTS
$F_{1}$ and $F_{2}$ which are of the same size, i.e. each can be
order injected into the other, but which are not order isomorphic,
e.g. see Proposition \ref{prop8.8} below.
\vspace{.5cm}

Recall that a LOTS $X$ is $\R$-bounded when it admits an order injection into $\R_{\d}$ for some countable ordinal $\d$. Proposition \ref{propnew04a}
implies that if $X$ is an $\R$-bounded HLOTS, then  for any countable tail-like ordinal $\b$, $\widehat{X_\b}$
 is an $\R$-bounded CHLOTS.

\begin{cor}\label{cor8.3a} If $X$ is an $\R$-bounded CHLOTS, then for
$\a$ a sufficiently large countable tail-like ordinal $X_{\a} \cong \R_{\a}$. \end{cor}

\begin{proof} $\R$ injects into any CHLOTS and so if $X$ is an $\R$-bounded CHLOTS, it has size between that of $\R$
and that of $\R_{\d}$ for some countable $\d$.  The result follows from Theorem \ref{theo8.3} (b).

\end{proof}\vspace{.5cm}

If $X$ is an IHLOTS with completion $F$, then the size of the
completion of $X_{\omega}$ lies between $F$ and $F_{\omega}$.  The
first conjecture might be that it is order isomorphic with
$F_{\omega}$ but we will now see that this is rarely true. We will
require some preliminary study of the completion of the elements
of the tower over $X$.

Let $X$ be  an IHLOTS with completion the CHLOTS $F$ and let $Y$
be the complementary IHLOTS, i.e. $ Y = F \setminus X$.  Assume
$-1 < +1$ in $X$ and that $J$ is the interval $[-1,+1] \subset F$.
Let $\alpha$ be a limit ordinal. We can regard $X_{\alpha}$ as the
branch space of the subtree $T$ of the simple tree on $X,\alpha $
with $S_{0} = X$ but $S_{p} \cong  J \cap X$ instead of all of
$X$. When we apply Proposition \ref{prop5.6} we can identify the completion
with the branch space of the tree completion $\hat{T}$.  So we
will use
\begin{equation}\label{eq8.23}\begin{split}
\widehat{X_{\alpha}} \quad = \quad X_{\alpha} \ \cup \ Y^1 \ \cup \hspace{5cm}\\
\bigcup _{0 <\beta < \alpha} \{ x \in F^{\beta +1} : x(i) \in \begin{cases} \  \  \ X  \quad \mbox{for} \ i = 0,\\
  \ J \cap X   \ \mbox{for} \ 0 < i < \beta, \\  \ J \cap Y   \ \mbox{for} \ i = \b \end{cases} \}.
\end{split} \end{equation}

With these identifications the canonical projections for $0 <
\beta \leq \alpha$
\begin{equation}\label{eq8.24}
\begin{split}
\pi^{\a}_{\b} : X_{\alpha} \rightarrow X_{\beta}, \hspace{3cm} \\
\hat{\pi}^{\a}_{\b} : \widehat{X_{\alpha}} \rightarrow
\widehat{X_{\beta}} \hspace{3cm}
\end{split}
\end{equation}
are the obvious restriction maps.  With $\beta = 1$ we identify
the completion with $F$ and so define $\hat{\pi}^{\alpha} :
\widehat{X_{\alpha}} \rightarrow F $ which we write as $\hat{\pi}$ when the
subscript is unambiguous. We use it to define, for $y \in F$
\begin{equation}\label{eq8.25}
P(y) \quad = \quad \{ z \in \widehat{X_{\alpha}} :  z(0) = y
\}\ = \ (\hat{\pi})^{-1}(y).
\end{equation}\index{$P(y)$}
Since $\hat{\pi}$ is a continuous order surjection, $P(y)$ is a
nonempty compact interval in $\widehat{X_{\alpha}}$ for all $y \in F$.

We carry over the tree concept of \emph{height}\index{height}, defining for  $ x
\in \widehat{X_{\alpha}}$
\begin{equation}\label{eq8.26}
h(x) \ = \ \begin{cases} \alpha \quad \mbox{for} \ x \in
X_{\alpha} \\ \beta + 1 \quad \mbox{for} \ x \in F^{\beta
+1}.\end{cases}
\end{equation}

We identify $Y$ with the elements of height $1$ by $x \mapsto x(0)$, so that $Y \subset \widehat{X_{\alpha}}$.

If $h(x) \geq \epsilon$, then we will write $x|\epsilon$ for
$\hat{\pi}^{\a}_{\ep}(x)$, i.e. the restriction of the
map $x$ to the subset $\epsilon$ of its domain.  Clearly,
\begin{equation}\label{eq8.27}
h(x) > \epsilon \qquad \Longleftrightarrow \qquad
\hat{\pi}_{\epsilon}^{\alpha}(x) \in X_{\epsilon}.
\end{equation}

Now suppose $ 0 < \epsilon < \alpha $ and $ w \in X_{\epsilon} $.
We define the compact subinterval $J_{w} \subset
\widehat{X_{\alpha}}$ and the map $\hat{\pi}_{w} : J_{w}
\rightarrow J$ by
\begin{equation}\label{eq8.28}
\begin{split}
J_{w} \quad = \quad   (\hat{\pi}^{\a}_{\ep})^{-1}(w) \ = \
 \{ z \in \widehat{X_{\alpha}} : z|\epsilon = w \} \ = \ [w-,w+], \\
\hat{\pi}_{w}(z) \quad = \quad z(\epsilon). \hspace{4cm}
\end{split}
\end{equation}
It is clear from description (\ref{eq8.23}) that $\hat{\pi}_{w}$ is an
order surjection for all $w \in X_{\epsilon}$. In particular,
$J_{w}$ is nontrivial. For  $w \in X_{\epsilon}, y \in J$ we
extend definition (\ref{eq8.25})
\begin{equation}\label{eq8.29}
P_{w}(y) \quad = \quad \{ z \in J_{w} :  z(\epsilon) = y \}\
= \ (\hat{\pi}_{w})^{-1}(y).
\end{equation}
As before each $P_{w}(y)$ is a nonempty compact subinterval of
$J_{w} \subset \widehat{X_{\alpha}}$.

Now assume that $I = [x,y]$ is a nontrivial, closed interval in
$\widehat{X_{\alpha}}$, so that $x < y$.  We denote by
$\epsilon(I)$ the equality level of the pair $x,y$, i.e. $min \
\{i : x_{i} \not= y_{i} \} $.  With $\epsilon = \epsilon(I)$ we
have, for all $z \in I$
\begin{equation}\label{eq8.30}
\begin{split}
h(z) > \epsilon, \hspace{5cm}\\
x|\epsilon \ = \ z|\epsilon \ = \ y|\epsilon,\qquad \mbox{and}
\hspace{4cm}\\  x(\epsilon) \ \leq \ z(\epsilon) \ \leq \
y(\epsilon), \hspace{4cm}
\end{split}
\end{equation}
with at least one of the latter inequalities strict.  We call the
common element of $F^{\epsilon}$ the \emph{stem of $I$}\index{stem of $I$}\index{$stem(I)$}.  From
(\ref{eq8.27}) we have
\begin{equation}\label{eq8.31}
stem(I) \ \in \ X_{\epsilon(I)}. \hspace{2.5cm}
\end{equation}

Finally, we define
\begin{equation}\label{eq8.32}
\begin{split}
span(I) \quad = \quad [x(\epsilon),y(\epsilon)] \hspace{2cm}\\
span^{\circ}(I) \quad = \quad (x(\epsilon),y(\epsilon)).
\hspace{1.5cm}
\end{split}
\end{equation}
So that $span(I)$\index{$span(I)$} is a nontrivial, compact subinterval of $F$ and
is contained in $ J $ when $ \epsilon > 0$. Its interior
$span^{\circ}(I)$ is nonempty.

The most important special case occurs when $\epsilon(I) = 0 $ in
which case $stem(I) = \emptyset $.  Clearly, with $I = [x,y]$
\begin{equation}\label{eq8.33}
\begin{split}
\epsilon(I) = 0 \quad \Longleftrightarrow \quad x(0) < y(0) \\ \mbox{in which case} \hspace{5cm} \\
span(I) \quad = \quad [x(0),y(0)] \quad = \quad \hat{\pi}(I)
\end{split}
\end{equation}
Notice that (\ref{eq8.30}) implies that
\begin{equation}\label{eq8.34}
Y \cap I \ \not= \ \emptyset \qquad \Rightarrow \qquad \epsilon(I)
= 0.
\end{equation}

\begin{lem}\label{lem8.4}   Let $I, I_{1}, I_{2}$ be nontrivial, closed subintervals of $\widehat{X_{\alpha}}$.
\begin{enumerate}
\item[(a)]  Given $y \in F$ the compact interval  $P(y)$ is trivial, i.e. is a
singleton,  if and only if  $y \in Y$.  Given  $\epsilon > 0$,  $w \in X_{\epsilon}$ and
$y \in J$ the compact interval  $P_{w}(y)$ is trivial, i.e. is a singleton,  if and only if  $y \in Y \cap J$.
\item[(b)]  If $z \in \widehat{X_{\alpha}} $  with  $ \hat{\pi}_{\epsilon(I)}^{\alpha}(z) = stem(I) $
 , then $h(z) > \epsilon(I)$.  If, in addition, $z(\epsilon(I)) \in span^{\circ}(I) $ , then  $ z \in I$.
\item[(c)]  If $I_{1}$ and $I_{2}$ are disjoint subintervals such that
$\epsilon(I_{1}) = \epsilon(I_{2})$ and $stem(I_{1}) = stem(I_{2})$, then
$span(I_{1})$ and $span(I_{2})$ are non-overlapping subintervals of F, i.e.
\begin{equation}\label{eq8.35}
span^{\circ}(I_{1}) \ \cap \ span^{\circ}(I_{2}) \quad = \quad
\emptyset  \hspace{1.5cm}
\end{equation}
\end{enumerate}
\end{lem}

\begin{proof}  (a):   This is obvious from the
identification (\ref{eq8.23}).

(b):    Let $I = [x,y]$.  By (\ref{eq8.31}) and (\ref{eq8.27}) $
\hat{\pi}_{\epsilon(I)}^{\alpha}(z) = stem(I) $ implies $h(z) >
\epsilon(I)$ and $x|\epsilon  =  z|\epsilon  =  y|\epsilon $.
If, in addition,     $z \in span^{\circ} (I)$, then $x(\epsilon) <
z(\epsilon) < y(\epsilon)$ and so $z \in (x,y)$. \vspace{.2cm}

(c):   Obvious from (a) and (b).

\end{proof} \vspace{.5cm}

\begin{theo}\label{theo8.5}  Let $X$ be an IHLOTS with completion the CHLOTS $F$ and let $Y$
be the complementary IHLOTS $ Y = F \setminus X$.  If for any pair $\alpha$ and
$\beta$ of countable limit ordinals the completion $\widehat{X_{\alpha}}$ is order
isomorphic to $F_{\beta}$, then $X$ is a first category subset of $F$.  That is, $Y$
contains a dense $G_{\delta}$ subset of $F$.
\end{theo}

\begin{proof}  Assume that $f :
\widehat{X_{\alpha}}  \rightarrow  F_{\beta} $ is an order
surjection. By Proposition \ref{prop2.1}(a) $f$ is continuous and
topologically proper.  We will show that if $Y$ does not contain a
particular dense $G_{\delta}$ set which we will construct, then $f$ is
not injective.

Using the projection $\pi_{i}^{\b} : F_{\beta} \rightarrow
F_{i} $ for $0 < i < \beta$ we define, for each $y \in Y \subset
\widehat{X_{\a}}$
\begin{equation}\label{eq8.36}
Q(y,i) \quad = \quad (\pi_{i}^{\b}  \circ
f)^{-1}(\pi_{i}^{\b} (f(y))) \ = \ f^{-1}(J_{f(y)|i}),
\end{equation}
where the latter equation uses a definition analogous to (\ref{eq8.28}). For any  $w \in
F_{i}$,  $J_{w} = (\pi_{i}^{\b})^{-1} (w)$ is a nontrivial, compact interval in $F_{\beta}$
and so $Q(y,i)$ is a nontrivial, compact interval in
$\widehat{X_{\alpha}}$.

Clearly, for $y_{1}, y_{2} \in Y$ and $i < \beta$
\begin{equation}\label{eq8.37}\begin{split}
Q(y_{1},i) \cap Q(y_{2},i) \  \not= \  \emptyset \quad \Rightarrow \hspace{2.5cm} \\
\quad f(y_{1})|i = f(y_{2})|i \quad \Rightarrow \hspace{3cm} \\ \quad Q(y_{1},i) \
= \  Q(y_{2},i). \hspace{3cm} \end{split}
\end{equation}

Because $y \in Q(y,i)$, (\ref{eq8.34}) implies $\epsilon(Q(y,i)) = 0$ for
all $y \in Y$ and $0 < i < \beta $.  It follows from (\ref{eq8.33}),
(\ref{eq8.37}) and Lemma \ref{lem8.4}(c) that distinct members of the set of
intervals
\begin{equation}\label{eq8.38}
\mathcal{Q}_{i}\quad = \quad \{ span(Q(y,i)) : y \in Y \} \
= \ \{ \pi(Q(y,i)) : y \in Y \}
\end{equation}
are non-overlapping.  Since $y \in span(Q(y,i))$ and $Y$ is dense
in X, $\mathcal{Q}_{i}$ has a dense union for each positive $i <
\beta$.  As each member of $\mathcal{Q}_{i}$ is nontrivial the
open set
\begin{equation}\label{eq8.39}
O_{i} \quad = \quad  \bigcup \ \{ span^{\circ}(Q(y,i)) :y \in Y \}
\hspace{1cm}
\end{equation}
is dense in $F$ for each positive $i < \beta$.

Since $F$ is locally compact and $\beta$ is countable, the Baire
Category Theorem implies that $D = \cap O_{i}$ is a dense
$G_{\delta}$ subset of $F$. If there exists $ t \in  X \cap D $
then  by Lemma \ref{lem8.4}(a) the interval $P(t)$ is nontrivial.  We will show
that $f$ is constant on $P = P(t)$.  It suffices to show that
$\pi_{i}^{\beta} \circ f$ is constant on $P$ for each $i < \beta $
.

Since $t \in O_{i}$ there exists $y \in Y$ such that $t \in
span^{\circ}(Q(y,i))$.  With $I = Q(y,i), \epsilon(I) = 0$, and
$stem(I) = \emptyset$. If $x \in P$, then $ x(\ep(I)) = \hat \pi (x) = t \in span^{\circ}(I)$.
So Lemma \ref{lem8.4}(b) implies that $x \in I$. That is,
 $P \subset Q(y,i)$, and so $\pi_{i}^{\beta}(f(x)) =
\pi_{i}^{\beta}(f(y))$ for all $x \in P$.

Thus, $f$ can be injective only when $D \subset Y$.

\end{proof}\vspace{.5cm}

\subsection{IHLOTS in $\R$}

Now we consider examples which are constructed from IHLOTS in
$\R$. First we will see that there are many such.

\begin{prop}\label{prop8.6}  Let $\mathcal{G}$\index{$\mathcal{G}$} be the countable group of positive,affine
 transformations of $\R$ with rational coefficients, i.e. of the
 form $t \mapsto at + b $  with $a,b \in \Q$ and $a > 0$.  If $X$
 is a nonempty, proper subset of $\R$ which is invariant with respect
 to the action of $\mathcal{G}$, then $X$ is an IHLOTS.  In particular, any proper
  subfield of $\R$ is an IHLOTS as is its complement.
\end{prop}

\begin{proof}  Assume first that $X$ contains
some rational number and so that $\Q \subset X$. Apply
Lemma \ref{lem5.19} with $W = \Q$ to see that $X$ is a HLOTS. If
$\Q \cap X = \emptyset $, then the same result show that
$\R \setminus X $  is a HLOTS and so by Proposition \ref{prop3.4}(d)
$X$ itself is a HLOTS.  Since $X$ and its complement are clearly
dense they are both IHLOTS with completion $\R$. A
subfield of $\R$ contains $\Q$ and so is
invariant under $\mathcal{G}$.

Notice that this is a special case of Proposition \ref{propauto6}.

\end{proof}\vspace{.5cm}

We will use a bit of classical topology. A topological space $X$ is called a \emph{Polish space}\index{Polish space} when it is a s
second countable space which admits a complete metric. so, of course, the complete metric space $\R$ is Polish. It is a classical result of
Alexandroff and Hausdorff, see \cite{K} page 208, that a
 $G_{\delta}$ subset of a complete metric space admits a complete metric. Hence, a $G_{\d}$ subset of a Polish space is Polish.

A \emph{Cantor set}\index{Cantor set} is a non-empty, zero-dimensional, compact, perfect metric space.
Any Cantor set is homeomorphic to the classical Cantor Set in $[0,1]$.

\begin{lem}\label{lem8.7}  Let $X$ be a nonempty Polish space.
\begin{enumerate}
\item[(a)]  If $X$ has no isolated points, i.e. $X$ is \emph{perfect}, then $X$ contains
a Cantor set and so is uncountable.
\item[(b)]  If $X_{0} = \cup  \{ U : U $ is a countable, open subset of $X  \}$, then the
open set  $X_{0}$ is a countable Polish space. The complementary closed set,  $X_{1} = X \setminus X_{0}$,
 is a Polish space with no isolated points.
 \item[(c)] The space $X$ is perfect  if and only if  every nonempty open subset is uncountable.
\end{enumerate}
\end{lem}

\begin{proof} (a) Choose a complete metric and use the usual dichotomy
procedure.  With $A_{0} = X$, define for each word $x \in
\{-1,+1\}^{n+1}$ a closed set $A_{x}$ with a nonempty interior, of
diameter at most $2^{-n}$ such that $A_{x,\pm1} \subset Int \
A_{x}$.  From the Cantor Intersection Theorem we obtain a
topological embedding of $\{-1,+1\}^{\omega}$ into $X$.

(b)  Because $X$ has a countable base it follows that $X_{0}$ is
countable.  If $x \in X_{1}$, then any neighborhood $U$ of $x$ in
$X$ is uncountable and so $U \cap X_{1} = U \setminus X_{0}$ is
uncountable.  Any $G_{\delta}$ subset of a Polish space is Polish.

(c) If $X$ is a perfect Polish space, then every nonempty open subset is a perfect Polish space which therefore
contains a Cantor set by (a). Thus, every nonempty open subset is uncountable.  The converse is obvious.

\end{proof}\vspace{.5cm}

\begin{prop}\label{prop8.8}  Let $X$ be an IHLOTS in $\R$ and let $\alpha$ be an
infinite tail-like ordinal.  If $X$ contains a Cantor set, then $X_{\alpha}$,
its completion $\widehat{X_{\alpha}}$ and $\R_{\alpha}$ all have the same size.
  If, in addition, $X$ is not a first category subset of $\R$, then no two of
   these HLOTS are
   homeomorphic and so not order isomorphic.  In particular, if $X$ is the IHLOTS of irrational
   numbers, i.e. $X = \I$, then $\widehat{X_{\alpha}}$
      and $\R_{\alpha}$ are CHLOTS of the same size which are not homeomorphic and so not order isomorphic.
\end{prop}

\begin{proof} $X_{\alpha}$ is a subset of
$\R_{\alpha}$ and so  by Proposition \ref{prop4.6}(b)
$\widehat{X_{\a}}$ injects into  $\R_{\alpha}$.  We
can assume that the Cantor set $C$ is contained in the
distinguished interval $J$ of $X$. The order isomorphism between
$\Q$ and the set of left end-points in $C$, excluding
$max \ C$, shows that $\Q$ injects into $C$ and so by
Proposition \ref{prop4.6}(b) again $\R$ injects into $C$. In fact, if $\tilde C$ is the classical Cantor Set $C$ with the right end-points
 and the min removed, then the usual Cantor map from $C$ onto $[0,1]$ restricts to an order isomorphism from
 $\tilde C$ onto $(0,1)$.

By Proposition \ref{prop4.6}(c) $\R_{\alpha} \subset \R^{\alpha} $ injects into
$C^{\alpha} \subset X_{\alpha}$.

If $X$ is not of first category,
then by Theorem \ref{theo8.5} the CHLOTS $\widehat{X_{\alpha}}$ and
$\R_{\alpha}$ are not order isomorphic.
If  $h : \widehat{X_{\alpha}} \to \R_{\alpha}$ were a homeomorphism, then by Lemma \ref{lem3.1}
it be either an order isomorphism or an order$^*$ isomorphism.  Since $\R$ is symmetric, an order isomorphism
would then exist.

The IHLOTS $X_{\alpha}$ has dense holes and so is not locally compact as are $\widehat{X_{\alpha}}$ and $ \R_{\alpha}$. So it
is not homeomorphic to either of them.

The set of irrationals is not of first category by the Baire
Category Theorem and it contains a Cantor set by Lemma \ref{lem8.7}(a).

\end{proof} \vspace{.5cm}

We call a subset $A$ of a Polish space $X$ a \emph{Mycielski set}\index{Mycielski set}
if it is a countably infinite union of Cantor sets in $X$, see, e.g. \cite{A2}.

\begin{lem}\label{lem8.9} Let $X$ be a Polish space.
\begin{enumerate}
\item[(a)] A countable union of Mycielski sets in $X$ is a Mycielski set.
\item[(b)] The nonempty intersection of a Mycielski set and an open set is a Mycielski set when it is nonempty.
\item[(c)] If $A$ is a Mycielski set and $B$ is a countable subset of the closure
of $A$ in $X$, then $A \cup B$  is a Mycielski set.
\end{enumerate}
\end{lem}

\begin{proof}  (a): Obvious.

(b): If $A = \bigcup_{n} C_{n}$ with $C_{n}$ a Cantor set, and
$U$ is open, then $A \cap U = \bigcup_{n} (C_{n} \cap U)$.
Since $C_{n} \cap U$ is open in $C_{n}$ it is, when nonempty, the countable union
of nonempty clopen subsets of $C_{n}$ each of which is a Cantor set.

(c): By (a) it suffices to show that $ A \cup \{x\} $ is a
Mycielski set when $x$ is a limit point of $A$.  Let $U_{n}$ be
the open ball of radius $2^{-n}$ centered at $x$.  By (b), $A \cap
U_{n}$ is a Mycielski set and since $x \in \overline{A}$ each is
nonempty. If $C_{n}$ is a Cantor set in $A \cap U_{n}$, then $
C_{x} \ = \  \{x\} \cup (\bigcup_{n}  C_{n})$ is a Cantor
set since it is closed, zero-dimensional and without isolated points. So $A \cup \{x\} = A \cup C_{x}$ is a Mycielski set.

\end{proof}\vspace{.5cm}

\noindent{\bfseries Remark.} If $A$ is a Mycielski set and $D$ is a
countable set, then $A \setminus D$ need not be a Mycielski set.
If $X$ is a Polish space and $D$, dense in $X$, is a countable
union of closed nowhere dense subsets of $X$, then by the Baire
Category Theorem $X \setminus D$ is a $G_{\delta}$ subset but not
an $F_{\sigma}$ subset of $X$.  In particular, if $A$ is a Cantor
set and $D$ is a countable, dense subset of $A$, then $A \setminus D$
is not  $\sigma$-compact and so is not a Mycielski set.

\begin{theo}\label{theo8.10} Let $X$ be a dense, proper subset of $\R$ and let $Y = \R \setminus X $.
\begin{enumerate}
\item[(a)]  The following conditions are equivalent.
\begin{itemize}
\item[(1)]  $X$ is a Mycielski set in $\R$.
\item[(2)]  $Y$ is a dense subset of $\R$, $X$ is an $F_{\sigma}$
subset and every nonempty open subset of $X$ is uncountable.
\item[(3)]  $X$ is an uncountable $F_{\sigma}$ subset of $\R$ which is a HLOTS.
\item[(4)]  $Y$ is a dense subset of $\R$ which is order isomorphic to $\Q^{\omega}$.
\end{itemize}
\item[(b)]  If $X_{1}$ is a dense Mycielski subset of $\R$, then there
exists $f \in H_{+}(\R)$ such that $f(X) = f(X_{1})$.
\end{enumerate}
\end{theo}

\begin{proof} (a) (2) $\Rightarrow$  (1):  Let $X
= \bigcup_{n} A_{n}$ with each $A_{n}$ closed in ${\R}$.
Let $\mathcal{B}$ be a countable base for $\R$.
For each $n$ and $U \in \mathcal{B}$ we use Lemma \ref{lem8.7} to decompose
$U \cap A_{n} = B(U,n) \cup C(U,n)$ with $B(U,n)$ countable and
$C(U,n)$ perfect or empty. The closure $\bar{C}(U,n) \subset
A_{n}$ is perfect and is nowhere dense because $Y$ is dense.
Thus, each nonempty $\bar{C}(U,n)$ is a Cantor set.  For each $U$
$U \cap X$ is uncountable and so some $\bar{C}(U,n)$ is nonempty.
Hence, $\tilde{X} = \bigcup_{U,n} \{ \bar{C}(U,n) \}$ is a dense
Mycielski set. Hence, $X = \tilde{X} \cup (\bigcup_{U,n}  B(U,n))$
 is a Mycielski set by Lemma \ref{lem8.9}(c).

(1) $\Rightarrow$ (4):  If $C$ is a Cantor set in $\R$
with $min = a$ and $max = b$, then we will call the components of
the open set $[a,b] \setminus C$ the \emph{complementary intervals
for $C$}. The LOTS of complementary intervals for $C$ is order
isomorphic with $\Q$.

We will repeatedly use the following:

{\bfseries Fact } If $\epsilon > 0$ and $a < b \ \in \  X$, then
there exists a Cantor set $C \subset X$ with $a = min \ C,b = max
\ C$, and such that the diameter of each complementary subinterval
is less than $\epsilon$. We will call such a $C$ an \emph{$X$
Cantor set for $[a,b]$ with mesh less than $\epsilon$}.

\begin{proof}:  Choose $f : {\Z} \rightarrow (a,b)$ an order
injection with image $\pm$cofinal and such that for all
$i$ $f(i+1) - f(i) < \epsilon /2$. For each $i$ use Lemma \ref{lem8.9}(b) to
choose a Cantor set $C(i) \subset X \cap (f(i),f(i+1))$. Let $C =
\{a,b\} \cup (\bigcup_{i} C(i) )$.

\end{proof}

Now write $X$ as the countable union of Cantor sets $C(n)$ and
proceed inductively.

Because $X$ is dense we can choose an order injection $ f :
{\Z} \rightarrow X$ with image $\pm$cofinal in
$\R$ and $f(i+1) - f(i) < 1$ for all $i$.  For each $i$
choose an $X$ Cantor set for $[f(i),f(i+1)]$ which contains
$[f(i),f(i+1)] \cap C(0)$.  Let $A(0)$ be the closed set which is
their union. Choose an order isomorphism $q \mapsto J(q)$ from
$\Q$ to the set of intervals complementary to $A(0)$ in
$\R$.  We will call a closed subset $A$ of $\R$
an \emph{extended Cantor set} if $A \cap [f(i),f(i+1)]$ is a
Cantor set for all $i \in {\Z}$.  If $A$ is a perfect,
closed subset of $\R$ with $A(0) \subset A \subset X$
then $A$ is an extended Cantor set because $Y$ is dense.

Assume that we have defined for $i = 0,...,n$ an extended Cantor
set $A(i)$ and an order isomorphism from $\Q^{i}$ to the
set of complementary intervals for $A(i)$ such that for $i =
1,...,n$ and $q_{0}...q_{i} \in \Q^{i}$
\begin{equation}\label{eq8.40}
A(i-1) \cup C(i) \ \subset  \ A(i)\qquad \mbox{and} \qquad
\overline{J(q_{0}...q_{i})} \subset J(q_{0}...q_{i-1}).
\end{equation}
Furthermore, the complementary intervals for $A(i)$ have diameter
at most $2^{-i+1}$

For the next step, choose for each $A(n)$-complementary interval
$J(q_{0}...q_{n}) = (a,b)$ an $X$ Cantor set for $[a,b]$ which
contains $[a,b] \cap C(n+1)$ and which has mesh less than
$2^{-n}$. Choose an order isomorphism $q \mapsto
J(q_{0}...q_{n}q)$ from $\Q$ to the set of complementary
intervals. Let $A(n+1)$ be the union of $A(n)$ together with these
newly constructed Cantor sets.

From the construction, $X = \bigcup _{n} A(n) $ and for $x \in
\Q^{\omega}$ the intersection
\begin{equation}\label{eq8.41}
\bigcap \{ J(x(0)...x(n)) : n \in \omega \} \quad = \quad \bigcap
\{ \overline{J(x(0)...x(n))} : n \in \omega \}
\end{equation}
is a single point of $Y$.  If we  denote this point $g(x)$, then $g
: \Q^{\omega} \rightarrow Y $ is an order isomorphism.

(4) $\Rightarrow$ (3) and (b): By Corollary \ref{cor6.2} $\Q^{\omega}$ is a HLOTS.
From (4) it follows that $Y$ is a HLOTS
and since $Y$ is dense, it has completion $\R$.  By
Proposition \ref{prop3.4}(d), its complement $X$ is a HLOTS as well.

If $X_{1}$ as well as $X$ satisfy condition (4), then there is an
order isomorphism between the complements $Y_{1}$ and $Y$ because
both are isomorphic to $\Q^{\omega}$. The extension to
the completion $\R$ restricts to an isomorphism
between $X_{1}$ and $X$.

Since (1) implies (4) we can start with any dense Mycielski set
$X_{1}$.  The isomorphism shows that $X$ is an uncountable $F_{\s}$.

(3) $\Rightarrow$ (2):  Since $X$ is dense it has completion
$\R$ and the complementary IHLOTS $Y$ is dense as well
because $X$ is a proper subset of $\R$.  Any open
interval in $X$ is order isomorphic to $X$ itself and so is
uncountable.

\end{proof} \vspace{.5cm}

Beginning with any Mycielski set or Cantor set in $\R$ we can close it up under the action of the rational affine group described in
Proposition \ref{prop8.6}. The resulting union is a dense Mycielski set and it
is an IHLOTS by Proposition \ref{prop8.6}. Once we know that any two dense Mycielski sets
in $\R$ are order isomorphic, as in Theorem \ref{theo8.10} above,
it becomes clear that they are all IHLOTS.

\begin{cor}\label{cor8.11} If $X$ is a dense, Mycielski subset of $\R$, then $X$ is an IHLOTS
and the completion $\widehat{X_{\omega}}$ is order isomorphic to $\R_{\omega}$.
\end{cor}

\begin{proof}   By Theorem \ref{theo8.10}
$\Q^{\omega}$ is isomorphic to $ Y $, the complement of $X$ in $\R$ and with completion
$\R$.  That is, $\widehat{\Q^{\omega}} \setminus \Q^{\omega} \cong X$.

Lemma \ref{lem8.1} implies that $(\widehat{\Q^{\omega}} \setminus \Q^{\omega})^{\omega} $ has
completion isomorphic to $ \R_{\omega}$.  So  $X^{\omega}$ has completion
isomorphic to $\R_{\omega}$.  By
 Proposition \ref{prop6.3} $X^{\omega}$ and $X_{\omega}$ have isomorphic completions.

\end{proof} \vspace{.5cm}

  We saw in Proposition \ref{propnew04a} that a tree of $\Q$ type with countable height has an $\R$-bounded branch space.
  In Corollary \ref{corcon07} we saw that an Aronszajn tree does not have an $\R$-bounded branch space.  Nonetheless, it can
  happen that a tree with height $\Omega$  has an $\R$-bounded branch space.

  \begin{ex}\label{exOmegaRbounded} There exists an $\Omega$-bounded normal tree of $\Q$ type and of height $\Omega$ whose branch space is $\R$-bounded.\end{ex}

  \begin{proof} Let $A$ be an Aronszajn tree of $\Q$ type, see Corollary \ref{cor6.16}. If $\a$ is a countable limit ordinal
  then by Proposition \ref{prop5.9} (a), the branch space $X(A^{\a})$ has completion isomorphic to $\R$ and so we can choose
  $i_{\a}$ an order injection from $X(A^{\a})$ to $(-1,1) = J^{\circ}$.

  Let $T^{\om + 1}$ be the simple tree on $\Q, \om + 1$ so that $L_{\om} = \Q^{\om}$.

  Choose a surjection $p \to \a(p)$ from the uncountable set $\Q^{\om}$ to the set of infinite limit ordinals less than $\Omega$.

  We construct $T$ so that for each $p \in \Q^{\om}$, $T_p \cong A^{\a(p)}$.

  It follows that $X(T)$ is isomorphic to the sum $\sum_{p \in \Q^{\om}} \ X(A^{\a(p)})$. Using the sum map $\sum_{p \in \Q^{\om}} i_{\a(p)}$
  we obtain an injection from $X(T)$ into $\Q^{\om} \times J$. From Theorem \ref{theo8.10} we can embed $\Q^{\om}$ in $\R$ and so
  obtain an order injection from $X(T)$ into $\R \times J = \R_2$ which injects into $\R_{\om}$.

  It is clear that $T$ has branches of arbitrarily large countable height, but no uncountable branches.  Thus, $T$ is $\Omega$-bounded and
  with height $\Omega$.

  \end{proof}   \vspace{1cm}

\subsection{The Hart-van Mill Construction}

We conclude by describing the results of Hart and van Mill.

Let $F$ be a perfect Polish space, i.e. one with no isolated points, and so,
by Lemma \ref{lem8.7}(c), every nonempty open subset of $F$
is uncountable.  Following Hart and
van Mill, we call $Y \subset F$ a \emph{Bernstein subset},\index{Bernstein subset}\index{subset! Bernstein}\index{B set}
hereafter a \emph{B set}, if $Y$ meets every Cantor set in $F$, or,
equivalently, if its complement $X = F \setminus Y$ does not
contain a Cantor set.  (If, as in \cite{A}, we consider the Furstenberg family
generated by the Cantor sets of $F$, then the Bernstein sets are
the members of the dual family.) We call $Y$ a \emph{Bi-Bernstein
subset}, hereafter a \emph{BB set}, if both $Y$ and its complement
$X$ are Bernstein subsets.\index{Bi-Bernstein subset}\index{subset! Bi-Bernstein}\index{BB set}

We denote by ${\mathbf c}$ the cardinal number of  $\R$
and so of every Cantor set.  Recall that a
$G_{\delta}$ subset of $\R$ is a Polish space and so, if it is uncountable, it
 contains a Cantor set by Lemma \ref{lem8.7}.

\begin{lem}\label{lem8.12}  Let $F$ be a perfect Polish space.
\begin{enumerate}
\item[(a)] If $U$ is an open subset of $F$ and $x \in U$, then there exists a Cantor set
 $C$ such that $x \in C$ and $C \subset U$.
\item[(b)] Assume that $Y$ is a B set in $F$.  If $A$ is any uncountable, $G_{\delta}$
subset of $F$, then $Y \cap A$ has cardinality ${\mathbf c}$.  In particular, if $A$ is a
dense $G_{\delta}$ subset of $F$, then  $Y \cap A $
is dense in $F$.
\item[(c)] If $Y$ is a BB set in $F$, then the complement $F \setminus Y$ is a BB set in $F$.
\item[(d)] If $Y_{1}, Y_{2}$ are BB sets in $F$ and $Y_{1} \subset Y \subset Y_{2}$, then $Y$ is a BB set in $F$.

\end{enumerate}
\end{lem}

\begin{proof}(a): Let $U_{1}, U_{2},...$ be a
sequence of open subsets of $U$, each containing $x$ and with
diameter tending to zero.  By Lemma \ref{lem8.7}(a) each $U_{n}$ contains a
Cantor set $C_{n}$. Let $C = \{x\} \cup (\bigcup \{C_{n}\})$.

(b):  Let $C$ be a Cantor subset of $A$. There exists a
homeomorphism $s :C \times C \rightarrow C$ where $C \times C$ has
the usual product topology, ignoring the order structure.  For
each $x \in C$, $s(C \times \{x\})$ is a Cantor set which meets
$Y$. As $x$ varies over $C$ we obtain a pairwise disjoint family
of cardinality ${\mathbf c}$ which consists of nonempty subsets of
$Y \cap A$.  In particular, if $A$ is a dense $G_{\delta}$ and $U$
is a nonempty open set, then  $Y$  meets
$A \cap U$.

(c) and (d) are obvious.

\end{proof}\vspace{.5cm}

Of course, for us it is the special order results which are of
importance.  For $A \subset \R$ we define a subset of the
AS double $\R'$
\begin{equation}\label{eq8.42}
\R \vee A' \quad = \quad \R \times \{-1\}
\ \cup \ A \times \{+1\}.
\end{equation}
The gap pairs of $\R \vee A'$ are $x^{-} < x^{+}$  for $x
\in A$.

\begin{lem}\label{lem8.13}  Let $X$ be a LOTS, $Y$ be a B set in $\R$ and $A \subset \R$  disjoint from $Y$.
\begin{enumerate}
\item[(a)]  If $f : Y \rightarrow X$ is a order map, then the image $f(Y)$
 is countable or has cardinality ${\mathbf c}$.
\item[(b)]  If $f : \R \vee A' \rightarrow X$ is an order map such
that $f(Y \times \{-1\})$ is countable, then $A_{f} = \{ t \in A : f(t^{-}) < f(t^{+}) \}$
 is countable and the image $f(\R \vee A')$ is countable.
\end{enumerate}
\end{lem}

\begin{proof}(a):  Let $B = \{ x \in X : f^{-1}(x)
$ contains more than one point$\}$.  For $x \in B$, $f^{-1}(x)$
is a nontrivial interval in $Y$ and so $B$ is countable because
$Y$ is separable.  For each $x \in B$ let $I(x)$ be the smallest
closed interval in $\R$ which contains $f^{-1}(x)$, i.e.
the convex hull of the closure in $\R$.  Let $E_{1}$ be
the countable set of endpoints of the intervals $I(x)$ and let $E$
be the complement in $\R$ of union of the intervals
$I(x)$. Thus, $E$ is a $G_{\delta}$ subset of $\R$.  If
$E$ is countable, then $f(Y) = B \cup f(E_{1} \cap Y) \cup f(E \cap
Y)$ is countable.  If $E$ is uncountable, then by Lemma \ref{lem8.12}(b) $E \cap Y$
has cardinality ${\mathbf c}$.  Since $f$ is injective on $E \cap
Y$ the image has cardinality ${\mathbf c}$.

(b):    Identify $\R \setminus A$ with the subset
$(\R \setminus A) \times \{-1\} \subset \R \vee
A'$ so that $Y \subset \R\vee A'$.  For each $x \in f(Y)$
$Y \cap f^{-1}(x)$ is a nonempty interval in $Y$.
Since $Y$ is dense in $\R$, the closure in $\R$, $\overline{Y \cap f^{-1}(x)}$, is a
nonempty interval in $\R$ and its $\R$ interior is mapped by $f$ to $x$.
Notice that here we use that $f$ is order preserving rather than continuity of $f$, which is not assumed.

Let $F = \bigcup \{ \overline{Y \cap f^{-1}(x)} : x \in f(Y) \} $, where the closure is again
taken in  $\R$. Thus, $F$ is a countable union of
closed intervals in $\R$. Let $F_{1}$ be the countable
collection of endpoints of these intervals and let $F_{2} =
\R \setminus F$ which is countable because it is a
$G_{\delta}$ set disjoint from the B set $Y$. If $t \in
\overline{Y \cap f^{-1}(x)} \setminus F_{1}$ for some $ x \in
f(Y)$, then $f(t^{-}) =  x $ and if, in addition, $t \in A$, then
$f(t^{+}) =  x $.  Hence, $A_f \cup \{ t \in \R :
f(t^{-}) \not\in f(Y) \} \cup \{ t \in A : f(t^{+}) \not\in f(Y)
\} \  \subset \ F_{1} \cup F_{2}$ and so is countable.  Thus,
$A_f$ and the image of $f$ are countable.

\end{proof} \vspace{.5cm}

\begin{ex}\label{example3.3x} Products $\Z \times X$ and $\Q \times X$ with $X$ an IHLOTS. \end{ex}

If $X$ is any IHLOTS, we can choose an order injection $z : \Z \to \hat X \setminus X$ which is $\pm$cofinal in $\hat X$.
Using this we see that $\Z \times X \cong X$.

On the other hand, if $X$ is an IHLOTS dense in $\R$, then $\Q \times X$ is isomorphic to $\tilde X = ([0,1] \setminus C) \cap X$ with
$C$ the Cantor Set in $[0,1]$.

If $X$ is a dense Mycielski set, then $([0,1] \setminus C) \cap X$ is a Mycielski set dense in $[0,1] \setminus C$ and hence in $(0,1)$.
It then follows from Theorem \ref{theo8.10} that $X \cong \tilde X \cong \Q \times X$ and so $\Q \times X$ is an IHLOTS.

On the other hand, if $X$ is a BB-set, and the pair $a < b$ in $X$ is contained in a component of $[0,1] \setminus C$, then the interval
$(a,b) \cap \tilde X$ is a BB-set in the real interval $(a,b)$. If $a$ and $b$ lie in different components of $[0,1] \setminus C$, then
 the intersection of $C$ and the real interval $(a,b)$ contains a Cantor set.  Hence, $(a,b) \cap \tilde X$ is not a BB-set in $(a,b)$.
  It follows that $\Q \times X \cong \tilde X$ is not doubly transitive.  \vspace{.5cm}

\begin{df}\label{df8.14} Let $V \subset \R$ and $\mathcal{H}$ be a nonempty set
of subsets of $\R$. We say that $\mathcal{H}$ is a
\emph{Hart-van Mill collection with base set $V$}\index{Hart-van Mill collection}\index{Hart-van Mill collection!base set $V$}
when the following conditions hold.
\begin{enumerate}
\item[(i)]  $\Q \ \subset \ V$.
\item[(ii)]  Each $Y \in \mathcal{H} \cup \{V\} $ is a BB set, that is
$\mathcal{G}$ invariant, where $\mathcal{G}$ is the countable group of
positive, affine transformations of $\R$ with rational coefficients.
\item[(iii)] The elements of $\mathcal{H} \cup \{V\} $ are pairwise disjoint.
\item[(iv)]  If  $Y \in \mathcal{H}$, then $-Y = \{ -x : x \in Y \} $
is an element of $\mathcal{H}$ distinct from $Y$.
\item[(v)]   If  $f : \R \rightarrow \R $ is an order map
and $Y \in \mathcal{H}$ is such that the cardinality of $f(Y) \setminus Y $
is ${\mathbf c} $, then the cardinality of  $f(Y) \cap V$  is ${\mathbf c}$.
\end{enumerate}
By replacing $V$ by $V \cup -V$ we may assume that $V = -V$.
\end{df}
\vspace{.5cm}

\noindent{\bfseries Remark.} If $f : Y \rightarrow {\mathbb
R}$ is an order map with $Y \in \mathcal{H}$, then we can extend
$f$ to $\R$ by defining $\hat f(t) =  sup \ f((-\infty,t] \cap
Y)$.  The  order map $\hat f$ on $\R$ extends $f$ and so $\hat f(Y) = f(Y)$.
 So we can apply condition (v) even when $f$ is
only defined on $Y$. \vspace{.5cm}

If $Y \in \mathcal{H}$, then we let $X(Y) = \R \setminus
Y$.  More generally, if $\mathcal{J}$ is a nonempty subset of
$\mathcal{H}$ we call $X(\mathcal{J}) = \R \setminus \cup
\mathcal{J}$ \emph{the associated IHLOTS} for $\mathcal{J}$.  By
Proposition \ref{prop8.6} each $X(\mathcal{J})$ is an IHLOTS containing
$\Q$ with completion $\R$.
In addition, $V \subset X(\mathcal{J}) \subset \R \setminus Y$ for
$Y \in \mathcal{J}$ implies that $X(\mathcal{J})$ is a BB set by Lemma \ref{lem8.12} (c) and (d). We use $[-1,+1]$ as
the distinguished interval in each $X(\mathcal{J})$.\vspace{.5cm}

We will do some preliminary setup work which will be used twice.
It extends the notation of the proof of Theorem \ref{theo8.5}.

Assume that $ \mathcal{J} ,\mathcal{J}_{1} $ are subsets of a
Hart-van Mill collection $\mathcal{H}$ with base set $V$, that $Y
\in \mathcal{J} \setminus \mathcal{J}_{1}$ and that $\alpha, \beta
$ are infinite ordinals. Let $X = X(\mathcal{J}), X_{1} =
X(\mathcal{J}_{1})$ be the associated IHLOTS so that
\begin{equation}\label{eq8.43}
Y \quad \subset \quad  X_{1} \setminus X \qquad \mbox{and} \qquad
V \quad \subset \quad X \cap X_{1}.
\end{equation}
Assume that  $f: \widehat{X_{\alpha}} \rightarrow
\widehat{(X_{1})_{\beta}}$ is a continuous order map.

Since $Y$ is contained in the complement of $X$,  $Y \subset
\widehat{X_{\alpha}}$ consisting of elements of height $1$. Define for each $i < \beta$
\begin{equation}\label{eq8.44}
Y_{i} \quad = \quad \begin{cases} Y \qquad \mbox{for} \ i =
0 \\  \{ y \in Y : \hat{\pi}_{i}^{\b}(f(y)) \in (X_{1})_{i} \}
\quad \mbox{for} \ i > 0  \end{cases}
\end{equation}
When $0 < i $ and  $y \in Y_{i}$ let
\begin{equation}\label{eq8.45}
Q(y,i) \quad = \quad (\hat{\pi}_{i}^{\b} \circ
f)^{-1}(\hat{\pi}_{i}^{\b}(f(y))) = f^{-1}(J_{f(y)|i}) \quad
\subset \widehat{X_{\alpha}},
\end{equation}
where the latter equation uses definition (\ref{eq8.28}) above, which we now recall:

For $w \in (X_1)_i$, we define the compact subinterval $J_{w} \subset
\widehat{(X_1)_{\b}}$ and the map $\hat{\pi}_{w} : J_{w}
\rightarrow J$ by
\begin{equation}\label{eq8.28a}
\begin{split}
J_{w} \quad = \quad   (\hat{\pi}^{\b}_{\i})^{-1}(w) \ = \
 \{ z \in \widehat{(X_1)_{\b}} : z|i = w \} \ = \ [w-,w+]. \\
\hat{\pi}_{w}(z) \quad = \quad z(i). \hspace{5cm}
\end{split}
\end{equation}

As in proof of Theorem \ref{theo8.5}, for $y_{1}, y_{2} \in Y_{i}$
\begin{equation}\label{eq8.46}\begin{split}
Q(y_{1},i) \cap Q(y_{2},i) \  \not= \  \emptyset \quad \Rightarrow \hspace{2.5cm}\\
 f(y_{1})|i = f(y_{2})|i \quad \Rightarrow \hspace{3cm} \\ Q(y_{1},i) \
= \  Q(y_{2},i).\hspace{3cm}\end{split}
\end{equation}

Because $y \in Q(y,i)$, (\ref{eq8.34}) implies $\epsilon(Q(y,i)) = 0$ for
all $y \in Y_{i}$ and $0 < i < \beta $.  It follows from (\ref{eq8.33}),
(\ref{eq8.46}) and Lemma \ref{lem8.4}(c) that distinct members of the set of
intervals
\begin{equation}\label{eq8.47}
\mathcal{Q}_{i}\quad = \quad \{ span(Q(y,i)) : y \in Y_{i}
\} \ = \ \{ \hat{\pi}(Q(y,i)) : y \in Y_{i} \}
\end{equation}
are non-overlapping.  Since $y \in span(Q(y,i))$, we have
\begin{equation}\label{eq8.48}
\begin{split}
Y_{i} \cup O_{i} \ \subset \ \bigcup \mathcal{Q}_{i} \qquad \mbox{with}\\
O_{i} \quad = \quad  \bigcup \ \{ span^{\circ}(Q(y,i)) :y \in Y_i \}
\hspace{1cm}
\end{split}
\end{equation}
an open subset of $\R$.

\begin{lem}\label{lem8.15} Let $ 0 < i < \beta$.
\begin{enumerate}
\item[(a)] $Y \cap O_{i} \subset Y_{i} $.
\item[(b)] For all $ y  \in Y_{i} $ the closed interval $Q(y,i)$ is nontrivial.
\item[(c)] The set $\mathcal{Q}_{i}$ is countable.
\item[(d)] $O_{i}$ is a dense subset of $\bigcup \mathcal{Q}_{i}$.
\item[(e)] $Y_{i}$ is dense in $\R$  if and only if  $O_{i}$ is.
\item[(f)] The image projection $ \hat{\pi}(f(Y))$ is a countable subset of $\R$.
\item[(g)] If there exists $D$ a dense subset of $\widehat{X_{\alpha}}$ such
that $ \hat{\pi}(f(D))$ is a countable subset of $\R$, then
$ \hat{\pi} \circ f$ is  constant on $\widehat{X_{\alpha}}$.
\end{enumerate}
\end{lem}

\begin{proof}(a):  If $t \in Y \cap
span^{\circ}Q(y,i) $, then by Lemma \ref{lem8.4}(b) $\hat{\pi}^{-1}(t) = P(t) \subset Q(y,i) $. In
particular, by (\ref{eq8.45}) $f(t) \in J_{f(y)|i}$ and so by (\ref{eq8.28a})
$f(t)|i = f(y)|i \in (X_{1})_{i} $ which says that $t \in Y_{i}$.

(b):   Because $f$ is continuous and $\widehat{X_{\alpha}}$ is
connected, its image $F = f(\widehat{X_{\alpha}})$ is connected
and so is convex in $\widehat{(X_{1})_{\beta}}$. Let $y \in Y_{i}$
so that $f(y)|i \in (X_{1})_{i}$.\vspace{.25cm}

\textbf{Case 1:} Assume that there exist $a,b \in F$ such that $
\hat{\pi}_{i}^{\b}(a) <  f(y)|i  < \hat{\pi}_{i}^{\b}(b) $.
As defined by (\ref{eq8.28a}) the compact interval $J_{f(y)|i} =
(\hat{\pi}_{i}^{\b})^{-1}(f(y)|i) $ is nontrivial.  It is
entirely contained in the interval $(a,b)$ and so in the convex
set $F$.  Hence, the preimage $Q(y,i)$ is nontrivial. \vspace{.25cm}

\textbf{Case 2:} Assume $ f(y)|i \ = \  max \ \hat{\pi}_{i}^{\b}(F) $.
Let $ t \in Y \cap [y,\infty) $ which is an infinite set because
$Y$ is unbounded.  Since $\hat{\pi}_{i}^{\b} \circ f$ is an order map $f(y)|i =
\hat{\pi}_{i}^{\b}(f(t)) $ and so $t \in Q(y,i)$.  By a similar
argument $Q(y,i)$ is nontrivial when  $ f(y)|i \ = \  min \
\hat{\pi}_{i}^{\b}(F) $.\vspace{.25cm}

(c):   The intervals in $\mathcal{Q}_{i}$ are non-overlapping and by
(b) they are nontrivial. Since $\R$ is separable, the set
of intervals is countable.

(d):  The interior of a nontrivial interval is dense in the
interval.

(e):  Since $Y$ is a BB set it is dense in $\R$ and so $Y
\cap O_{i}$ is dense in $O_{i}$.  From part (a) it follows that if
$O_{i}$ is dense in $\R$, then $Y_{i}$ is.  On the other
hand, $Y_{i} \subset \bigcup \mathcal{Q}_{i}$ and so the union is
dense in $\R$ when $Y_{i}$ is.  From part (d) it then
follows that $O_{i}$ is dense in $\R$.

(f):  Define the order map $\tilde{f} = \hat{\pi} \circ f :
Y \rightarrow \R$. Notice that $Y_{1} =
(\tilde{f})^{-1}(X_{1})$ and for $y \in Y_{1}$ \ $Q(y,1) =
(\tilde{f})^{-1}(\tilde{f}(y))$.  By (c) the set $\mathcal{Q}_{1}$
is countable and so $X_{1} \cap \tilde{f}(Y)$ is countable. By
(\ref{eq8.43}) this set contains $(Y \cup V) \cap \tilde{f}(Y)$.  By
condition (v) of the Hart-van Mill collection, adjusted by the
remark after Definition \ref{df8.14}, it follows that $\tilde{f}(Y)
\setminus Y$ has cardinality less than ${\mathbf c}$.  Since $Y
\cap \tilde{f}(Y)$ is countable, the image $\tilde{f}(Y)$ has
cardinality less than ${\mathbf c}$.  So by Lemma \ref{lem8.13}(a),
$\tilde{f}(Y)$ is countable.

(g):  Assume that $a, b \in D$ with $ \hat{\pi}(f(a)) <
\hat{\pi}(f(b))$. We can choose a point $x \in  (X_{1}
\cap (\hat{\pi}(f(a)),\hat{\pi}(f(b))) \setminus  \hat{\pi}(f(D))$ since the image is
countable. The open interval $(x-,x+) \subset
\widehat{(X_{1})_{\beta}}$ is disjoint from $f(D)$.  Because $f$ is
continuous and $D$ is dense, $(x-,x+)$ is disjoint from
$f(\widehat{X_{\alpha}})$ which is connected.  This contradicts
$f(a) < x- < x+ < f(b)$. Hence, $ \hat{\pi} \circ f$ is  constant on $D$ and so on $\widehat{X_{\alpha}}$.

\end{proof}\vspace{.5cm}

\noindent{\bfseries Remark.} It is, of course, part (f) which really uses
the Hart-van Mill properties. \vspace{.5cm}

In preparation for what follows we define $\F $\index{$\F $} to be the set of order maps from $\R$ to $\R$, a pointwise closed semigroup of
real functions.

\begin{lem}\label{lem8.16a} (a) Any $f \in \F $ has only countably many discontinuities.

(b) The set $\F $ has cardinality ${\mathbf c} $. \end{lem}

\begin{proof} Let $D$ be a dense subset of $\R$. For an order map $g : D \to \R$, define $g_+,g_- \in \F$ by
\begin{equation}\label{eq8.48a}
\begin{split}
g_-(t) \ = \  sup \{ g(d) : d \in (-\infty,t)\cap D \}, \hspace{2cm} \\
g_+(t) \ = \ inf \{ g(d) : d \in (t,\infty)\cap D \}, \hspace{2cm}
\end{split}
\end{equation}
for all $t \in \R$.

For any $f : \R \to \R$, $f \in \F$  if and only if  $g = f|D$ is an order map and $g_- \leq f \leq g_+$. Furthermore, $f$ is discontinuous
at $t$  if and only if  $g_-(t) < g_+(t)$. By separability the family of nonempty intervals $(g_-(t),g_+(t)) $ is countable and so $f$ has at most
countably many discontinuities. If $f$ is continuous, then it is uniquely determined by $g$. If not, then there are ${\mathbf c} $
choices of $f$ between $g_-$ and $g_+$.

If $D$ is countable, then there are ${\mathbf c} $ maps $g$ from $D$ to $\R$.

It follows that the cardinality of $\F$ is ${\mathbf c} $.

\end{proof} \vspace{.5cm}

The amazing result of Hart and van Mill (1985) is the following.

\begin{theo}\label{theo8.16}
\begin{enumerate}
\item[(a)] There exists a Hart-van Mill collection $\mathcal{H}$ of cardinality
${\mathbf c}$.  \item[(b)] Let $\mathcal{H}$ be a Hart-van Mill collection with
base set $V$. Let $\mathcal{J}$ and $\mathcal{J}_{1}$ be subsets of  $\mathcal{H}$
with associated complementary IHLOTS $X$ and $X_{1}$, respectively. If $\mathcal{J}
\subset \mathcal{J}_{1}$, then $X_{1} \subset X$ and so for every ordinal $\alpha$,
$\widehat{(X_{1})_{\alpha}}$ injects into $\widehat{X_{\alpha}}$. If $\mathcal{J}$
is not a subset of  $\mathcal{J}_{1}$, then $\widehat{(X_{1})_{\omega}}$ does not
inject into $\widehat{X_{\omega}}$.  So if neither $\mathcal{J}$ nor $\mathcal{J}_{1}$
includes the other, then the CHLOTS $\widehat{X_{\omega}}$ and $\widehat{(X_{1})_{\omega}}$
 are not comparable with respect to size, i.e. neither injects into the other.
 In particular, if for some $Y \in \mathcal{J}$ we have $-Y \not\in \mathcal{J}$
, then the CHLOTS  $\widehat{X_{\omega}}$ is not even comparable in size with its reverse CHLOTS.
\end{enumerate}
\end{theo}

\begin{proof} (a): Let $\tilde{\mathcal{G}}$ \index{$\tilde{\mathcal{G}}$}
denote the group of nonconstant affine transformations of
$\R$ with rational coefficients, i.e. of the form $t
\mapsto at + b$ with $a \not= 0$ and $a,b$ rational. It contains
$\mathcal{G}$ as a subgroup of index two.  If $t \in {\mathbb I}$,
the set of irrationals, then $g \mapsto g(t)$ is an injective map
from $\tilde{\mathcal{G}}$ to $\R$. In particular, if $t
\in {\mathbb I}$, then the orbit sets $\mathcal{G}t$ and
$-\mathcal{G}t = \mathcal{G}(-t)$ are disjoint.

For $f \in \F$ let $S(f,\mathcal{A}) = \{ t \in
\R : f(t) \not\in \mathcal{A}t \}$ for $\mathcal{A} = \mathcal{G}$
or $\tilde{\mathcal{G}}$. Hart and van Mill call $f$
\emph{singular} if the set $f(S(f, \mathcal{G}))$ has cardinality
${\mathbf c}$. Observe that for each $g \in \mathcal{G}$ the
equation
\begin{equation}\label{eq8.49}
f(t) \quad = \quad -g(t)
\end{equation}
has at most one solution since $f$ is an order map and $-g$ is a
decreasing function.  It follows that the set $S(f, \mathcal{G})
\setminus S(f, \tilde{\mathcal{G}})$ is countable. So if $f$ is
singular the set $f(S(f,\tilde{\mathcal{G}}))$ has cardinality
${\mathbf c}$.

Let $\F' $ denote the set of singular functions in
$\F$. For example all translations by elements of
${\mathbb I}$ lie in $\F'$.  So from Lemma \ref{lem8.16a} it follows that the
cardinality of $\mathcal{F}'$ is ${\mathbf c}$.

For each  $f \in \F'$ choose a set $K(f) \subset S(f,
\tilde{\mathcal{G}})$ such that $f|K(f)$ is injective and $f(K(f))
= f(S(f,\tilde{\mathcal{G}}))$.  That is,
for each $z \in f(S(f,\tilde{\mathcal{G}}))$, we choose one point of $(f|S(f,\tilde{\mathcal{G}}))^{-1}(z)$.

 Let ${\mathbf c}$
denote the cardinal, i.e. the first ordinal with cardinality
${\mathbf c}$ and let $ \{f_{ij} : i,j < {\mathbf c}\} $ be a
listing of $\F'$ so that each function occurs ${\mathbf
c}$ times in each row and let $ \{ C_{ij} : i,j < {\mathbf c} \}$
be a similar listing of the Cantor sets in $\R$.

We will find points
$x(ij,0,\pm),x(ij,1,\pm),y(ij,0,\pm),y(ij,1,\pm)$ in
${\mathbb I}$ such that:
\begin{enumerate}
\item[(1)] $ x(ij,0,+), \  -x(ij,0,-) \ \in K(f_{ij}) \ $  and \\ $ y(ij,0,+) = f(x(ij,0,+))$,
$\quad y(ij,0,-) = f(-x(ij,0,-))$\\ so that $ y(ij,0,\pm) \in f(K(f_{ij})) = S(f_{ij},\tilde{\mathcal{G}}) $.
\item[(2)] $ x(ij,1,+), \  -x(ij,1,-), \ y(ij,1,\pm) \in C_{ij} $.
\item[(3)]  If $ (i,j,\alpha,\epsilon) \not= (i',j',\alpha',\epsilon')$ \\
with $i,j,i',j' \in {\mathbf c}; \alpha,\alpha' \in 2, \epsilon,\epsilon' \in \{+,-\}$,\\
then the four orbit sets $\tilde{\mathcal{G}}(x(ij,\alpha,\epsilon)),
\tilde{\mathcal{G}}(x(i'j',\alpha',\epsilon'))$,\\$\tilde{\mathcal{G}}(y(ij,\alpha,\epsilon)),
\tilde{\mathcal{G}}(y(i'j',\alpha',\epsilon'))$ are pairwise disjoint.
\end{enumerate}

To construct these points we choose a bijection from the index set
$ {\mathbf c} \times {\mathbf c}$ to ${\mathbf c}$ itself and so
well-order the index set with order type the cardinal ${\mathbf
c}$.
For index $ij$ we let $H(ij)$ be the set of rationals
together with the $\tilde{\mathcal{G}}$ orbit of all of the $x$
and $y$ points with index $i'j'$ preceding $ij$.  Thus, $H(ij)$
has cardinality less than ${\mathbf c}$. (Notice that for this reason
the lexicographic ordering on the product won't work for our purposes).

Since $C_{ij}$ has
cardinality ${\mathbf c}$ we can choose $ x(ij,1,+),-x(ij,1,-),$\\
$y(ij,1,\pm) \in C_{ij} \setminus H(ij) $ with distinct
$\tilde{\mathcal{G}}$ orbits and adjoin these four orbits to
$H(ij)$ to define $H(ij,1)$. Since $f_{ij}$ is singular the set
$K(f_{ij})$ is defined and has cardinality ${\mathbf c}$ and on it
$f_{ij}$ is injective. So we can choose $x(ij,0,+) \in   K(f_{ij})
\setminus (H(ij,1) \cup (f_{ij})^{-1}(H(ij,1))$.  By definition
the $\tilde{\mathcal{G}}$ orbits of $x(ij,0,+)$ and  $ y(ij,0,+) =
f(x(ij,0,+))$ are distinct.  Adjoin these two orbits to $H(ij,1)$
to define $H(ij,2)$.  Finally, choose $-x(ij,0,-) \in  K(f_{ij})
\setminus (H(ij,2) \cup (f_{ij})^{-1}(H(ij,2))$ and let $y(ij,0,-)
= f(-x(ij,0,-))$.  Notice that $x(ij,0,-)$ is in the
$\tilde{\mathcal{G}}$ orbit of $-x(ij,0,-)$ but not in the
$\mathcal{G}$ orbit since the point is irrational.

The members of family $\mathcal{H} =\{ Y_{i}^{\epsilon} : i
\in \mathbf c,\, \epsilon \in \{+,-,\} \}$ and the base set $V$
are given by
\begin{equation}\label{eq8.50}
\begin{split}
Y_{i}^{+} \ \  = \ \ \mathcal{G}\cdot \{x(ij,\alpha,\epsilon) :
j \in {\mathbf c}, \alpha \in 2, \epsilon \in \{+,-\} \}  \ \quad \mbox{for} \ i \in {\mathbf c}, \\
Y_{i}^{-} \ \  =  \ \mathcal{G} \cdot \{-x(ij,\alpha,\epsilon) : j \in {\mathbf c},
\alpha \in 2, \epsilon \in \{+,-\} \} \quad \mbox{for} \ i \in {\mathbf c},\\
V  \ \ \ =   \   \mathcal{G} \cdot (\{1\} \cup
\{y(ij,\alpha,\epsilon) : i, j \in {\mathbf c}, \alpha \in 2,
\epsilon \in \{+,-\} \}).\hspace{.55cm}
\end{split}
\end{equation}
Clearly, $Y_{i}^{-} = -Y_{i}^{+}$ and by condition (3) the indexed
family $\mathcal{H} \cup \{ V \} $ is pairwise disjoint.  By
condition (2) each member is a B set and so by disjointness  each
is a BB set.

Finally, suppose that $f \in \F$ and that for some
$Y \in \mathcal{H}$ the set $f(Y) \setminus Y$ has cardinality
${\mathbf c}$. Since $Y$ is $\mathcal{G}$ invariant, $ \{ t \in Y
: f(t) \not\in Y \} \subset S(f,\mathcal{G})$.  Hence, $f$ is
singular and so for each $i \in {\mathbf c}$ the set $Z(f,i)
= \{ j \in {\mathbf c} : f = f_{ij} \}$ has cardinality
${\mathbf c}$. If $ Y = Y_{i}^{\pm} $, then $ \{ y(ij,0,\pm) : j
\in Z(f,i) \}$ is a subset of $V \cap f(Y)$ of cardinality
${\mathbf c}$.  Thus, $\mathcal{H}$ is a Hart-van Mill collection
with base set $V$.

(b): If $\mathcal{J} \subset \mathcal{J}_{1}$, then by definition
$X_{1} \subset X$. So for every ordinal $\alpha$,
$\widehat{(X_{1})_{\alpha}}$ injects into $\widehat{X_{\alpha}}$
by Proposition \ref{prop4.6}(b),(c).

Now assume that $Y \in \mathcal{J} \setminus \mathcal{J}_{1}$.

Assuming that $ f : \widehat{X_{\omega}} \rightarrow
\widehat{(X_{1})_{\omega}}$ is a continuous order map, we will
show, following Hart and van Mill, that $f$ is a constant. By
Proposition \ref{prop4.4} this exactly says that $\widehat{(X_{1})_{\omega}}$ does
not inject into $\widehat{X_{\omega}}$.  We will apply the
preliminaries leading up to Lemma \ref{lem8.15} with $\alpha = \beta =
\omega$.

It suffices to show that each coordinate function $f_{n}$ is a
constant.  We begin with $f_{0} = \hat{\pi} \circ f :
\widehat{X_{\omega}} \rightarrow \R$.

For $0 < n < \omega$ let
\begin{equation}\label{eq8.51}
A_{n} \quad = \quad \{ w \in X_{n} : f_{0}(w-) < f_{0}(w+)
\}. \hspace{2cm}
\end{equation}
Observe that $\pi^{n}_{m}$ maps $A_{n}$ into $A_{m}$ if $0 < m
\leq n$, i.e. $w \in A_{n}$ implies $w|m \in A_{m}$, because $
(w|m)- \leq w- < w+ \leq (w|m)+ $.

Since $Y$ is dense in $\R$ it is easy to see from (\ref{eq8.23})
that the closure $\bar{Y}$ in $\widehat{X_{\omega}}$ satisfies
\begin{equation}\label{eq8.52}
\bar{Y} \quad \cong \quad \R \vee X' \hspace{4cm}
\end{equation}
in the notation of (\ref{eq8.42}).

From Lemma \ref{lem8.15}(f), $f_{0}(Y)$ is countable and so from Lemma \ref{lem8.13}(b)
$f_{0}(\overline{Y})$ and $A_{1}$ are countable sets.

Now for $0 < n < \omega$ and $w \in X_{n}$ let
\begin{equation}\label{eq8.53}
Y(w) \quad = \quad \{z \in J_{w} : z(n) \in Y \cap J \}.
\end{equation}

Recall that $J_{w}$ is the interval $ [w-,w+] $ in
$\widehat{X_{\omega}}$. We define the retraction $r : \R
\rightarrow J $ by mapping $(-\infty,-1]$ to $-1$ and
$[+1,\infty)$ to $+1$. Then define the order surjection  $r_{w} :
\widehat{X_{\omega}} \rightarrow J_{w} $ by
\begin{equation}\label{eq8.54}
r_{w}(z)_{i} \quad = \quad \begin{cases}  w_{i} \qquad \mbox{for} \ i < n  \\
r(z_{0}) \qquad \mbox{for} \ i = n  \\  z_{i-n} \qquad \mbox{for}
\ i > n. \end{cases}
\end{equation}

Apply Lemma \ref{lem8.15}(f) and then Lemma \ref{lem8.13}(b) to $f_{0} \circ r_{w}$ and
conclude that $f_{0}(\overline{Y(w)})$ and $A_{n+1} \cap J_{w}$
are countable sets.  Since $A_{n+1}$ projects to $A_{n}$ it
follows that $A_{n+1} = \bigcup \{ A_{n+1} \cap J_{w} : w \in
A_{n} \} $ and so by induction $A_{n}$ is countable for all $0 < n
< \omega$. Define
\begin{equation}\label{eq8.55}
\begin{split}
Y(\omega) \quad = \quad Y \cup \bigcup \{ Y(w) : w \in X_{n} \ \mbox{for some } 0 < n < \omega \}  \\
Z \quad = \quad \overline{Y} \cup \bigcup \{ \overline{Y(w)}
: w \in A_{n} \ \mbox{for some } 0 < n < \omega \}
\end{split}
\end{equation}
where the closures are taken in $\widehat{X_{\omega}}$.

The image
$f_{0}(Z)$ is countable. We now show that $f_{0}(Y(\omega))$ is a
subset of $f_{0}(Z)$ and so it is countable as well.

We show that $f_{0}(x) \in f_{0}(Z)$  if $x \in Y(w)$ with $w \in
X_{n}$.  If $w \in A_{n}$, then $x \in Z$.  Let $m = min \{ i : 0 <
i \leq n $ and $ w|i \not\in A_{i} \}$. If $i = 1$, then
$(w|i)-,(w|i)+  \in \overline{Y}$ and if $i > 1$, then $w|(i-1) \in
A_{i-1}$ and $(w|i)-,(w|i)+  \in \overline{Y(w|i-1)} $.
Furthermore, $f_{0}$ is constant on the interval $[(w|i)-,(w|i)+]$
which contains $x$. Thus, $f_{0}(x) = f_{0}((w|i)-) \in f(Z) $.

Because $Y(\omega)$ is dense in $\widehat{X_{\omega}}$ it follows
from Lemma \ref{lem8.15}(g) that $f_{0}$ is constant.

We complete the proof by using induction to show that $f_{n}$ is
constant for all $n$. If $f_{i}$ is a constant for all $i < n$
then define $\tilde{f} :  \widehat{X_{\omega}} \rightarrow
\widehat{(X_{1})_{\omega}}   $ by $\tilde{f}(x)_{j} =
f(x)_{j+n}$, i.e. just forget the first $n$ coordinates. Because
$f$ is constant on the first $n$ coordinates, this is a continuous
order map.  By the above initial step result the $\tilde{f}_{0} =
f_{n}$ is constant.

\end{proof} \vspace{.5cm}

\noindent{\bfseries Remark.} For every IHLOTS $X \subset \R$ the
CHLOTS $F = \widehat{X_{\omega}}$ has size between $\R$
and $\R_{\omega} $. It follows from Theorem \ref{theo8.3}(b) that
with $\alpha$ a  tail-like ordinal with $\alpha \geq
\omega^{\omega \cdot 2}$, then
\begin{equation}\label{eq8.56}
F_{\alpha} \quad \cong \quad \R_{\alpha}.  \hspace{2cm}
\end{equation}
In particular, for $F = \widehat{X(Y)_{\omega}}$ with  $Y  \in
\mathcal{H}$, $F$ is not symmetric but $F_{\alpha}$ is symmetric.
\vspace{.5cm}

Use the Axiom of Choice to select a subset $\mathcal{H}_{+}$ of
$\mathcal{H}$ so that for all $Y \in \mathcal{H}$ exactly one
member of the pair $\{Y,-Y\}$ lies in $\mathcal{H}_{+}$ and let
$\mathcal{H}_{-}$ be the complement. For every  $\mathcal{A}
\subset \mathcal{H}_{+}$  define $\mathcal{J}(\mathcal{A}) =
\mathcal{A} \cup (\mathcal{H}_{-} \setminus \{ -Y : Y \in
\mathcal{A} \}) $. If $\mathcal{A}_{1} \not= \mathcal{A}_{2}$, then
neither of the two sets
$\mathcal{J}(\mathcal{A}_{1}),\mathcal{J}(\mathcal{A}_{2})$
contains the other. So from the Hart-van Mill Theorem \ref{theo8.16} we
obtain a family of CHLOTS of cardinality $2^{{\mathbf c}}$ each with
size between  $\R$ and $\R_{\omega}$ no two of
which are comparable with respect to size.

Our final result combines the arguments of Hart and van Mill with
Theorem \ref{theo8.5}.

\begin{theo}\label{theo8.17} Assume that $\mathcal{H}$ is a Hart-van Mill collection
of subsets of $\R$. For distinct, nonempty subsets
$\mathcal{J}, \mathcal{J}_{1}$ of $\mathcal{H}$ let
$X = X(\mathcal{J})$ and $X_{1} = X(\mathcal{J}_{1})$ be the associated IHLOTS.
If $\alpha$ and $\beta$ are countable limit ordinals, then $\widehat{X_{\alpha}}$ is not
order isomorphic to $\widehat{(X_{1})_{\beta}}$. In particular, if for some
$Y \in \mathcal{J}$ we have $-Y \not\in \mathcal{J}$, then   $\widehat{X_{\alpha}}$ is not symmetric.
\end{theo}

\begin{proof}  Since the two subsets are
distinct, we can assume that $Y \in \mathcal{J}  \setminus
\mathcal{J}_{1}$.

We will assume that $f : \widehat{X_{\alpha}} \rightarrow
\widehat{(X_{1})_{\beta}}$ is an order isomorphism and derive a
contradiction by showing that $f$ is not injective. As before we
apply the preliminaries leading up to Lemma \ref{lem8.15}.

We prove by induction on $i < \beta$ that the open set $O_{i}$ is
dense in $ \R$. By Lemma \ref{lem8.15}(e) this is equivalent to
$Y_{i}$ being dense in $ \R$.\vspace{.25cm}

\textbf{Case 1:} For the initial step, $i = 1$, we apply Lemma \ref{lem8.15}(f) to see
that $\hat{\pi}(f(Y))$ is a countable subset of $\R$.
On $\R\setminus X_{1} $ the map $\hat{\pi}$
is injective and we have assumed that $f$ is injective. Hence,
$B_{\emptyset} = \{ y \in Y : f(y)(0) \not\in X_{1} \} $ is
countable. By the Baire Category Theorem and Lemma \ref{lem8.12}(b), $Y_{1} =
Y \setminus B_{\emptyset}$ is dense in $\R$.

For the case $ i = j+1 $ with $j > 0$ we first fix $y$ in the set
$Y_{j}$ and prove that $Y_{i} \cap span^{\circ}Q(y,j)$ is dense in
$span^{\circ}Q(y,j) $. It will, then follow that $Y_{i}$ is dense
in $O_{j}$ which is dense in $\R$ by inductive
hypothesis. Thus, $Y_{i}$ is dense in $\R$ in this case.

To analyze $Y_{j+1} \cap span^{\circ}Q(y,j)$ we use a variation of
the initial argument.  Let $r : \widehat{X_{\alpha}} \rightarrow
Q(y,j) $ be the canonical retraction.  That is, if $Q(y,j) =
[a,b]$ map $(-\infty,a]$ to $a$ and $[b.\infty)$ to $b$. Define
the order injection $s_{w} : J_{f(y)|j} \rightarrow
\widehat{(X_{1})_{\beta \setminus j}}$ by forgetting the
coordinates in $j$.  We apply Lemma \ref{lem8.15}(f) to $s_{w} \circ f \circ
r $. The analogue of $\hat{\pi} \circ f$ becomes, in this
case, $\hat{\pi}_{f(y)|j} \circ f \circ r$. Here $\hat{\pi}_{f(y)|j} = \hat{\pi}_{w}$ of (\ref{eq8.28a}) with $w = f(y)|j$.

It follows that $\hat{\pi}_{f(y)|j}(f(Y \cap span^{\circ}Q(y,j)))$ is a countable
subset of $J = [-1,+1]$. On  $ \{ z \in J_{f(y)|j} : z(j) \not\in
X_{1} \} $ the map $ \hat{\pi}_{f(y)|j} $ is injective and we have
assumed that $f$ is injective.  Hence, $B_{f(y)|j} = \{ t
\in Y \cap span^{\circ}Q(y,j): f(t)(j) \not\in X_{1} \} $ is
countable. For all $t \in Y \cap spanQ(y,j)$ we have $f(t)|j =
f(y)|j $.  Hence, $Y_{j+1} \cap span^{\circ}Q(y,j) = Y \cap
(span^{\circ}Q(y,j) \setminus B_{f(y)|j}) $. As before the Baire
Category Theorem and Lemma \ref{lem8.12}(b) imply that $Y_{j+1} \cap
span^{\circ}Q(y,j)$ is dense in $span^{\circ}Q(y,j)$.\vspace{.25cm}

\textbf{Case 2:} When $i$ is a limit ordinal Lemma \ref{lem8.15}(a) implies that
\begin{equation}\label{eq8.57}
Y_{i} \quad = \quad \bigcap \{ Y_{j} : j < i \} \quad \supset
\quad (\bigcap \{ O_{j} : j < i \}) \cap Y.
\end{equation}
Because $\beta$ is countable, $\bigcap \{ O_{j} : j < i \}$ is a
$G_{\delta}$ set which is dense by induction hypothesis and so the
intersection with $Y$ is dense by Lemma \ref{lem8.12}(b).  Thus, $Y_{i}$ is
dense in this case as well. \vspace{.25cm}

Having completed the induction we see, as in Case 1, that
$(\bigcap \{ O_{j} : j < \beta \}) \cap X$ is dense in $\R$.  As in the final portion of the proof of Theorem \ref{theo8.5}, if $y$
is in this set, then $P(y)$ is nontrivial and $f$ is constant on
$P(y)$.  Hence, $f$ is not injective, contradicting our initial
assumption.

\end{proof}\vspace{.5cm}

 It follows from Corollary \ref{cor8.3a} that the tower $\{ F_{\om^{\g}} \}$ above each
$F = \widehat{ X(\mathcal{J})_{\a} }$ coincides at a sufficiently high level with the tower above $\R$.

On the other hand,
from Corollary \ref{cornew03} we immediately obtain the following.

\begin{theo}\label{theo8.18} Assume that $\mathcal{H}$ is a Hart-van Mill collection
of subsets of $\R$. For not necessarily distinct, nonempty subsets
$\mathcal{J}, \mathcal{J}_{1}$ of $\mathcal{H}$ let
$X = X(\mathcal{J})$ and $X_{1} = X(\mathcal{J}_{1})$ be the associated IHLOTS.
If $\a, \b$ are infinite tail-like ordinals with $\a > \om \cdot \b$, then
there does not exist an order injection from $\widehat{(X_1)_{\a}}$ into
$\widehat{(X_2)_{\b}}$. \end{theo} \vspace{.5cm}

\begin{cor}\label{cor8.19} Assume that $\mathcal{H}$ is a Hart-van Mill collection
of subsets of $\R$. For a subset
$\mathcal{J}$ of $\mathcal{H}$ let
$X = X(\mathcal{J})$. The transfinite sequence of CHLOTS $\{ \widehat{(X)_{\om^{\g}}}, \ 0 < \g < \Omega \}$
 and the transfinite sequence of CHLOTS Cantor Spaces $\{ C(\widehat{(X)_{\om^{\g}}}), \ 0 < \g < \Omega \}$
are nondecreasing in size. . Furthermore, if $\g_1 > 1 + \g_2$, then
$\widehat{(X)_{\om^{\g_1}}}$ is strictly bigger in size than
$\widehat{(X)_{\om^{\g_2}}}$ and $C(\widehat{(X)_{\om^{\g_1}}})$ is strictly bigger in size than
$C(\widehat{(X)_{\om^{\g_2}}})$. \end{cor}

\begin{proof} This follows from Theorem \ref{theonew04} with $\d = 1$ and so with $\g_0 = 0$.

\end{proof} \vspace{.5cm}

\noindent{\bfseries Remark:} If $\g_2 \geq \g_3 = \om$, then $1 + \g_2 = \g_2$. Hence,
after the $\om^{\om}$ level, the members of the towers are strictly increasing in size. \vspace{.5cm}

Using a Hart-van Mill collection of cardinality ${\mathbf c}$ we
obtain a collection of cardinality $2^{{\mathbf c}}$ consisting of
CHLOTS with size between $\R$ and $\R_{\omega}$ each forming the base of
a tower of CHLOTS  of height $\Omega$. The separate towers do not intersect by Theorem \ref{theo8.17}. These
towers are disjoint from the tower $\{ \R_{\a} \}$ over $\R$  as well by Theorem \ref{theo8.5}.
The corresponding Cantor Space towers do not intersect either by (\ref{eq3.37a}). \vspace{1cm}

\section{\textbf{Zero Dimensional LOTS}}

If $J$ is a nonempty, bounded, clopen convex subset of a complete LOTS $X$, then $y- = sup \ J$ is an element of $J$
and $y+ = inf \ (y-, \infty) \not\in J$. So $y- < y+$ is a gap pair in $X$. Similarly, there is a gap pair $x- < x+$
with $x- \not\in J$ and $x+ =  \in J$. Hence, $J$ is the clopen interval $[x+,y-]$.
It follows that $X$ is zero-dimensional, i.e. the clopen convex sets form a base for $X$, if and only if $X$ is
\emph{gap pair dense} \index{gap pair dense} \index{LOTS!gap pair dense} that is, for every $x < y \in X$ there
exists a gap pair $z- < z+ \in X$ with $x \le z- < z+ \le y$.

Assume that $X$ is a perfect, complete LOTS, i.e. it has no isolated points. We
obtain the quotient space $F$ by identifying
each gap pair with a single point.  That is, we use the equivalence relation
\begin{equation}\label{eqZD00}
\{ (x,y) \in X \times X : x = y \ \text{ or} \ \{ x, y \} \ \text{ is a gap pair in } \ X \}.
\end{equation}
 On $F$ there is a unique ordering so that $\pi : X \to F$
is an order surjection which is continuous and
topologically proper by Proposition \ref{prop2.1}(a). $F$ is complete because
$X$ is and if $\pi(x) < \pi(y)$ in $F$, then
$\{ x, y \}$ is not a gap pair and so the interval $(x,y)$ is infinite.
It follows that $F$ is order dense and so is connected.
Let $A \subset F$ consist of the classes of the gap pairs. We call $F$ the
\emph{connected quotient} \index{connected quotient} of $X$
and $A$ the \emph{gap pair set} \index{gap pair set}. If $F$ has extrema,
we let $F^{\circ}$ denote $F$ with the max and min removed
if they exist.

Conversely, if $F$ is a connected LOTS and  $A \subset F$ we extend \ref{eq8.42}
by defining a subset of the AS double $F'$
\begin{equation}\label{eqZD01}
F \vee A' \quad = \quad F \times \{-1\}
\ \cup \ A \times \{+1\}.
\end{equation}
\index{$F \vee A'$}The gap pairs of $F \vee A'$ are $a- < a+$  for $a \in A$.
It is clear that we can identify $F$ with the connected quotient of
$F \vee A'$ and $A$ with its gap pair set.
On the other hand, if $X$ is a perfect, complete LOTS with connected quotient
$F$ and gap pair set $A$, then $X$ is isomorphic to $F \vee A'$.

In the case of the AS double $F'$ itself the gap pair set $A$ is all of $F$.

\begin{prop}\label{propZD01} Let $X$ be a perfect, complete LOTS with connected quotient $F$ and gap pair set $A$.
\begin{itemize}
\item[(a)] $X$ is compact or first countable if and only if $F$ satisfies the corresponding property.

\item[(b)] The LOTS $F^{\circ}$ is unbounded.

\item[(c)]  $X$ is gap pair dense if and only if  $A$ is dense in $F$. In that case
the inclusion of the dense set $A$ into $F$ is an
embedding and $A$ is unbounded. The induced map from the completion $\hat A$ to $F$
is an isomorphism onto $F^{\circ}$.

\item[(d)] If $f$ is an order automorphism of $X$, then there is a unique order
automorphism $g$ of $F$ such that $\pi \circ f = g \circ \pi$.
The induced automorphism satisfies $g(A) = A$. Conversely, if $g$ is an order
automorphism of $F$ such that $g(A) = A$, then it is induced
from a unique automorphism $f$ of $X$.
\end{itemize} \end{prop}

\begin{proof}  (a): The compactness result is clear because $\pi$ is topologically
proper. Every point of $F$ (except the maximum if any) is the limit
of a decreasing sequence if and only if every point of $X$ which is not a left
end-point is the limit of a decreasing sequence.  Similarly
for increasing sequences, and these two observations yield the first countability equivalence.

(b):  No extreme point of $X$ can be part of a gap pair because $X$ has no isolated points.
So if we remove the max and min, if any, from
$X$ and $F$, the resulting LOTS are unbounded.

(c): The density equivalence is clear and the embedding result follows from
Proposition \ref{prop2.1}(b).

An extreme point of $X$ cannot be part of a gap pair
because $X$ has no isolated points. It follows that
the dense set $A$ is contained in $F^{\circ}$ and so is unbounded. Because $F^{\circ}$
is connected and unbounded, it is the completion of
the dense subset $A$.

(d): An automorphism $f$ on $X$ maps gap pairs to gap pairs and so induces $g$.
Conversely, given $g$, $f$ can be regarded as the restriction to
$F \vee A'$ of the automorphism $g'$ of $F'$.

\end{proof} \vspace{.5cm}

\begin{df}\label{dfZD02} Let $X$ be a perfect, complete, zero-dimensional LOTS. We say that $X$
satisfies the \emph{clopen interval condition}
\index{clopen interval condition} when any two bounded clopen intervals are isomorphic. \end{df} \vspace{.5cm}

\begin{theo}\label{theoZD03} Let $X$ be a perfect, complete, zero-dimensional LOTS.
\begin{itemize}
\item[(a)] $X$ satisfies the clopen interval condition if and only if $A$ is doubly transitive.

\item[(b)] If $X$ is first countable and $\s$-bounded, then it satisfies the clopen interval
condition if and only if $A$ is a
HLOTS
in which case $F^{\circ}$ is a CHLOTS. In addition, $X$ is then topologically homogeneous.

\item[(c)] If $X$ is  topologically homogeneous, then it is first countable.
\end{itemize}
\end{theo}

\begin{proof} (a): Let $a < b, c < d$ in the gap pair set $A$. Because $A$ is unbounded we
can choose $e_1 < min(a,c) < max(b,d) <e_2$ in $A$.

Assume $X$ satisfies the clopen interval condition. Combine isomorphisms
\begin{equation}\label{eqZD02}
[e_1+,a-] \cong [e_1+,c-], [a+,b-] \cong [c+,d-], [b+,e_2-] \cong [d+,e_2-]
\end{equation}
with the identity on $(-\infty,e_1-] \cup [e_2+,\infty)$ to obtain an automorphism $f$ of $X$.
By Proposition \ref{propZD01}(d) $f$
 induces $g$ on $F$ which preserves $A$ and maps
$[a,b]$  to $[c,d]$. Hence, $A$ is doubly transitive.

Conversely, assume that $A$ is doubly transitive.  There exists an automorphism $h$ of $A$
which maps $[a,b]$ in $A$ to $[c,d]$ in A.
The completion $\hat h$ is an automorphism of $F^{\circ}$ which extends to an
automorphism $g$ of $F$ such that $g(A) = A$.
By Proposition \ref{propZD01}(d) again $g$ is induced by an automorphism $f$ on $X$
which clearly maps $[a+,b-]$ to $[c+,d-]$.

(b): If $X$ is first countable and $\s$-bounded as well as complete, then it is of countable type by
Proposition \ref{prop2.5}(d)  and so $F$ and $A$ are of
countable type by Proposition \ref{prop2.5} (f) and (g). So $A$ is an IHLOTS if and
only if it is doubly transitive by
Proposition \ref{prop3.4}(a) in which case the completion $F^{\circ}$ is a
CHLOTS by Proposition \ref{prop3.4}(c).

The topological homogeneity follows from the clopen interval condition with
various cases. We provide just a sample.

Suppose $a \in A, b \in F^{\circ} \setminus A$.
Let $e_1, e_2 \in A$ with $e_1 < a, b < e_2$. Choose  decreasing sequences
$\{a_n \}$ in $A \cap (a,e_2)$ converging to $a$,
$\{b_n \}$ in $A \cap (b,e_2)$ converging to $b$ and an increasing sequence
$\{ c_n \}$ in $A \cap (e_1,b)$ converging to $b$.
Define
\begin{itemize}
\item $J_0 = [e_1+,a-], J_1 = [a_1+,e_2-]$ and \\ $J_k = [a_k+,a_{k-1}-]$ for $k \ge 2$.
\item $K_0 = [e_1+,c_1-]$ and $K_{2k} = [c_{k}+,c_{k+1}-]$ for $k \ge 1$.
\item $K_1 = [b_1+,e_2-]$ and $K_{2k+1} = [b_{k+1}+,b_{k}-]$ for $k \ge 1$.
\end{itemize}
Choose isomorphisms $J_k \cong K_k$ for $k \ge 0$. Together these extend to a
homeomorphism of $[e_1+,e_2-]$ which maps
$a+$ to $b$. Use the identity on the complementary set.

(c): A bounded sequence of distinct points in $X$ converges to a point by completeness.  Topological homogeneity then implies that
every point $x \in X$ is the limit of some sequence in $X \setminus \{ x \}$. If $x = a-$ is a left endpoint, then such a sequence must
consist of points below $x$. By going to a subsequence we can assume the sequence $\{ y_n \}$ converging to $x$ is increasing.
Thus, $\{ (y_n, a+) \}$ is a neighborhood base for $x$. By topological homogeneity it follows that every point has a
countable neighborhood base, i.e. $X$ is first countable.

\end{proof} \vspace{.5cm}

\noindent {\bfseries Remark:} The topological homogeneity argument in (b) is due to Maurice \cite{M2}. \vspace{.5cm}

\begin{cor}\label{corZD03a} If $A$ is any doubly transitive LOTS with completion $F = \hat A$, then $X = F \vee A'$
is a perfect, complete, zero-dimensional LOTS which satisfies the clopen interval condition. Furthermore, if
$a < b \in A$, then the interval $[a+, b-]$, clopen in $X$, is a perfect, compact, zero-dimensional LOTS
which satisfies the clopen interval condition. \end{cor}

\begin{proof}  This is clear from Theorem \ref{theoZD03} (a).

\end{proof} \vspace{.5cm}

We turn now to trees of type $2 = \{ 0, 1 \}$. Any tree $T$ of type $2$ is isomorphic with a subtree of the simple tree on
$2, \a$ where $\a = h(T)$ and so we will restrict attention to such subtrees so that $p \in T$ with $o(p) = \b$ is an element
of $2^{\b}$. In particular, if $T$ is normal with $h(T)$ a limit ordinal, then the height of every branch is a
limit ordinal and as in \ref{eq6.28} we can identify the branch space
$X(T)$ with
\begin{equation}\label{eqZD03}
\begin{split}
\{x \in 2^{\beta} : \beta \  \mbox{is an infinite limit ordinal},
x \not\in T, \\ \mbox{and} \ x|\epsilon \in T \ \mbox{for all} \
\epsilon < \beta \}. \hspace{2cm}
\end{split}\end{equation}
For $x, y \in X(T)$ we have $x < y$ if there exists $\ep = \ep(x,y) < h(x),h(y)$ such that $x_i = y_i $ for all $i < \ep$,
$x_{\ep} = 0 $ and $y_{\ep} = 1$.

\begin{lem}\label{lemZD04} Let $T$ be a normal tree of $2$ type with height a limit ordinal. The pair $x < y \in X(T)$ is a
gap pair if and only if there exists $\ep < h(x), h(y)$ such that
\begin{align}\label{eqZD04}
\begin{split}
x_i = y_i \ \text{for all} \ i < \ep, \ \ &x_{\ep} = 0, y_{\ep} = 1, \ \text{and} \\
x_j = 1 \ \text{for all} \  \ep < j < h(x), \ \ &y_j = 0 \ \text{for all} \  \ep < j < h(y).
\end{split}\end{align}
\end{lem}

\begin{proof} It is clear that if $x, y$ satisfy (\ref{eqZD04}), then the interval $(x,y)$ is empty.
On the other hand, if for some $ \ep < j < h(x)$ $x_j = 0$, then
there exists $z \in X(T)$ with $z_i = x_i$ for $i < j$ and $z_j = 1$.  So $x < z < y$.  Similarly, if
for some $ \ep < j < h(y)$ $y_j = 1$ there exists $z$ with $x < z < y$.

\end{proof} \vspace{.5cm}

We say that a branch $x$ \emph{eventually equals $0$}\index{eventually equals $0$} if there exists $\b < h(x)$ such that
$x_i = 0$ for all $i \ge \b$. In that case, we let $\b^*(x)$ be the minimum of such ordinals $\b$. Similarly, we say that
$x$ \emph{eventually equals $1$} if there exists $\b < h(x)$ such that
$x_i = 1$ for all $i \ge \b$. In that case, we let $\b^*(x)$ be the minimum of such ordinals $\b$.  From Lemma \ref{lemZD04}
we see that $x < y$ is a gap pair if and only if $x$  eventually equals $1$, $y$ eventually equals $0$ and $\b^*(x) = \b^*(y) = \ep + 1$.
Thus, $x$ is a member of a gap pair if and only if it  eventually equals $0$ or eventually equals $1$ and, in addition,
 $\b^*(x)$ is a successor ordinal.

 In particular, if we define $\bar 0, \bar 1 \in 2^{h(T)}$ by $\bar 0_i = 0, \bar 1_i = 1$ for all $i < h(T)$, then
 by normality and (\ref{eqZD03}) there are unique ordinals $\g_0, \g_1 \le h(T)$ such that $\bar 0|\g_0, \bar 1|\g_1 \in X(T)$.
 Somewhat abusively, we will denote these branches $\bar 0$ and $\bar 1$ so that $h(\bar 0) = \g_0, h(\bar 1) = \g_1$.
 It is clear that $\bar 0$ is the minimum element of $X(T)$ and $\bar 1$ is the maximum.

\begin{theo}\label{theoZD05} If $T$ is a normal tree of type $2$ with height a limit ordinal, then the branch space
$X(T)$ is a perfect, compact, zero-dimensional LOTS. If $T$ is $\Omega$-bounded, then $X(T)$ is first countable. \end{theo}

\begin{proof} By Proposition \ref{prop5.5a} $X(T)$ is complete. Since $2$ is bounded, $X(T)$ is bounded and so is compact.
By Proposition \ref{prop5.9} it is of countable type, and so is first countable, if $T$ is $\Omega$-bounded. From (\ref{eqZD04})
it is clear that no left end-point is a right end-point and so there are no isolated points.

Now suppose that $x < z < y$ in $X(T)$. With $\ep = \ep(x,y)$ we have that $z_i = x_i = y_i$ for all $i < \ep$. Suppose that
$z_{\ep} = 0$. Since $x < z$, there exists $k$ with $\ep < k < h(x), h(z)$ such that $x_k = 0$ and $z_k = 1$. Let
$w \in 2^{h(T)}$ with $w_i = z_i$ for all $i \le k$ and $w_j = 0$ for all $j > k$. By normality and (\ref{eqZD03})
there is a limit ordinal $\g$ with $ k < \g$ such that $w|\g \in X(T)$. Thus, $w|\g$ is a right end-point which lies
between $x$ and $y$ as does its associated left end-point. We proceed similarly if $z_{\ep} = 1$.

Thus, $X(T)$ is gap pair dense and so is zero-dimensional.

\end{proof} \vspace{.5cm}

\begin{cor}\label{theoZD05a} If $\a$ is a limit ordinal, then
$2^{\a}$ is a perfect, compact, zero-dimensional LOTS. If $\a$ is countable, then $2^{\a}$ is first countable. \end{cor}

 \begin{proof} Apply Theorem \ref{theoZD05} to the $2, \a$ simple tree.

 \end{proof}\vspace{.5cm}

\begin{lem}\label{lemZD06} Let $T$ be an additive tree of type $2$ with height the tail-like ordinal $\a$.
Assume that the branches $\bar 0, \bar 1 \in X(T)$ have height $\a$.
\begin{itemize}
\item[(i)] If $x \in X(T)$ eventually equals $0$, then there is an isomorphism from $[x,\bar 1]$ to $X(T) = [\bar 0,\bar 1]$.
\item[(ii)] If $x \in X(T)$ eventually equals $1$, then there is an isomorphism from $[\bar 0, x]$ to $X(T) = [\bar 0,\bar 1]$.
\end{itemize}
\end{lem}

\begin{proof} Observe first that additivity implies that every branch which eventually equals $0$ or eventually equals $1$
has height $\a$.

Assume $x$ eventually equals $0$ and $\b = \b^*(x)$. Let $K = \{ \b \} \cup \{ k : x_k = 0 $ and $ k < \b \}$ and let
$r$ be an isomorphism from the well-ordered set $K \subset \a$ onto an ordinal $\g + 1  \le \b + 1 < \a$, so that $r(\b) = \g$.

Define $p^{\b} = x|\b$ and for $k \in K$ with $k < \b$, $p^k_i = x_i$ for $i < k$ and $p^k_k = 1$. Thus, $p^k \in T$ for all $k \in K$
with $o(p^{\b}) = \b$, $o(p^k) = k+1$ for $k < \b$. Notice that if $\b = \g + 1$, then $x_{\g} = 1$ by definition of $\b^*$.

Define $q^{\b} = \bar 0|r(\b)$ and for $k \in K$ with $k < \b$, $q^k_i = 0$ for $i < r(k)$ and $q^k_{r(k)} = 1$. Notice that $r(\b) < \a = h(\bar 0)$.
This is where we need $\a = h(\bar 0)$. If, for example, it happened that $h(\bar 0) < r(\b)$, then we could not define all the elements $q^k$.

We show that
\begin{equation}\label{eqZD05}
[x,\bar 1] = \bigcup_{k \in K} X(T_{p^k}), \quad \text{and} \quad [\bar 0, \bar 1] = \bigcup_{k \in K} X(T_{q^k}).
\end{equation}

Note first that $x \in X(T_{p^{\b}})$. Now assume $x < y$, and let $\ep = \ep(x,y)$.

If $\ep \ge \b$, then again  $y \in X(T_{p^{\b}})$. If $\ep < \b$, then $x_{\ep} = 0$ and $y_{\ep} = 1$. So $\ep \in K$ and $y \in T_{p^{\ep}}$.

Next $\bar 0 \in X(T_{q^{\b}})$. Now assume $\bar 0 < y$ and let $\ep = \ep(\bar 0,y)$.

If $\ep \ge r(\b)$, then again $ y \in X(T_{q^{\b}})$. If $\ep < r(\b)$, then $y_{\ep} = 1$ and $\ep = r(k)$ for some $k \in K$ with $k < \b$.
So $y \in X(T_{q^{k}})$.

Each of the unions in (\ref{eqZD05}) is disjoint. The isomorphism is defined, using additivity, by:
\begin{equation}\label{eqZD06}
p^k + z \mapsto q^k + z \quad \text{for all} \ \ k \in K, z \in X(T).
\end{equation}

The proof of (ii) is completely analogous.

\end{proof} \vspace{.5cm}

\begin{theo}\label{theoZD07} If $T$ is an additive tree of type $2$ with height the tail-like ordinal $\a$ and
 the branches $\bar 0, \bar 1 \in X(T)$ have height $\a$, then $X(T)$ is a compact, perfect, zero-dimensional LOTS which
 satisfies the clopen interval condition. \end{theo}

 \begin{proof} Assume $x < y \in X(T)$ with $x$ eventually equal to $0$ and $y$ eventually equal to $1$. We prove that the interval
 $[x,y]$ is isomorphic to $X(T)$. This includes the case when $[x,y]$ is a clopen interval.

 Let $\ep = \ep(x,y)$ so that $x_{\ep} = 0$ and $y_{\ep} = 1$. Let $a = x|(\ep + 1), b = y|(\ep + 1)$.

 The truncations $\t_{\ep +1}(x), \t_{\ep +1}(y)$ are, respectively, eventually equal to $0$ and to $1$. By Lemma \ref{lemZD06} there are
 isomorphisms $f_0 : [\t_{\ep +1}(x),\bar 1] \to [\bar 0, \bar 1]$ and $f_1 : [\bar 0, \t_{\ep +1}y)] \to [\bar 0, \bar 1]$.

 Let $0, 1$ denote
 the elements of level 1 of the tree. Define $f : [x,y] \to X(T)$ by
 \begin{equation}\label{eqZD07}
 f(z) = \begin{cases} 0 + f_0(\t_{\ep +1}(z))  \ if \ z_{\ep} = 0,\\
 1 + f_1(\t_{\ep +1}(z)) \ if \ z_{\ep} = 1. \end{cases}
 \end{equation}
 Notice that for $z \in [x,y]$ if $z_{\ep} = 0$, then $z|(\ep +1) = a$ and if $z_{\ep} = 1$ then $z|(\ep +1) = b$.
 The result follows because $T$ is the disjoint union of $T_0$ and $T_1$.

 \end{proof} \vspace{.5cm}

  As a corollary we obtain the following extension of a theorem of Maurice \cite{M}, who proved the result
  when $\a$ is countable.

  \begin{cor}\label{corZD08} If $\a$ is an infinite tail-like ordinal,
   then $2^{\a}$ is a compact, perfect, symmetric, zero-dimensional LOTS which
 satisfies the clopen interval condition. \end{cor}

 \begin{proof} Apply Theorem \ref{theoZD07} to the $2, \a$ simple tree. Interchanging $0$ and $1$, we obtain an order* automorphism
 of $2^{\a}$ and so it is symmetric.

 \end{proof}\vspace{.5cm}

 In general if $\a$ is a limit ordinal, let $\pi_{\a} : 2^{\a} \to F_{\a}$ be the projection to the
 connected quotient of the zero-dimensional LOTS $2^{\a}$. As in Corollary \ref{corZD08} $2^{\a}$ is
 symmetric and so the quotient $F_{\a}$ is symmetric as well.

 \begin{theo}\label{theoZD09} If $\a > \b$ are limit ordinals, then $2^{\a}$ is bigger than $2^{\b}$ and
 $F_{\a}$ is bigger than $F_{\b}$. In particular, $2^{\a}$ and $2^{\b}$ are not isomorphic.
 The connected quotients $F_{\a}$ and $F_{\b}$ are not homeomorphic. \end{theo}

 \begin{proof} From the order surjection $\pi_{\a}$ we obtain an order injection $i_{\a} : F_{\a} \to 2^{\a}$.
 Explicitly we map each $a \in A$ to the left end-point $a-$ of the pair. This induces the injection
 $i_{\a}' : F_{\a}' \to (2^{\a})' = 2^{\a + 1}$.

 Assume that $h : 2^{\b} \to 2^{\a}$ is an order injection. If $a < b $ in $2^{\b}$ and
 $h(a) < h(b)$ is a gap pair in $2^{\a}$, then  $a < b$ is a gap pair in $2^{\b}$. Otherwise, there exists $c$ with $a < c < b$
 and so $h(a) < h(c) < h(b)$. It follows that the order map $\pi_{\a} \circ h \circ i_{\b} : F_{\b} \to F_{\a}$ is injective.

 Now choose $z \in 2^{\a \setminus \b}$ with $z$ not eventually equal to $0$ or eventually equal to $1$. Let
 $z_0 = 0 + z, z_1 = 1 + z$, with $0, 1$ here regarded as elements of order $1$ in the simple tree.
 Since $\a > \b$, $\a \setminus \b$ is a limit ordinal and so $z_0, z_1 \in 2^{\a \setminus \b}$
 with neither eventually equal to $0$ or $1$.

 Define the order injection $f : (2^{\b})' \to 2^{\a}$ by
  \begin{equation}\label{eqZD08}
  f(a-) \ = \ a + z_0, \quad \text{and} \quad f(a+) = a + z_1.
  \end{equation}
  By the choice of $z$ it follows that the image of $f$ is disjoint from all of the gap pairs in $2^{\a}$.
  Hence, $\pi_{\a} \circ f : (2^{\b})' \to F_{\a}$ is injective.

  Because $2^{\b}$ and $F_{\b}$ are compact LOTS both are order simple by Corollary \ref{cor4.10}.

  From the projection $\pi_{\b}^{\a} : 2^{\a} \to 2^{\b}$ we obtain an order injection $h : 2^{\b} \to 2^{\a}$ and
  as noted above, $\pi_{\a} \circ h \circ i_{\b} : F_{\b} \to F_{\a}$ is an order injection. Thus, $2^{\a}$ is at least as big as
  $2^{\b}$ and $F_{\a}$ is at least as big as   $F_{\b}$.

 If there were an order injection $q : 2^{\a} \to 2^{\b}$, then $q \circ f : (2^{\b})' \to 2^{\b}$ would be an order injection,
 violating the order simplicity of $2^{\b}$.

 If there were an order injection $q : F_{\a} \to F_{\b}$, then  $q \circ (\pi_{\a}\circ f) \circ i_{\b}' : (F_{\b})' \to F_{\b}$
 would be an order injection,  violating the order simplicity of $F_{\b}$.

 In particular, we see that $2^{\a}$ is not isomorphic to  $2^{\b}$ and
 $F_{\a}$ not isomorphic to  $F_{\b}$. By Lemma \ref{lem3.1}  a homeomorphism between connected LOTS is either an order isomorphism or
 an order* isomorphism. Since $F_{\a}$ is symmetric, an order* isomorphism would yield an order isomorphism. Thus $F_{\a}$ not homeomorphic to  $F_{\b}$.

 \end{proof} \vspace{.5cm}

 \noindent {\bfseries Remarks:} For $\a$ and $\b$ countable, Maurice \cite{M} proved that $2^{\a}$ is not isomorphic to  $2^{\b}$.

 The zero-dimensional LOTS $2^{\om}$ is isomorphic to the Cantor Set and the connected quotient
 $F_{\om}$ is isomorphic to the unit interval in $\R$.  Hence, $2^{\om}$ embeds in $F_{\om}$ and, in particular, $2^{\om}$ and $F_{\om}$
 have the same size.

\vspace{1cm}

\section{\textbf{Appendix: Treybig's Homogeneity Theorem}}

In this Appendix we present the proof of Treybig's Homogeneity Theorem \ref{theotreybig}. \vspace{.5cm}

 Let $X$ be a LOTS and $H_+(X)$ be the group of order automorphisms.

 \begin{lem}\label{lemA01} If $f_1, f_2 \in H_+(X)$ and for some $a \in X$, $f_1(a) = f_2(a)$ then
 $f_3 \in H_+(X)$ where $f_3(x) = \begin{cases} f_1(x) \ \text{for} \ x \le a,\\ f_2(x) \ \text{for} \ x \ge a. \end{cases} $\end{lem}

 \begin{proof} Let $ g = f_2^{-1} \circ f_1 \in H_+$ so that $g(a) = a$. Hence $x \ge a \ \Leftrightarrow \ g(x) \ge a$. So
 $g'$ defined to be $g$ on $(-\infty,a]$ and the identity on $[a,\infty)$ is in $H_+$ and $f_3 = f_2 \circ g'$.

 \end{proof} \vspace{.5cm}

 \begin{lem}\label{lemA02} If $f_1, f_2 \in H_+(X)$ and for some $a \in X$, $f_1(a) > f_2(a)$, then there exists $f_3 \in H_+(X)$
 which equals $f_1$ on an open interval containing $a$ and such that $f_3(x) \ge f_2(x)$ for all $x$. \end{lem}

  \begin{proof} As before let $ g = f_2^{-1} \circ f_1 \in H_+$ so that $g(a) > a$. Inductively, $k \mapsto g^k(a)$ defines an increasing
  bi-infinite sequence for $k \in \Z$ and $g(x) > x$ for all
  $x \in J = \bigcup_{k \in \Z} [g^k(a),g^{k+1}(a)] = \bigcup_{k \in \Z} (g^{k-1}(a),g^{k+1}(a)) $.  So $J$ is an open convex set
  which contains $a$ and $g(J) = J$. Points above $J$ are mapped above $J$ and the points below $J$ are mapped below $J$.
  So $g'$ defined to be $g$ on $J$ and the identity otherwise is in $H_+$. Let $f_3 = f_2 \circ g'$.

 \end{proof} \vspace{.5cm}

 \begin{lem}\label{lemA03} Assume $X$ is complete. If there exists a sequence $\{ f_n \} \subset H_+$ which converges pointwise to the identity on $X$
 and such that $f_n(x) < x$ for all $x \in X$ and $n \in \N$, then
 for each $u \in X$ $\{ f_n^k(u) : n \in \N, k \in \Z \}$ is dense in $X$. In particular, $X$ is separable. \end{lem}

 \begin{proof} If $\{ f_n^k(u) \}$ were bounded below for some $n$ then the infimum would be a fixed point of $f_n$, contradicting the assumption
 that $f_n(x) < x$ for all $x$. Similarly, the sequence is not bounded above for any $n$.

 Let $a < b \in X$. By pointwise convergence
 there exists $n \in \N$ such that $a < f_n(b) < b$.

 Because the orbits are unbounded the set $K = \{ f_n^j(u) : j \in \Z$  and $f_n^j(u) < b \}$ is nonempty. Let $z$ be its supremum so that $z \le b$.
There exists $j$ such that $f_n(z) < f_n^j(u) \le z$. If $z \le a$ then $z < f_n^{j-1}(u) \le f_n^{-1}(z) \le f_n^{-1}(a) < b$. That is,
$f_n^{j-1}(u) \in K$ violating the definition of $z$.

Since $a < z$ there exists $f_n^k(u) \in K$ with $a < f_n^k(u)$.  That is, $f_n^k(u)$ is in the open interval $(a,b)$. Density follows.

\end{proof} \vspace{.5cm}

\noindent {\bfseries Remark:}  Completeness is needed. Let $X = \R \times \R$ which is not complete and not separable.
Let $f_n(x_1,x_2) = (x_1,x_2 - \frac{1}{n})$. \vspace{.5cm}

For the transitive LOTS $\Z$ for each pair $a, b \in \Z$ there is a unique element $f \in H_+(\Z)$ such that $f(a) = b$.
This \emph{translational uniqueness } does not occur for any nontrivial connected  LOTS.

\begin{prop}\label{propA04} If $X$ is a nontrivial connected, transitive LOTS, then it is not true that
for every $a, b \in X$ there is a unique element $f \in H_+(X)$ such that $f(a) = b$.\end{prop}

 \begin{proof} Assume instead that translational uniqueness holds for $X$. This implies that if $f, g \in H_+$ and $f(a) = g(a)$ for some $a$
 then $f(x) = g(x)$ for all $x$. It follows that if $f(a) > g(a)$, then $f(x) > g(x)$ for all $x$, because if $f(b) \le g(b)$ for some $b$
 then $f(c) = g(c)$ for some $c $ between $a$ and $b$, because $X$ is connected.

 Let $\{ a_n \}$ be an increasing sequence in $X$ converging to $a$. Let $f_n \in H_+$ with $f_n(a) = a_n$. Hence, for all $x$
 $f_n(x) < f_{n+1}(x) < x = 1(x)$. Where $1 \in H_+$ is the identity map.

 Claim: For all $x$ $\{ f_n(x) \}$ converges to $x$.

If not, then there is some $b$ such that $c = sup \{ f_n(b) \} < b$. Fix $g \in H_+$ with $g(b) = c$.
  Since for all $n$, $f_n(b) < c = g(b) < 1(b)$,
 it follows that $f_n(x) < g(x) < x$ for all $x$.  Applied with $x = a$ this contradicts the convergence of $\{ a_n = f_n(a) \}$ to $a$.

 Thus, $\{ f_n \}$ converges to $1$ pointwise.

 From Lemma \ref{lemA03} it follows that $X$ is separable and so, since it is connected and unbounded, it is isomorphic to $\R$.  However,
 translational uniqueness clearly does not hold for $\R$.

\end{proof} \vspace{.5cm}

\begin{cor}\label{corA05} If $X$ is a nontrivial connected, transitive LOTS, then one of the following holds.
\begin{itemize}
\item[(i)] For all $a \in X$ there exists $f \in H_+(X)$ and $b \in X$ with $a < b$ such that $f(x) = x$ for all $x \le a$ and $f(x) < x$ for
all $x \in (a,b]$.
 \item[(ii)] For all $a \in X$ there exists $f \in H_+(X)$ and $b \in X$ with $b < a$ such that $f(x) = x$ for all $x \ge a$ and $f(x) > x$ for
all $x \in [b,a)$.
\end{itemize}\end{cor}

\begin{proof} From Proposition \ref{propA04} there exist $f_1, f_2 \in H_+$ and $s, t \in X$ such that $f_1(t) = f_2(t)$ and $f_1(s) > f_2(s)$.
Letting $g = f_2^{-1} \circ f_1$ we have $g(t) = t$ and $g(s) > s$.

Assume $t < s$. Let $J$ be the connected component containing $s$ of the open set $\{ x  : g(x) > x \}$ and let $t_1$ be its infimum, i.e. its
left endpoint. Hence,  $t \le t_1$ and $g(t_1) = t_1 $.  Thus, $g(x) > x$ for all $x \in (t_1,s] \subset J$.

For any $a \in X$, let $h \in H_+$ such that $h(a) = t_1$ and let $r = h^{-1} \circ g \circ h$ so that $r(a) = a$ and $r(x) > x $ for all $x \in (a,b]$ with
$b = h^{-1}(s)$. We can apply Lemma \ref{lemA01} to replace $r$ by $r'$ which agrees with $r$ on $[a,\infty)$ and with the identity on $(-\infty,a]$.

Let $f = r'^{-1}$.

This is case (i).  We similarly obtain (ii) if $s < t$.

\end{proof} \vspace{.5cm}

\begin{lem}\label{lemA06} Let $X$ be a nontrivial connected, transitive LOTS and let $a, b, c \in X$ with $a < b$. For all $x > c$ there exists
$y \in (c,x)$ and $f \in H_+$ mapping $[a,b]$ to $[c,y]$. \end{lem}

\begin{proof} Let $U$ be the set of $x > c$ such that $[a,b]$ can be mapped onto $[c,x]$ by a member of $H_+$ and let
$V$ be the set of $x < c$ such that $[a,b]$ can be similarly mapped  onto $[x,c]$. Observe that by transitivity $U$ and $V$ are nonempty.
Let $u = inf \ U$ and $v = sup \ V$. Our goal is to prove $u = c$.
Suppose instead that $c < u$.\vspace{.25cm}

{\bfseries Case 1:} Assume that (i) of Corollary \ref{corA05} holds. Choose $w_1 \in X, g \in H_+$ so that $c < w_1 < u$ and $g(x) = x$ for
$x \in (-\infty,w_1]$ and $g(x) < x$ for $x \in (w_1,w_2]$ with $w_1 < w_2 < u$. Choose $w \in (w_1,w_2)$ and $k \in H_+$ such that
$k(w) = u$. By Lemma \ref{lemA02} we can choose $k$ so that $k(x) \ge x$ for all $x \in X$. Define $u_2 = k(w_2) > k(w) = u$ and
$u_1 = k(w_1) \ge w_1 > c$.  Hence, $c < u_1 < u < u_2$. Let $f = k \circ g \circ k^{-1}$ so that $f(x) = x$ for $x \in (-\infty,u_1]$ and
$f(x) < x$ for $x \in (u_1,u_2]$. It follows that for any $y \in (u_1,u_2]$ the sequence $f^n(y)$ is decreasing with limit $u_1$.

There exists $x_0 \in (u,u_2)$ and $f_0 \in H_+$ which maps $[a,b]$ onto $[c,x_0]$. Clearly, $f^n  \circ f_0$ maps $[a,b]$ onto $[c,f^n(x_0)]$.
Because the sequence $\{ f^n(x_0) \}$ converges to $u_1$, it follows that for sufficiently large $n$, $f^n(x_0) < u$. This contradicts the
definition of $u$.\vspace{.25cm}

 {\bfseries Case 2:} Assume that (ii) of Corollary \ref{corA05} holds. From an argument similar to the one for Case 1, it follows that $v = c$.
 Choose $u_1 \in (c,u)$ and $g \in H_+$ so that $u_1 = g(c)$. Hence, $g^{-1}(c) < c$. Because $v = c$, there exists $h \in H_+$ mapping
 $[a,b]$ onto $[z,c]$ with $z$ a point of $(g^{-1}(c),c)$. Let $u_2 = g(z)$ so that $c < u_2 < u_1 < u$. $g \circ h$ maps $[a,b]$ onto $[u_2,u_1]$.
 Let $k \in H_+$ with $k(a) = c$. By definition of $u$, $k(b) \ge u$. Since $k(a) < (g \circ h)(a)$ and $k(b) > (g \circ h)(b)$,
 there exists $t \in (a,b)$ such that  $k(t) = (g \circ h)(t)$.  Apply Lemma \ref{lemA01} to define $f$ to equal $ k$ on $(-\infty,t]$ and equal to
 $g \circ h$ on $[t,\infty)$. Since $f$ maps $[a,b]$ onto $[c,u_1]$ we again contradict the definition of $u$.

\end{proof} \vspace{.5cm}

Now we complete the proof of Treybig's Theorem.

\begin{theo}\label{theoA07} If $X$ is a nontrivial connected, transitive LOTS, then $X$ is doubly transitive. \end{theo}

\begin{proof} Given $a < b$ and $c < d$. Choose $g \in H_+$ with $g(b) = d$.\vspace{.25cm}

{\bfseries Case 1:} $g(a) \le c$. There exists $h$ mapping $[a,b] $ to $[c,e]$ with $c < e < d$. Because $g(a) \le c = h(a)$ and
$g(b) = d > e = h(b)$, there exists $t \in [a,b]$ such that
$g(t) = h(t)$. If $f = h$ on $(-\infty,t]$ and $f = g$ on $[t,\infty)$ then $f$ maps $[a,b]$ to $[c,d]$. \vspace{.25cm}

{\bfseries Case 2:} $g(a) > c$. There exists $h^{-1}$ mapping $[c,d]$ onto $[a,e]$ with $a < e < b$. So $h$ maps $[a,e]$ onto
$[c,d]$. Because $h(a) = c < g(a)$ and $h(e) = d = g(b) > g(e)$, there exists $t \in (a,e)$ such that $h(t) = g(t)$. For $f$ use
$h$ on $(-\infty,t]$ and $g$ on $[t,\infty)$ to map $[a,b]$ to $[c,d]$.

\end{proof} \vspace{1cm}

 \newpage

\bibliographystyle{amsplain}

\printindex

\end{document}